\documentclass[oneside,10pt]{article}%
\usepackage{amsfonts}
\usepackage{amsmath}
\usepackage{amssymb}
\usepackage{amsthm}
\usepackage{graphicx}
\usepackage[colorlinks=true,linkcolor=RawSienna,citecolor=RawSienna,urlcolor=RawSienna,bookmarksopen=true,pdfstartview=FitB]%
{hyperref}
\usepackage[left=1.35in,right=1.35in,top=1in,bottom=1in]{geometry}
\usepackage{natbib}
\usepackage[dvipsnames]{xcolor}
\usepackage[singlespacing]{setspace}
\usepackage{float}
\usepackage{xr-hyper}%
\providecommand{\U}[1]{\protect\rule{.1in}{.1in}}
\providecommand{\U}[1]{\protect\rule{.1in}{.1in}}
\allowdisplaybreaks
\newtheorem{theorem}{Theorem}

\newtheorem{definition}{Definition}[section]
\newtheorem{example}{Example}
\newtheorem{lemma}{Lemma}

\DeclareMathOperator{\sech}{sech}
\begin{document}

\title{\vspace{-1.2cm}{Uniformly consistent proportion estimation for composite
hypotheses via integral equations: ``the case of location-shift
families''}}
\author{Xiongzhi Chen\thanks{Department of Mathematics and Statistics, Washington
State University, Pullman, WA 99164, USA; E-mail:
\texttt{xiongzhi.chen@wsu.edu}}}
\date{}
\maketitle

\begin{abstract}
We consider estimating the proportion of random variables for two types of
composite null hypotheses: (i) the means or medians of the random variables
belonging to a non-empty, bounded interval; (ii) the means or medians of the
random variables belonging to an unbounded interval that is not the whole real
line. For each type of composite null hypotheses, uniformly consistent
estimators of the proportion of false null hypotheses are constructed
for random variables whose distributions are members of a Type I
location-shift family. Further, uniformly
consistent estimators of certain functions of a bounded null on the means or
medians are provided for the random variables mentioned earlier;
these functions are continuous and of bounded variation.
The estimators are constructed via solutions to Lebesgue-Stieltjes integral
equations and harmonic analysis, do not rely on a concept of p-value, and have various applications.

\end{abstract}

\tableofcontents

\section{Introduction}

\label{sec1intro}

This is a sequel to \cite{Xiongzhi2025} on using Lebesgue-Stieltjes integral
equations to uniformly consistently estimate the proportion of false null
hypotheses and its dual, the proportion of true null hypotheses
for two commonly encountered composite null hypotheses, referred to as
``one-sided null'', i.e., a one-sided, unbounded interval to which the null
parameter belongs, and ``bounded null'', i.e., a non-empty, bounded interval
to which the null parameter belongs. The motivations behind using such
composite nulls, the important roles of these proportions in statistical
modeling and inference, the current status of works on estimating these
proportions, and the advantages of using Lebesgue-Stieltjes integral equations
over existing methods to estimate these proportions have been discussed in
detail in \cite{Chen:2018a} and \cite{Xiongzhi2025}. So, we refer readers to
these two works instead of repeating them here.

Let $\Theta_{0}$ be a subset of $\mathbb{R}$ that has a non-empty interior and
non-empty complement $\Theta_{1}=\mathbb{R}\setminus\Theta_{0}$. For each
$i\in\left\{  1,\ldots,m\right\}  $, let $z_{i}$ be a random variable with
mean or median $\mu_{i}$, such that, for some integer $m_{0}$ between $0$ and
$m$, $\mu_{i}\in\Theta_{0}$ for each $i\in I_{0,m}=\left\{  1,\ldots
,m_{0}\right\}  $ and $\mu_{i}\in\Theta_{1}$ for each $i\in I_{1,m}=\left\{
m_{0}+1,\ldots,m\right\}  $. Consider simultaneously testing the null
hypothesis $H_{i0}:\mu_{i}\in\Theta_{0}$ versus the alternative hypothesis
$H_{i1}:\mu_{i}\in\Theta_{1}$ for all $i\in\left\{  1,\ldots,m\right\}  $.
Then the proportion of true null hypothesis (``null proportion'' for short) is
defined as $\pi_{0,m}=m^{-1} \left\vert I_{0,m}\right\vert =m^{-1}\sum
_{i=1}^{m}1_{\Theta_{0}}\left(  \mu_{i}\right)  $, and the proportion of false
null hypotheses (``alternative proportion'' for short) $\pi_{1,m}=1-\pi_{0,m}%
$. In other words,%
\begin{equation}
\pi_{1,m}=m^{-1} \left\vert I_{1,m}\right\vert = m^{-1}\sum_{i=1}^{m}%
1_{\Theta_{1}}\left(  \mu_{i}\right)  . \label{defPi}%
\end{equation}
Our target is to consistently estimate $\pi_{1,m}$ as $m\rightarrow\infty$,
and we will focus on the ``bounded null'' $\Theta_{0}=\left(  a,b\right)  $
for some fixed, finite $a,b\in U$ with $a<b$ and the ``one-sided null''
$\Theta_{0}=\left(  -\infty,b\right)  $, both of which are composite nulls.
However, the strategy to be introduced next to achieve this target applies to
general $\Theta_{0}$ (and hence general $\Theta_{1}$).

\subsection{The strategy via solutions to Lebesgue-Stieltjes integral
equations}

\label{SecStrategy}

Now let us describe the strategy of estimation via solutions to
Lebesgue-Stieltjes integral equations that was inspired by the work of
\cite{Jin:2008}, proposed by \cite{Chen:2018a}, and further developed
\cite{Xiongzhi2025}. Let $\mathbf{z}=\left(  z_{1},\ldots,z_{m}\right)  ^{\top}$ and
$\boldsymbol{\mu}=\left(  \mu_{1},...,\mu_{m}\right)  ^{\top}$. Denote by
$F_{\mu_{i}}$ the CDF of $z_{i}$ for $i\in\left\{  1,\ldots,m\right\}  $ and
suppose each $F_{\mu_{i}}$ is a member of a set $\mathcal{F}$ of CDFs such
that $\mathcal{F}=\left\{  F_{\mu}:\mu\in U\right\}  $ for some non-empty,
known $U$ in $\mathbb{R}$. For the rest of the paper, we assume that each
$F_{\mu}$ is uniquely determined by $\mu$ and that $U$ has a non-empty
interior. Recall the definition of $\pi_{1,m}$ in (\ref{defPi}). The strategy
to estimate $\pi_{1,m}$ via Lebesgue-Stieltjes integral equations approximates
each indicator function $1_{\Theta_{1}}\left(  \mu_{i}\right)  $, and is
stated below.

Suppose for each fixed $\mu\in U$, we can approximate the indicator function
$1_{\Theta_{0}}\left(  \mu\right)  $ by

\begin{itemize}
\item[C1)] A ``discriminant function'' $\psi\left(  t,\mu\right)  $ satisfying
$\lim_{t\rightarrow\infty}\psi\left(  t,\mu\right)  =1_{\Theta_{0}}\left(
\mu\right)  $, and

\item[C2)] A ``matching function'' $K:\mathbb{R}^{2}\rightarrow\mathbb{R}$
that does not depend on any $\mu\in\Theta_{1}$ and satisfies the
Lebesgue-Stieltjes integral equation%
\begin{equation}
\psi\left(  t,\mu\right)  =\int K\left(  t,x\right)  dF_{\mu}\left(  x\right)
\textcolor{blue}{ = \mathbb{E}_{Z \sim F_{\mu}}\left[K\left(t,Z\right)\right]}
, \forall\mu\in U. \label{eq3}%
\end{equation}

\end{itemize}

\noindent Then the ``average discriminant function''
\begin{equation}
\varphi_{m}\left(  t,\boldsymbol{\mu}\right)  =m^{-1}\sum_{i=1}^{m}\left\{
1-\psi\left(  t,\mu_{i}\right)  \right\}  \label{eq4a}%
\end{equation}
satisfies $\lim_{t\rightarrow\infty}\varphi_{m}\left(  t,\boldsymbol{\mu
}\right)  =\pi_{1,m}$ for any fixed $m$ and $\boldsymbol{\mu}$. Further, the
``empirical matching function''
\begin{equation}
\hat{\varphi}_{m}\left(  t,\mathbf{z}\right)  =m^{-1}\sum_{i=1}^{m}\left\{
1-K\left(  t,z_{i}\right)  \right\}  \label{eq4b}%
\end{equation}
satisfies $\mathbb{E}\left\{  \hat{\varphi}_{m}\left(  t,\mathbf{z}\right)
\right\}  =\varphi_{m}\left(  t,\boldsymbol{\mu}\right)  $ for any fixed $m,t$
and $\boldsymbol{\mu}$. Namely, $\hat{\varphi}_{m}\left(  t,\mathbf{z}\right)
$ is an unbiased estimator of $\varphi_{m}\left(  t,\boldsymbol{\mu}\right)
$. We will reserve the notation $\left(  \psi,K\right)  $ for a pair of
discriminant function and matching function and the notations $\varphi
_{m}\left(  t,\boldsymbol{\mu}\right)  $ and $\hat{\varphi}_{m}\left(
t,\mathbf{z}\right)  $ as per (\ref{eq4a}) and (\ref{eq4b}) unless otherwise
noted. The pair $\left(  \psi,K\right)  $ converts proportion estimation into
solving a specific Lebesgue-Stieltjes integral equation.

When the difference
\begin{equation}
e_{m}\left(  t\right)  =\hat{\varphi}_{m}\left(  t,\mathbf{z}\right)
-\varphi_{m}\left(  t,\boldsymbol{\mu}\right)  \label{eq2d}%
\end{equation}
is small for large $t$, $\hat{\varphi}_{m}\left(  t,\mathbf{z}\right)  $ will
accurately estimate $\pi_{1,m}$. Since $\varphi_{m}\left(  t,\boldsymbol{\mu
}\right)  =\pi_{1,m}$ or $\hat{\varphi}_{m}\left(  t,\mathbf{z}\right)
=\pi_{1,m}$ rarely happens, $\hat{\varphi}_{m}\left(  t,\mathbf{z}\right)  $
usually employs an increasing sequence $\left\{  t_{m}\right\}  _{m\geq1}$
with $\lim_{m\rightarrow\infty}t_{m}=\infty$ in order to achieve consistency
in the sense that%
\begin{equation}
\left\vert \pi_{1,m}^{-1}\hat{\varphi}_{m}\left(  t_{m},\mathbf{z}\right)
-1\right\vert \rightsquigarrow0\text{ as }m\rightarrow\infty,
\label{defConsistency}%
\end{equation}
where ``$\rightsquigarrow$'' denotes ``convergence in probability''. Following
the convention set by \cite{Chen:2018a}, we refer to $t_{m}$ as the ``speed of
convergence'' of $\hat{\varphi}_{m}\left(  t_{m},\mathbf{z}\right)  $.
Throughout the paper, consistency of a proportion estimator is defined via
(\ref{defConsistency}) to accommodate the scenario $\lim_{m\rightarrow\infty
}\pi_{1,m}=0$. Further, the accuracy of $\hat{\varphi}_{m}\left(
t_{m},\mathbf{z}\right)  $ in terms of estimating $\pi_{1,m}$ and its speed of
convergence depend on how fast $\pi_{1,m}^{-1}e_{m}\left(  t_{m}\right)  $
converges to $0$ and how fast $\pi_{1,m}^{-1}\varphi_{m}\left(  t_{m}%
,\boldsymbol{\mu}\right)  $ converges to $1$. This general principle also
applies to the works of \cite{Jin:2008,Jin:2007} and \cite{Chen:2018a}.

By duality, $\psi_{m}\left(  t,\boldsymbol{\mu}\right)
=1-\varphi_{m}\left(  t,\boldsymbol{\mu}\right)  $ satisfies $\pi_{0,m}%
=\lim_{t\rightarrow\infty}\psi_{m}\left(  t,\boldsymbol{\mu}\right)
$ for any fixed $m$ and $\boldsymbol{\mu}$, and $\hat{\psi}_{m}
\left(  t,\mathbf{z}\right)  =1-\hat{\varphi}_{m}\left(  t,\mathbf{z}\right)
$ satisfies $\mathbb{E}[  \hat{\psi}_{m}\left(  t,\mathbf{z}%
\right)  ]  =\psi_{m}\left(  t,\boldsymbol{\mu}\right)  $ for
any fixed $m,t$ and $\boldsymbol{\mu}$. Moreover, $\hat{\psi}_{m}
\left(  t,\mathbf{z}\right)  $ will accurately estimate $\pi_{0,m}$ when
$e_{m}\left(  t\right)  $ is suitably small for large $t$, and the stochastic
oscillations of $\hat{\psi}_{m}\left(  t,\mathbf{z}\right)  $ and
$\hat{\varphi}_{m}\left(  t,\mathbf{z}\right)  $ are the same and is
quantified by $e_{m}\left(  t\right)  $.

\subsection{Main contributions and summary of results}

In this work, we deal with the estimation problem when the random variables
have distribution functions as members of a Type I location-shift family (to
be introduced by \autoref{DefTLS} in \autoref{SecBackground}). Our main
contributions are summarized as follows:

\begin{itemize}
\item ``Construction I'': Construction of proportion estimators for testing a
bounded null on the means or medians of random variables whose distributions
are members of a Type I location-shift family, and under independence between
these random variables the speeds of convergence and uniform consistency
classes of these estimators.

\item ``Construction II'': Construction of proportion estimators for testing a
one-sided null on the means or medians of random variables whose distributions
are members of a Type I location-shift family, and under independence between
these random variables the speeds of convergence and uniform consistency
classes of these estimators.

\item Extension of Construction I to estimate the ``proportions'' induced by a
function of the parameter that is continuous and of bounded variation (see
\autoref{SecExtensions} for details), and under independence between these
random variables the speeds of convergence and uniform consistency classes of
these estimators. This considerably extends the constructions in Section 6 of
\cite{Jin:2008} and strengthens Theorem 13 there.
\end{itemize}

Note here ``consistency'' is defined via the ``ratio'' (see the definition in
(\ref{defConsistency}) to be introduced later) rather than the ``difference''
between an estimator and the alternative proportion, in order to account for a
diminishing alternative proportion that converges to $0$ as the number $m$ of
hypotheses tends to $\infty$, since for such a proportion, an estimate that
converges to zero in probability as $m \to\infty$ is consistent in the classic
definition of ``consistency'' that is based on such difference, and
consistency in terms of (\ref{defConsistency}) implies consistency in the
classic sense.

\textcolor{blue}{As a side product, we have derived concentration inequalities for three complicated empirical processes that are induced by
unbounded or non-Lipschitz functions, and these inequalities have independent interests; see Section B.3 in the supplementary material.}

As stated in \cite{Xiongzhi2025}, we have found similar things: (i) for
estimating the proportion of false null hypotheses under independence for both
bounded and one-sided nulls, the maximal speeds of convergence of our new
estimators are of the same order as those for the proportion estimators for a
point null in \cite{Chen:2018a}; (ii) for a proportion estimator that is
constructed via a solution to some Lebesgue-Stieltjes equation as an
approximator to the indicator of (a transform of) the parameter set under the
alternative hypothesis, the sparsest alternative proportion such an estimator
is able to consistently estimate can never be of larger order than the maximal
speed of convergence of the solution to its targeted indicator function. Our
simulation studies show that the new estimators often perform much better than
the MR estimator, Storey's estimator, and proportion estimators of
\cite{Hoang:2022a,Hoang:2022b} that are based on randomized p-values. The
computer codes for implementing the proposed proportion estimators
\textcolor{blue}{are written in the \textrm{R} language and} can be found on
the author's website at \url{http://archive.math.wsu.edu/faculty/xchen/welcome.php} or \url{https://xiongzhichen.github.io/}.

\subsection{Notations and conventions}

The notations and conventions we will use throughout are stated as follows:
$C$ denotes a generic, positive constant whose values may differ at different
occurrences; $O\left(  \cdot\right)  $ and $o\left(  \cdot\right)  $ are
respectively Landau's big O and small o notations; $\mathbb{E}$ and
$\mathbb{V}$ are respectively the expectation and variance with respect to the
probability measure $\Pr$; $\mathbb{R}$ and $\mathbb{C}$ are respectively the
set of real and complex numbers; $\Re$ denotes the real part of a complex
number; $\mathbb{N}$ denotes the set of non-negative integers, and
$\mathbb{N}_{+}=\mathbb{N}\setminus\left\{  0\right\}  $; $\nu$ the Lebesgue
measure, and when no confusion arises, the usual notation $d\cdot$ for
differential will be used in place of $\nu\left(  d\cdot\right)  $; for a
real-valued (measurable) function $f$ defined on some (measurable)
$A\subseteq\mathbb{R}$, $\left\Vert f\right\Vert _{p}=\left\{  \int%
_{A}\left\vert f\left(  x\right)  \right\vert ^{p}\nu\left(  dx\right)
\right\}  ^{1/p}$ and $L^{p}\left(  A\right)  =\left\{  f:\left\Vert
f\right\Vert _{p}<\infty\right\}  $ for $1\leq p < \infty$, $\left\Vert
f\right\Vert _{\infty}$ is its essential supremum, and $\left\Vert
f\right\Vert _{\mathrm{TV}}$ is the total variation of $f$ on $A$
\textcolor{blue}{when $A$ is a non-empty closed interval of $\mathbb{R}$}; for
a set $A\subseteq\mathbb{R}^{d}$, $\left\vert A\right\vert $ is the
cardinality of $A$, and $1_{A}$ the indicator of $A$; $\partial_{\cdot}$
denotes the derivative with respect to the subscript;
\textcolor{blue}{for a set $\mathcal{S}$, $\mathcal{S}^{\mathbb{N}}$ is the
$\aleph$-Cartesian product of $\mathcal{S}$, where $\aleph$ is the cardinality
of $\mathbb{N}$, i.e., an element of $\mathcal{S}^{\mathbb{N}}$ is a countably infinite, ordered sequence of elements of $\mathcal{S}$};
$\mathcal{N}_{m}\left(  \mathbf{u},\mathbf{S}\right)  $ denotes the density
and distribution of the $m$-dimensional Gaussian random vector with mean
vector $\mathbf{u}$ and covariance matrix $\mathbf{S}$, and $\Phi$ the
cumulative distribution function (CDF) of $X\sim\mathcal{N}_{1}\left(
0,1\right)  $; \textcolor{blue}{for a Lebesgue
measurable, bivariate function $g$ defined on $\mathbb{R}^{2}$, the integral
$\int_{\mathbb{R}}dy\int_{\mathbb{R}}g\left(  x,y\right)  dx$ is understood as the
iterated integral $\int_{\mathbb{R}}dy\left(  \int_{\mathbb{R}}g\left(
x,y\right)  dx\right)  $}.

\subsection{Organization of article}

The rest of the article is organized as follows. We formulate in
\autoref{SecPre} the problem of proportion estimation, provide the needed
background, and give an overview of the main constructions used or proposed
here. We develop in \autoref{SecTypeI} uniformly consistent proportion
estimators for multiple testing the means or medians of random variables from
a Type I location-shift family, and extend in \autoref{SecExtensions} one such
estimator to estimate ``proportions'' induced by a continuous function of a
bounded null that is of bounded variation. We
\textcolor{blue}{provide in \autoref{SecNumericalStudies} a simulation study and}
end the article with a discussion in \autoref{SecConcAndDisc}. \textcolor{blue}{Some auxiliary results are provided in the appendix, and technical proofs are given in the supplementary material.}

\section{Preliminaries}

\label{SecPre}

We provide in \autoref{SecBackground} a very brief background on
location-shift families, in \autoref{secResultsPoint} the key results of
\cite{Chen:2018a} on proportion estimation for the setting of a point null
that are needed here for the constructions of proportion estimators for the
setting of composite nulls, \textcolor{blue}{and in \autoref{SecIllustration} an overview and illustration of the constructions}.

\subsection{Type I location-shift family}

\label{SecBackground}

Recall $\mathcal{F}=\left\{  F_{\mu}:\mu\in U\right\}  $ and that the CDF
$F_{\mu_{i}}$ of $z_{i}$ for each $i$ is a member of $\mathcal{F}$. First, let
us discuss the setting where $\mathcal{F}$ is a location-shift family since it
is widely used. Recall the definition of location-shift family, i.e.,
$\mathcal{F}$ is a location-shift family if and only if $z+\mu^{\prime}$ has
CDF $F_{\mu+\mu^{\prime}}$ whenever $z$ has CDF $F_{\mu}$ for $\mu,\mu
+\mu^{\prime}\in U$. Let $\hat{F}_{\mu}\left(  t\right)  =\int e^{\iota
tx}dF_{\mu}\left(  x\right)  $ be the characteristic function (CF) of $F_{\mu
}$ where $\iota=\sqrt{-1}$. We can write $\hat{F}_{\mu}=r_{\mu}e^{\iota
h_{\mu}}$, where $r_{\mu}$ is the modulus of $\hat{F}_{\mu}$ and $h_{\mu}$ is
the argument of $\hat{F}_{\mu}$ (to be determined case-wise). If $\mathcal{F}$
is a location-shift family, then $\hat{F}_{\mu}\left(  t\right)  =\hat{F}%
_{\mu_{0}}\left(  t\right)  \exp\left\{  \iota t\left(  \mu-\mu_{0}\right)
\right\}  $ for all $\mu,\mu_{0}\in U$ and $r_{\mu}$ does not depend on $\mu$.
In order to construct proportion estimators for a one-sided or bounded null or
a function of the null proportion (see \autoref{SecExtensions} for details)
using the strategy in \autoref{SecStrategy} when $\mathcal{F}$ is a
location-shift family, we need:

\begin{definition}
\label{DefTLS}$\mathcal{F}$ is a ``Type I location-shift family'' if
$\mathcal{F}$ is a location-shift family for which $\hat{F}_{0}$ has no real
zeros and $\hat{F}_{0}=r_{0}$.
\end{definition}

With respect to \autoref{DefTLS}, we have a few remarks ready. First, $\hat
{F}_{0}=r_{0}$ holds if $h_{0}\equiv0$, and \textcolor{blue}{this} implies
$h_{\mu}\left(  t\right)  =\mu t$ and $\hat{F}_{\mu}\left(  t\right)  =\hat
{F}_{0}\left(  t\right)  \exp\left(  \iota t\mu\right)  $ for all $\mu\in U$
and $t\in\mathbb{R}$. \textcolor{blue}{Second,} if $F_{0}$ has a density
function $f_{0}$, then $f_{0}$ is an even function. In particular, $\hat
{F}_{0}=r_{0}$ if and only if the CF of $F_{0}$ is real, even and
non-negative. \textcolor{blue}{Third}, the requirement ``$\hat{F}_{0}$ has no
real zeros'' facilitates division by $\hat{F}_{0}$ to construct matching
functions $K$. Note that if $\hat{F}_{0}$ has no real zeros, then $\hat
{F}_{\mu}$ has no real zeros for all $\mu\in U$. Since $r_{\mu}\equiv r_{0}$
for all $\mu\in U$ for any location-shift family, then $r_{\mu}(t)/r_{0}(t)=1$
for all $\mu\in U$ and all $t\in\mathbb{R}$ whenever $r_{0}$ has no real
zeros. The condition ``$\sup_{t \in\mathbb{R}}r_{\mu}(t)/r_{0}(t) < \infty$
for each $\mu\in U$'' ensures the validity of the generalized Riemann-Lebesgue
lemma as Theorem 3 of \cite{Costin:2016} to construct $K$. Namely,
\autoref{DefTLS} is a special case of Definition 1 of \cite{Chen:2018a} that
states ``$\mathcal{F}$ is a location-shift family with Riemann-Lebesgue type
characteristic functions (RL-CFs)'' and that is needed to construct matching
functions for a point null there. Third, the requirement ``$\hat{F}_{0}$ has
no real zeros and $\hat{F}_{0}=r_{0}$'' is not restrictive since it is
satisfied by many location-shift families used in practice, which include
Gaussian family and four other families to be provided in
\autoref{secExamples}. Since enumerating all Type I location-shift families is
equivalent to the open problem of classifying real, even functions on
$\mathbb{R}$ whose Fourier transforms are positive (which will not be pursued
here; see, e.g., \cite{Giraud:2014,Tuck:2006}), we will only provide in
\autoref{ConstructLocShift} two methods to construct these families.

As revealed by \cite{Jin:2008,Chen:2018a} and will be seen later, the key
appeal of a Type I location-shift family in constructing proportion estimators
for both a one-sided null and a bounded null respectively is that its
location-shift property utilizes the additive structure of the group $\left(
\mathbb{R},+\right)  $, i.e., $\mathbb{R}$ with addition ``$+$'', and couples
very well with convolution on an additive group, Fourier transform, and
Dirichlet integrals to provide solutions to the Lebesgue-Stieltjes integral
equation (\ref{eq3}).

\subsection{Constructions of proportion estimators for a point null}

\label{secResultsPoint}

Since the bounded null $\Theta_{0}=\left(  a,b\right)  $ and the one-sided
null $\Theta_{0}=\left(  -\infty,b\right)  $ have two finite boundary points
$a$ and $b$, in order to consistently estimate $\pi_{1,m}=m^{-1}\sum_{i=1}%
^{m}1_{\Theta_{1}}\left(  \mu_{i}\right)  $ using the strategy in
\autoref{SecStrategy} for these two composite hypotheses, we can use the
proportion estimators of \cite{Chen:2018a} for a point null to specifically
account for the proportion of $\mu_{i}$'s that are equal to $a$ or $b$ when
$\mathcal{F}$ is a location-shift family. Consider a point null $\Theta
_{0}=\left\{  \mu_{0}\right\}  $ for a fixed $\mu_{0}\in U$. We can state the
constructions of \cite{Chen:2018a} of uniformly consistent estimators of
$\pi_{1,m}$ for the point null when $\mathcal{F}$ is a location-shift family
with RL-CFs as follows in terms of a discriminant function
and a matching function that satisfy conditions C1) and C2):

\begin{theorem}
\label{ThmPoinNull}Let $\omega$ be an even, bounded, probability density
function on $\left[  -1,1\right]  $ and $\mathcal{F}$ be a Type I
location-shift family. For $\mu^{\prime}\in U$, define
\[
K_{1,0}\left(  t,x;\mu^{\prime}\right)  ={\int_{\left[  -1,1\right]  }}%
\dfrac{\omega\left(  s\right)  \cos\left\{  ts\left(  x-\mu^{\prime}\right)
\right\}  }{r_{\mu^{\prime}}\left(  ts\right)  }ds
\]
and let%
\[
\psi_{1,0}\left(  t,\mu;\mu^{\prime}\right)  =\int K_{1,0}\left(
t,x;\mu^{\prime}\right)  dF_{\mu}\left(  x\right)  .
\]
Then
\[
\psi_{1,0}\left(  t,\mu;\mu^{\prime}\right)  ={\int_{\left[  -1,1\right]  }%
}\omega\left(  s\right)  \cos\left\{  ts\left(  \mu-\mu^{\prime}\right)
\right\}  ds.
\]
For the point null\ $\Theta_{0}=\left\{  \mu_{0}\right\}  $ with $\mu_{0}\in
U$, $\left(  \psi,K\right)  =\left(  \psi_{1,0}\left(  t,\mu;\mu_{0}\right)
,K_{1,0}\left(  t,x;\mu_{0}\right)  \right)  $. In particular, $\psi_{1,0}\left(  t,\mu_{0};\mu
_{0}\right)  =1$ for all $t$.
\end{theorem}

In this work, we will further assume that $\omega$ is of bounded variation
unless otherwise noted. For example, the triangular density $\omega\left(
s\right)  =\left(  1-\vert s\vert\right)  1_{\left[  -1,1\right]  }\left(
s\right)  $ or the uniform density $\omega\left(  s\right)  =0.5
\times1_{\left[  -1,1\right]  }\left(  s\right)  $ can be used.
\autoref{ThmPoinNull} will be used by the constructions to be introduced in
\autoref{SecTypeI} and \autoref{SecExtensions}.

\subsection{Overview and illustration of constructions}

\label{SecIllustration}

We give an overview on the constructions provided in \autoref{secResultsPoint}
and to be introduced in later sections. Even though conceptually these
constructions can be more elegantly stated and better understood in terms of
complex analysis (as their proofs reveal), we describe them using terms of
real analysis wherever feasible. Let $\ast$ denote the additive convolution
with respect to $\left(  \mathbb{R},+\right)  $. For two functions $\tilde{f}$
and $\tilde{g}$, recall their additive convolution as $( \tilde{f}\ast
\tilde{g}) \left(  z\right)  =\int\tilde{g}\left(  y\right)  \tilde{f}\left(
z-y\right)  dy$; see, e.g., Section 1.1.6 of \cite{rudin2017fourier} for such
a definition. There are several conventions on the definition of Fourier
transform $\mathcal{H}$ on the group $\left(  \mathbb{R},+\right)  $, and here
we use $\left(  \mathcal{H}f\right)  \left(  t\right)  =\int e^{-\iota
tx}f\left(  x\right)  dx$ and define $\mathcal{H}_{1}$ as $\left(
\mathcal{H}_{1}f\right)  \left(  t\right)  =\left(  2\pi\right)  ^{-1}\int
f\left(  x\right)  e^{\iota tx}dx$ for $f\in L^{1}\left(  \mathbb{R}\right)
$; see, e.g., the definition of Fourier transform in Section 1.2.3 and the
examples on page 13 of \cite{rudin2017fourier}. If needed, for $\mathcal{H}f$
or $\mathcal{H}_{1}f$, we will write out the argument of $f$ that is
integrated out by $\mathcal{H}$ or $\mathcal{H}_{1}$.

Also, we need two Dirichlet integrals%

\[
\mathcal{D}_{1}\left(  t,\mu;a,b\right)  =\frac{1}{\pi}\int_{\left(
\mu-b\right)  t}^{\left(  \mu-a\right)  t}\frac{\sin y}{y}dy\text{ \ \ and
\ }\mathcal{D}_{2}\left(  t,\mu;b\right)  =\frac{1}{\pi}\int_{0}^{t}\frac
{\sin\left\{  \left(  \mu-b\right)  y\right\}  }{y}dy
\]
and%

\begin{equation}
\mathcal{D}_{1,\infty}\left(  \mu;a,b\right)  =\lim_{t\rightarrow\infty
}\mathcal{D}_{1}\left(  t,\mu;a,b\right)  =\left\{
\begin{array}
[c]{lll}%
1 & \text{if} & a<\mu<b\\
2^{-1} & \text{if} & \mu=a\text{ or }\mu=b\\
0 & \text{if} & \mu<a\text{ or }\mu>b
\end{array}
\right.  , \label{EqDirichlet1}%
\end{equation}
and
\begin{equation}
\mathcal{D}_{2,\infty}\left(  \mu;b\right)  =\lim_{t\rightarrow\infty
}\mathcal{D}_{2}\left(  t,\mu;b\right)  =\left\{
\begin{array}
[c]{lll}%
2^{-1} & \text{if} & \mu>b\\
0 & \text{if} & \mu=b\\
-2^{-1} & \text{if} & \mu<b
\end{array}
\right.  . \label{EqDB1}%
\end{equation}
The identities (\ref{EqDirichlet1}) and (\ref{EqDB1})
\textcolor{blue}{are a consequence of} \autoref{lm:Dirichlet} in
the supplementary material.

Let $f_{\mu}$ be the density of $F_{\mu}$ with respect to the Lebesgue
measure. Since $\mathcal{F}$ is a Type I location-shift family, then $f_{0}$
is an even function, $\hat{F}_{0} \equiv r_{0}$, and $\hat{F}_{\mu}$ has no
real zeros for all $\mu\in U$. For the point null $\Theta=\left\{  \mu
_{0}\right\}  $, the pair $\left(  \psi_{1,0},K_{1,0}\right)  $ given by
\autoref{ThmPoinNull} satisfies the following two identities: $\psi
_{1,0}\left(  t,\mu;\mu_{0}\right)  =\left(  \mathcal{H}\omega\right)  \left(
t\left(  \mu-\mu_{0}\right)  \right)  $ since $\omega$ is even, i.e.,
$\psi_{1,0}\left(  t,\mu;\mu_{0}\right)  $ is the Fourier transform of
$\omega$ evaluated at $t\left(  \mu-\mu_{0}\right)  $, and%
\begin{align}
\psi_{1,0}\left(  t,\mu;\mu_{0}\right)   &  =\int K_{1,0}\left(  t,x;\mu
_{0}\right)  dF_{\mu}\left(  x\right)  =\int K_{1,0}\left(  t,x;\mu
_{0}\right)  f_{0}\left(  x- \mu\right)  dx,
\label{ExplainPointNullLocationShift}%
\end{align}
where the second identity only uses the location-shift property of
$\mathcal{F}$. Since $f_{0}$ is even, $f_{0}\left(  y-\mu\right)
=f_{0}\left(  \mu-y\right)  $, and (\ref{ExplainPointNullLocationShift}) is
just
\textcolor{blue}{$\psi_{1,0}\left(  t,\mu;\mu_{0}\right)=\left( K_{1,0}\left(t,\cdot;\mu_0\right)\ast f_{0}\right)\left(\mu\right)$},
and $K_{1,0}$ is the inverse Fourier transform of $\mathcal{H}\psi
_{1,0}/\mathcal{H}f_{0}$.

For the bounded null $\Theta=\left(  a,b\right)  $, we first let $\tilde
{g}\left(  t\right)  =\int_{a}^{b}\exp\left(  -\iota yt\right)  dy$, i.e.,
$\tilde{g}$ is the Fourier transform of the indicator function $1_{\left(
a,b\right)  }\left(  y\right)  $ of $\left(  a,b\right)  $, and then find
$K_{1}\left(  t,x\right)  $ such that%
\begin{equation}
\begin{aligned} \psi_{1}\left( t,\mu\right) & =\int K_{1}\left( t,x\right) dF_{\mu }\left( x\right) =\int K_{1}\left( t,y\right) f_{0}\left( y-\mu\right) dy\\ & = \mathcal{D}_{1}\left( t,\mu;a,b\right) =\frac{1}{2\pi}\int_{-t}^{t}\exp\left( \iota\mu s\right) \tilde{g}\left( s\right) ds, \label{ExplainLocationShiftBoundedNull}\end{aligned}
\end{equation}
where the second identity uses the location-shift property of $\mathcal{F}$,
and the last \textcolor{blue}{is a consequence of} \autoref{LmDirichlet} in
the supplementary material. Again since $f_{0}$ is even,
(\ref{ExplainLocationShiftBoundedNull}) is just
\[
\psi_{1}=K_{1}\ast f_{0}= \mathcal{H}_{1}\left(  1_{\left[  -t,t\right]  }
\times\mathcal{H}1_{\left(  a,b\right)  } \right)  ,
\]
and $K_{1}$ is the inverse Fourier transform of $\mathcal{H}\psi
_{1}/\mathcal{H}f_{0}$. Now, if we set $\psi\left(  t,\mu\right)  =\psi
_{1}\left(  t,\mu\right)  -2^{-1}\left\{  \psi_{1,0}\left(  t,\mu;a\right)
+\psi_{1,0}\left(  t,\mu;b\right)  \right\}  $, then $\lim_{t\rightarrow
\infty}\psi\left(  t,\mu\right)  =1_{\left(  a,b\right)  }\left(  \mu\right)
$ due to (\ref{EqDirichlet1}). So, if we set $K\left(  t,x\right)
=K_{1}\left(  t,x\right)  -2^{-1}\left\{  K_{1,0}\left(  t,x;a\right)
+K_{1,0}\left(  t,x;b\right)  \right\}  $, then the pair $\left(
\psi,K\right)  $ solves the Lebesgue-Stieltjes equation (\ref{eq3}).

For the one-sided null $\Theta=\left(  -\infty,b\right)  $ with some fixed,
finite $b$, we can assume $b=0$ without loss of generality. Then we need to
find $K_{1}\left(  t,x\right)  $ such that%
\begin{align}
\psi_{1}\left(  t,\mu\right)   &  =\int K_{1}\left(  t,x\right)  dF_{\mu
}\left(  x\right)  =\int K_{1}\left(  t,x\right)  f_{0}\left(  x-\mu\right)
dx= \mathcal{D}_{2}\left(  t, \mu;0\right) \label{EPLocationShiftOneSide}\\
&  =\frac{1}{2\pi}\int_{0}^{t}dy\int_{-1}^{1}\mu\exp\left(  \iota
ys\mu\right)  ds ,\nonumber
\end{align}
where the second identity uses the location-shift property of $\mathcal{F}$,
and the fourth is \textcolor{blue}{a consequence of} \autoref{LmDirichlet} in
the supplementary material. Since $f_{0}$ is even, then%
\begin{equation}
\psi_{1}\left(  t,\mu\right)  =\left(  K_{1}\ast f_{0}\right)  \left(
t,\mu\right)  =\int_{0}^{t}\mu\times\left(  \mathcal{H}_{1} 1_{\left[
-1,1\right]  }\right)  \left(  y\mu\right)  dy. \label{EPLocShiftOneSided}%
\end{equation}
However, (\ref{EqDB1}), (\ref{EPLocationShiftOneSide}) and
\autoref{lm:Dirichlet} in the supplementary material imply that there exists some
constant $t_{\ast}\geq2$ such that $\psi_{1}\left(  t,\mu\right)  \geq1/4$
uniformly for all sufficiently large $t\geq t_{\ast}$ and all $\mu\geq1$,
i.e., $\psi_{1}\left(  t,\cdot\right)  \notin L^{1}\left(  \mathbb{R}\right)
$ for each such $t$. Hence the Fourier transform of $\psi_{1}\left(
t,\mu\right)  $ in (\ref{EPLocShiftOneSided}) in the argument $\mu$ is
undefined for each fixed, sufficiently large $t$, and we cannot apply Fourier
transform to $\psi_{1}\left(  t,\mu\right)  $ in the argument $\mu$ to obtain
the Fourier transform of $K_{1}$ (and then invert the Fourier transform of
$K_{1}$ to obtain $K_{1}$). This is the key difference between constructing
$K$ for a one-sided null and those for a point null or bounded null, and
requires additional innovation. Instead, finding $K_{1}$ for a one-sided null
employs a different method, and $K_{1}$ is the real part of%
\[
K_{1}^{\dagger}\left(  t,x\right)  =\frac{1}{2\pi}\int_{0}^{1}dy\int_{-1}%
^{1}\frac{1}{\iota y}\left(  \frac{d}{ds}\frac{\exp\left(  \iota tysx\right)
}{r_{0}\left(  tys\right)  }\right)  ds;
\]
see details in the proof of \autoref{ThmTypeIOneSided}. If we set $\psi\left(
t,\mu\right)  =2^{-1}-\psi_{1}\left(  t,\mu\right)  -2^{-1}\psi_{1,0}\left(
t,\mu;0\right)  $, then $\lim_{t\rightarrow\infty}\psi\left(  t,\mu\right)
=1_{\left(  -\infty,0\right)  }\left(  \mu\right)  $ due to (\ref{EqDB1}). So,
if we set $K\left(  t,x\right)  =2^{-1}-K_{1}\left(  t,x\right)
-2^{-1}K_{1,0}\left(  t,x;0\right)  $, then the pair $\left(  \psi,K\right)  $
solves the Lebesgue-Stieltjes equation (\ref{eq3}).

\begin{figure}[th]
\centering
\includegraphics[height=0.35\textheight,width=1\textwidth]{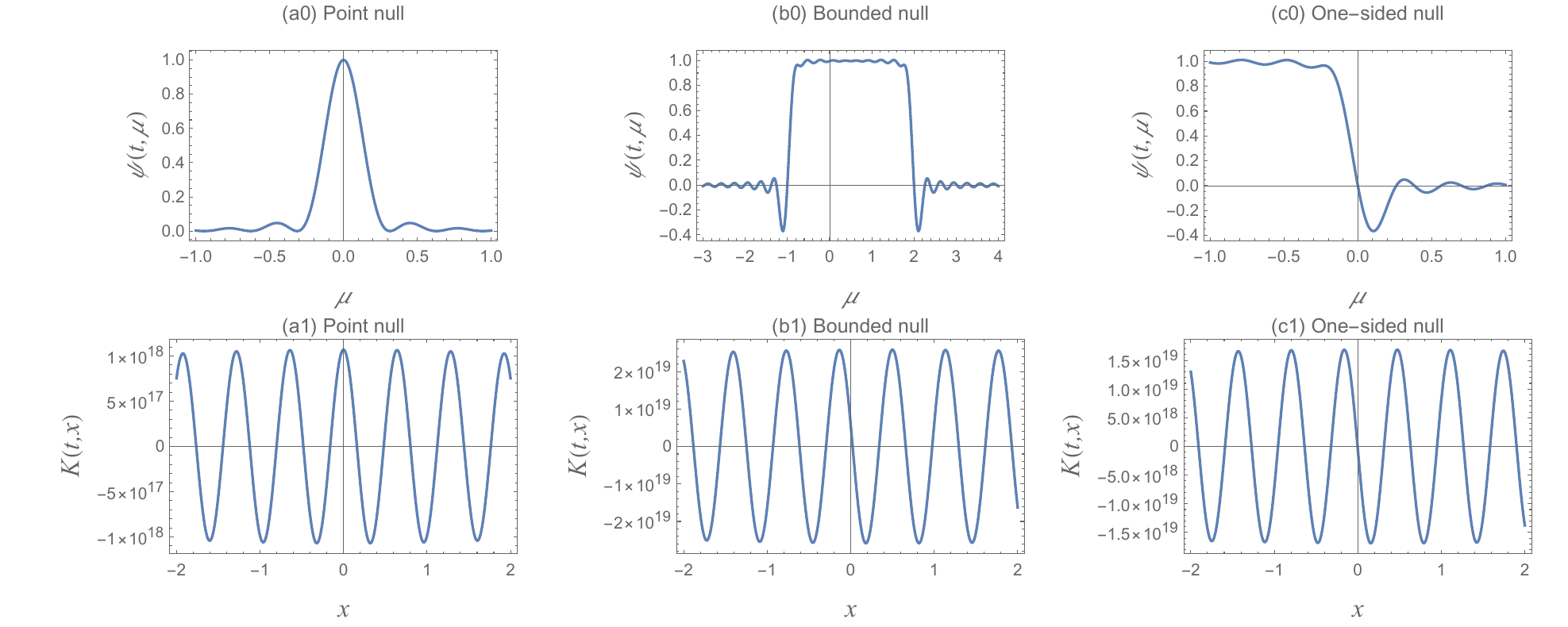}
\vspace{-0.4cm} \caption{We have set $t=20$ in each discriminant function
$\psi\left(  t,\mu\right)  $, and $t=10$ in each matching function $K\left(
t,x\right)  $ since the magnitude of $K$ is quite large when $t$ is large. For
$t=20$, $\psi\left(  t,\mu\right)  $ already approximates relatively well
$1_{\Theta_{0}}\left(  \mu\right)  $ for each of the three types of nulls.}%
\label{figIllustrations}%
\end{figure}

In \autoref{figIllustrations}, we provide visualizations of the pair $\left(
\psi,K\right)  $ for the point null $\Theta_{0}=\left\{  0\right\}  $, bounded
null $\Theta_{0}=\left(  -1,2\right)  $, and one-sided null $\Theta
_{0}=\left(  -\infty,0\right)  $, respectively, when each $F_{\mu}$ is the CDF
of $X_{\mu}\sim\mathcal{N}_{1}\left(  \mu,1\right)  $ and $\omega\left(
s\right)  =\left(  1-\left\vert s\right\vert \right)  1_{\left[  -1,1\right]
}\left(  s\right)  $, i.e., the triangular density on $\left[  -1,1\right]  $
is used. Details on how we numerically compute $\left(  \psi,K\right)  $ are
given in the beginning of \autoref{SecNumericalStudies} and the appendix. We see that $K$ is an
oscillator in each construction.

\section{Constructions for Type I location-shift families}

\label{SecTypeI}

Recall $\mathcal{F}=\left\{  F_{\mu}:\mu\in U\right\}  $. We will refer to as
\textquotedblleft Construction I\textquotedblright\ the construction of
estimators of $\pi_{1,m}$ for the bounded null $\Theta_{0}=\left(  a,b\right)
\cap U$ for fixed, finite $a,b\in U$ with $a<b$, and as \textquotedblleft
Construction II\textquotedblright\ the construction of estimators of
$\pi_{1,m}$ for the one-sided null $\Theta_{0}=\left(  -\infty,b\right)  \cap
U$ for a fixed, finite $b\in U$. By default, each $F_{\mu}\in\mathcal{F}$ is
uniquely determined by $\mu$. So, all $F_{\mu},\mu\in U$ have the same scale
parameter, if any, when $\mathcal{F}$ is a location-shift family.

\subsection{The case of a bounded null}

\label{SecConstructionIV}

Construction I for the bounded null for a Type I location-shift family is
provided by

\begin{theorem}
\label{ThmIV}Assume $\mathcal{F}$ is a Type I location-shift family. Set
\begin{equation}
K_{1}\left(  t,x\right)  =\frac{t}{2\pi}\int_{a}^{b}dy\int_{\left[
-1,1\right]  }\frac{\cos\left\{  ts\left(  x-y\right)  \right\}  }
{r_{0}\left(  ts\right)  }ds. \label{eq13e}%
\end{equation}
Then
\begin{equation}
\psi_{1}\left(  t,\mu\right)  =\int K_{1}\left(  t,x\right)  dF_{\mu}\left(
x\right)  =\frac{1}{\pi}\int_{\left(  \mu-b\right)  t}^{\left(  \mu-a\right)
t}\frac{\sin y}{y}dy, \label{GaussianBoundedpsiFunc}%
\end{equation}
and the desired $\left(  \psi,K\right)  $ is
\begin{equation}
\left\{
\begin{array}
[c]{l}%
K\left(  t,x\right)  =K_{1}\left(  t,x\right)  -2^{-1}\left\{  K_{1,0}\left(
t,x;a\right)  +K_{1,0}\left(  t,x;b\right)  \right\} \\
\psi\left(  t,\mu\right)  =\psi_{1}\left(  t,\mu\right)  -2^{-1}\left\{
\psi_{1,0}\left(  t,\mu;a\right)  +\psi_{1,0}\left(  t,\mu;b\right)  \right\}
\end{array}
\right.  . \label{IV-a}%
\end{equation}

\end{theorem}

Note $\lim_{t\rightarrow\infty}\psi\left(  t,\mu\right)  =1_{\left(
a,b\right)  }\left(  \mu\right)  $ for $\psi\left(  t,\mu\right)  $ in
(\ref{IV-a}). Define%
\begin{equation}
g\left(  t,\mu\right)  =\int_{\left[  -1,1\right]  }\frac{1}{r_{\mu}\left(
ts\right)  }ds\text{ for}\ \mu\in U,t\in\mathbb{R} \label{eq3b}%
\end{equation}
and%
\begin{equation}
u_{m}=\min_{\tau\in\left\{  a,b\right\}  }\min_{\left\{  j:\mu_{j}\neq
\tau\right\}  }\left\vert \mu_{j}-\tau\right\vert . \label{eq10b}%
\end{equation}
Then $g$ measures the average reciprocal modulus of the CF $\hat{F}_{\mu}$ of
$F_{\mu}$ on $\left[  -1,1\right]  $. As already shown by \cite{Chen:2018a},
$g$ plays a critical role in bounding the oscillation of $e_{m}\left(
t\right)  $ for the estimator $\hat{\varphi}_{m}\left(  t,\mathbf{z}\right)
$, and a larger $g$ usually gives a smaller maximal speed of convergence for
$\hat{\varphi}_{m}\left(  t,\mathbf{z}\right)  $ to achieve consistency. In
contrast, $u_{m}$ measures the minimal distance from $\mu_{j}$ to the boundary
points $a$ and $b$ of the bounded null $\Theta_{0}$, and a suitable magnitude
for $u_{m}$ is needed for the estimator induced by $K_{1,0}\left(
t,x;a\right)  $ and $K_{1,0}\left(  t,x;b\right)  $ in (\ref{IV-a}) to
consistently estimate the proportion of $\mu_{j}$'s that are equal to $a$ or
$b$; see \autoref{ThmIVa} below and Theorems 2 and 3 of \cite{Chen:2018a}.

For the concept of uniform consistency, we adopt the following definition and
interpretation from \cite{Xiongzhi2025}:

\begin{definition}
\label{DefUniformConsistency}Given a family $\mathcal{F}$, consider a
non-empty set $\mathcal{Q}\left(  \mathcal{F}\right)
\textcolor{blue}{=\mathcal{Q}_1\left(  \mathcal{F}\right) \times \mathcal{Q}_2\left(  \mathcal{F}\right)}
\subseteq\mathbb{R}^{\mathbb{N}}\times
\textcolor{blue}{\left(\mathcal{P}\left(\mathbb{R}\right)\right)^{\mathbb{N}}}$
and its coordinate projection $\mathcal{Q}_{m}\left(  \mathcal{F}\right)
=\mathcal{Q}_{m,1}\left(  \mathcal{F}\right)  \times\mathcal{Q}_{m,2}\left(
\mathcal{F}\right)  $ onto $\mathbb{R}^{m}\times
\textcolor{blue}{\mathcal{P}\left(\mathbb{R}\right)}$ for each $m\in
\mathbb{N}_{+}$,
\textcolor{blue}{where $\mathcal{P}\left(\mathbb{R}\right)$ is the power set of $\mathbb{R}$, $\mathcal{Q}_{m,1}\left(  \mathcal{F}\right) $ the first $m$ coordinates of $\mathcal{Q}_1\left(  \mathcal{F}\right)$, and $\mathcal{Q}_{m,2}\left(
\mathcal{F}\right) $ the $m$-th coordinate of $\mathcal{Q}_2\left(  \mathcal{F}\right)$}.
Then $\mathcal{Q}\left(  \mathcal{F}\right)  $ is said to be a
\textquotedblleft uniform consistency class\textquotedblright\ for
$\hat{\varphi}_{m}\left(  t,\mathbf{z}\right)  $ and $\hat{\varphi}_{m}\left(
t,\mathbf{z}\right)  $ is said to be \textquotedblleft uniformly
consistent\textquotedblright\ on $\mathcal{Q}\left(  \mathcal{F}\right)  $ if
\begin{equation}
\sup\nolimits_{\boldsymbol{\mu}\mathcal{\in Q}_{m,1}\left(  \mathcal{F}%
\right)  }\left\vert \pi_{1,m}^{-1}\sup\nolimits_{t\in\mathcal{Q}_{m,2}\left(
\mathcal{F}\right)  }\hat{\varphi}_{m}\left(  t,\mathbf{z}\right)
-1\right\vert \rightsquigarrow0\,\ \text{as }m\rightarrow\infty.
\label{defUCC}%
\end{equation}

\end{definition}

Namely, a uniform consistency class for an estimator $\hat{\varphi}_{m}\left(
t,\mathbf{z}\right)  $ is the asymptotic setting (as $m\rightarrow\infty$) for
the family $\mathcal{F}$ under which $\hat{\varphi}_{m}\left(  t,\mathbf{z}%
\right)  $ at its maximal speed of convergence (that is indicated by
$\sup\nolimits_{t\in\mathcal{Q}_{m,2}\left(  \mathcal{F}\right)  }$ in
\autoref{DefUniformConsistency}) still maintains consistency uniformly over
all settings of all $F_{\mu_{i}},i=1,\ldots,m$ (that is indicated by
$\sup\nolimits_{\boldsymbol{\mu}\mathcal{\in Q}_{m}\left(  \mathcal{F}\right)
}$ in \autoref{DefUniformConsistency}) and hence uniformly over all settings
of the alternative proportion $\pi_{1,m}$ that are induced by such $F_{\mu
_{i}},i=1,\ldots,m$.

For a constant $\rho\geq0$, define
\begin{equation}
\mathcal{B}_{\textcolor{blue}{j},m}\left(  \rho\right)  =\left\{
\boldsymbol{\mu}\in\mathbb{R}^{m}:m^{-1}\sum\nolimits_{i=1}^{m}\left\vert
\mu_{i}\right\vert ^{\textcolor{blue}{j}}\leq\rho\right\}  . \label{eq19f}%
\end{equation}
Even when $\rho>0$, the set $\mathcal{B}_{1,m}\left(  \rho_{m}\right)  $
allows $\lim_{m\rightarrow\infty}\max_{1\leq j\leq m}\left\vert \mu
_{j}\right\vert =\infty$ and was used by \cite{Jin:2008} and others in various
estimation or testing problems. For $\mu\in U$, let $X_{\left(  \mu\right)  }$
have CDF\ $F_{\mu}$. The following theorem provides uncertainty quantification
of the estimation error $e_{m}\left(  t\right)  $ of $\hat{\varphi}_{m}\left(
t,\mathbf{z}\right)  $ and a uniform consistency class for $\hat{\varphi}%
_{m}\left(  t,\mathbf{z}\right)  $.

\begin{theorem}
\label{ThmIVa}Assume that $\mathcal{F}$ is a Type I location-shift family for
which $\int\left\vert x\right\vert ^{2}dF_{\mu}\left(  x\right)  <\infty$ for
each $\mu\in U$, and that $\left\{  z_{i}\right\}  _{i=1}^{m}$ are
independent. Then%
\begin{equation}
\mathbb{V}\left\{  e_{m}\left(  t\right)  \right\}  \leq2m^{-1}g^{2}\left(
t,0\right)  \left(  \left\Vert \omega\right\Vert _{\infty}^{2}+\pi^{-2}\left(
b-a\right)  ^{2}t^{2}\right)  , \label{BndNullVarianceBnd}%
\end{equation}
and there are positive sequences $\left\{  \tau_{m}\right\}  _{m\geq
1},\left\{  \gamma_{m}\right\}  _{m\geq1}$ and $\left\{  \rho_{m}\right\}
_{m\geq1}$ such that, for all sufficiently large $m$, with probability at
least $1-o\left(  1\right)  $
\begin{equation}
\sup_{\boldsymbol{\mu}\in\mathcal{B}_{1,m}\left(  \rho_{m}\right)  }\sup
_{t\in\left[  0,\tau_{m}\right]  }\left\vert e_{m}\left(  t\right)
\right\vert \leq\left\{  \left(  b-a\right)  \frac{\tau_{m}}{2\pi}+\left\Vert
\omega\right\Vert _{\infty}\right\}  \Upsilon\left(  \gamma_{m}m^{-1/2}%
,\tau_{m}\right)  \label{BndeBoundedNullLocationShift}%
\end{equation}
where%
\begin{equation}
\Upsilon\left(  \lambda,\tilde{\tau}\right)  =2\lambda\sup_{t\in\left[
0,\tilde{\tau}\right]  }\int_{\left[  0,1\right]  }\frac{ds}{r_{0}\left(
ts\right)  }\text{ for }\lambda,\tilde{\tau}\geq0. \label{UpsilonSimple}%
\end{equation}
Further, a uniform consistency class for the estimator $\hat{\varphi}%
_{m}\left(  t,\mathbf{z}\right)  $ is%
\begin{equation}
\mathcal{Q}\left(  \mathcal{F}\right)  =\left\{
\begin{array}
[c]{l}%
\left(  2\tau_{m}+R_{1}\left(  \rho_{m}\right)  \right)  ^{-1}m^{\vartheta
-1/2}\gamma_{m}^{-1}\ln\tau_{m}\rightarrow\infty\\
\tau_{m}^{2}m^{2\vartheta}\exp\left(  -2^{-1}\gamma_{m}^{2}\right)
+m^{-2\vartheta}\gamma_{m}^{2}\ln^{-2}\tau_{m}=o\left(  1\right) \\
\tau_{m}^{-1}\sup_{\boldsymbol{\mu}\in\Omega_{m}\left(  \rho_{m}\right)  }%
\pi_{1,m}^{-1}\left(  1+u_{m}^{-1}\right)  =o\left(  1\right) \\
\tau_{m}\Upsilon\left(  \gamma_{m}m^{-1/2},\tau_{m}\right)  \sup
_{\boldsymbol{\mu}\in\Omega_{m}\left(  \rho_{m}\right)  }\pi_{1,m}%
^{-1}=o\left(  1\right)
\end{array}
\right\}  , \label{LocationBoundedNullUCS}%
\end{equation}
where $\vartheta>0$ is a constant,
\begin{equation}
\Omega_{m}\left(  \rho_{m}\right)  \subseteq\mathcal{B}_{1,m}\left(  \rho
_{m}\right)  \text{ \ and }R_{1}\left(  \rho_{m}\right)  =2\mathbb{E}\left[
\left\vert X_{(0)}\right\vert \right]  +2\max\left\{  \left\vert a\right\vert
,\left\vert b\right\vert \right\}  +2\rho_{m}. \label{UCCBoundedNullSetAndPho}%
\end{equation}
In particular, for the Gaussian family with variance $\sigma^{2}>0$, the
$\mathcal{Q}\left(  \mathcal{F}\right)  $ in (\ref{LocationBoundedNullUCS})
can be set as
\begin{equation}
\mathcal{Q}\left(  \mathcal{F}\right)  =\left\{
\begin{array}
[c]{l}%
0<\gamma<2^{-1}\sigma^{-2};\tau_{m}=\sqrt{2\gamma\ln m}\\
q>2\vartheta;0\leq\vartheta^{\prime}<\vartheta-2^{-1};R_{1}\left(  \rho
_{m}\right)  =O\left(  m^{\vartheta^{\prime}}\right) \\
\pi_{1,m}^{-1}=O\left(  \sqrt{\ln m}\right)  ;u_{m}\geq\left(  \tau
_{m}\right)  ^{-1}\ln\ln m
\end{array}
\right\}  \label{MainClass}%
\end{equation}
for constants $q$, $\gamma$, $\vartheta$ and $\vartheta^{\prime}$.
\end{theorem}

Note that $\Upsilon\left(  \lambda,\tilde{\tau}\right)  =\lambda\sup
_{t\in\left[  0,\tilde{\tau}\right]  }$ $g\left(  t,0\right)  $.
\autoref{ThmIVa} has the following two implications: (i) the variance of the
error $e_{m}\left(  t\right)  $ of the estimator $\hat{\varphi}_{m}\left(
t,\mathbf{z}\right)  $ is characterized by the average reciprocal modulus $g$
at the generating measure $F_{0}$ of the location-shift family $\mathcal{F}$,
the supremum norm $\Vert\omega\Vert_{\infty}$ of the employed density $\omega
$, the \textquotedblleft size\textquotedblright\ $b-a$ of the bounded null,
the number $m$ of hypotheses to test, and the magnitude of the tuning
parameter $t$, and (ii) the uniform, absolute deviation of $e_{m}\left(
t\right)  $ is characterized similarly as is the variance of $e_{m}\left(
t\right)  $, except that a uniformity requirement is imposed on the average
reciprocal modulus $g$ via $\Upsilon\left(  \lambda,\tilde{\tau}\right)  $ in
(\ref{UpsilonSimple}). An example of (\ref{MainClass}) is%

\begin{equation}
\mathcal{U}_{m}=\left\{  \boldsymbol{\mu}\in\mathcal{B}_{1,m}\left(  \rho
_{m}\right)  :%
\begin{array}
[c]{l}%
\vartheta=3/4;q=7/4;\tau_{m}=\sqrt{2^{-1}\sigma^{-2}\ln m}\\
0\leq\vartheta^{\prime}<1/4;R_{1}\left(  \rho_{m}\right)  \leq\tilde
{C}m^{\vartheta^{\prime}}\\
\pi_{1,m}^{-1}\leq\tilde{C}\sqrt{\ln m};u_{m}\geq\left(  \tau_{m}\right)
^{-1}\ln\ln m
\end{array}
\right\}  \label{MainClassEg}%
\end{equation}
for a constant $\tilde{C}>0$, for which $\sup\nolimits_{\boldsymbol{\mu}%
\in\mathcal{U}_{m}}\left(  \left\vert \pi_{1,m}^{-1}\sup\nolimits_{t\in
\left\{  \tau_{m}\right\}  }\hat{\varphi}_{m}\left(  t,\mathbf{z}\right)
-1\right\vert \right)  \rightsquigarrow0$ as\ $m\rightarrow\infty$. We remark
that for the Gaussian family with variance $\sigma^{2}>0$, the inequality
(\ref{BndNullVarianceBnd}) becomes%
\[
\mathbb{V}\left\{  e_{m}\left(  t\right)  \right\}  \leq\frac{32\exp\left(
t^{2}\sigma^{2}\right)  }{mt^{4}\sigma^{4}}\left(  \left\Vert \omega
\right\Vert _{\infty}^{2}+\pi^{-2}\left(  b-a\right)  ^{2}t^{2}\right)  ,
\]
and with $\tau_{m}=\sqrt{2\gamma\ln m}$ and $\gamma_{m}=\sqrt{2q\ln m}$,
(\ref{BndeBoundedNullLocationShift}) becomes
\[
\sup_{\boldsymbol{\mu}\in\mathcal{B}_{1,m}\left(  \rho_{m}\right)  }\sup
_{t\in\left[  0,\sqrt{2\gamma\ln m}\right]  }\left\vert e_{m}\left(  t\right)
\right\vert \leq\left\{  \left(  b-a\right)  \frac{\sqrt{2\gamma\ln m}}{2\pi
}+\left\Vert \omega\right\Vert _{\infty}\right\}  \frac{2\sqrt{2q}}%
{\gamma\sigma^{2}\sqrt{\ln m}}m^{\gamma\sigma^{2}-1/2}%
\]
with probability at least $1-o\left(  1\right)  $, as will be shown by the
proof of \autoref{ThmIVa}.

The uniform consistency class (\ref{LocationBoundedNullUCS}) is a general and
informative result, and is applicable whenever each CDF $F_{\mu}$ of a Type I
location-shift family has finite variance. First, it characterizes the
performance of the estimator $\hat{\varphi}_{m}\left(  t,\mathbf{z}\right)  $
at its maximal speed $\tau_{m}$ of convergence uniformly over the set
$\Omega_{m}\left(  \rho_{m}\right)  $ for $\boldsymbol{\mu}$, where
$\lim_{m\rightarrow\infty}\rho_{m}=\infty$ is allowed rather than restricting
$\rho_{m}=\rho$ for all $m\geq1$. Second, through its conditions $\tau
_{m}\Upsilon\left(  \gamma_{m}m^{-1/2},\tau_{m}\right)  \sup_{\boldsymbol{\mu
}\in\Omega_{m}\left(  \rho_{m}\right)  }\pi_{1,m}^{-1}=o\left(  1\right)  $,
we we can get a sense on the magnitude of the class of alternative proportion
$\pi_{1,m}$ for which $\hat{\varphi}_{m}\left(  t,\mathbf{z}\right)  $ is
uniformly consistent at its maximal speed $\tau_{m}$. Third, in
(\ref{LocationBoundedNullUCS}), the maximal speed of $\hat{\varphi}_{m}\left(
t,\mathbf{z}\right)  $ is set as $t=\tau_{m}$. While the condition $\tau
_{m}^{-1}\pi_{1,m}^{-1}\left(  1+u_{m}^{-1}\right)  =o\left(  1\right)  $
implies that $\tau_{m}$ should be larger if we want to be able to consistently
estimate smaller $\pi_{1,m}$, the condition $\tau_{m}\Upsilon\left(
\gamma_{m}m^{-1/2},\tau_{m}\right)  \pi_{1,m}^{-1}=o\left(  1\right)  $
implies that $\tau_{m}$ should be smaller if we want to be able to
consistently estimate smaller $\pi_{1,m}$. So, these two conditions together
balance the magnitudes of $\tau_{m}$ and $\pi_{1,m}$ in order for
$\hat{\varphi}_{m}\left(  \tau_{m},\mathbf{z}\right)  $ to achieve
consistency. Fourth, the first two rows of conditions of $\mathcal{Q}\left(
\mathcal{F}\right)  $ in (\ref{LocationBoundedNullUCS}) ensure that the
concentration inequality (\ref{BndeBoundedNullLocationShift}) for
$e_{m}\left(  t\right)  $ holds asymptotically with probability $1$, the third
ensures $\left\vert \pi_{1,m}^{-1}\hat{\varphi}_{m}\left(  t,\mathbf{z}%
\right)  -1\right\vert \rightarrow0$, and the last ensures that the upper
bound in (\ref{BndeBoundedNullLocationShift}) converges to zero, all as
$m\rightarrow\infty$.

It is informative to compare the set $\mathcal{Q}$ based on (\ref{MainClass}%
)\ for the bounded null and the uniform consistency class $\mathcal{\tilde{Q}%
}$ for a point null given by Theorem 3 of \cite{Chen:2018a}, both for Gaussian
random variables $\left\{  z_{i}\right\}  _{i=1}^{m}$ each with mean $\mu_{i}$
and standard deviation $1$. For this scenario, we set $t=\tau_{m}%
=\sqrt{2\gamma\ln m}$ in the estimators $\hat{\varphi}_{m}\left(
t,\mathbf{z}\right)  $, and $\gamma_{m}=\ln m$ for $\mathcal{\tilde{Q}}$.
Then
\begin{equation}
\tau_{m}^{-1}\left(  1+u_{m}^{-1}\right)  =o\left(  \pi_{1,m}\right)  \text{
and }\pi_{1,m}^{-1}\tau_{m}\Upsilon\left(  \gamma_{m}m^{-1/2},\tau_{m}\right)
\leq Cm^{\gamma-0.5}\pi_{1,m}^{-1}\text{ \ for }\mathcal{Q} \label{uccNow}%
\end{equation}
and%
\begin{equation}
\tau_{m}u_{m,0}\rightarrow\infty\text{ and }\pi_{1,m}^{-1}\Upsilon\left(
\gamma,\tau_{m},\gamma_{m},r_{\mu_{0}}\right)  \leq C\left(  \sqrt{2\gamma\ln
m}\right)  ^{-1}m^{\gamma-0.5}\pi_{1,m}^{-1}\text{ for }\mathcal{\tilde{Q}},
\label{ucc2018a}%
\end{equation}
where $u_{m,0}=\min\left\{  \left\vert \mu_{j}-\mu_{0}\right\vert :\mu_{j}%
\neq\mu_{0}\right\}  $. Here the derivation of (\ref{uccNow}) can be found in
the proof of \autoref{ThmIVa}, and the derivation of (\ref{ucc2018a}) can be
found in Corollary 3 of \cite{Chen:2018a}. So, the maximal speeds of
convergence of corresponding estimators have the same order as $\sqrt{\ln m}$
due to $t_{m}=\sqrt{2\gamma\ln m}$ for both $\mathcal{Q}$ and $\mathcal{\tilde
{Q}}$, even though for $\mathcal{\tilde{Q}}$ the maximal speed is $\sqrt{\ln
m}$ which is achieved when $\gamma=2^{-1}$ and $\liminf_{m\rightarrow\infty
}\pi_{1,m}>0$. However, since $\tau_{m}^{-1}\left(  1+u_{m}^{-1}\right)
=o\left(  \pi_{1,m}\right)  $ is required for $\mathcal{Q}$ as demanded by the
speed of convergence of $\varphi_{m}\left(  t,\boldsymbol{\mu}\right)  $ to
$\pi_{1,m}$, the sparsest $\pi_{1,m}$ contained in $\mathcal{Q}$ is usually
larger in order than $\tau_{m}^{-1}$ when $Cm^{\gamma-0.5}\pi_{1,m}%
^{-1}=o\left(  1\right)  $. In contrast, for $\mathcal{\tilde{Q}}$ the speed
of convergence of the corresponding $\varphi_{m}\left(  t,\boldsymbol{\mu
}\right)  $ to $\pi_{1,m}$ does not depend on $\pi_{1,m}$ but depends only on
$u_{m,0}$ via $\tau_{m}u_{m,0}\rightarrow\infty$, and thus the sparsest
$\pi_{1,m}$ contained in $\mathcal{\tilde{Q}}$ can be of order arbitrarily
close to (even though not equal to) $m^{-0.5}$ when $C\left(  \sqrt{2\gamma\ln
m}\right)  ^{-1}m^{\gamma-0.5}\pi_{1,m}^{-1}=o\left(  1\right)  $.

For two Type I location-shift families $\mathcal{F}_{1}$ and $\mathcal{F}_{2}%
$, their corresponding uniform consistency classes $\mathcal{Q}\left(
\mathcal{F}_{1}\right)  $ and $\mathcal{Q}\left(  \mathcal{F}_{2}\right)  $ of
the form (\ref{LocationBoundedNullUCS}), made to have the same $\vartheta$,
$\left\{  \tau_{m}\right\}  _{m\geq1}$, $\left\{  \gamma_{m}\right\}
_{m\geq1}$, $\left\{  \rho_{m}\right\}  _{m\geq1}$ and $\Omega_{m}\left(
\rho_{m}\right)  $, satisfy $\mathcal{Q}\left(  \mathcal{F}_{1}\right)
\subseteq\mathcal{Q}\left(  \mathcal{F}_{2}\right)  $ when $\mathcal{F}%
_{1}\succeq\mathcal{F}_{2}$, where the ordering $\succeq$ means that
$r_{1,\mu}\left(  t\right)  \geq r_{2,\mu}\left(  t\right)  $ for all $\mu\in
U$ and $t\in\mathbb{R}$ and $r_{i,\mu}$ is the modulus of the CF\ $\hat
{F}_{i,\mu}$ of an $F_{i,\mu}\in\mathcal{F}_{i}$ for $i\in\left\{
1,2\right\}  $. This was also shown by the discussion below Theorem 3 in
\cite{Chen:2018a} for the setting of a point null hypothesis. Roughly
speaking, the larger the moduli of the CF's are, the larger the sparsest
alternative proportion an estimator $\hat{\varphi}_{m}\left(  t,\mathbf{z}%
\right)  $ is able to consistently estimate, and the more likely it has a
slower maximal speed of convergence to achieve consistency.

\subsection{The case of a one-sided null}

\label{secConstructionV}

When $\mathcal{F}$ is a location-shift family, it suffices to set $b=0$ for
the one-sided null $\Theta_{0}=\left(  -\infty,b\right)  $. Construction II
for the one-sided null for a Type I location-shift family is provided by:

\begin{theorem}
\label{ThmTypeIOneSided}Suppose $\Theta_{0}=\left(  -\infty,0\right)  $.
Assume $\mathcal{F}$ is a Type I location-shift family such that $F_{0}$ is
differentiable, $\int\left\vert x\right\vert dF_{\mu}\left(  x\right)
<\infty$ for each $\mu\in U$ and
\begin{equation}
\int_{0}^{t}\frac{1}{y}dy\int_{-1}^{1}\left\vert \frac{d}{ds}\frac{1}
{r_{0}\left(  ys\right)  }\right\vert ds<\infty\text{ \ for each }t>0.
\label{CondVc1}%
\end{equation}
Set
\[
K_{1}^{\dagger}\left(  t,x\right)  =\frac{1}{2\pi}\int_{0}^{1}dy\int_{-1}
^{1}\frac{1}{\iota y}\frac{d}{ds}\frac{\exp\left(  \iota tysx\right)  }
{r_{0}\left(  tys\right)  }ds.
\]
Then%
\begin{equation}
\psi_{1}\left(  t,\mu\right)  =\int K_{1}^{\dagger}\left(  t,x\right)
dF_{\mu}\left(  x\right)  =\frac{1}{\pi}\int_{0}^{t}\frac{\sin\left(  \mu
y\right)  }{y}dy, \label{GaussianOnesidepsi}%
\end{equation}
and the desired $\left(  \psi,K\right)  $ is
\begin{equation}
\left\{
\begin{array}
[c]{l}%
K\left(  t,x\right)  =2^{-1}-\Re\left\{  K_{1}^{\dagger}\left(  t,x\right)
\right\}  -2^{-1}K_{1,0}\left(  t,x;0\right) \\
\psi\left(  t,\mu\right)  =2^{-1}-\psi_{1}\left(  t,\mu\right)  -2^{-1}
\psi_{1,0}\left(  t,\mu;0\right)
\end{array}
\right.  . \label{V-a}%
\end{equation}
Set $K_{1}\left(  t,x\right)  =\Re\left\{  K_{1}^{\dagger}\left(  t,x\right)
\right\}  $. Then
\[
K_{1}\left(  t,x\right)  =\frac{1}{2\pi}\int_{0}^{1}dy\int_{-1}^{1}\left[
\frac{\sin\left(  ytsx\right)  }{y}\left\{  \frac{d}{ds}\frac{1}{r_{0}\left(
tys\right)  }\right\}  +\frac{tx\cos\left(  tysx\right)  }{r_{0}\left(
tys\right)  }\right]  ds.
\]

\end{theorem}

Note that $\lim_{t\rightarrow\infty}\psi\left(  t,\mu\right)  =1_{\left(
-\infty,0\right)  }\left(  \mu\right)  $ for $\psi\left(  t,\mu\right)  $ in
(\ref{V-a}). \autoref{ThmTypeIOneSided} cannot be applied to the Cauchy family
since none of its members has finite first-order absolute moment. However, it
is applicable to the Gaussian family (as given by the example below) and three
other families given in \autoref{secExamples} (and others that are not
presented in this work).

\begin{example}
Gaussian family $\mathcal{N}\left(  \mu,\sigma^{2}\right)  $ with mean $\mu$
and standard deviation $\sigma>0$, for which
\[
\frac{dF_{\mu}}{d\nu}\left(  x\right)  =f_{\mu}\left(  x\right)  =\left(
\sqrt{2\pi}\sigma\right)  ^{-1}\exp\left\{  -2^{-1}\sigma^{-2}\left(
x-\mu\right)  ^{2}\right\}  .
\]
The CF of $f_{\mu}$ is $\hat{f}_{\mu}\left(  t\right)  =\exp\left(  \iota
t\mu\right)  \exp\left(  -2^{-1}t^{2}\sigma^{2}\right)  $. So, $r_{\mu}
^{-1}\left(  t\right)  =\exp\left(  2^{-1}t^{2}\sigma^{2}\right)  $ and
$\hat{f}_{0}=r_{0}$. Further,
\[
\frac{1}{y}\frac{d}{ds}\frac{1}{r_{0}\left(  ys\right)  }=\frac{1}{y}%
\sigma^{2}sy^{2}\exp\left(  2^{-1}y^{2}s^{2}\sigma^{2}\right)  =\sigma
^{2}sy\exp\left(  2^{-1}y^{2}s^{2}\sigma^{2}\right)
\]
and condition (\ref{CondVc1}) is satisfied.
\end{example}

Recall $\left\Vert \boldsymbol{\mu}\right\Vert _{\infty}=\max_{1\leq i\leq
m}\left\vert \mu_{i}\right\vert $ and $\Upsilon\left(  \lambda,\tilde{\tau
}\right)  $ in (\ref{UpsilonSimple}). For the estimator $\hat{\varphi}%
_{m}\left(  t,\mathbf{z}\right)  $, the concentration property of the
difference $e_{m}\left(  t\right)  $ (defined by (\ref{eq2d})) depends
critically on $\partial_{t}\left\{  1/r_{0}\left(  t\right)  \right\}  $, as
can be inferred from \autoref{ThmTypeIOneSided}. Even though properties of
$\partial_{t}\left\{  1/r_{0}\left(  t\right)  \right\}  $ can be quite
different for different distribution families (as illustrated by the examples
in \autoref{secExamples}), we still provide a general treatment on the
oscillation of $e_{m}\left(  t\right)  $ and on the consistency of
$\hat{\varphi}_{m}\left(  t,\mathbf{z}\right)  $ at the cost of cumbersome
notations. In contrast, a uniform consistency class of the estimator
$\hat{\varphi}_{m}\left(  t,\mathbf{z}\right)  $ for the Gaussian family has a
much simpler formulation.

\begin{theorem}
\label{ThmTypeI-V}Set $\tilde{u}_{m}=\min_{\left\{  j:\mu_{j}\neq0\right\}
}\left\vert \mu_{j}\right\vert $. Under the hypotheses of
\autoref{ThmTypeIOneSided} and assume $\int\left\vert x\right\vert ^{2}%
dF_{\mu}\left(  x\right)  <\infty$ for each $\mu\in U$. Suppose that $\left\{
z_{i}\right\}  _{i=1}^{m}$ are independent. Then%
\begin{equation}
\mathbb{V}\left\{  e_{m}\left(  t\right)  \right\}  \leq\frac{4}{\pi^{2}%
m}\left[  \bar{r}_{0}^{2}\left(  t\right)  +t^{2}\check{r}_{0}^{2}\left(
t\right)  \tilde{D}_{m}\right]  +\frac{\left\Vert \omega\right\Vert _{\infty
}^{2}}{2m}g^{2}\left(  t,0\right)  , \label{eq15ex}%
\end{equation}
where $\tilde{D}_{m}=m^{-1}\sum_{i=1}^{m}\left(  \sigma_{i}^{2}+\mu_{i}%
^{2}\right)  $, $\sigma_{i}^{2}=\mathbb{V}\left[  z_{i}\right]  $,%
\begin{equation}
\bar{r}_{0}\left(  t\right)  =\sup_{\left(  y,s\right)  \in\left[  0,1\right]
\times\left[  -1,1\right]  }{\left\vert y^{-1}\partial_{s}\left\{
1/r_{0}\left(  tys\right)  \right\}  \right\vert }\text{ and }\check{r}%
_{0}\left(  t\right)  =\sup_{\left(  y,s\right)  \in\left[  0,1\right]
\times\left[  -1,1\right]  }\frac{1}{r_{0}\left(  tys\right)  }.
\label{CFInducedOneSide}%
\end{equation}
Further, there are positive sequences $\left\{  \gamma_{m}\right\}  _{m\geq1}%
$, $\left\{  \tau_{m}\right\}  _{m\geq1}$ and $\left\{  \tilde{\rho}%
_{m}\right\}  _{m\geq1}$ and a positive constant $\vartheta_{1}$, such that
for all $m$ large enough, with probability at least $1-o\left(  1\right)  ,$%
\begin{align}
\sup_{\boldsymbol{\mu}\in\mathcal{\tilde{B}}_{m}\left(  \tilde{\rho}%
_{m}\right)  }\sup_{t\in\left[  0,\tau_{m}\right]  }\left\vert e_{m}\left(
t\right)  \right\vert  &  \leq2^{-1}\left\Vert \omega\right\Vert _{\infty
}\Upsilon\left(  \gamma_{m}m^{-1/2},\tau_{m}\right) \nonumber\\
&  +\Upsilon_{1}\left(  \gamma_{m}m^{-1/2},\tau_{m}\right)  +\tau_{m}%
\Upsilon_{2}\left(  \gamma_{m}m^{-\vartheta_{1}},\tau_{m}\right)  ,
\label{BoundSupWholeErrorX1}%
\end{align}
where $\mathcal{\tilde{B}}_{m}\left(  \tilde{\rho}_{m}\right)  =\mathcal{B}%
_{1,m}\left(  \sqrt{\tilde{\rho}_{m}}\right)  \cap\mathcal{B}_{2,m}\left(
\tilde{\rho}_{m}\right)  $ and%
\begin{equation}
\left\{
\begin{array}
[c]{l}%
\Upsilon_{1}\left(  \lambda,\tilde{\tau}\right)  =\dfrac{\lambda}{\pi}%
\sup_{t\in\left[  0,\tilde{\tau}\right]  }%
{\displaystyle\int_{0}^{1}}
dy%
{\displaystyle\int_{0}^{1}}
\left\vert y^{-1}\partial_{s}\left\{  1/r_{0}\left(  tys\right)  \right\}
\right\vert ds\smallskip\\
\Upsilon_{2}\left(  \lambda,\tilde{\tau}\right)  =\dfrac{\lambda}{\pi}%
\sup_{t\in\left[  0,\tilde{\tau}\right]  }%
{\displaystyle\int_{0}^{1}}
dy%
{\displaystyle\int_{0}^{1}}
\dfrac{ds}{r_{0}\left(  tys\right)  }%
\end{array}
\right.  \text{ for }\lambda,\tilde{\tau}\geq0. \label{Upsilon1and2Oneside}%
\end{equation}
Further, a uniform consistency class for $\hat{\varphi}_{m}\left(
t,\mathbf{z}\right)  $ is%
\begin{equation}
\mathcal{Q}\left(  \mathcal{F}\right)  =\left\{
\begin{array}
[c]{l}%
m^{-2\vartheta}\gamma_{m}^{2}\ln^{-2}\tau_{m}=o\left(  1\right)  ;\tau_{m}%
^{2}m^{\vartheta}\exp\left(  -2^{-1}\gamma_{m}^{2}\right)  =o\left(  1\right)
\\
\left(  \gamma_{m}R_{1}\left(  \sqrt{\tilde{\rho}_{m}}\right)  \right)
^{-1}m^{\vartheta-1/2}\ln\tau_{m}\rightarrow\infty\\
\left(  \gamma_{m}R_{2}\left(  \tilde{\rho}_{m}\right)  \right)
^{-1}m^{\vartheta-\vartheta_{1}}\ln\tau_{m}\rightarrow\infty;m^{-1-2\vartheta
+2\vartheta_{1}}\gamma_{m}^{2}\ln^{-2}\tau_{m}=o\left(  1\right) \\
\tau_{m}m^{\vartheta}\exp\left(  -2^{-1}\kappa^{-2}m^{1-2\vartheta_{1}}\left(
\gamma_{m}^{2}-2\ln\tau_{m}\right)  \right)  =o\left(  1\right) \\
\left[  1-2\Pr\left\{  X_{\left(  0\right)  }>\kappa-\left\Vert
\boldsymbol{\mu}\right\Vert _{\infty}\right\}  \right]  ^{m}=1+o\left(
1\right) \\
\tau_{m}^{-1}\sup_{\boldsymbol{\mu}\in\Omega_{m}\left(  \tilde{\rho}%
_{m}\right)  }\left(  1+\tilde{u}_{m}^{-1}\right)  \pi_{1,m}^{-1}=o\left(
\pi_{1,m}\right) \\
\left[  \Upsilon\left(  \gamma_{m}m^{-1/2},\tau_{m}\right)  +\Upsilon
_{1}\left(  \gamma_{m}m^{-1/2},\tau_{m}\right)  \right]  \sup_{\boldsymbol{\mu
}\in\Omega_{m}\left(  \tilde{\rho}_{m}\right)  }\pi_{1,m}^{-1}=o\left(
1\right) \\
\tau_{m}\Upsilon_{2}\left(  \gamma_{m}m^{-\vartheta_{1}},\tau_{m}\right)
\sup_{\boldsymbol{\mu}\in\Omega_{m}\left(  \tilde{\rho}_{m}\right)  }\pi
_{1,m}^{-1}=o\left(  1\right)
\end{array}
\right\}  \label{OneSideUCC}%
\end{equation}
for some positive constant $\vartheta$ and positive sequence $\left\{
\kappa=\kappa_{m}\right\}  _{m\geq1}$, where%
\begin{equation}
R_{1}\left(  \sqrt{\tilde{\rho}_{m}}\right)  =2\mathbb{E}\left[  \left\vert
X_{(0)}\right\vert \right]  +2\sqrt{\tilde{\rho}_{m}}\text{ and }R_{2}\left(
\tilde{\rho}_{m}\right)  =4\mathbb{E}\left[  X_{\left(  0\right)  }%
^{2}\right]  +4\tilde{\rho}_{m} \label{L1L2Rhos}%
\end{equation}
amd $\Omega_{m}\left(  \rho_{m}\right)  \subseteq$ $\mathcal{\tilde{B}}%
_{m}\left(  \tilde{\rho}_{m}\right)  $. In particular, for the Gaussian family
with variance $\sigma^{2}>0$, a uniform consistency class for $\hat{\varphi
}_{m}\left(  t,\mathbf{z}\right)  $ is
\begin{equation}
\mathcal{Q}\left(  \mathcal{F}\right)  =\left\{
\begin{array}
[c]{l}%
0<\vartheta_{1}<2^{-1};q>\vartheta;0\leq\vartheta^{\prime\prime}%
<\vartheta-2^{-1};R_{2}\left(  \tilde{\rho}_{m}\right)  =O\left(
m^{\vartheta^{\prime\prime}}\right) \\
0<\gamma<\vartheta_{1}\sigma^{-2};\tau_{m}=\sqrt{2\gamma\ln m};\vartheta
_{2}\geq1;\kappa-\left\Vert \boldsymbol{\mu}\right\Vert _{\infty}%
=\sqrt{2\sigma^{2}\vartheta_{2}\ln m}\\
\tau_{m}^{-1}\sup_{\boldsymbol{\mu}\in\Omega_{m}\left(  \tilde{\rho}%
_{m}\right)  }\pi_{1,m}^{-1}\left(  1+\tilde{u}_{m}^{-1}\right)  =o\left(
1\right) \\
m^{\sigma^{2}\gamma-\vartheta_{1}}\ln m\sup_{\boldsymbol{\mu}\in\Omega
_{m}\left(  \tilde{\rho}_{m}\right)  }\pi_{1,m}^{-1}=o\left(  1\right)
\end{array}
\right\}  . \label{eq12d}%
\end{equation}
where $q,\gamma,\vartheta_{2}$ and $\vartheta^{\prime\prime}$ are also constants.
\end{theorem}

We remarked that $\tilde{u}_{m}$ in \autoref{ThmTypeI-V} should be set as
$\tilde{u}_{m}=\min_{\left\{  j:\mu_{j}\neq b\right\}  }\left\vert \mu
_{j}-b\right\vert $ when $b\neq0$ in the one-sided null $\Theta_{0}=\left(
-\infty,b\right)  $. Similar to the $u_{m}$ in (\ref{eq10b}), a suitable
magnitude for $\tilde{u}_{m}$ is needed for the estimator induced by
$K_{1,0}\left(  t,x;0\right)  $ in (\ref{V-a}) to consistently estimate the
proportion of $\mu_{j}$'s that are equal to $0$; see Theorems 2 and 3 of
\cite{Chen:2018a}. Similar to the bounds given by \autoref{ThmIVa},
\autoref{ThmTypeI-V} states that the variance of the error $e_{m}\left(
t\right)  $ of the estimator $\hat{\varphi}_{m}\left(  t,\mathbf{z}\right)  $
is characterized by the average reciprocal modulus $g$ at the generating
measure $F_{0}$ of the location-shift family $\mathcal{F}$, the supremum norm
$\Vert\omega\Vert_{\infty}$ of the employed density $\omega$, the
\textquotedblleft average squared effective size\textquotedblright%
\ $m^{-1}\sum_{i=1}^{m}\mu_{i}^{2}$ of the one-sided null, the average
variance $m^{-1}\sum_{i=1}^{m}\sigma_{i}^{2}$ of the random variables, the
number $m$ of hypotheses to test, the magnitude of the tuning parameter $t$,
and the two quantities in (\ref{CFInducedOneSide}) that are induced by the
reciprocal of the modulus $r_{0}\left(  \cdot\right)  $ of the CF of $F_{0}$.
Further, the lager $g$ is, the more likely the estimator $\hat{\varphi}%
_{m}\left(  t,\mathbf{z}\right)  $ will have a slower maximal speed of
convergence to achieve consistency. In contrast, the uniform concentration of
$e_{m}\left(  t\right)  $ is mainly characterized by certain uniformity of
three \textquotedblleft averages\textquotedblright\ induced by the reciprocal
of the modulus of the CF of $F_{0}$ that are defined by (\ref{UpsilonSimple})
and (\ref{Upsilon1and2Oneside}). Roughly speaking, the large these averages
are, the less $\left\vert e_{m}\left(  t\right)  \right\vert $ is concentrated.

The uniform consistency class $\mathcal{Q}\left(  \mathcal{F}\right)  $ in
(\ref{OneSideUCC}) looks quite complicated since to obtain it we have to deal
with three empirical processes simultaneously, one of which is induced by
unbounded and non-Lipschitz functions (as can been from the function
$K_{1}\left(  t,x\right)  $ in \autoref{ThmTypeIOneSided} with a $z_{i}$ in
place of $x$ and the proof of \autoref{ThmTypeI-V}). Firstly, the first five
rows of conditions in $\mathcal{Q}\left(  \mathcal{F}\right)  $ ensure that
the concentration inequality (\ref{BoundSupWholeErrorX1}) for $e_{m}\left(
t\right)  $ holds asymptotically with probability $1$, the sixth row ensures
$\left\vert \pi_{1,m}^{-1}\hat{\varphi}_{m}\left(  t,\mathbf{z}\right)
-1\right\vert \rightarrow0$, and the last two ensure that the upper bound in
(\ref{BoundSupWholeErrorX1}) converges to $0$, all as $m\rightarrow\infty$. In
particular, the condition $\left[  1-2\Pr\left\{  X_{\left(  0\right)
}>\kappa-\left\Vert \boldsymbol{\mu}\right\Vert _{\infty}\right\}  \right]
^{m}=1+o\left(  1\right)  $ is the consequence of emplying the truncation
technique to deal with an unbounded empirical processes. Second, the quantity
$R_{2}\left(  \tilde{\rho}_{m}\right)  $ in $\mathcal{Q}\left(  \mathcal{F}%
\right)  $ already involves a control on $m^{-1}\sum_{i=1}^{m}\mu_{i}^{2}$,
i.e., on the size of $\mathcal{B}_{2,m}\left(  \tilde{\rho}_{m}\right)  $.
This is more involved than that in \autoref{ThmIVa} where only $R_{1}\left(
\rho_{m}\right)  $ is used to control the size of $\mathcal{B}_{1,m}\left(
\rho_{m}\right)  $. Third, similar to the ordering via subsetting for uniform
consistency classes for the estimator $\hat{\varphi}_{m}\left(  t,\mathbf{z}%
\right)  $ for a bounded null for different Type I\ location-shift families
via the ordering of the magnitude of corresponding moduli of CF's when other
things are set the same, which has been discussed at the end of
\autoref{SecConstructionIV}, such an ordering via subsetting can be
established for uniform consistency classes for the estimator $\hat{\varphi
}_{m}\left(  t,\mathbf{z}\right)  $ for a one-sided null for different Type
I\ location-shift families via the ordering of the magnitudes of the moduli of
the CF's and of the partial derivatives of the reciprocals of these moduli
when other things are set equal, as revealed by the upper bound in
(\ref{BoundSupWholeErrorX1}) and the representation of the uniform consistency
class $\mathcal{Q}\left(  \mathcal{F}\right)  $ in (\ref{OneSideUCC}).
However, we will not provide details on this here.

Despite the complexity of $\mathcal{Q}\left(  \mathcal{F}\right)  $ for a
general Type I location-shift family, a uniform consistency class for the
Gaussian family with variance $\sigma^{2}>0$ has a much simpler representation
as in (\ref{eq12d}), due to its simple form of CF's. Specifically, an example
of (\ref{eq12d}) is%
\[
\mathcal{U}_{m}=\left\{  \boldsymbol{\mu}\in\mathcal{\tilde{B}}_{m}\left(
\rho_{m}\right)  :%
\begin{array}
[c]{l}%
0\leq\vartheta^{\prime\prime}<1/8;R_{2}\left(  \tilde{\rho}_{m}\right)
\leq\tilde{C}m^{\vartheta^{\prime\prime}}\\
\tau_{m}=\sqrt{\left(  1/4\right)  \sigma^{-2}\ln m};\kappa-\left\Vert
\boldsymbol{\mu}\right\Vert _{\infty}=\sqrt{2\sigma^{2}\ln m}\\
\pi_{1,m}^{-1}\leq\tilde{C}\sqrt{\ln m};\tilde{u}_{m}\geq\left(  \tau
_{m}\right)  ^{-1}\ln\ln m
\end{array}
\right\}
\]
for a constant $\tilde{C}>0$, for which $\sup\nolimits_{\boldsymbol{\mu}%
\in\mathcal{U}_{m}}\left(  \left\vert \pi_{1,m}^{-1}\sup\nolimits_{t\in
\left\{  \tau_{m}\right\}  }\hat{\varphi}_{m}\left(  t,\mathbf{z}\right)
-1\right\vert \right)  \rightsquigarrow0$ as\ $m\rightarrow\infty$. Further,
for such a family, the variance upper bound in (\ref{eq15ex}) becomes
\[
\mathbb{V}\left\{  \hat{\varphi}_{m}\left(  t,\mathbf{z}\right)  \right\}
\leq\frac{4t^{2}\exp\left(  t^{2}\sigma^{2}\right)  }{\pi^{2}m}\left(
t^{2}\sigma^{4}+\sigma^{2}+m^{-1}\sum_{i=1}^{m}\mu_{i}^{2}\right)
+\frac{\left\Vert \omega\right\Vert _{\infty}^{2}}{m}\frac{8\exp\left(
t^{2}\sigma^{2}\right)  }{t^{4}\sigma^{4}},
\]
and with $\tau_{m}=\sqrt{2\gamma\ln m}$ and $\gamma_{m}=\sqrt{2q\ln m}$, the
concentration inequality (\ref{BoundSupWholeErrorX1}) becomes%
\[
\sup_{\boldsymbol{\mu}\in\mathcal{\tilde{B}}_{m}\left(  \tilde{\rho}%
_{m}\right)  }\sup_{t\in\left[  0,\sqrt{2\gamma\ln m}\right]  }\left\vert
e_{m}\left(  t\right)  \right\vert \leq\bar{\Upsilon}\left(  q,\gamma
,\vartheta_{1}\right)
\]
with probability at least $1-o\left(  1\right)  $, where%
\[
\bar{\Upsilon}\left(  q,\gamma,\vartheta_{1}\right)  =\frac{\left\Vert
\omega\right\Vert _{\infty}\sqrt{2q\ln m}}{\sigma^{2}\gamma\ln m}m^{\sigma
^{2}\gamma-1/2}+\dfrac{\sqrt{2q\ln m}}{\pi}m^{\gamma\sigma^{2}-1/2}%
+\dfrac{2\sqrt{\gamma q}\ln m}{\pi}m^{\gamma\sigma^{2}-\vartheta_{1}},
\]
as will be shown by the proof of \autoref{ThmTypeI-V}.

\subsection{Additional examples of Type I location-shift families}

\label{secExamples}

We provide four additional Type I location-shift families, all of which were
discussed in Section 3.1 of \cite{Chen:2018a}. \autoref{ThmTypeIOneSided}
applies to each of the Laplace, Logistic and Hyperbolic Secant families (given
below), and for each such family, a uniform consistency class for the
estimator $\hat{\varphi}_{m}\left(  t,\mathbf{z}\right)  $ can be obtained for
the case of a bounded null using the same techniques in
\autoref{SecConstructionIV} and for the case of a one-sided null using the
same techniques in \autoref{secConstructionV}, respectively. However, we will
not pursue this here.

\begin{example}
\label{egLaplace}Laplace family $\mathsf{Laplace}\left(  \mu,2\sigma
^{2}\right)  $ with mean $\mu$ and standard deviation $\sqrt{2}\sigma>0$ for
which%
\[
\frac{dF_{\mu}}{d\nu}\left(  x\right)  =f_{\mu}\left(  x\right)  =\frac
{1}{2\sigma}\exp\left(  -\sigma^{-1}\left\vert x-\mu\right\vert \right)
\]
and the CF of $f_{\mu}$ is $\textcolor{blue}{\hat{F}}_{\mu}\left(  t\right)
=\left(  1+\sigma^{2}t^{2}\right)  ^{-1}\exp\left(  \iota t\mu\right)  $. So,
$r_{\mu}^{-1}\left(  t\right)  =1+\sigma^{2}t^{2}$ and
$\textcolor{blue}{\hat{F}}_{0}=r_{0}$. Further,
\[
\frac{1}{y}\frac{d}{ds}\frac{1}{r_{0}\left(  ys\right)  }=\frac{1}{y}
2\sigma^{2}sy^{2}=2\sigma^{2}sy,
\]
condition (\ref{CondVc1}) is satisfied, and \autoref{ThmTypeIOneSided} applies
to this example.
\end{example}

\begin{example}
Logistic family $\mathsf{Logistic}\left(  \mu,\sigma\right)  $ with mean $\mu$
and scale parameter $\sigma>0$, for which
\[
\frac{dF_{\mu}}{d\nu}\left(  x\right)  =f_{\mu}\left(  x\right)  =\frac
{1}{4\sigma}{\sech}^{2}\left(  \frac{x-\mu}{2\sigma}\right)
\]
and the CF of $f_{\mu}$ is $\textcolor{blue}{\hat{F}}_{\mu}\left(  t\right)
=\frac{\pi\sigma t}{\sinh\left(  \pi\sigma t\right)  }\exp\left(  \iota
t\mu\right)  $. So, $r_{\mu}^{-1}\left(  t\right)  =\frac{\sinh\left(
\pi\sigma t\right)  } {\pi\sigma t}$ and $\textcolor{blue}{\hat{F}}_{0}=r_{0}%
$. Further,
\[
\left\vert \frac{1}{y}\frac{d}{ds}\frac{1}{r_{0}\left(  ys\right)
}\right\vert =o\left(  \left\vert ys\right\vert \right)  \text{ \ as
\ }y\rightarrow0+\text{ for each fixed }s,
\]
condition (\ref{CondVc1}) is satisfied, and \autoref{ThmTypeIOneSided} applies
to this example.
\end{example}

\begin{example}
Cauchy family $\mathsf{Cauchy}\left(  \mu,\sigma\right)  $ with median $\mu$
and scale parameter $\sigma>0$, for which
\[
\frac{dF_{\mu}}{d\nu}\left(  x\right)  =f_{\mu}\left(  x\right)  =\frac{1}
{\pi\sigma}\frac{\sigma^{2}}{\left(  x-\mu\right)  ^{2}+\sigma^{2}}%
\]
and the CF of $f_{\mu}$ is $\hat{F}_{\mu}\left(  t\right)  =\exp\left(
-\sigma\left\vert t\right\vert \right)  \exp\left(  \iota t\mu\right)  $.
Since $\int\left\vert x\right\vert dF_{\mu}\left(  x\right)  =\infty$,
\autoref{ThmTypeIOneSided} cannot be applied to the case of one-sided null.
\end{example}

\begin{example}
Hyperbolic Secant family $\mathsf{HSecant}\left(  \mu,\sigma\right)  $ with
mean $\mu$ and scale parameter $\sigma>0$, for which%
\[
\frac{dF_{\mu}}{d\nu}\left(  x\right)  =f_{\mu}\left(  x\right)  =\frac
{1}{2\sigma}\frac{1}{\cosh\left(  \pi\frac{x-\mu}{\sigma}\right)  };
\]
see, e.g., Chapter 1 of \cite{Fischer:2014}. The identity%
\[
\int_{-\infty}^{+\infty}e^{\iota tx}\frac{dx}{\pi\cosh\left(  x\right)
}=\cosh\left(  2^{-1}\pi t\right)  ,
\]
implies $\hat{F}_{\mu}\left(  t\right)  =\sigma^{-1}\exp\left(  -\iota
t\mu\sigma^{-1}\right)  {\sech}\left(  t\sigma^{-1}\right)  $. So, $r_{\mu
}^{-1}\left(  t\right)  =\sigma\cosh\left(  t\sigma^{-1}\right)  $ and
$\hat{F}_{0}=r_{0}$. Further,%
\[
\left\vert \frac{1}{y}\frac{d}{ds}\frac{1}{r_{0}\left(  ys\right)
}\right\vert =O\left(  \left\vert s+o\left(  ys\right)  \right\vert \right)
\text{ when }y\rightarrow0\text{ for each fixed }s,
\]
condition (\ref{CondVc1}) is satisfied, and \autoref{ThmTypeIOneSided} applies
to this example.
\end{example}

\section{Extension of constructions for a bounded null}

\label{SecExtensions}

The settings for extending the construction for a bounded null $\Theta
_{0}=\left(  a,b\right)  $ are identical to those provided by
\cite{Xiongzhi2025}, as will be quoted below. Specifically, we consider estimating
\begin{equation}
\check{\pi}_{0,m}=m^{-1}\sum\nolimits_{\left\{  i\in\left\{  1,\ldots
,m\right\}  :\mu_{i}\in\Theta_{0}\right\}  }\phi\left(  \mu_{i}\right)
\label{eq4c}%
\end{equation}
for a suitable function $\phi$, where $\check{\pi}_{0,m}\in\left[  0,1\right]
$ does not necessarily hold. Two examples of $\phi$ are $\phi(x)=|x|^{p}$ for
some $p>0$ and $\phi(x;c_{\ast})=\min\{|x|^{p},c_{\ast}^{p}\}$ for some
$p,c_{\ast}>0$, which were discussed by \cite{Jin:2008} and \cite{Xiongzhi2025}.
To formulate the extension, we quote two definitions, a convention and a lemma
from \cite{Xiongzhi2025} as follows. First, for a $\phi\in L^{1}\left(  \left[
a,b\right]  \right)  $ (with finite $a$ and $b$ such that $a<b$), define
\[
\mathcal{D}_{\phi}\left(  t,\mu;a,b\right)  =\frac{1}{\pi}\int_{a}^{b}%
\frac{\sin\left\{  \left(  \mu-y\right)  t\right\}  }{\mu-y}\phi\left(
y\right)  dy\text{ \ for \ }t,\mu\in\mathbb{R}.
\]
Second, we quote:

\begin{definition}
\label{defDiffratio}
\textcolor{blue}{Let $\phi$ be a function defined on $U.$} For each $\mu\in U
$ and $z\in\left[  \mu-b,\mu-a\right]  $, define
\begin{equation}
W_{\mu}\left(  z\right)  =z^{-1}\left(  \phi\left(  \mu-z\right)  -\phi\left(
\mu\right)  \right)  \text{ for }\ z\neq0 \label{eqWmu}%
\end{equation}
and%
\begin{equation}
W_{\mu}^{-}\left(  z\right)  =\left\{
\begin{array}
[c]{lll}%
W_{\mu}\left(  z\right)  & \text{for} & \mu-a\geq z>0\\
\lim_{z\rightarrow0+}W_{\mu}\left(  z\right)  & \text{for} & z=0\text{ if the
limit exists}%
\end{array}
\right.  \label{eqWmuA}%
\end{equation}
and%
\begin{equation}
W_{\mu}^{+}\left(  z\right)  =\left\{
\begin{array}
[c]{lll}%
-W_{\mu}\left(  z\right)  & \text{for} & \mu-b\leq z<0\\
-\lim_{z\rightarrow0-}W_{\mu}\left(  z\right)  & \text{for} & z=0\text{ if the
limit exists}%
\end{array}
\right.  . \label{eqWmuB}%
\end{equation}

\end{definition}

\textcolor{blue}{
Third, for the rest of this section, we assume that $U$ is an interval
containing $\left[  a,b\right]  $ as a subset, and take the convention that
the total variation of a function defined on $U$ is the supremum of the total
variation of this function on each compact interval contained in $U$ (to avoid
unneeded complications of defining total variation of a function on a
measurable set).
}Fourth, we quote:

\begin{lemma}
\label{ThmDirichletInterval}If $\phi\in L^{1}\left(  \left[  a,b\right]
\right)  $, then setting $\hat{\phi}\left(  s\right)  =\int_{a}^{b}\phi\left(
y\right)  \exp\left(  -\iota ys\right)  dy$ gives
\begin{equation}
\mathcal{D}_{\phi}\left(  t,\mu;a,b\right)  =\frac{t}{2\pi}\int_{-1}^{1}
\hat{\phi}\left(  ts\right)  \exp\left(  \iota\mu ts\right)  ds. \label{eq4d}%
\end{equation}
On the other hand, if $\phi$ is continuous and of bounded variation on
$\left[  a,b\right]  $, then
\begin{equation}
\mathcal{D}_{\phi,\infty}\left(  \mu;a,b\right)  := \lim_{t\rightarrow\infty
}\mathcal{D}_{\phi}\left(  t,\mu;a,b\right)  =\left\{
\begin{array}
[c]{lll}%
\phi\left(  \mu\right)  & \text{if} & a<\mu<b\\
2^{-1}\phi\left(  \mu\right)  & \text{if} & \mu=a\text{ or }\mu=b\\
0 & \text{if} & \mu<a\text{ or }\mu>b
\end{array}
. \right.  \label{eqA1}%
\end{equation}
If in addition both $W_{\mu}^{+}$ and $W_{\mu}^{-}$ are well-defined and of
bounded variation for each fixed $\mu\in U$, then
\begin{equation}
\left\vert \mathcal{D}_{\phi}\left(  t,\mu;a,b\right)  -\mathcal{D}%
_{\phi,\infty}\left(  \mu;a,b\right)  \right\vert \leq\frac{4C_{\mu}\left(
\phi\right)  }{\pi t}+\frac{4\left\Vert \phi\right\Vert _{\infty}}{t}\left(
\frac{1}{b-a}+\frac{3}{\delta_{\mu,a,b}}\right)  , \label{eqA4}%
\end{equation}
when $\min\left\{  t\delta_{\mu,a,b},t\left(  b-a\right)  \right\}  \geq2$,
where $\delta_{\mu,a,b}=\min_{\mu\notin\left\{  a,b\right\}  }\left\{
\left\vert \mu-a\right\vert ,\left\vert \mu-b\right\vert \right\}  $ and
\begin{equation}
C_{\mu}\left(  \phi\right)  =\left\Vert W_{\mu}^{-}\right\Vert _{\infty
}+\left\Vert W_{\mu}^{-}\right\Vert _{\mathrm{TV}}+\left\Vert W_{\mu}%
^{+}\right\Vert _{\infty}+\left\Vert W_{\mu}^{+}\right\Vert _{\mathrm{TV}}.
\label{defphiBnd}%
\end{equation}

\end{lemma}

Inequality (\ref{eqA4}) gives the speed of convergence of $\mathcal{D}_{\phi
}\left(  t,\mu;a,b\right)  $ and helps determine the speed of convergence of
the estimators to be constructed below:

\begin{theorem}
\label{ThmExtA}Let $\phi$ be continuous and of bounded variation on $\left[
a,b\right]  $. Assume $\mathcal{F}$ is a Type I location-shift family and set
\begin{equation}
K_{1}\left(  t,x\right)  =\frac{t}{2\pi}\int_{a}^{b}\phi\left(  y\right)
dy\int_{\left[  -1,1\right]  }\frac{\cos\left\{  ts\left(  x-y\right)
\right\}  }{r_{0}\left(  ts\right)  }ds. \label{eq13a}%
\end{equation}
Then
\begin{equation}
\psi_{1}\left(  t,\mu\right)  =\int K_{1}\left(  t,x\right)  dF_{\mu}\left(
x\right)  =\mathcal{D}_{\phi}\left(  t,\mu;a,b\right)  \label{eq13aaa}%
\end{equation}
and the desired $\left(  \psi,K\right)  $ for estimating $\check{\pi}_{0,m}$
is
\begin{equation}
\left\{
\begin{array}
[c]{l}%
K\left(  t,x\right)  =K_{1}\left(  t,x\right)  -2^{-1}\left\{  \phi\left(
a\right)  K_{1,0}\left(  t,x;a\right)  +\phi\left(  b\right)  K_{1,0}\left(
t,x;b\right)  \right\} \\
\psi\left(  t,\mu\right)  =\psi_{1}\left(  t,\mu\right)  -2^{-1}\left\{
\phi\left(  a\right)  \psi_{1,0}\left(  t,\mu;a\right)  +\phi\left(  b\right)
\psi_{1,0}\left(  t,\mu;b\right)  \right\}
\end{array}
\right.  . \label{eq2b}%
\end{equation}

\end{theorem}

Note that $\lim_{t\rightarrow\infty}\psi\left(  t,\mu\right)  =1_{\left(
a,b\right)  }\left(  \mu\right)  \phi\left(  \mu\right)  $ for $\psi\left(
t,\mu\right)  $ in (\ref{eq2b}). The constructions in \autoref{ThmExtA} can be
easily modified to estimate any linear function of $\check{\pi}_{0,m}$, which
will not be discussed here. In particular, if we set $a=-b$ with $b>0$ and
take $K_{1}$ in (\ref{eq13a}) and $\psi_{1}$ in (\ref{eq13aaa}), then the
construction $\left(  \psi_{1},K_{1}\right)  $ reduces to those for Gaussian
family in Section 6 of \cite{Jin:2008}. Moreover, when $\phi$ is the constant
function $1$ on $\left[  a,b\right]  $, i.e., $\phi\left(  t\right)  \equiv1$,
(\ref{eq13a}) reduces to (\ref{eq13e}).

Define%
\begin{equation}
\hat{\varphi}_{m}\left(  t,\mathbf{z}\right)  =m^{-1}\sum_{i=1}^{m}K\left(
t,z_{i}\right)  \text{ \ and \ }\varphi_{m}\left(  t,\boldsymbol{\mu}\right)
=m^{-1}\sum_{i=1}^{m}\mathbb{E}\left\{  K\left(  t,z_{i}\right)  \right\}
\label{eq21b}%
\end{equation}
with $K$ in (\ref{eq2b}) and set $e_{m}\left(  t\right)  =\hat{\varphi}%
_{m}\left(  t,\mathbf{z}\right)  -\varphi_{m}\left(  t,\boldsymbol{\mu
}\right)  $. Then $\hat{\varphi}_{m}\left(  t,\mathbf{z}\right)  $ defined by
(\ref{eq21b}) estimates $\check{\pi}_{0,m}$ defined by (\ref{eq4c}).
Consistency of the estimator $\hat{\varphi}_{m}\left(  t,\mathbf{z}\right)  $
given by (\ref{eq21b}) can be obtained for independent $\left\{
z_{i}\right\}  _{i=1}^{m}$ via almost identical arguments as those for the
proof of \autoref{ThmIVa}; see \autoref{ThmFinal} below.

For the rest of this section, we assume that $\phi$ is continuous and of
bounded variation on $\left[  a,b\right]  $. Recall $u_{m}\ $in (\ref{eq10b})
and $C_{\mu}\left(  \phi\right)  $ defined by (\ref{defphiBnd}). We quote from
\cite{Xiongzhi2025} the following condition:

\begin{description}
\item[C3)] $\phi$ has finite left and right derivatives at each interior point
of $U$, a finite right derivative at the left boundary point of $U$, and a
finite left derivative at the right boundary point of $U$, and $\left\Vert
\phi\right\Vert _{1,\infty}:=\sup_{\mu\in U}C_{\mu}\left(  \phi\right)
<\infty$.
\end{description}

With these preparations, we have the uniform consistency of $\hat{\varphi}%
_{m}\left(  t,\mathbf{z}\right)  $ as:

\begin{theorem}
\label{ThmFinal}Assume $\left\{  z_{i}\right\}  _{i=1}^{m}$ are independent
whose CDFs are members of a Type I\ location-shift family $\mathcal{F}$. Then
the following hold:

\begin{enumerate}
\item For the estimator $\hat{\varphi}_{m}\left(  t,\mathbf{z}\right)  $ in
(\ref{eq21b}) of $\check{\pi}_{0,m}=m^{-1}\sum\nolimits_{\left\{  i\in\left\{
1,\ldots,m\right\}  :\mu_{i}\in\left(  a,b\right)  \right\}  }\phi\left(
\mu_{i}\right)  $,
\[
\mathbb{V}\left\{  e_{m}\left(  t\right)  \right\}  \leq2\left\Vert
\phi\right\Vert _{\infty}^{2}m^{-1}g^{2}\left(  t,0\right)  \left(  \pi
^{-2}\left(  b-a\right)  ^{2}t^{2}+\left\Vert \omega\right\Vert _{\infty}%
^{2}\right)  ,
\]
and, there are positive sequences $\left\{  \gamma_{m}\right\}  _{m\geq1}$,
$\left\{  \tau_{m}\right\}  _{m\geq1}$ and $\left\{  \rho_{m}\right\}
_{m\geq1}$, such that, for all sufficiently larget $m$, with probability at
least $1-o\left(  1\right)  $
\[
\sup_{\boldsymbol{\mu}\in\mathcal{B}_{1,m}\left(  \rho_{m}\right)  }\sup
_{t\in\left[  0,\tau_{m}\right]  }\left\vert e_{m}\left(  t\right)
\right\vert \leq\left\{  \frac{\left(  b-a\right)  }{2\pi}\tau_{m}+\left\Vert
\omega\right\Vert _{\infty}\right\}  \left\Vert \phi\right\Vert _{\infty
}\Upsilon\left(  \gamma_{m}m^{-1/2},\tau_{m}\right)  .
\]

\item For the estimator $\hat{\varphi}_{1,m}\left(  t,\mathbf{z}\right)
=m^{-1}\sum_{i=1}^{m}K_{1}\left(  t,z_{i}\right)  $ with $K_{1}$ in
(\ref{eq13a}) that estimates
\begin{equation}
\tilde{\pi}_{0,m}=m^{-1}\sum\nolimits_{\left\{  i\in\left\{  1,\ldots
,m\right\}  :\mu_{i}\in\left[  a,b\right]  \right\}  }\phi\left(  \mu
_{i}\right)  , \label{eq21d}%
\end{equation}
then
\[
\mathbb{V}\left\{  e_{1,m}\left(  t\right)  \right\}  \leq\pi^{-2}m^{-1}%
t^{2}\left(  b-a\right)  ^{2}\left\Vert \phi\right\Vert _{\infty}^{2}%
g^{2}\left(  t,0\right)  ,
\]
and there are positive sequences $\left\{  \gamma_{m}\right\}  _{m\geq1}$,
$\left\{  \tau_{m}\right\}  _{m\geq1}$ and $\left\{  \rho_{m}\right\}
_{m\geq1}$, such that, for all $m$ large enough, with probability at least
$1-o\left(  1\right)  $,%
\[
\sup_{\boldsymbol{\mu}\in\mathcal{B}_{1,m}\left(  \rho_{m}\right)  }\sup
_{t\in\left[  0,\tau_{m}\right]  }\left\vert e_{1,m}\left(  t\right)
\right\vert \leq\frac{\left(  b-a\right)  \tau_{m}\left\Vert \phi\right\Vert
_{\infty}}{2\pi}\Upsilon\left(  \tau_{m}m^{-1/2},\tau_{m}\right)  .
\]

\end{enumerate}

If in addition condition C3) holds, then a uniform consistency class for
either $\hat{\varphi}_{m}\left(  t,\mathbf{z}\right)  $ of $\check{\pi}_{0,m}$
or $\hat{\varphi}_{1,m}\left(  t,\mathbf{z}\right)  $ of $\tilde{\pi}_{0,m}$
is the $\mathcal{Q}\left(  \mathcal{F}\right)  $ in
(\ref{LocationBoundedNullUCS}), with $\check{\pi}_{0,m}$ or $\tilde{\pi}%
_{0,m}$ in place of $\pi_{1,m}$, for which $\mathcal{Q}\left(  \mathcal{F}%
\right)  $ in (\ref{LocationBoundedNullUCS}) takes for the form of
(\ref{MainClass}) with an example given by (\ref{MainClassEg}), when
$\mathcal{F}$ is the Gaussian family with variance $\sigma^{2}>0$, again with
$\check{\pi}_{0,m}$ or $\tilde{\pi}_{0,m}$ in place of $\pi_{1,m}$.
\end{theorem}

The \textquotedblleft uniform consistency class\textquotedblright\ in
\autoref{ThmFinal} bears the meaning of \autoref{DefUniformConsistency} but
with $\check{\pi}_{0,m}$ or $\tilde{\pi}_{0,m}$ in place of $\pi_{1,m}$. The
assertion of \autoref{ThmFinal} on estimating $\tilde{\pi}_{0,m}$ for the
Gaussian family complements and strengthens much Theorem 13 of \cite{Jin:2008}%
, since the latter work in our notations requires $\phi$ to be absolutely
continuous, deals with the case $\left[  a,b\right]  $ being a symmetric
interval, and only shows $\sup_{\mathcal{\tilde{B}}_{1,m}\left(  \rho\right)
}\left\vert \hat{\varphi}_{1,m}\left(  t,\mathbf{z}\right)  -\tilde{\pi}%
_{0,m}\right\vert =o\left(  1\right)  $ for a subset $\mathcal{\tilde{B}%
}_{1,m}\left(  \rho\right)  $ of $\mathcal{B}_{1,m}\left(  \rho_{m}\right)  $
(that is defined by (\ref{eq19f}) and allows $\rho_{m}\rightarrow\infty$).
Since the construction in \autoref{ThmExtA} reduces to that in \autoref{ThmIV}
when $\phi\left(  t\right)  \equiv1$, $t\in U$, there is no surprise that
variance bounds and concentration bounds in \autoref{ThmFinal} differ from
those in \autoref{ThmIVa} by a factor of $\left\Vert \phi\right\Vert _{\infty
}^{2}$ and $\left\Vert \phi\right\Vert _{\infty}$ respectively, and that the
uniform consitency classes in \autoref{ThmFinal} and \autoref{ThmIVa} are the
same, except with $\check{\pi}_{0,m}$ or $\tilde{\pi}_{0,m}$ in place of
$\pi_{1,m}$, since asymptotically a factor of $\left\Vert \phi\right\Vert
_{\infty}^{2}$ or $\left\Vert \phi\right\Vert _{\infty}$ does not affect the
consistency of the corresponding estimator $\hat{\varphi}_{m}\left(
t,\mathbf{z}\right)  $ or $\hat{\varphi}_{1,m}\left(  t,\mathbf{z}\right)  $.


\section{Simulation study}

\label{SecNumericalStudies}

We will present a simulation study on the proposed estimators, with a
comparison to the ``MR'' estimator of \cite{Meinshausen:2006} or Storey's
estimator of \cite{Storey:2004} for the case of a one-sided. For one-sided
null $\Theta_{0}=\left(  -\infty,b\right)  \cap U$, when $X_{0}$ is an
observation from a random variable $X$ with CDF $F_{\mu}$, $\mu\in U$, its
one-sided p-value is computed as $1-F_{b}\left(  X_{0}\right)  $. For one-side
null on the means of Gaussian random variables, we will also compare our
estimator with one minus the null proportion estimator of Section 4.3 of
\cite{Hoang:2022b} that has its first tuning parameter (for thresholding
randomized p-values) as $0.5$ and second tuning parameter $\tilde{c}_{0}$ (for
constructing randomized p-values), which leads to the ``HD'' estimator with
$\tilde{c}_{0} =0.5$ and the ``HD1'' estimator with the value for $\tilde
{c}_{0}$ that is optimally determined from data from a grid of $50$ equally
spaced candidate values for $\tilde{c}_{0}$.

We numerically implement the solution $\left(  \psi,K\right)  $ in two cases
as follows: (a) if $\psi$ or $K$ is defined by a univariate integral, then the
univariate integral is approximated by a Riemann sum based on an equally
spaced partition with norm $0.01$ of the corresponding domain of integration;
(b) if $\psi$ or $K$ is defined by a double integral, then the double integral
is computed as an iterated integral, for which each univariate integral is
computed as if it were case (a). We choose norm $0.01$ for a partition so as
to reduce a bit the computational complexity of the proposed estimators when
the number of hypotheses to test is very large. However, we will not explore
here how much more accurate these estimators can be when finer partitions are
used to obtain the Riemman sums, or explore here which density function
$\omega(s)$ on $[-1,1]$ should be used to give the best performances to the
proposed estimators among all continuous densities on $[-1,1]$ that are of
bounded variation. By default, we will choose the triangular density
$\omega\left(  s\right)  =\left(  1-\vert s\vert\right)  1_{\left[
-1,1\right]  }\left(  s\right)  $, since numerical evidence in
\cite{Jin:2008,Chen:2018a} shows that this density leads to good performances
of the proposed estimators for the setting of a point null.

The MR estimator (designed for continuous p-values) is implemented as follows:
let the ascendingly ordered p-values be $p_{\left(  1\right)  }<p_{\left(
2\right)  }<\cdots<p_{\left(  m\right)  }$ for $m>4$, set $b_{m}^{\ast
}=m^{-1/2}\sqrt{2\ln\ln m}$, and define
\[
q_{i}^{\ast}=\left(  1-p_{\left(  i\right)  }\right)  ^{-1}\left\{
im^{-1}-p_{\left(  i\right)  }-b_{m}^{\ast}\sqrt{p_{\left(  i\right)  }\left(
1-p_{\left(  i\right)  }\right)  }\right\}  ;
\]
then $\hat{\pi}_{1,m}^{\mathsf{MR}}=\min\left\{  1,\max\left\{  0,\max_{2\leq
i\leq m-2}q_{i}^{\ast}\right\}  \right\}  $ is the MR estimator. Storey's
estimator will be implemented by the \textsf{qvalue} package (version 2.14.1)
via the \textrm{`pi0.method=smoother'} option. All simulations will be done
with \textsf{R} version 3.5.0.

For an estimator $\hat{\pi}_{1,m}$ of $\pi_{1,m}$ or an estimator $\hat{\pi
}_{0,m}$ of $\tilde{\pi}_{0,m}$, its accuracy is measured by the excess
$\tilde{\delta}_{m}=\hat{\pi}_{1,m}\pi_{1,m}^{-1}-1$ or $\tilde{\delta}%
_{m}=\hat{\pi}_{0,m}\tilde{\pi}_{0,m}^{-1}-1$. For each experiment, the mean
$\mu_{m}^{\ast}$ and standard deviation $\sigma_{m}^{\ast}$ of $\tilde{\delta
}_{m}$ is estimated from independent realizations of the experiment. Among two
estimators, the one that has smaller $\sigma_{m}^{\ast}$ is taken to be more
stable, and the one that has both smaller $\sigma_{m}^{\ast}$ and smaller
$\left\vert \mu_{m}^{\ast}\right\vert $ is better. \textcolor{blue}{
In each boxplot in each figure of simulation results to be presented later, the horizontal bar has been programmed to represent the mean of the quantity being plotted and the black dots represent the outliers from the quantity being plotted.}

\subsection{Simulation design and results}

\label{simDesign}

We will simulate $\mathbf{z}\sim\mathcal{N}_{m}\left(  \boldsymbol{\mu
},\boldsymbol{\Sigma}\right)  $ with $\boldsymbol{\Sigma}$ as the identity
matrix. For $a<b$, let $\mathsf{U}\left(  a,b\right)  $ be
\textcolor{blue}{the uniform
distribution} on the closed interval $\left[  a,b\right]  $. We consider $6$
values for $m=10^{3}$, $5\times10^{3}$, $10^{4}$, $5\times10^{4}$, $10^{5}$ or
$5\times10^{5}$, and $2$ sparsity levels $\pi_{1,m}=0.2$ (indicating the dense
regime) or $\left(  \ln\ln m\right)  ^{-1}$ (indicating the moderately sparse
regime). The speed of the proposed estimators $t_{m}=\sqrt{0.99\ln m}$ (i.e.,
$t_{m}$ has tuning parameter $\gamma=0.495$) and $u_{m}=\tilde{u}_{m}=\left(
\ln\ln m\right)  ^{-1}$, where $u_{m}$ and $\tilde{u}_{m}$ are respectively
defined by (\ref{eq10b}) and \autoref{ThmTypeI-V}. This ensures $t_{m}%
^{-1}\left(  1+\max\left\{  u_{m}^{-1},\tilde{u}_{m}^{-1}\right\}  \right)
=o\left(  \pi_{1,m}\right)  $ and the consistency of the proposed estimator as
per \autoref{ThmIVa} and \autoref{ThmTypeI-V}. The simulated data are
generated as follows:

\begin{itemize}
\item Scenario 1 ``estimating $\pi_{1,m}$ for a bounded null'': set $a=-1$ and
$b=2$; generate $m_{0}$ $\mu_{i}$'s independently from $\mathsf{U}\left(
a+u_{m},b-u_{m}\right)  $, $m_{11}$ $\mu_{i}$'s independently from
$\mathsf{U}\left(  b+u_{m},b+6\right)  $, and $m_{11}$ $\mu_{i}$'s
independently from $\mathsf{U}\left(  a-4,a-u_{m}\right)  $, where
$m_{11}=\max\left\{  1,\lfloor  0.5m_{1}\rfloor  -\lfloor  m/\ln\ln m\rfloor
\right\}  $ and $\textcolor{blue}{\lfloor x \rfloor}$ is the integer part of $x \in \mathbb{R}$; set half of the remaining $m-m_{0}-2m_{11}$ $\mu_{i}$'s to be
$a$, and the rest to be $b$.

\item Scenario 2 ``estimating $\pi_{1,m}$ for a one-sided null'': set $b=0$;
generate $m_{0}$ $\mu_{i}$'s independently from $\mathsf{U}\left(
-4,b-u_{m}\right)  $, and $\lfloor  0.9m_{1}\rfloor  $ $\mu_{i}$'s
independently from $\mathsf{U}\left(  b+u_{m},b+6\right)  $; set the rest
$\mu_{i}$'s to be $b$.

\item Scenario 3 \textquotedblleft estimating average, truncated
2-norm\textquotedblright, i.e., estimating $\tilde{\pi}_{0,m}$ in
(\ref{eq21d}) with $\phi\left(  t\right)  =\left\vert t\right\vert
^{2}1_{\left\{  \left\vert t\right\vert \leq b\right\}  }\left(  t\right)  $
for a fixed $b>0$: set $b=2$ and $a=-2$; generate $m_{0}$ $\mu_{i}$'s
independently from $\mathsf{U}\left(  a,b\right)  $, $\lfloor  0.5m_{1}\rfloor
$ $\mu_{i}$'s independently from $\mathsf{U}\left(  b+u_{m},b+6\right)  $, and
the rest $\mu_{i}$'s independently from $\mathsf{U}\left(  b-4,b-u_{m}\right)
$. In this setting, $C^{-1}\pi_{1,m}\leq\tilde{\pi}_{1,m}\leq C\pi_{1,m}$
holds for some constant $C>0$ and $t_{m}^{-1}=o\left(  \tilde{\pi}%
_{1,m}\right)  $ holds, ensuring the consistency of the proposed estimator
$\hat{\varphi}_{1,m}$ as per \autoref{ThmFinal}.
\end{itemize}

Scenario 1 models the setting that when testing a bounded null in practice, it
is unlikely that there is always a positive proportion of means or medians
that are equal to either of the two boundary points, and Scenario 2 takes into
account that when testing a one-sided null with $0$ as the boundary point, it
is likely that there is a positive or diminishing proportion of means or
medians that are equal to $0$, as in differential gene expression studies.
Each triple of $\left(  m,\pi_{1,m},\Theta_{0}\right)  $ or $\left(
m,\tilde{\pi}_{0,m},\Theta_{0}\right)  $ determines an experiment, and there
are $36$ experiments in total. Each experiment is repeated independently $200$ times.
\textcolor{blue}{
As mentioned in the beginning of this section, the new estimators are implemented by numerical approximation, which we call ``numerical versions'',
and this numerical implementation/approximation causes numerical error. Further, the simulations are for the ``numerical versions'' of the new estimators rather than the new estimators themselves.}

\autoref{fig1} visualizes the simulation results,
\textcolor{blue}{and Table 1 provides numerical summaries that complement \autoref{fig1}. Please note again that these results are for the numerical implementation, i.e., numerical approximation, of the new estimators, rather than the new estimators themselves, even though the interpretations of the results will be for the new estimators.}
In this figure, Storey's estimator
is not shown since it is always $0$ for all experiments in Scenario 2. Such a
strange behavior of Storey's estimator has not been reported before and is
worth investigation but is not our focus here. A plausible explanation for
this is that Storey's estimator excessively over-estimates $\pi_{0,m}$ when no
p-value is uniformly distributed under the null. The following six
observations can be made: (i) for estimating the alternative proportion for a
one-sided null, the proposed estimator is more accurate than the MR estimator,
and it shows a trend of convergence towards consistency in the dense regime.
(ii) for estimating the alternative proportion for a bounded null, the
proposed estimator is accurate, and it shows a trend of convergence
towards consistency in the dense regime. (iii) the proposed estimator very
accurately estimates the average, truncated $2$-norm, with a trend of
convergence towards consistency. (iv) The MR estimator does not seem to
actively capture the changes in the number of alternative hypotheses as the
number of hypotheses varies because this estimator was quite inaccurate but
varied little as the alternative proportion changes. (v) For testing one-sided
null on the means of Gaussian random variables, in terms of accuracy, the HD1
estimator performs worse than all of the HD, MR and proposed estimators and
hence is not shown, whereas the HD estimator performs better than the MR
estimator but worse than the proposed estimator.
\textcolor{blue}{(vi) non-asymptotically the new estimator $\hat{\varphi}_m\left(t,\boldsymbol{\mu}\right)$ of the proportion of false nulls $\pi_{1,m}$ often over-estimates
$\pi_{1,m}$, meaning that its dual $\hat{\psi}_m\left(t,\boldsymbol{\mu}\right)$, which estimates the proportion of true nulls $\pi_{0,m}=1-\pi_{1,m}$, usually
under-estimates $\pi_{0,m}$. In terms of false discovery rate (FDR) control in nonasymptotic settings, an adaptive FDR procedure that uses the new estimators $\hat{\psi}_m\left(t,\boldsymbol{\mu}\right)$ may fail to maintain a prespecified nominal FDR, even though such a procedure may have larger power compared to its non-adaptive counterparts.
}
We remark that the accuracy
and speed of convergence of the proposed estimators can be improved by
employing more accurate Riemann sums for the integrals than currently used.

\textcolor{blue}{
Almost identical to what have been numerically observed \cite{Xiongzhi2025}, let us explain how the numerical implementation and numerical errors affect the observed trend of convergence for the new estimators.
Let $\hat{\pi}_{1,m}^{\textrm{New}}$ denote our ``New" estimators and $\hat{\pi}^{\dagger,\textrm{New}}_{1,m}$ their numerical implementations. Then,
the numerical error of the ``New" estimators $\hat{\pi}_{1,m}^{\textrm{new}}$ is $\tilde{e}_m^{\textrm{New}} = \hat{\pi}_{1,m}^{\textrm{New}} - \hat{\pi}^{\dagger,\textrm{New}}_{1,m}$. Recall $\tilde{\delta}_m = \hat{\pi}_{1,m}/\pi_{1,m} -1$, where $\hat{\pi}_{1,m}$ is an estimate of $\pi_{1,m}$, and that $\tilde{\delta}_m$ converges to $0$ as $m \to \infty$ is equivalent to the consistency of the estimator $\hat{\pi}_{1,m}$. Due to our numerical approximation, for our ``New" estimator, $\tilde{\delta}_m$ is actually computed as
\begin{equation*}
  \tilde{\delta}_m = \frac{\hat{\pi}^{\dagger,\textrm{New}}_{1,m}}{\pi_{1,m}} -1
  = \frac{\hat{\pi}^{\textrm{New}}_{1,m}}{\pi_{1,m}} - \frac{\tilde{e}_m^{\textrm{New}}}{\pi_{1,m}}  -1.
\end{equation*}
Since our theory has rigorously proved that $\hat{\pi}^{\textrm{New}}_{1,m}/\pi_{1,m} -1$ converges to $0$ in probability as $m \to \infty$, we see the actual $\tilde{\delta}_m$ computed for our ``New" estimator $\hat{\pi}^{\textrm{New}}_{1,m}$, as given above, satisfies
\begin{equation*}
  \tilde{\delta}_m  \approx \frac{\tilde{e}_m^{\textrm{New}}}{\pi_{1,m}} \quad \text{with high probability for large } m.
\end{equation*}
However, the numerical error $\tilde{e}_m^{\textrm{New}}$ may not converge to $0$ as $m \to \infty$.
So, in the dense regime when $\pi_{1,m} = 0.2$, we will see a trend of convergence for $\tilde{\delta}_m $ as $m$ increases. But such a convergence may stall if $\tilde{e}_m^{\textrm{New}}$ does not decrease with $m$. This is exactly what happened
for the dense regime in \autoref{fig1}. In contrast, in the moderately sparse regime when $\pi_{1,m}=1/\ln{(\ln{m})}$, $\pi_{1,m}$ monotonically decreases as $m$ increases and $\pi_{1,m}$ converges to $0$ as $m \to \infty$. So, when $\tilde{e}_m^{\textrm{New}}$ is not of smaller order than $\pi_{1,m}=1/\ln{(\ln{m})}$ as $m$ increases, we may see the actual $\tilde{\delta}_m$ on average increases with $m$. This is exactly what happened
for the moderately sparse regime in \autoref{fig1}.
Unless we increase the numerical precision or equivalently reduce the numerical error $\tilde{e}_m^{\textrm{New}}$ (dynamically also with respect to $m$), increasing $m$ but keeping the current numerical approximation scheme as described earlier will not allow us to see a clear trend of convergence of our
``New" estimators in the moderately sparse regime where $\lim_{m \to \infty}\pi_{1,m}=0$. Unfortunately, we have to defer to another manuscript how to systematically improve the numerical accuracy of implementing the new estimators. Nevertheless, for practical applications where we do not need to repeat an experiment many times as is done in the simulations here, we recommend keeping as many terms and using as fine partitions as one's computational recourses allow when respectively truncating the power series and forming Riemann sums that numerically approximate the definitions of the new estimators.
}

\section{Discussion}

\label{SecConcAndDisc}

For multiple testing a bounded or one-sided null on the means or medians of
random variables whose CDFs are members of a Type I location-shift family, we
have constructed uniformly consistent estimators of the corresponding
proportion of false null hypotheses via solutions to Lebesgue-Stieltjes
integral equations, for which consistency, speeds of convergence, and uniform
consistency classes have been obtained under independence between these random
variables. The strategy proposed in the Discussion section of
\cite{Chen:2018a} or that in Section 2.3 of \cite{Jin:2008}, i.e., choosing a
speed of convergence $t_{m}=\sqrt{2\gamma\ln m}$ with a suitable $\gamma$ that
controls the variance of the error term $e_{m}(t_{m})$ of a proportion
estimator for a finite $m$, can be used to adaptively determine the speed of
convergence (and hence the tuning parameter $\gamma$) for the proposed
estimators. These estimators can be used to develop adaptive versions of the
FDR procedure of \cite{Chen:2018d}, the \textquotedblleft BH\textquotedblright%
\ procedure of \cite{Benjamini:2001} under the conditions of their Theorem
5.2, the FDR and FNR procedures of \cite{Sarkar:2006}, or any other
conservative FDR or FNR procedure that is applicable to multiple testing
composite null hypotheses.

Our settings and results can be extended in several aspects as follows. First,
the constructions and uniform consistency of the proportion estimators
provided here can be easily extended to the setting where the null parameter
set belongs to the algebra generated by bounded, one-sided and point nulls.
Here the term \textquotedblleft algebra\textquotedblright\ refers to the
family of sets generated by applying any finite combination of set union,
intersection or complement to these three types of nulls. This covers null
sets $\left(  -\infty,b\right]  $ and $\left[  a,b\right]  $ that are more
conventional in classic textbooks such as \cite{Lehmann:2005}, even though we
have chosen null sets $\left(  -\infty,b\right)  $ and $\left(  a,b\right)  $
mainly due to how Dirichlet integrals given in \autoref{SecIllustration}
converge to their discontinuous, piecewise linear limiting functions;
see \autoref{secClosdeNulls} for details on dealing with null sets $\left(
-\infty,b\right]  $ and $\left[  a,b\right]  $. Second, the speed of
convergence and uniform consistency class for Construction I and Construction
II can be obtained for each non-Gaussian Type I location-shift family given in
\autoref{secExamples} to which our theory applies. Third, it is possible to
establish the uniform consistency of the proportion estimators in
\autoref{SecExtensions} for weakly dependent random variables that are
bivariate Gaussian via their associated Hermite polynomials. Fourth, for the
settings of one-sided and bounded nulls and suitable functionals of the
bounded null, respectively, following the principles in Section 3 of
\cite{Jin:2008}, all constructions we have provided can be applied to
consistently estimate the mixing proportions for two-component mixture models
at least one of whose components follows a Gaussian distribution.

Our proportion estimators may lead to consistent estimators of or tests on the
``sparsity level'' of regression coefficients in high-dimensional, sparse
Gaussian linear models, as explained below. Consider the model $\mathbf{y}%
=\mathbf{X}\boldsymbol{\beta}+\boldsymbol{\varepsilon}$, where $\mathbf{y}%
\in\mathbb{R}^{n}$ is the vector of observed response, $\mathbf{X}%
\in\mathbb{R}^{n\times p}$ a known design matrix, $\boldsymbol{\beta}%
\in\mathbb{R}^{p}$ the vector of unknown coefficients, and
$\boldsymbol{\varepsilon}\sim\mathcal{N}_{n}\left(  0,\boldsymbol{\Sigma
}\right)  $ the vector of random errors. Model selection via penalized least
squares to consistently estimate $\boldsymbol{\beta}$ in the setting $p\gg n$
and $\min\left\{  n,p\right\}  \rightarrow\infty$ often assumes an order for
the ``sparsity level'' $\varpi$, i.e., the number of nonzero entries of
$\boldsymbol{\beta}$, relative to $n$ and $p$; see, e.g.,
\cite{Zhao:2006,Lv:2009,Zhang:2010,Su16,JM2018}. However, there does not seem
to exist a consistent estimator of or test on $\varpi$ when $\mathbf{X}$ is
not a diagonal matrix. If we set $\boldsymbol{\mu}=\mathbf{X}\boldsymbol{\beta
}=\left(  \mu_{1},\ldots,\mu_{n}\right)  ^{\intercal}$ and assume
$\boldsymbol{\Sigma}=\sigma^{2}\boldsymbol{I}$ for some $\sigma^{2}>0$ where
$\boldsymbol{I}$ is the identity matrix, then we can look at $\ddot{\pi}%
_{1,n}=n^{-1}\sum_{i=1}^{n}\tilde{\phi}\left(  \mu_{i}\right)  $ for suitable
functions $\tilde{\phi}$ as test functions on properties of $\boldsymbol{\mu}$
and hence on $\boldsymbol{\beta}$, and consistent estimators of $\ddot{\pi
}_{1,n}$ for various $\tilde{\phi}$ (which include $\phi$ discussed in
\autoref{SecExtensions}) may lead to good lower or upper bounds on or
consistent estimators of $\varpi$ for nontrivial $\mathbf{X}$ and
$\boldsymbol{\beta}$.

Finally, our constructions essentially utilize suitable group structures of
the domain of the parameters and arguments of the CDFs and then apply special
transforms that are adapted to such group structures to obtain solutions to
the Lebesgue-Stieltjes integral equation. This general principle is also
applicable to probability distributions on Lie groups (which contains, e.g.,
the unit sphere a special case), and tools of harmonic analysis on Lie groups
can be used to construction proportion estimators for these distributions.
This has applications to modeling and analysis of data on Lie groups. We will
report on all these in other articles.

\appendix{}

\section{Type I location-shift families and closed or half-closed nulls}

We provide two methods to construct Type I location-shift families in \autoref{ConstructLocShift}, adaptations of our methods to the settings of closed or
half-closed nulls in \autoref{secClosdeNulls}, and implementations of our methods and estimators for Gaussian family in \autoref{Implementation}.

\subsection{Constructing Type I location-shift families}

\label{ConstructLocShift}

Recall \autoref{DefTLS} of Type I location-shift family in the main text. We
first provide two examples that are not Type I location-shift families:

\begin{example}
Let%
\[
\tilde{f}_{0}\left(  x\right)  =\left\{
\begin{array}
[c]{ccc}%
0.5e^{-x} & \text{if} & x\geq0\\
\frac{1}{\sqrt{2\pi}\sigma}\exp\left(  \frac{-x^{2}}{2\sigma^{2}}\right)  &
\text{if} & x\leq0
\end{array}
\right.  ,
\]
where $\sigma=\sqrt{2/\pi}$. Then $\tilde{f}_{0}$ is a density function on
$\mathbb{R}$. Further, $\int_{0}^{\infty}x\tilde{f}_{0}\left(  x\right)
dx=0.5$ and $\int_{-\infty}^{0}x\tilde{f}_{0}\left(  x\right)  dx=-1/\pi$. Let
$X$ have density $\tilde{f}_{0}\left(  x\right)  $ and $Y=X-\left(
0.5-1/\pi\right)  $. Denote the density of $Y$ by $f_{0}$. Then $\mathbb{E}
\left(  Y\right)  =0$ and $F_{\mu}\left(  x\right)  =\int_{-\infty}^{x}
f_{0}\left(  y-\mu\right)  dy,\mu\in\mathbb{R}$ form a location-shift family
$\mathcal{F}$. But this $\mathcal{F}$ is not a Type I location-shift family
since $f_{0}$ is not an even function.
\end{example}

\begin{example}
Consider a P\'{o}lya-type characteristic function (CF) $\hat{F}\left(
t\right)  =\left\{  1-\left\vert t\right\vert \right\}  1_{\left\{  \left\vert
t\right\vert \leq1\right\}  }$ (see, e.g., page 85 of \cite{Lukacs:1970}).
Then the density function corresponding to $\hat{F}\left(  t\right)  $ is
$f\left(  x\right)  =\pi^{-1}x^{-2}\left(  1-\cos x\right)  $ for
$x\in\mathbb{R}$, and $F_{\mu}\left(  x\right)  =\int_{-\infty}^{x}%
f_{0}\left(  y-\mu\right)  dy,\mu\in\mathbb{R}$ form a location-shift family
$\mathcal{F}$. However, $\hat{F}\left(  t\right)  =0$ for all $\left\vert
t\right\vert \geq1$. So, this $\mathcal{F}$ is not a Type I location-shift family.
\end{example}

Now let us describe two methods to construct Type I location-shift families,
for which we need the following ``P\'{o}lya's condition'' (see, e.g., page 83
of \cite{Lukacs:1970}):

\begin{lemma}
[G. P\'{o}lya]\label{lemmaPolyaType}Let $h:\mathbb{R}\rightarrow\mathbb{R}$ be
a continuous function. If $h$ is even, $h\left(  0\right)  =1$, $h\left(
t\right)  $ is convex for $t>0$, and $\lim_{t\rightarrow\infty}h\left(
t\right)  =0$, then $h$ is the CF of an absolutely continuous CDF.
\end{lemma}

A CF satisfying the conditions in \autoref{lemmaPolyaType} is called a
``P\'{o}lya-type CF''. Since the density function $f$ of an absolutely
continuous CDF $F$ is even if and only if the CF $\hat{F}$ of $F$ is real, we
see that a random variable $X$ whose CDF $F_{0}$ has a real CF $\hat{F}_{0}$
must have zero expectation. So, \autoref{lemmaPolyaType} tells us that, if we
pick a P\'{o}lya-type CF $h_{0}$ that is nowhere $0$ on $\mathbb{R}$ and is
Lebesgue integrable on $\mathbb{R}$, then its Fourier inverse $\hat{h}_{0}$
generates a Type I location-shift family with members CDFs as $F_{\mu}\left(
x\right)  =\int_{-\infty}^{x}\hat{h}_{0}\left(  y-\mu\right)  dy,\mu
\in\mathbb{R}$.

A second method to construct a Type I location-shift family is based on a key
result of \cite{Tuck:2006}, which can also be derived by the arguments in the
proof of Theorem 4.3.1 of \cite{Lukacs:1970}. This result is

\begin{theorem}
[\cite{Tuck:2006}]\label{ThmTuck}If a function $f$ defined on $\mathbb{R}%
_{>0}$ is Lebesgue integrable on $\mathbb{R}_{>0}$ and strictly convex on
$\mathbb{R}_{>0}$, then the function $g\left(  t\right)  =\int_{0}^{\infty
}f\left(  x\right)  \cos\left(  tx\right)  dx$ is positive for all $t>0$.
\end{theorem}

\autoref{ThmTuck} says that a convex function has a positive \textquotedblleft
Fourier-cosine\textquotedblright\ transform, and it shows us a second way to
construct Type I location-shift families, i.e., take an even density function
$f_{0}$ on $\mathbb{R\ }$such that $f_{0}$ is strictly convex on
$\mathbb{R}_{>0}$, then $F_{\mu}\left(  x\right)  =\int_{-\infty}^{x}%
f_{0}\left(  y-\mu\right)  dy,\mu\in\mathbb{R}$ form a Type I location-shift family.

Finally, the complex analytic techniques of \cite{Tuck:2006} can be used to
construct Type I location-shift families from an even density function on
$\mathbb{R}$ that is bell-shaped. These techniques mainly use analytic
continuation, contour deformation, and residue calculus, which we will not
detail here.

\subsection{Estimators for closed or half-closed nulls}

\label{secClosdeNulls}

Let us show how to adapt the constructions, the estimators, their
concentration inequalities, and their consistency results to estimating the
proportion $\pi_{1,m}$ when the null hypotheses are closed or half-closed
intervals; see \autoref{SecSummaryBndClosedNull}, \autoref{SecSummaryOneSidedClosedNull} and \autoref{secSummaryExtension}, respectively.

\subsubsection{The case of a bounded null}

\label{SecSummaryBndClosedNull}

When $\Theta_{0}=\left[  a,b\right]  $, we can just set%
\begin{equation}
\left\{
\begin{array}
[c]{l}%
K\left(  t,x\right)  =K_{1}\left(  t,x\right)  +2^{-1}\left\{  K_{1,0}\left(
t,x;a\right)  +K_{1,0}\left(  t,x;b\right)  \right\} \\
\psi\left(  t,\mu\right)  =\psi_{1}\left(  t,\mu\right)  +2^{-1}\left\{
\psi_{1,0}\left(  t,\mu;a\right)  +\psi_{1,0}\left(  t,\mu;b\right)  \right\}
\end{array}
\right.  \label{BndClosedNull}%
\end{equation}
in comparison to the construction when $\Theta_{0}=\left(  a,b\right)  $ as%
\begin{equation}
\left\{
\begin{array}
[c]{l}%
K\left(  t,x\right)  =K_{1}\left(  t,x\right)  -2^{-1}\left\{  K_{1,0}\left(
t,x;a\right)  +K_{1,0}\left(  t,x;b\right)  \right\} \\
\psi\left(  t,\mu\right)  =\psi_{1}\left(  t,\mu\right)  -2^{-1}\left\{
\psi_{1,0}\left(  t,\mu;a\right)  +\psi_{1,0}\left(  t,\mu;b\right)  \right\}
\end{array}
\right.  . \label{BndOpenNull}%
\end{equation}
The definitions of the estimator and its expectation for either $\Theta
_{0}=\left(  a,b\right)  $ or $\Theta_{0}=\left[  a,b\right]  $ remain
identical as%
\[
\hat{\varphi}_{m}\left(  t,\mathbf{z}\right)  =m^{-1}\sum_{i=1}^{m}\left\{
1-K\left(  t,z_{i}\right)  \right\}  \text{ \ and\ }\varphi_{m}\left(
t,\boldsymbol{\mu}\right)  =m^{-1}\sum_{i=1}^{m}\left\{  1-\psi\left(
t,\mu_{i}\right)  \right\}  .
\]

When $\Theta_{0}=\left(  a,b\right)  $, in the proofs for the estimator
$\hat{\varphi}_{m}\left(  t,\mathbf{z}\right)  $, we have used $e_{1,m}\left(
t\right)  =\hat{\varphi}_{1,m}\left(  t,\mathbf{z}\right)  -\varphi
_{1,m}\left(  t,\boldsymbol{\mu}\right)  $, where%
\[
\hat{\varphi}_{1,m}\left(  t,\mathbf{z}\right)  =m^{-1}\sum_{i=1}^{m}%
K_{1}\left(  t,z_{i}\right)  \ \text{and }\varphi_{1,m}\left(
t,\boldsymbol{\mu}\right)  =\mathbb{E}\left\{  \hat{\varphi}_{1,m}\left(
t,\mathbf{z}\right)  \right\}  ,
\]
$e_{1,0,m}\left(  t,\tau\right)  =\hat{\varphi}_{1,0,m}\left(  t,\mathbf{z}%
;\tau\right)  -$ $\varphi_{1,0,m}\left(  t,\boldsymbol{\mu};\tau\right)  $,
$\tau\in\left\{  a,b\right\}  $, where%
\[
\hat{\varphi}_{1,0,m}\left(  t,\mathbf{z};\tau\right)  =m^{-1}\sum_{i=1}%
^{m}K_{1,0}\left(  t,z_{i};\tau\right)  \text{ and }\varphi_{1,0,m}\left(
t,\boldsymbol{\mu};\tau\right)  =\mathbb{E}\left\{  \hat{\varphi}%
_{1,0,m}\left(  t,\mathbf{z};\tau\right)  \right\}  ,
\]
and%
\begin{equation}
e_{m}\left(  t\right)  =\hat{\varphi}_{m}\left(  t,\mathbf{z}\right)
-\varphi_{m}\left(  t,\boldsymbol{\mu}\right)  =-e_{1,m}\left(  t\right)
+2^{-1}e_{1,0,m}\left(  t,a\right)  +2^{-1}e_{1,0,m}\left(  t,b\right)  .
\label{ErrorBndOpenNull}%
\end{equation}
When $\Theta_{0}=\left(  a,b\right)  $, to upper bound the variance of
$e_{m}\left(  t\right)  $, we have upper bounded the variances of
$e_{1,m}\left(  t\right)  $, $e_{1,0,m}\left(  t,a\right)  $ and
$e_{1,0,m}\left(  t,b\right)  $ individually, and then directly replaced each
variance in each summand on the right-hand side of the inequality%
\begin{equation}
\mathbb{V}\left[  e_{m}\left(  t\right)  \right]  \leq2\mathbb{V}\left\{
e_{1,m}\left(  t\right)  \right\}  +\mathbb{V}\left[  e_{1,0,m}\left(
t,a\right)  \right]  +\mathbb{V}\left[  e_{1,0,m}\left(  t,b\right)  \right]
\label{VarianceBndOpenNull}%
\end{equation}
with these individual variance upper bounds. In addition, when $\Theta
_{0}=\left(  a,b\right)  $, to upper bound the deviation of $\left\vert
e_{m}\left(  t\right)  \right\vert $, we have upper bounded the deviation of
each of $\left\vert e_{1,m}\left(  t\right)  \right\vert $, $\left\vert
e_{1,0,m}\left(  t,a\right)  \right\vert $ and $\left\vert e_{1,0,m}\left(
t,b\right)  \right\vert $ individually, and then directly replaced each of
$\left\vert e_{1,m}\left(  t\right)  \right\vert $, $\left\vert e_{1,0,m}%
\left(  t,a\right)  \right\vert $ and $\left\vert e_{1,0,m}\left(  t,b\right)
\right\vert $ in the right-hand side of the inequality%
\begin{equation}
\left\vert e_{m}\left(  t\right)  \right\vert \leq\left\vert e_{1,m}\left(
t\right)  \right\vert +2^{-1}\left\vert e_{1,0,m}\left(  t,a\right)
\right\vert +2^{-1}\left\vert e_{1,0,m}\left(  t,b\right)  \right\vert
\label{DeviationBndOpenNull}%
\end{equation}
by these individual upper bounds.

In the setting of the closed null $\Theta_{0}=\left[  a,b\right]  $,
(\ref{BndClosedNull}) implies that (\ref{ErrorBndOpenNull}) becomes%
\begin{equation}
e_{m}\left(  t\right)  =\hat{\varphi}_{m}\left(  t,\mathbf{z}\right)
-\varphi_{m}\left(  t,\boldsymbol{\mu}\right)  =-e_{1,m}\left(  t\right)
-2^{-1}e_{1,0,m}\left(  t,a\right)  -2^{-1}e_{1,0,m}\left(  t,b\right)
\label{ErrorBndClosedNull}%
\end{equation}
However, (\ref{VarianceBndOpenNull}) and (\ref{DeviationBndOpenNull}) remain
valid for $e_{m}\left(  t\right)  $ in (\ref{ErrorBndClosedNull}).

When $\Theta_{0}=\left(  a,b\right)  $,%
\begin{equation}
\varphi_{m}\left(  t,\boldsymbol{\mu}\right)  =1-\varphi_{1,m}\left(
t,\boldsymbol{\mu}\right)  +2^{-1}\varphi_{1,0,m}\left(  t,\boldsymbol{\mu
};a\right)  +2^{-1}\varphi_{1,0,m}\left(  t,\boldsymbol{\mu};b\right)
=\sum_{i=1}^{5}\widetilde{d}_{1,m} \label{OracleBnd}%
\end{equation}
and%
\begin{equation}
\pi_{1,m}^{-1}{\varphi_{m}\left(  t,\boldsymbol{\mu}\right)  -1}=\pi
_{1,m}^{-1}\widetilde{d}_{1,m}-1+\pi_{1,m}^{-1}\widetilde{d}_{2,m}+\pi
_{1,m}^{-1}\widetilde{d}_{3,m}+\pi_{1,m}^{-1}\widetilde{d}_{4,m}+\pi
_{1,m}^{-1}\widetilde{d}_{5,m}, \label{OracleBndAdjusted}%
\end{equation}
where%
\[
\left\{
\begin{array}
[c]{l}%
\widetilde{d}_{1,m}=1-m^{-1}\sum\nolimits_{\left\{  j:\mu_{j}\in\left(
a,b\right)  \right\}  }\psi_{1}\left(  t,\mu_{j}\right) \\
\widetilde{d}_{2,m}=-m^{-1}\sum\nolimits_{\left\{  j:\mu_{j}=a\right\}  }%
\psi_{1}\left(  t,\mu_{j}\right)  +2^{-1}m^{-1}\sum\nolimits_{\left\{
j:\mu_{j}=a\right\}  }{\tilde{\psi}}_{1,0}\left(  t,\mu_{j};a\right) \\
\widetilde{d}_{3,m}=-m^{-1}\sum\nolimits_{\left\{  j:\mu_{j}=b\right\}  }%
\psi_{1}\left(  t,\mu_{j}\right)  +2^{-1}m^{-1}\sum\nolimits_{\left\{
j:\mu_{j}=b\right\}  }{\tilde{\psi}}_{1,0}\left(  t,\mu_{j};b\right) \\
\widetilde{d}_{4,m}=2^{-1}m^{-1}\sum\nolimits_{\left\{  j:\mu_{j}\neq
a\right\}  }{\tilde{\psi}}_{1,0}\left(  t,\mu_{j};a\right)  +2^{-1}m^{-1}%
\sum\nolimits_{\left\{  j:\mu_{j}\neq b\right\}  }{\tilde{\psi}}_{1,0}\left(
t,\mu_{j};b\right) \\
\widetilde{d}_{5,m}=-m^{-1}\sum\nolimits_{\left\{  j:\mu_{j}<a\right\}  }%
\psi_{1}\left(  t,\mu_{j}\right)  -m^{-1}\sum\nolimits_{\left\{  j:\mu
_{j}>b\right\}  }\psi_{1}\left(  t,\mu_{j}\right)
\end{array}
\right.  .
\]
Further, when $\Theta_{0}=\left(  a,b\right)  $, to upper bound $\left\vert
\pi_{1,m}^{-1}{\varphi_{m}\left(  t,\boldsymbol{\mu}\right)  -1}\right\vert $,
we have replaced each $\left\vert \widetilde{d}_{j,m}\right\vert ,2\leq
j\leq5$ by its upper bound $\hat{d}_{j,m},2\leq j\leq5$ and replaced
$\left\vert \pi_{1,m}^{-1}\widetilde{d}_{1,m}-1\right\vert $ by its upper
bound $\hat{d}_{0,m}$ in the inequality%
\begin{equation}
\left\vert \pi_{1,m}^{-1}{\varphi_{m}\left(  t,\boldsymbol{\mu}\right)
-1}\right\vert \leq\left\vert \pi_{1,m}^{-1}\widetilde{d}_{1,m}-1\right\vert
+\pi_{1,m}^{-1}\left\vert \widetilde{d}_{2,m}\right\vert +\pi_{1,m}%
^{-1}\left\vert \widetilde{d}_{3,m}\right\vert +\pi_{1,m}^{-1}\left\vert
\widetilde{d}_{4,m}\right\vert +\pi_{1,m}^{-1}\left\vert \widetilde{d}%
_{5,m}\right\vert . \label{OracleBndUpperBound}%
\end{equation}

In case $\Theta_{0}=\left[  a,b\right]  $, (\ref{OracleBnd}) becomes%
\begin{align*}
\varphi_{m}\left(  t,\boldsymbol{\mu}\right)   &  =1-\varphi_{1,m}\left(
t,\boldsymbol{\mu}\right)  -2^{-1}\varphi_{1,0,m}\left(  t,\boldsymbol{\mu
};a\right)  -2^{-1}\varphi_{1,0,m}\left(  t,\boldsymbol{\mu};b\right) \\
&  =\widetilde{d}_{1,m}+\widetilde{d}_{2,m}^{\ast}+\widetilde{d}_{3,m}^{\ast
}-\widetilde{d}_{4,m}+\widetilde{d}_{5,m}%
\end{align*}
and%
\begin{equation}
\left\vert \pi_{1,m}^{-1}{\varphi_{m}\left(  t,\boldsymbol{\mu}\right)
-1}\right\vert \leq\left\vert \pi_{1,m}^{-1}\left(  \widetilde{d}%
_{1,m}+\widetilde{d}_{2,m}^{\ast}+\widetilde{d}_{3,m}^{\ast}\right)
-1\right\vert +\pi_{1,m}^{-1}\left\vert \widetilde{d}_{4,m}\right\vert
+\pi_{1,m}^{-1}\left\vert \widetilde{d}_{5,m}\right\vert
\label{OracleBndClosed}%
\end{equation}
where%
\[
\left\{
\begin{array}
[c]{l}%
\widetilde{d}_{2,m}^{\ast}=-m^{-1}\sum\nolimits_{\left\{  j:\mu_{j}=a\right\}
}\psi_{1}\left(  t,\mu_{j}\right)  -2^{-1}m^{-1}\sum\nolimits_{\left\{
j:\mu_{j}=a\right\}  }{\tilde{\psi}}_{1,0}\left(  t,\mu_{j};a\right) \\
\widetilde{d}_{3,m}^{\ast}=-m^{-1}\sum\nolimits_{\left\{  j:\mu_{j}=b\right\}
}\psi_{1}\left(  t,\mu_{j}\right)  -2^{-1}m^{-1}\sum\nolimits_{\left\{
j:\mu_{j}=b\right\}  }{\tilde{\psi}}_{1,0}\left(  t,\mu_{j};b\right)
\end{array}
\right.  .
\]
However,
\[
\pi_{1,m}^{-1}\left\vert \widetilde{d}_{1,m}+\widetilde{d}_{2,m}^{\ast
}+\widetilde{d}_{3,m}^{\ast}-1\right\vert \leq\hat{d}_{0,m}+\pi_{1,m}^{-1}%
\hat{d}_{2,m}+\pi_{1,m}^{-1}\hat{d}_{3,m}.
\]
So, the upper bound for $\left\vert \pi_{1,m}^{-1}{\varphi_{m}\left(
t,\boldsymbol{\mu}\right)  -1}\right\vert $ in (\ref{OracleBndUpperBound})
when $\Theta_{0}=\left(  a,b\right)  $ is also an upper bound for $\left\vert
\pi_{1,m}^{-1}{\varphi_{m}\left(  t,\boldsymbol{\mu}\right)  -1}\right\vert $
in (\ref{OracleBndClosed}) when $\Theta_{0}=\left[  a,b\right]  $.

In summary, all results we have derived for the construction
(\ref{BndOpenNull}) when $\Theta_{0}=\left(  a,b\right)  $ remain valid for
the construction (\ref{BndOpenNull}) when $\Theta_{0}=\left[  a,b\right]  $.

\subsubsection{The case of a one-sided null}

\label{SecSummaryOneSidedClosedNull}

When $\Theta_{0}=(-\infty,0]$, we can just set%
\begin{equation}
\left\{
\begin{array}
[c]{l}%
K\left(  t,x\right)  =2^{-1}-\Re\left\{  K_{1}^{\dagger}\left(  t,x\right)
\right\}  +2^{-1}K_{1,0}\left(  t,x;0\right) \\
\psi\left(  t,\mu\right)  =2^{-1}-\psi_{1}\left(  t,\mu\right)  +2^{-1}%
\psi_{1,0}\left(  t,\mu;0\right)
\end{array}
\right.  . \label{OneSideClosedNull}%
\end{equation}
in comparison to the construction when $\Theta_{0}=\left(  -\infty,0\right)  $
as%
\begin{equation}
\left\{
\begin{array}
[c]{l}%
K\left(  t,x\right)  =2^{-1}-\Re\left\{  K_{1}^{\dagger}\left(  t,x\right)
\right\}  -2^{-1}K_{1,0}\left(  t,x;0\right) \\
\psi\left(  t,\mu\right)  =2^{-1}-\psi_{1}\left(  t,\mu\right)  -2^{-1}%
\psi_{1,0}\left(  t,\mu;0\right)
\end{array}
\right.  . \label{OneSideOpenNull}%
\end{equation}
The definitions of the estimator and its expectation for either $\Theta
_{0}=(-\infty,0]$ or $\Theta_{0}=\left(  -\infty,0\right)  $ remain identical
as%
\[
\hat{\varphi}_{m}\left(  t,\mathbf{z}\right)  =m^{-1}\sum_{i=1}^{m}\left\{
1-K\left(  t,z_{i}\right)  \right\}  \text{ \ and\ }\varphi_{m}\left(
t,\boldsymbol{\mu}\right)  =m^{-1}\sum_{i=1}^{m}\left\{  1-\psi\left(
t,\mu_{i}\right)  \right\}  .
\]
Then $e_{m}\left(  t\right)  =e_{1,m}\left(  t\right)  +2^{-1}e_{1,m,0}\left(
t,0\right)  $ when $\Theta_{0}=\left(  -\infty,0\right)  $ becomes%
\begin{equation}
e_{m}\left(  t\right)  =e_{1,m}\left(  t\right)  -2^{-1}e_{1,0,m}\left(
t,0\right)  \text{ when }\Theta_{0}=(-\infty,0].
\label{ErrorOneSideClosedNull}%
\end{equation}

When $\Theta_{0}=(-\infty,0]$, to upper bound the variance of $e_{m}\left(
t\right)  $, we have upper bounded the variances of $e_{1,m}\left(  t\right)
$ and $e_{1,0,m}\left(  t,0\right)  $ individually, and then directly replaced
each variance in each summand on the right-hand side of the inequality%
\begin{equation}
\mathbb{V}\left[  e_{m}\left(  t\right)  \right]  \leq2\mathbb{V}\left\{
e_{1,m}\left(  t\right)  \right\}  +2^{-1}\mathbb{V}\left[  e_{1,0,m}\left(
t,0\right)  \right]  \label{VBndOneSideNull}%
\end{equation}
with these individual variance upper bounds. In addition, when $\Theta
_{0}=(-\infty,0]$, to upper bound the deviation of $\left\vert e_{m}\left(
t\right)  \right\vert $, we have upper bounded the deviation of each of
$\left\vert e_{1,m}\left(  t\right)  \right\vert $ and $\left\vert
e_{1,0,m}\left(  t,0\right)  \right\vert $ individually, and then directly
replaced each of $\left\vert e_{1,m}\left(  t\right)  \right\vert $ and
$\left\vert e_{1,0,m}\left(  t,0\right)  \right\vert $ in the right-hand side
of the inequality%
\begin{equation}
\left\vert e_{m}\left(  t\right)  \right\vert \leq\left\vert e_{1,m}\left(
t\right)  \right\vert +2^{-1}\left\vert e_{1,0,m}\left(  t,0\right)
\right\vert \label{ErrorBndOneSideNull}%
\end{equation}
by these individual upper bounds. However, (\ref{VBndOneSideNull}) and
(\ref{ErrorBndOneSideNull}) remain valid for $e_{m}\left(  t\right)  $ in
(\ref{ErrorOneSideClosedNull}) when $\Theta_{0}=(-\infty,0]$.

When $\Theta_{0}=\left(  -\infty,0\right)  $, we have%
\[
\varphi_{m}\left(  t,\boldsymbol{\mu}\right)  =2^{-1}+\varphi_{1,m}\left(
t,\boldsymbol{\mu}\right)  +2^{-1}\varphi_{1,0,m}\left(  t,\boldsymbol{\mu
};{0}\right)  =\bar{d}_{1,m}+\bar{d}_{2,m}+\bar{d}_{3,m}+\bar{d}_{4,m},
\]
where%
\[
\left\{
\begin{array}
[c]{l}%
\bar{d}_{1,m}=m^{-1}\sum\nolimits_{\left\{  i:\mu_{i}>{0}\right\}  }\left(
2^{-1}+\psi_{1}\left(  t,\mu_{i}\right)  \right) \\
\bar{d}_{2,m}=m^{-1}\sum\nolimits_{\left\{  i:\mu_{i}={0}\right\}  }\left(
2^{-1}+\psi_{1}\left(  t,\mu_{i}\right)  +2^{-1}\tilde{\psi}_{1,0}\left(
t,\mu_{i};{0}\right)  \right) \\
\bar{d}_{3,m}=m^{-1}\sum\nolimits_{\left\{  i:\mu_{i}<{0}\right\}  }\left(
2^{-1}+\psi_{1}\left(  t,\mu_{i}\right)  \right) \\
\bar{d}_{4,m}=2^{-1}m^{-1}\sum\nolimits_{\left\{  i:\mu_{i}\neq{0}\right\}
}\tilde{\psi}_{1,0}\left(  t,\mu_{i};{0}\right)
\end{array}
\right.
\]
and specifically $\bar{d}_{2,m}=m^{-1}\sum\nolimits_{\left\{  i:\mu_{i}%
={0}\right\}  }1$. Further, when $\Theta_{0}=\left(  -\infty,0\right)  $, to
upper bound $\left\vert \pi_{1,m}^{-1}{\varphi_{m}\left(  t,\boldsymbol{\mu
}\right)  -1}\right\vert $, we have replaced each $\left\vert \bar{d}%
_{j,m}\right\vert ,3\leq j\leq4$ by its upper bound and replaced $\left\vert
\pi_{1,m}^{-1}\left(  \bar{d}_{1,m}+\bar{d}_{2,m}\right)  -1\right\vert $ by
its upper bound in the inequality%
\[
\left\vert \pi_{1,m}^{-1}{\varphi_{m}\left(  t,\boldsymbol{\mu}\right)
-1}\right\vert \leq\left\vert \pi_{1,m}^{-1}\left(  \bar{d}_{1,m}+\bar
{d}_{2,m}\right)  -1\right\vert +\pi_{1,m}^{-1}\left\vert \bar{d}%
_{3,m}\right\vert +\pi_{1,m}^{-1}\left\vert \bar{d}_{4,m}\right\vert ,
\]
where%
\begin{equation}
\pi_{1,m}^{-1}\left(  \bar{d}_{1,m}+\bar{d}_{2,m}\right)  -1=\pi_{1,m}%
^{-1}m^{-1}\sum\nolimits_{\left\{  i:\mu_{i}>{0}\right\}  }\left(  \psi
_{1}\left(  t,\mu_{i}\right)  -2^{-1}\right)  . \label{OracleTrickOneSide}%
\end{equation}
Specifically, the upper bound on $\left\vert \pi_{1,m}^{-1}\left(  \bar
{d}_{1,m}+\bar{d}_{2,m}\right)  -1\right\vert $ is based on the inequality
$\left\vert \psi_{1}\left(  t,\mu_{i}\right)  -2^{-1}\right\vert \leq2\left(
t\tilde{u}_{m}\right)  ^{-1}$ for $\mu_{i}>0$, where $\tilde{u}_{m}%
=\min_{\left\{  j:\mu_{j}\neq0\right\}  }\left\vert \mu_{j}\right\vert $.

In contrast, when $\Theta_{0}=(-\infty,0]$, we have%
\[
\varphi_{m}\left(  t,\boldsymbol{\mu}\right)  =2^{-1}+\varphi_{1,m}\left(
t,\boldsymbol{\mu}\right)  -2^{-1}\varphi_{1,0,m}\left(  t,\boldsymbol{\mu
};{0}\right)  =\bar{d}_{1,m}+\bar{d}_{2,m}^{\ast}+\bar{d}_{3,m}-\bar{d}%
_{4,m},
\]
and%
\[
\left\vert \pi_{1,m}^{-1}{\varphi_{m}\left(  t,\boldsymbol{\mu}\right)
-1}\right\vert \leq\left\vert \pi_{1,m}^{-1}\bar{d}_{1,m}-1\right\vert
+\pi_{1,m}^{-1}\left\vert \bar{d}_{3,m}\right\vert +\pi_{1,m}^{-1}\left\vert
\bar{d}_{4,m}\right\vert ,
\]
where%
\[
\bar{d}_{2,m}^{\ast}=m^{-1}\sum\nolimits_{\left\{  i:\mu_{i}={0}\right\}
}\left(  2^{-1}+\psi_{1}\left(  t,\mu_{i}\right)  -2^{-1}\tilde{\psi}%
_{1,0}\left(  t,\mu_{i};{0}\right)  \right)  =0.
\]
However, again%
\[
\pi_{1,m}^{-1}\bar{d}_{1,m}-1=\pi_{1,m}^{-1}m^{-1}\sum\nolimits_{\left\{
i:\mu_{i}>{0}\right\}  }\left(  \psi_{1}\left(  t,\mu_{i}\right)
-2^{-1}\right)  ,
\]
whose right-hand side is identical for that of (\ref{OracleTrickOneSide}).
Therefore, the upper bound we have derived for $\left\vert \pi_{1,m}%
^{-1}{\varphi_{m}\left(  t,\boldsymbol{\mu}\right)  -1}\right\vert $ when
$\Theta_{0}=\left(  -\infty,0\right)  $ is also an upper bound for $\left\vert
\pi_{1,m}^{-1}{\varphi_{m}\left(  t,\boldsymbol{\mu}\right)  -1}\right\vert $
when $\Theta_{0}=(-\infty,0]$.

In summary, all results we have derived for the setting $\Theta_{0}%
=(-\infty,0)$ for the construction (\ref{OneSideOpenNull}) remain valid for
the setting $\Theta_{0}=(-\infty,0]$ for the construction
(\ref{OneSideClosedNull}).

\subsubsection{The setting for the extensions}

\label{secSummaryExtension}

When $\Theta_{0}=\left[  a,b\right]  $, we can just set%
\begin{equation}
\left\{
\begin{array}
[c]{l}%
K\left(  t,x\right)  =K_{1}\left(  t,x\right)  +2^{-1}\left\{  \phi\left(
a\right)  K_{1,0}\left(  t,x;a\right)  +\phi\left(  b\right)  K_{1,0}\left(
t,x;b\right)  \right\} \\
\psi\left(  t,\mu\right)  =\psi_{1}\left(  t,\mu\right)  +2^{-1}\left\{
\phi\left(  a\right)  \psi_{1,0}\left(  t,\mu;a\right)  +\phi\left(  b\right)
\psi_{1,0}\left(  t,\mu;b\right)  \right\}
\end{array}
\right.  \label{ExtensionConstructionCloseNull}%
\end{equation}
in comparison to the construction for $\Theta_{0}=\left(  a,b\right)  $ as%
\begin{equation}
\left\{
\begin{array}
[c]{l}%
K\left(  t,x\right)  =K_{1}\left(  t,x\right)  -2^{-1}\left\{  \phi\left(
a\right)  K_{1,0}\left(  t,x;a\right)  +\phi\left(  b\right)  K_{1,0}\left(
t,x;b\right)  \right\} \\
\psi\left(  t,\mu\right)  =\psi_{1}\left(  t,\mu\right)  -2^{-1}\left\{
\phi\left(  a\right)  \psi_{1,0}\left(  t,\mu;a\right)  +\phi\left(  b\right)
\psi_{1,0}\left(  t,\mu;b\right)  \right\}
\end{array}
\right.  . \label{ExtOpenNull}%
\end{equation}
Again the definitions of the estimator and its expectation for either
$\Theta_{0}=\left(  a,b\right)  $ or $\Theta_{0}=\left[  a,b\right]  $ remain
identical as%
\[
\hat{\varphi}_{m}\left(  t,\mathbf{z}\right)  =m^{-1}\sum_{i=1}^{m}K\left(
t,z_{i}\right)  \text{ \ and\ }\varphi_{m}\left(  t,\boldsymbol{\mu}\right)
=m^{-1}\sum_{i=1}^{m}\psi\left(  t,\mu_{i}\right)  .
\]
For $\Theta_{0}=\left(  a,b\right)  $ we have%
\[
e_{m}\left(  t\right)  =\hat{\varphi}_{m}\left(  t,\mathbf{z}\right)
-\varphi_{m}\left(  t,\boldsymbol{\mu}\right)  =e_{1,m}\left(  t\right)
-2^{-1}\phi\left(  a\right)  e_{1,0,m}\left(  t,a\right)  -2^{-1}\phi\left(
a\right)  e_{1,0,m}\left(  t,b\right)  ,
\]
whereas for $\Theta_{0}=\left[  a,b\right]  $ we have%
\[
e_{m}\left(  t\right)  =\hat{\varphi}_{m}\left(  t,\mathbf{z}\right)
-\varphi_{m}\left(  t,\boldsymbol{\mu}\right)  =e_{1,m}\left(  t\right)
+2^{-1}\phi\left(  a\right)  e_{1,0,m}\left(  t,a\right)  +2^{-1}\phi\left(
a\right)  e_{1,0,m}\left(  t,b\right)  .
\]
Again, since have employed the inequalities%
\[
\left\{
\begin{array}
[c]{l}%
\mathbb{V}\left[  e_{m}\left(  t\right)  \right]  \leq2\mathbb{V}\left\{
e_{1,m}\left(  t\right)  \right\}  +\left\Vert \omega\right\Vert _{\infty}%
^{2}\mathbb{V}\left[  e_{1,0,m}\left(  t,a\right)  \right]  +\left\Vert
\omega\right\Vert _{\infty}^{2}\mathbb{V}\left[  e_{1,0,m}\left(  t,b\right)
\right] \\
\left\vert e_{m}\left(  t\right)  \right\vert \leq\left\vert e_{1,m}\left(
t\right)  \right\vert +2^{-1}\left\Vert \omega\right\Vert _{\infty}\left\vert
e_{1,0,m}\left(  t,a\right)  \right\vert +2^{-1}\left\Vert \omega\right\Vert
_{\infty}\left\vert e_{1,0,m}\left(  t,b\right)  \right\vert
\end{array}
\right.
\]
and then directly replaced each variance in each summand on the right-hand
side of the first inequality by their individual upper bounds and replaced
each of $\left\vert e_{1,m}\left(  t\right)  \right\vert $, $\left\vert
e_{1,0,m}\left(  t,a\right)  \right\vert $ and $\left\vert e_{1,0,m}\left(
t,b\right)  \right\vert $ by their individual upper bounds in the second
inequality, all variance upper bounds and deviation upper bounds we have
derived for $\left\vert e_{m}\left(  t\right)  \right\vert $ the setting
$\Theta_{0}=\left(  a,b\right)  $ remain valid for $\left\vert e_{m}\left(
t\right)  \right\vert $ in the setting $\Theta_{0}=\left[  a,b\right]  $ for
the construction (\ref{ExtensionConstructionCloseNull}).

In addition, when $\Theta_{0}=\left(  a,b\right)  $, we have%
\[
\varphi_{m}\left(  t,\boldsymbol{\mu}\right)  =m^{-1}\sum_{i=1}^{m}\left[
\psi_{1}\left(  t,\mu_{i}\right)  -2^{-1}\left\{  \phi\left(  a\right)
\psi_{1,0}\left(  t,\mu;a\right)  +\phi\left(  b\right)  {\tilde{\psi}}%
_{1,0}\left(  t,\mu;b\right)  \right\}  \right]  ,
\]
$\varphi_{m}\left(  t,\boldsymbol{\mu}\right)  =\sum_{j=1}^{5}d_{\phi
,j}\left(  t,\boldsymbol{\mu}\right)  $ and%
\begin{equation}
\left\vert \check{\pi}_{0,m}^{-1}\varphi_{m}\left(  t,\boldsymbol{\mu}\right)
-1\right\vert \leq\left\vert \check{\pi}_{0,m}^{-1}d_{\phi,1}\left(
t,\boldsymbol{\mu}\right)  -1\right\vert +\sum_{j=2}^{5}\check{\pi}_{0,m}%
^{-1}\left\vert d_{\phi,j}\left(  t,\boldsymbol{\mu}\right)  \right\vert ,
\label{OracleExtensionInequ}%
\end{equation}
where%
\[
\left\{
\begin{array}
[c]{l}%
d_{\phi,1}\left(  t,\boldsymbol{\mu}\right)  =m^{-1}\sum\nolimits_{\left\{
i:\mu_{i}\in\left(  a,b\right)  \right\}  }\psi_{1}\left(  t,\mu_{i}\right) \\
d_{\phi,2}\left(  t,\boldsymbol{\mu}\right)  =m^{-1}\sum\nolimits_{\left\{
i:\mu_{i}=a\right\}  }\left(  \psi_{1}\left(  t,\mu_{i}\right)  -2^{-1}%
\phi\left(  a\right)  {\tilde{\psi}}_{1,0}\left(  t,\mu_{i};a\right)  \right)
\\
d_{\phi,3}\left(  t,\boldsymbol{\mu}\right)  =m^{-1}\sum\nolimits_{\left\{
i:\mu_{i}=b\right\}  }\left(  \psi_{1}\left(  t,\mu_{i}\right)  -2^{-1}%
\phi\left(  b\right)  {\tilde{\psi}}_{1,0}\left(  t,\mu_{i};b\right)  \right)
\\
d_{\phi,4}\left(  t,\boldsymbol{\mu}\right)  =-m^{-1}\left(  \sum_{\left\{
i:\mu_{i}\neq a\right\}  }+\sum_{\left\{  i:\mu_{i}\neq b\right\}  }\right)
\left\{  2^{-1}\left[  \phi\left(  a\right)  {\tilde{\psi}}_{1,0}\left(
t,\mu_{i};a\right)  +\phi\left(  b\right)  {\tilde{\psi}}_{1,0}\left(
t,\mu_{i};b\right)  \right]  \right\} \\
d_{\phi,5}\left(  t,\boldsymbol{\mu}\right)  =m^{-1}\sum_{\left\{  i:\mu
_{i}<b\right\}  }\psi_{1}\left(  t,\mu_{i}\right)  +m^{-1}\sum_{\left\{
i:\mu_{i}>b\right\}  }\psi_{1}\left(  t,\mu_{i}\right)
\end{array}
\right.  .
\]
Further, when $\Theta_{0}=\left(  a,b\right)  $, to upper bound $\left\vert
\check{\pi}_{0,m}^{-1}\varphi_{m}\left(  t,\boldsymbol{\mu}\right)
-1\right\vert $, we have replaced each $\left\vert d_{\phi,j}\left(
t,\boldsymbol{\mu}\right)  \right\vert ,2\leq j\leq5$ by its upper bound
$\hat{d}_{\phi,j}\left(  t,\boldsymbol{\mu}\right)  ,2\leq j\leq5$ and
replaced $\left\vert \check{\pi}_{0,m}^{-1}d_{\phi,1}\left(  t,\boldsymbol{\mu
}\right)  -1\right\vert $ by its upper bound $\hat{d}_{\phi,0}\left(
t,\boldsymbol{\mu}\right)  $ directly in (\ref{OracleExtensionInequ}).

When $\Theta_{0}=\left[  a,b\right]  $, we have%
\[
\varphi_{m}\left(  t,\boldsymbol{\mu}\right)  =m^{-1}\sum_{i=1}^{m}\left[
\psi_{1}\left(  t,\mu_{i}\right)  +2^{-1}\left\{  \phi\left(  a\right)
\psi_{1,0}\left(  t,\mu;a\right)  +\phi\left(  b\right)  {\tilde{\psi}}%
_{1,0}\left(  t,\mu;b\right)  \right\}  \right]
\]
and%
\[
\varphi_{m}\left(  t,\boldsymbol{\mu}\right)  =d_{\phi,1}\left(
t,\boldsymbol{\mu}\right)  +d_{\phi,2}^{\ast}\left(  t,\boldsymbol{\mu
}\right)  +d_{\phi,3}^{\ast}\left(  t,\boldsymbol{\mu}\right)  -d_{\phi
,4}\left(  t,\boldsymbol{\mu}\right)  +d_{\phi,5}\left(  t,\boldsymbol{\mu
}\right)
\]
where%
\[
\left\{
\begin{array}
[c]{l}%
d_{\phi,2}^{\ast}\left(  t,\boldsymbol{\mu}\right)  =m^{-1}\sum
\nolimits_{\left\{  i:\mu_{i}=a\right\}  }\left(  \psi_{1}\left(  t,\mu
_{i}\right)  +2^{-1}\phi\left(  a\right)  {\tilde{\psi}}_{1,0}\left(
t,\mu_{i};a\right)  \right) \\
d_{\phi,3}^{\ast}\left(  t,\boldsymbol{\mu}\right)  =m^{-1}\sum
\nolimits_{\left\{  i:\mu_{i}=b\right\}  }\left(  \psi_{1}\left(  t,\mu
_{i}\right)  +2^{-1}\phi\left(  b\right)  {\tilde{\psi}}_{1,0}\left(
t,\mu_{i};b\right)  \right)
\end{array}
\right.  .
\]
Then%
\begin{align*}
\left\vert \check{\pi}_{0,m}^{-1}\varphi_{m}\left(  t,\boldsymbol{\mu}\right)
-1\right\vert  &  \leq\left\vert \check{\pi}_{0,m}^{-1}\left[  d_{\phi
,1}\left(  t,\boldsymbol{\mu}\right)  +d_{\phi,2}^{\ast}\left(
t,\boldsymbol{\mu}\right)  +d_{\phi,3}^{\ast}\left(  t,\boldsymbol{\mu
}\right)  \right]  -1\right\vert \\
&  +\check{\pi}_{0,m}^{-1}\left\vert d_{\phi,4}\left(  t,\boldsymbol{\mu
}\right)  \right\vert +\check{\pi}_{0,m}^{-1}\left\vert d_{\phi,5}\left(
t,\boldsymbol{\mu}\right)  \right\vert .
\end{align*}
However,%
\[
\left\vert \check{\pi}_{0,m}^{-1}\left[  d_{\phi,1}\left(  t,\boldsymbol{\mu
}\right)  +d_{\phi,2}^{\ast}\left(  t,\boldsymbol{\mu}\right)  +d_{\phi
,3}^{\ast}\left(  t,\boldsymbol{\mu}\right)  \right]  -1\right\vert \leq
\hat{d}_{\phi,0}\left(  t,\boldsymbol{\mu}\right)  +\check{\pi}_{0,m}^{-1}%
\hat{d}_{\phi,2}\left(  t,\boldsymbol{\mu}\right)  +\check{\pi}_{0,m}^{-1}%
\hat{d}_{\phi,3}\left(  t,\boldsymbol{\mu}\right)  .
\]
Therefore, the upper bound we have derived for $\left\vert \check{\pi}%
_{0,m}^{-1}\varphi_{m}\left(  t,\boldsymbol{\mu}\right)  -1\right\vert $ when
$\Theta_{0}=\left(  a,b\right)  $ is also an upper bound for $\left\vert
\check{\pi}_{0,m}^{-1}\varphi_{m}\left(  t,\boldsymbol{\mu}\right)
-1\right\vert $ when $\Theta_{0}=\left[  a,b\right]  $.

In summary, all results we have derived for the construction
(\ref{ExtOpenNull}) when $\Theta_{0}=\left(  a,b\right)  $ are valid for the
construction (\ref{ExtensionConstructionCloseNull}) when $\Theta_{0}=\left[
a,b\right]  $.

\subsection{Implementation of estimators for Gaussian family}

\label{Implementation}

For a location-shift family, it suffices to consider the point null
$\Theta_{0}=\left\{  0\right\}  $, bounded null $\Theta_{0}=\left(
a,b\right)  $ and one-sided null $\Theta_{0}=\left(  -\infty,0\right)  $. For
the Gaussian family with variance $\sigma^{2}>0$, we can always standard it to
have variance $1$. So, we only need to consider the setting where $\sigma
^{2}=1$. In this setting, we have the modulus of the CF's as $r_{0}\left(
t\right)  =\exp\left(  -2^{-1}t^{2}\right)  $ (independent of the mean
parameter $\mu$) and
\[
\frac{1}{y}\frac{d}{ds}\frac{1}{r_{0}\left(  ys\right)  }=\frac{1}{y}%
sy^{2}\exp\left(  2^{-1}y^{2}s^{2}\right)  =sy\exp\left(  2^{-1}y^{2}%
s^{2}\right)  .
\]
Let us use the triangular density $\omega\left(  s\right)  =\left(
1-|s|\right)  1_{\left[  -1,1\right]  }\left(  s\right)  $ for the
constructions. Recall $\mathbf{z}=\left(  z_{1},\ldots,z_{m}\right)
^{\intercal}$ and
\[
K_{1,0}\left(  t,x;\mu^{\prime}\right)  ={\int_{\left[  -1,1\right]  }}%
\dfrac{\omega\left(  s\right)  \cos\left\{  ts\left(  x-\mu^{\prime}\right)
\right\}  }{r_{\mu^{\prime}}\left(  ts\right)  }ds.
\]
Then we have the following:

\begin{itemize}
\item For the simple null $\Theta_{0}=\left\{  0\right\}  $, we have
\[
K_{1,0}\left(  t,x;0\right)  =2{\int_{0}^{1}}\left(  1-s\right)  \cos\left(
tsx\right)  \exp\left(  2^{-1}t^{2}s^{2}\right)  ds,
\]
and the estimator of $\pi_{1,m}=m^{-1}\sum_{i=1}^{m}1_{U\setminus\left\{
0\right\}  }\left(  \mu_{i}\right)  $ is%
\begin{align*}
\hat{\varphi}_{1,0,m}\left(  t,\mathbf{z}\right)   &  =m^{-1}\sum_{i=1}%
^{m}\left(  1-K_{1,0}\left(  t,z_{i};0\right)  \right) \\
&  =1-2m^{-1}\sum_{i=1}^{m}{\int_{0}^{1}}\left(  1-s\right)  \cos\left(
tsz_{i}\right)  \exp\left(  2^{-1}t^{2}s^{2}\right)  ds.
\end{align*}
This was provided by \cite{Chen:2018a}.

\item For the bounded null $\Theta_{0}=\left(  a,b\right)  $, we have%
\begin{align*}
K_{1}\left(  t,x\right)   &  =\frac{t}{2\pi}\int_{a}^{b}dy\int_{\left[
-1,1\right]  }\frac{\cos\left(  ts\left(  x-y\right)  \right)  }{r_{0}\left(
ts\right)  }ds\\
&  =\frac{t}{\pi}\int_{a}^{b}dy\int_{0}^{1}\cos\left(  ts\left(  x-y\right)
\right)  \exp\left(  2^{-1}t^{2}s^{2}\right)  ds,
\end{align*}
and%
\[
K_{1,0}\left(  t,x;a\right)  =2{\int_{0}^{1}}\left(  1-s\right)  \cos\left(
ts\left(  x-a\right)  \right)  \exp\left(  2^{-1}t^{2}s^{2}\right)  ds
\]
and%
\[
K_{1,0}\left(  t,x;b\right)  =2{\int_{0}^{1}}\left(  1-s\right)  \cos\left(
ts\left(  x-b\right)  \right)  \exp\left(  2^{-1}t^{2}s^{2}\right)  ds.
\]
Further, the estimator of $\pi_{1,m}=m^{-1}\sum_{i=1}^{m}1_{U\setminus\left(
a,b\right)  }\left(  \mu_{i}\right)  $ is%
\[
\hat{\varphi}_{m}\left(  t,\mathbf{z}\right)  =m^{-1}\sum_{i=1}^{m}\left[
1-K\left(  t,z_{i}\right)  \right]  ,
\]
where%
\[
K\left(  t,x\right)  =K_{1}\left(  t,x\right)  -2^{-1}\left[  K_{1,0}\left(
t,x;a\right)  +K_{1,0}\left(  t,x;b\right)  \right]  .
\]
Namely,%
\begin{align*}
\hat{\varphi}_{m}\left(  t,\mathbf{z}\right)   &  =1-m^{-1}\sum_{i=1}^{m}%
\frac{t}{\pi}\int_{a}^{b}dy\int_{0}^{1}\cos\left(  ts\left(  z_{i}-y\right)
\right)  \exp\left(  2^{-1}t^{2}s^{2}\right)  ds\\
&  \text{ \ }+m^{-1}\sum_{i=1}^{m}{\int_{0}^{1}}\left(  1-s\right)
\cos\left(  ts\left(  z_{i}-a\right)  \right)  \exp\left(  2^{-1}t^{2}%
s^{2}\right)  ds\\
&  \text{ \ }+m^{-1}\sum_{i=1}^{m}{\int_{0}^{1}}\left(  1-s\right)
\cos\left(  ts\left(  z_{i}-b\right)  \right)  \exp\left(  2^{-1}t^{2}%
s^{2}\right)  ds.
\end{align*}

\item In constrast, for the bounded null $\Theta_{0}=\left[  a,b\right]  $,
the estimator of $\pi_{1,m}=m^{-1}\sum_{i=1}^{m}1_{U\setminus\left[
a,b\right]  }\left(  \mu_{i}\right)  $ is%
\[
\hat{\varphi}_{m}\left(  t,\mathbf{z}\right)  =m^{-1}\sum_{i=1}^{m}\left[
1-K\left(  t,z_{i}\right)  \right]  ,
\]
where%
\[
K\left(  t,x\right)  =K_{1}\left(  t,x\right)  +2^{-1}\left[  K_{1,0}\left(
t,x;a\right)  +K_{1,0}\left(  t,x;b\right)  \right]  .
\]
Namely,%
\begin{align*}
\hat{\varphi}_{m}\left(  t,\mathbf{z}\right)   &  =1-m^{-1}\sum_{i=1}^{m}%
\frac{t}{\pi}\int_{a}^{b}dy\int_{0}^{1}\cos\left(  ts\left(  z_{i}-y\right)
\right)  \exp\left(  2^{-1}t^{2}s^{2}\right)  ds\\
&  \text{ \ }-m^{-1}\sum_{i=1}^{m}{\int_{0}^{1}}\left(  1-s\right)
\cos\left(  ts\left(  z_{i}-a\right)  \right)  \exp\left(  2^{-1}t^{2}%
s^{2}\right)  ds\\
&  \text{ \ }-m^{-1}\sum_{i=1}^{m}{\int_{0}^{1}}\left(  1-s\right)
\cos\left(  ts\left(  z_{i}-b\right)  \right)  \exp\left(  2^{-1}t^{2}%
s^{2}\right)  ds.
\end{align*}

\item For the one-sided null $\Theta_{0}=\left(  -\infty,0\right)  $, we have%
\begin{align*}
K_{1}\left(  t,x\right)   &  =\frac{1}{2\pi}\int_{0}^{1}dy\int_{-1}^{1}\left[
\frac{\sin\left(  ytsx\right)  }{y}\left\{  \frac{d}{ds}\frac{1}{r_{0}\left(
tys\right)  }\right\}  +\frac{tx\cos\left(  tysx\right)  }{r_{0}\left(
tys\right)  }\right]  ds\\
&  =\frac{1}{\pi}\int_{0}^{1}dy\int_{0}^{1}syt^{2}\sin\left(  ytsx\right)
\exp\left(  2^{-1}y^{2}t^{2}s^{2}\right)  ds\\
&  \text{ \ }+\frac{1}{\pi}\int_{0}^{1}dy\int_{0}^{1}tx\cos\left(
tysx\right)  \exp\left(  2^{-1}t^{2}y^{2}s^{2}\right)  ds
\end{align*}
and%
\[
K_{1,0}\left(  t,x;0\right)  =2{\int_{0}^{1}}\left(  1-s\right)  \cos\left(
tsx\right)  \exp\left(  2^{-1}t^{2}s^{2}\right)  ds.
\]
Further, the estimator of $\pi_{1,m}=m^{-1}\sum_{i=1}^{m}1_{U\setminus\left(
-\infty,0\right)  }\left(  \mu_{i}\right)  $ is%
\[
\hat{\varphi}_{m}\left(  t,\mathbf{z}\right)  =m^{-1}\sum_{i=1}^{m}\left[
1-K\left(  t,z_{i}\right)  \right]  ,
\]
where%
\[
K\left(  t,x\right)  =2^{-1}-K_{1}\left(  t,x\right)  -2^{-1}K_{1,0}\left(
t,x;0\right)  .
\]
Namely,%
\begin{align*}
\hat{\varphi}_{m}\left(  t,\mathbf{z}\right)   &  =2^{-1}+\frac{1}{\pi}%
m^{-1}\sum_{i=1}^{m}\int_{0}^{1}dy\int_{0}^{1}syt^{2}\sin\left(
ytsz_{i}\right)  \exp\left(  2^{-1}y^{2}t^{2}s^{2}\right)  ds\\
&  \text{ \ }+\frac{1}{\pi}m^{-1}\sum_{i=1}^{m}\int_{0}^{1}dy\int_{0}%
^{1}tx\cos\left(  tysz_{i}\right)  \exp\left(  2^{-1}t^{2}y^{2}s^{2}\right)
ds\\
&  \text{ \ }+m^{-1}\sum_{i=1}^{m}{\int_{0}^{1}}\left(  1-s\right)
\cos\left(  tsz_{i}\right)  \exp\left(  2^{-1}t^{2}s^{2}\right)  ds.
\end{align*}

\item In contrast, for the one-sided null $\Theta_{0}=(-\infty,0]$, the
estimator of $\pi_{1,m}=m^{-1}\sum_{i=1}^{m}1_{U\setminus(-\infty,0]}\left(
\mu_{i}\right)  $ is%
\[
\hat{\varphi}_{m}\left(  t,\mathbf{z}\right)  =m^{-1}\sum_{i=1}^{m}\left[
1-K\left(  t,z_{i}\right)  \right]  ,
\]
where%
\[
K\left(  t,x\right)  =2^{-1}-K_{1}\left(  t,x\right)  +2^{-1}K_{1,0}\left(
t,x;0\right)  .
\]
Namely,%
\begin{align*}
\hat{\varphi}_{m}\left(  t,\mathbf{z}\right)   &  =2^{-1}+\frac{1}{\pi}%
m^{-1}\sum_{i=1}^{m}\int_{0}^{1}dy\int_{0}^{1}syt^{2}\sin\left(
ytsz_{i}\right)  \exp\left(  2^{-1}y^{2}t^{2}s^{2}\right)  ds\\
&  \text{ \ }+\frac{1}{\pi}m^{-1}\sum_{i=1}^{m}\int_{0}^{1}dy\int_{0}%
^{1}tx\cos\left(  tysz_{i}\right)  \exp\left(  2^{-1}t^{2}y^{2}s^{2}\right)
ds\\
&  \text{ \ }-m^{-1}\sum_{i=1}^{m}{\int_{0}^{1}}\left(  1-s\right)
\cos\left(  tsz_{i}\right)  \exp\left(  2^{-1}t^{2}s^{2}\right)  ds.
\end{align*}

\item For the extension with respect to the bound null $\Theta_{0}=\left(
a,b\right)  $, the estimator of%
\[
\check{\pi}_{0,m}=m^{-1}\sum\nolimits_{\left\{  i\in\left\{  1,\ldots
,m\right\}  :\mu_{i}\in\left(  a,b\right)  \right\}  }\phi\left(  \mu
_{i}\right)
\]
is $\hat{\varphi}_{m}\left(  t,\mathbf{z}\right)  =m^{-1}\sum_{i=1}%
^{m}K\left(  t,z_{i}\right)  $ with%
\[
K\left(  t,x\right)  =K_{1}\left(  t,x\right)  -2^{-1}\left[  \phi\left(
a\right)  K_{1,0}\left(  t,x;a\right)  +\phi\left(  b\right)  K_{1,0}\left(
t,x;b\right)  \right]  .
\]
Namely,%
\begin{align*}
\hat{\varphi}_{m}\left(  t,\mathbf{z}\right)   &  =m^{-1}\sum_{i=1}^{m}%
\frac{t}{\pi}\int_{a}^{b}dy\int_{0}^{1}\cos\left(  ts\left(  z_{i}-y\right)
\right)  \exp\left(  2^{-1}t^{2}s^{2}\right)  ds\\
&  \text{ \ }-m^{-1}\phi\left(  a\right)  \sum_{i=1}^{m}{\int_{0}^{1}}\left(
1-s\right)  \cos\left(  ts\left(  z_{i}-a\right)  \right)  \exp\left(
2^{-1}t^{2}s^{2}\right)  ds\\
&  \text{ \ }-m^{-1}\phi\left(  b\right)  \sum_{i=1}^{m}{\int_{0}^{1}}\left(
1-s\right)  \cos\left(  ts\left(  z_{i}-b\right)  \right)  \exp\left(
2^{-1}t^{2}s^{2}\right)  ds.
\end{align*}

\item In contrast, for the extension with respect to the bound null
$\Theta_{0}=\left[  a,b\right]  $, the estimator of%
\[
\tilde{\pi}_{0,m}=m^{-1}\sum\nolimits_{\left\{  i\in\left\{  1,\ldots
,m\right\}  :\mu_{i}\in\left[  a,b\right]  \right\}  }\phi\left(  \mu
_{i}\right)
\]
is $\hat{\varphi}_{m}\left(  t,\mathbf{z}\right)  =m^{-1}\sum_{i=1}%
^{m}K\left(  t,z_{i}\right)  $ with%
\[
K\left(  t,x\right)  =K_{1}\left(  t,x\right)  +2^{-1}\left[  \phi\left(
a\right)  K_{1,0}\left(  t,x;a\right)  +\phi\left(  b\right)  K_{1,0}\left(
t,x;b\right)  \right]  .
\]
Namely,%
\begin{align*}
\hat{\varphi}_{m}\left(  t,\mathbf{z}\right)   &  =m^{-1}\sum_{i=1}^{m}%
\frac{t}{\pi}\int_{a}^{b}dy\int_{0}^{1}\cos\left(  ts\left(  z_{i}-y\right)
\right)  \exp\left(  2^{-1}t^{2}s^{2}\right)  ds\\
&  \text{ \ }+m^{-1}\phi\left(  a\right)  \sum_{i=1}^{m}{\int_{0}^{1}}\left(
1-s\right)  \cos\left(  ts\left(  z_{i}-a\right)  \right)  \exp\left(
2^{-1}t^{2}s^{2}\right)  ds\\
&  \text{ \ }+m^{-1}\phi\left(  b\right)  \sum_{i=1}^{m}{\int_{0}^{1}}\left(
1-s\right)  \cos\left(  ts\left(  z_{i}-b\right)  \right)  \exp\left(
2^{-1}t^{2}s^{2}\right)  ds.
\end{align*}

\end{itemize}

Here each integral is approximated by a Riemann sum. The iterated integrals
appear above satisfy Fubini's theorem for computing double integrals as
iterated intergrals. So, the order of computing a Riemann sum and then the
average with respect to $m$ can be interchanged in order to save a bit of
computational time.

\section*{Acknowledgements}

This research was funded by the New Faculty Seed Grant provided by Washington
State University. I am grateful to Prof. G\'{e}rard Letac for his constant
guidance and feedback on an earlier version of the manuscript and to Prof.
Jiashun Jin for his warm encouragements.
I am also very grateful to the anonymous reviewers for their very detailed, constructive
and critical comments that helped greatly improve the quality and presentation
of the manuscript.

\bibliographystyle{dcu}
\bibliography{estpidep}


\begin{table}[H]
\centering
\begin{tabular}{rllrr}
  \hline
  \hline
 $m$ & Sparsity & $\hat{E}(\tilde{\delta}_m)$  & $\hat{\sigma}(\tilde{\delta}_m)$ \\
  \hline
  \hline

    1000 & $\pi_{1,m}=0.2$ & -0.1421 & 0.2227 \\
   5000 & $\pi_{1,m}=0.2$ & -0.1557 & 0.1261 \\
   10000 & $\pi_{1,m}=0.2$ & -0.1465 & 0.0859 \\
   50000 & $\pi_{1,m}=0.2$ & -0.1473 & 0.0437 \\
   100000 & $\pi_{1,m}=0.2$ & -0.1451 & 0.0560 \\
   500000 & $\pi_{1,m}=0.2$ & -0.1404 & 0.0778 \\
   \hline
   1000 & $\pi_{1,m}=1/\ln{(\ln{m})}$ & -0.1557 & 0.3142 \\
   5000 & $\pi_{1,m}=1/\ln{(\ln{m})}$ & -0.1594 & 0.1653 \\
   10000 & $\pi_{1,m}=1/\ln{(\ln{m})}$ & -0.1373 & 0.1077 \\
   50000 & $\pi_{1,m}=1/\ln{(\ln{m})}$ & -0.1522 & 0.0609 \\
   100000 & $\pi_{1,m}=1/\ln{(\ln{m})}$ & -0.1516 & 0.0729 \\
   500000 & $\pi_{1,m}=1/\ln{(\ln{m})}$ & -0.1490 & 0.1053 \\

   \hline
   \hline
\end{tabular}
\label{TabL2Sim}
\caption{In the table, $\tilde{\delta}_m = \hat{\pi}_{1,m}/\pi_{1,m} -1$ (where $\hat{\pi}_{1,m}$ is an estimate of $\pi_{1,m}$),
 $\hat{E} (\tilde{\delta}_m )$ is the sample mean of $\tilde{\delta}_m$, and
$\hat{\sigma} (\tilde{\delta}_m )$ the sample standard deviation of $\tilde{\delta}_m$, where our new estimator $\hat{\pi}_{1,m}$ estimates the average truncated 2-norm
described by Scenario 3 in \autoref{simDesign}. When $\pi_{1,m}=0.2$, our proposed estimators
``New" show a clear trend of convergence to $0$ as $m$ increases. For $\pi_{1,m}=1/\ln{(\ln{m})}$ though, $\tilde{\delta}_m$ for our ``New'' estimators does not
show a clear trend of convergence to $0$ as $m$ increases. However, this is an artifact of the numerical error when implementing our ``New'' estimators,
as explained in \autoref{simDesign}.}
\end{table}

\begin{figure}[H]
\centering
\includegraphics[height=0.8\textheight,width=0.99\textwidth]{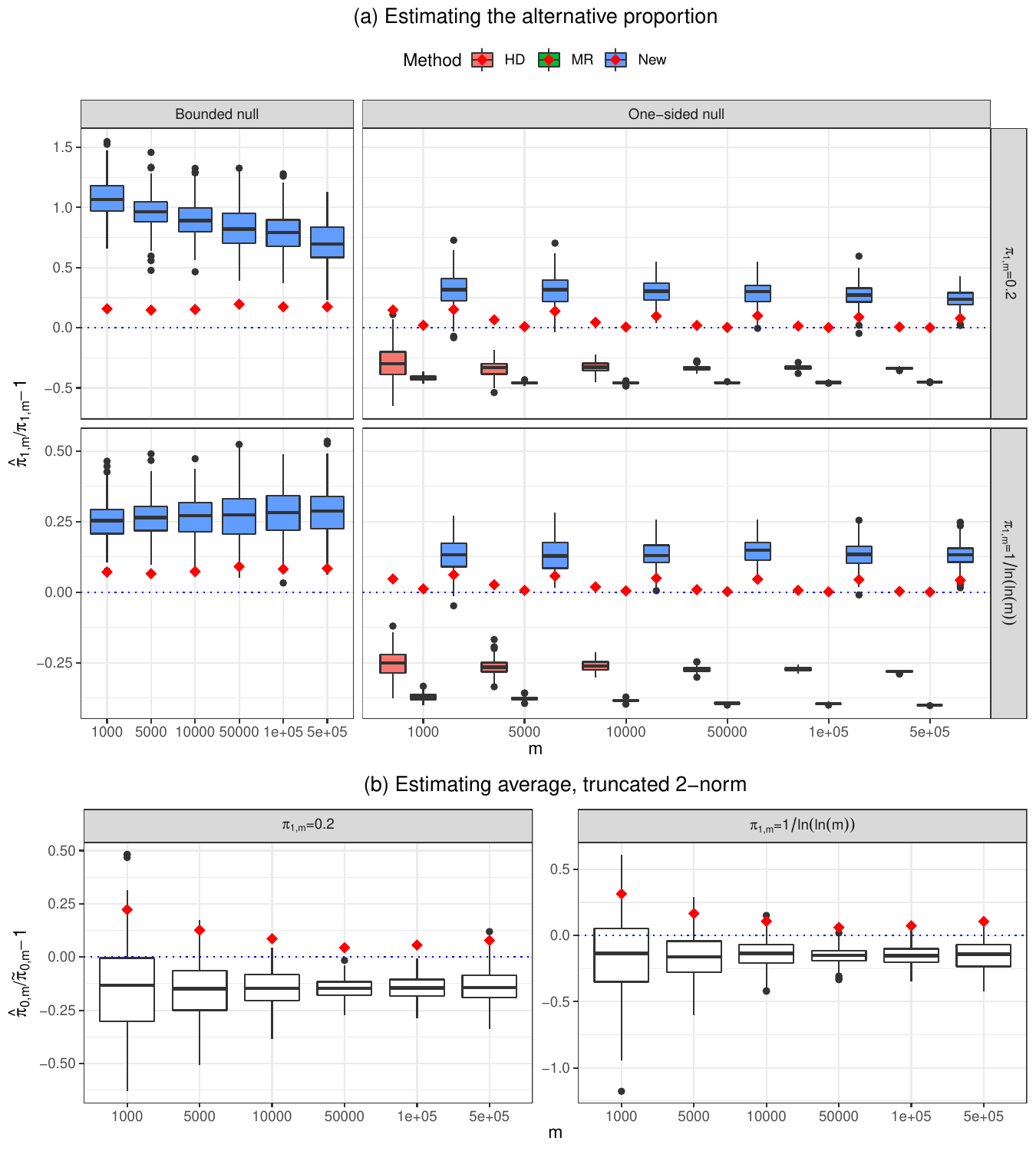}
\vspace{-0.4cm}\caption[Gaussian]{Boxplot of the excess $\tilde{\delta}_{m}$
(on the vertical axis) of an estimator $\hat{\pi}_{1,m}$ of ${\pi}_{1,m}$ in
panel (a)
\textcolor{blue}{as $\tilde{\delta}_{m}=\hat{\pi}_{1,m}\pi_{1,m}^{-1}-1$} (or
an estimator $\hat{\pi}_{0,m}$ of $\tilde{\pi}_{0,m}$ in panel (b)
\textcolor{blue}{as $\tilde{\delta}_{m}=\hat{\pi}_{0,m}\tilde{\pi}_{0,m}^{-1}-1$}).
The thick horizontal line and the diamond in each boxplot are respectively the
mean and standard deviation of $\tilde{\delta}_{m}$, and the dotted horizontal
line is the reference for $\tilde{\delta}_{m}=0$. Panel (a) is for Scenarios 1
and 2, and Panel (b) for Scenario 3, all described in \autoref{simDesign}. All
estimators were applied to Gaussian families, and the proposed estimator is
referred to as ``New''. For for each $m$ in each case of $\pi_{1,m}$ for a
one-side null, there are three boxplots, where the leftmost is for the ``HD''
estimator, the middle for the ``MR'' estimator, and the rightmost for ``New''.
No simulation was done for the ``MR'' or ``HD'' estimator for a bounded null
or estimating average, truncated 2-norm.}%
\label{fig1}%
\end{figure}

\begin{table}[H]
\centering
\begin{tabular}{rllllrr}
  \hline
  Null Type & Sparsity & $m$ & Method & $\hat{E}(\tilde{\delta}_m)$  & $\hat{\sigma}(\tilde{\delta}_m)$ \\
  \hline
 Bounded null & $\pi_{1,m}=0.2$ & 1000 & New & 1.0721 & 0.1571 \\
  Bounded null & $\pi_{1,m}=0.2$ & 5000 & New & 0.9573 & 0.1467 \\
  Bounded null & $\pi_{1,m}=0.2$ & 10000 & New & 0.8963 & 0.1519 \\
  Bounded null & $\pi_{1,m}=0.2$ & 50000 & New & 0.8255 & 0.1959 \\
  Bounded null & $\pi_{1,m}=0.2$ & 100000 & New & 0.7949 & 0.1740 \\
  Bounded null & $\pi_{1,m}=0.2$ & 500000 & New & 0.7051 & 0.1741 \\
  \hline
  Bounded null & $\pi_{1,m}=1/\ln{(\ln{m})}$ & 1000 & New & 0.2545 & 0.0712 \\
  Bounded null & $\pi_{1,m}=1/\ln{(\ln{m})}$ & 5000 & New & 0.2634 & 0.0654 \\
  Bounded null & $\pi_{1,m}=1/\ln{(\ln{m})}$ & 10000 & New & 0.2671 & 0.0730 \\
  Bounded null & $\pi_{1,m}=1/\ln{(\ln{m})}$ & 50000 & New & 0.2691 & 0.0907 \\
  Bounded null & $\pi_{1,m}=1/\ln{(\ln{m})}$ & 100000 & New & 0.2789 & 0.0820 \\
  Bounded null & $\pi_{1,m}=1/\ln{(\ln{m})}$ & 500000 & New & 0.2829 & 0.0839 \\
  \hline
  One-sided null & $\pi_{1,m}=0.2$ & 1000 & HD & -0.2951 & 0.1477 \\
  One-sided null & $\pi_{1,m}=0.2$ & 1000 & MR & -0.4192 & 0.0210 \\
  One-sided null & $\pi_{1,m}=0.2$ & 1000 & New & 0.3118 & 0.1517 \\
  One-sided null & $\pi_{1,m}=0.2$ & 5000 & HD & -0.3393 & 0.0661 \\
  One-sided null & $\pi_{1,m}=0.2$ & 5000 & MR & -0.4588 & 0.0092 \\
  One-sided null & $\pi_{1,m}=0.2$ & 5000 & New & 0.3069 & 0.1375 \\
  One-sided null & $\pi_{1,m}=0.2$ & 10000 & HD & -0.3288 & 0.0450 \\
  One-sided null & $\pi_{1,m}=0.2$ & 10000 & MR & -0.4582 & 0.0065 \\
  One-sided null & $\pi_{1,m}=0.2$ & 10000 & New & 0.3017 & 0.0965 \\
  One-sided null & $\pi_{1,m}=0.2$ & 50000 & HD & -0.3366 & 0.0193 \\
  One-sided null & $\pi_{1,m}=0.2$ & 50000 & MR & -0.4563 & 0.0029 \\
  One-sided null & $\pi_{1,m}=0.2$ & 50000 & New & 0.2878 & 0.1001 \\
  One-sided null & $\pi_{1,m}=0.2$ & 100000 & HD & -0.3317 & 0.0147 \\
  One-sided null & $\pi_{1,m}=0.2$ & 100000 & MR & -0.4539 & 0.0021 \\
  One-sided null & $\pi_{1,m}=0.2$ & 100000 & New & 0.2705 & 0.0884 \\
  One-sided null & $\pi_{1,m}=0.2$ & 500000 & HD & -0.3362 & 0.0066 \\
  One-sided null & $\pi_{1,m}=0.2$ & 500000 & MR & -0.4528 & 0.0009 \\
  One-sided null & $\pi_{1,m}=0.2$ & 500000 & New & 0.2385 & 0.0785 \\
  \hline
  One-sided null & $\pi_{1,m}=1/\ln{(\ln{m})}$ & 1000 & HD & -0.2543 & 0.0469 \\
  One-sided null & $\pi_{1,m}=1/\ln{(\ln{m})}$ & 1000 & MR & -0.3704 & 0.0120 \\
  One-sided null & $\pi_{1,m}=1/\ln{(\ln{m})}$ & 1000 & New & 0.1311 & 0.0617 \\
  One-sided null & $\pi_{1,m}=1/\ln{(\ln{m})}$ & 5000 & HD & -0.2648 & 0.0266 \\
  One-sided null & $\pi_{1,m}=1/\ln{(\ln{m})}$ & 5000 & MR & -0.3756 & 0.0064 \\
  One-sided null & $\pi_{1,m}=1/\ln{(\ln{m})}$ & 5000 & New & 0.1311 & 0.0573 \\
  One-sided null & $\pi_{1,m}=1/\ln{(\ln{m})}$ & 10000 & HD & -0.2605 & 0.0188 \\
  One-sided null & $\pi_{1,m}=1/\ln{(\ln{m})}$ & 10000 & MR & -0.3834 & 0.0048 \\
  One-sided null & $\pi_{1,m}=1/\ln{(\ln{m})}$ & 10000 & New & 0.1325 & 0.0498 \\
  One-sided null & $\pi_{1,m}=1/\ln{(\ln{m})}$ & 50000 & HD & -0.2731 & 0.0091 \\
  One-sided null & $\pi_{1,m}=1/\ln{(\ln{m})}$ & 50000 & MR & -0.3930 & 0.0021 \\
  One-sided null & $\pi_{1,m}=1/\ln{(\ln{m})}$ & 50000 & New & 0.1461 & 0.0459 \\
  One-sided null & $\pi_{1,m}=1/\ln{(\ln{m})}$ & 100000 & HD & -0.2732 & 0.0068 \\
  One-sided null & $\pi_{1,m}=1/\ln{(\ln{m})}$ & 100000 & MR & -0.3947 & 0.0015 \\
  One-sided null & $\pi_{1,m}=1/\ln{(\ln{m})}$ & 100000 & New & 0.1315 & 0.0446 \\
  One-sided null & $\pi_{1,m}=1/\ln{(\ln{m})}$ & 500000 & HD & -0.2807 & 0.0031 \\
  One-sided null & $\pi_{1,m}=1/\ln{(\ln{m})}$ & 500000 & MR & -0.3995 & 0.0007 \\
  One-sided null & $\pi_{1,m}=1/\ln{(\ln{m})}$ & 500000 & New & 0.1320 & 0.0428 \\
  \hline
\end{tabular}
\label{TabOriSim}
\caption{In the table, $\tilde{\delta}_m = \hat{\pi}_{1,m}/\pi_{1,m} -1$ (where $\hat{\pi}_{1,m}$ is an estimate of $\pi_{1,m}$),
 $\hat{E} (\tilde{\delta}_m )$ is the sample mean of $\tilde{\delta}_m$, and
$\hat{\sigma} (\tilde{\delta}_m )$ the sample standard deviation of $\tilde{\delta}_m$. When $\pi_{1,m}=0.2$, our proposed estimators
``New" show a clear trend of convergence to $0$ as $m$ increases. For $\pi_{1,m}=1/\ln{(\ln{m})}$ though, $\tilde{\delta}_m$ for our ``New'' estimators does not
show a clear trend of convergence to $0$ as $m$ increases. However, this is an artifact of the numerical error when implementing our ``New'' estimators,
as explained in \autoref{simDesign}.}
\end{table}

\newpage
\renewcommand*{\thefootnote}{\fnsymbol{footnote}}
\numberwithin{figure}{section}

\begin{center}
{\Large {Supplementary Material to ``Uniformly consistent proportion estimation for composite
hypotheses via integral equations: `the case of location-shift
families'"}\\
{\small Xiongzhi Chen (xiongzhi.chen@wsu.edu)}\\
{\small Department of Mathematics and Statistics, Washington State University}}
\bigskip
\end{center}

\autoref{AppProofA} contains some auxiliary results and their proofs,
\autoref{ProofsLocA} proofs related to Construction I, \autoref{AppProofsLocB} proofs related to Construction
II, \autoref{AppProofsExt} proofs related to
the extension of Construction I, and \autoref{SecConsistencyI} a correction to
the proof of Theorem 3 of \cite{Chen:2018a}.


\section{Some auxiliary results}

\label{AppProofA}

We quote in \autoref{AuxRes3} and \autoref{speedConv} six lemmas from \cite{Xiongzhi2025} and prove in
\autoref{SecConcentrations} results on the concentrations of three empirical
processes that are used in constructing the estimators.

\subsection{Dirichlet integral and Fourier transform}

\label{AuxRes3}

First, we quote the speed of convergence of Dirichlet integral as:

\begin{lemma}
\label{lm:Dirichlet}$\left\vert \int_{0}^{t}x^{-1}\sin xdx-2^{-1}%
\pi\right\vert \leq2\pi t^{-1}$ for $t\geq2$.
\end{lemma}

\noindent\autoref{lm:Dirichlet} implies the following identities (also
referred to as \textquotedblleft Dirichlet integral\textquotedblright), i.e.,
(\ref{EqDirichlet1}) and (\ref{EqDB1}) that have been presented in the main text.

Second, we quote the integral representations of solutions to the
Lebesgue-Stieltjes integral equation as:

\begin{lemma}
\label{LmDirichlet}For any $a,b,\mu,t\in\mathbb{R}$ with $a<b$,
\[
\int_{\left(  \mu-b\right)  t}^{\left(  \mu-a\right)  t}\frac{\sin v}%
{v}dv=\frac{1}{2}\int_{a}^{b}dy\int_{-1}^{1}t\exp\left\{  \iota\left(
\mu-y\right)  ts\right\}  ds.
\]
On the other hand, for any $b,\mu,t\in\mathbb{R}$,
\[
\frac{1}{\pi}\int_{0}^{t}\frac{\sin\left\{  \left(  \mu-b\right)  y\right\}
}{y}dy=\frac{1}{2\pi}\int_{0}^{t}dy\int_{-1}^{1}\left(  \mu-b\right)
\exp\left\{  \iota ys\left(  \mu-b\right)  \right\}  ds.
\]

\end{lemma}

Third, we quote the speed of convergence of the Riemann-Lebesgue lemma:

\begin{lemma}
\label{lm:OracleSpeed}Let $-\infty<a_{1}<b_{1}<\infty$. If $f:\left[
a_{1},b_{1}\right]  \rightarrow\mathbb{R}$ is of bounded variation, then%
\[
\left\vert \int_{\left[  a_{1},b_{1}\right]  }f\left(  s\right)  \cos\left(
ts\right)  ds\right\vert \leq4\left(  \left\Vert f\right\Vert _{\mathrm{TV}%
}+\left\Vert f\right\Vert _{\infty}\right)  \left\vert t\right\vert
^{-1}1_{\left\{  t\neq0\right\}  }\left(  t\right)
\]
and
\[
\left\vert \int_{\left[  a_{1},b_{1}\right]  }f\left(  s\right)  \sin\left(
ts\right)  ds\right\vert \leq4\left(  \left\Vert f\right\Vert _{\mathrm{TV}%
}+\left\Vert f\right\Vert _{\infty}\right)  \left\vert t\right\vert
^{-1}1_{\left\{  t\neq0\right\}  }\left(  t\right)  .
\]

\end{lemma}

\subsection{Speeds of convergence of discriminant functions}

\label{speedConv}

In order to quote and state speeds of convergence of discriminant functions, we will deliberately introduce definitions in the following 3 items that are independent of each other but are compatible with notations throughout the whole manuscript.

\begin{itemize}
\item First, we quote the speed of convergence of discriminant function for the bounded null. Recall $\psi_{1,0}\left(  t,\mu;\mu^{\prime}\right)  ={\int_{\left[  -1,1\right]  }%
}\omega\left(  s\right)  \cos\left\{  ts\left(  \mu-\mu^{\prime}\right)
\right\}  ds$ and define
\begin{equation}
{\tilde{\psi}_{1,0}\left(  t,\mu;\mu^{\prime}\right)  ={\int_{\left[  -1,1\right]  }%
}\omega\left(  s\right)  \cos\left\{  ts \sigma^{-1} \left(  \mu-\mu^{\prime}\right)
\right\}  ds}. \label{NewPsiPoint}
\end{equation}
Then setting $\sigma=1$ in $\tilde{\psi}_{1,0}\left(  t,\mu;\mu^{\prime
}\right)  $ implies that $\tilde{\psi}_{1,0}\left(  t,\mu;\mu^{\prime}\right)
=\psi_{1,0}\left(  t,\mu;\mu^{\prime}\right)  $.
Define
\[
\psi_{1}\left(  t,\mu\right)  =\frac{1}{\pi}%
\int_{\left(  \mu-b\right)  t}^{\left(  \mu-a\right)  t}\frac{\sin y}{y}dy
\]
and
\[
\varphi_{1,m}\left(  t,\boldsymbol{\mu}\right)  =m^{-1}\sum_{i=1}^{m}\psi
_{1}\left(  t,\mu_{i}\right)  \text{ \ and \ }\varphi_{1,0,m}\left(
t,\boldsymbol{\mu};\mu^{\prime}\right)  =m^{-1}\sum_{i=1}^{m}{\tilde{\psi}_{1,0}\left(
t,\mu_{i};\mu^{\prime}\right)}.
\]
We have:
\begin{lemma}
\label{SpeedOracleBoundedNull}
  Set $u_{m}=\min_{\tau\in\left\{  a,b\right\}  }\min_{\left\{
j:\mu_{j}\neq\tau\right\}  }\left\vert \mu_{j}-\tau\right\vert $ and%
\[
\varphi_{m}\left(  t,\boldsymbol{\mu}\right)  =1-\varphi_{1,m}\left(  t,\boldsymbol{\mu
}\right)  +2^{-1}\varphi_{1,0,m}\left(  t,\boldsymbol{\mu};a\right)
+2^{-1}\varphi_{1,0,m}\left(  t,\boldsymbol{\mu};b\right)  .
\]
Then, for positive $t$ such that $t\left(  b-a\right)  \geq2$
and $t u_{m}\geq2$,
\begin{equation}
\left\vert \pi_{1,m}^{-1} \varphi_{1,m}\left(  t,\boldsymbol{\mu}\right) -1\right\vert
\leq\frac{4{\sigma}\left(  \left\Vert \omega\right\Vert _{\mathrm{TV}}+\left\Vert
\omega\right\Vert _{\infty}\right)  +12}{\pi_{1,m}u_{m}t}+\frac{4}{\left(
b-a\right)  t\pi_{1,m}}. \label{eq10cx}%
\end{equation}
\end{lemma}

\item Second, we quote the speed of convergence of discriminant function for the one-side null.
Recall%
\[
{\tilde{\psi}_{1,0}\left(  t,\mu;\mu^{\prime}\right)  ={\int_{\left[  -1,1\right]  }%
}\omega\left(  s\right)  \cos\left\{  ts \sigma^{-1} \left(  \mu-\mu^{\prime}\right)
\right\}  ds}.
\]
Define
\[
\psi_{1}\left(  t,\mu\right)  =\frac{1}{\pi}\int%
_{0}^{t}\frac{\sin\left(  \left(  \mu-b\right)  y\right)  }{y}dy
\]
and%
\[
\varphi_{1,m}\left(  t,\boldsymbol{\mu}\right)  =m^{-1}\sum_{i=1}^{m}\psi
_{1}\left(  t,\mu_{i}\right)  \text{ \ and \ }\varphi_{1,0,m}\left(
t,\boldsymbol{\mu};b\right)  =m^{-1}\sum_{i=1}^{m}{\tilde{\psi}}_{1,0}\left(  t,\mu
_{i};b\right)  .
\]
Then we have:

\begin{lemma}
\label{SpeedOneSidedNull}

 Define $\tilde{u}_{m}=\min_{\left\{  j:\mu_{j}\neq b\right\}
}\left\vert \mu_{j}-b\right\vert $ and
\[
\varphi_{m}\left(  t,\boldsymbol{\mu}\right)  =2^{-1}+\varphi_{1,m}\left(
t,\boldsymbol{\mu}\right)  +2^{-1}\varphi_{1,0,m}\left(  t,\boldsymbol{\mu
};{b}\right)  .
\]
Then, when $t\tilde{u}_{m}\geq2$,%
\begin{equation}
\left\vert \pi_{1,m}^{-1}\varphi_{m}\left(  t,\boldsymbol{\mu}\right)
-1\right\vert \leq\frac{4+2{\sigma}\left(  \left\Vert \omega\right\Vert _{\infty
}+\left\Vert \omega\right\Vert _{\mathrm{TV}}\right)  }{t\tilde{u}_{m}%
\pi_{1,m}}.\label{eqdx4}%
\end{equation}
\end{lemma}

\item Third, we quote the speed of convergence of discriminant function for the extension.
Define%
\[
\psi_{1}\left(  t,\mu\right)  =\mathcal{D}_{\phi}\left(  t,\mu;a,b\right)
=\frac{1}{\pi}\int_{a}^{b}\frac{\sin\left\{  \left(  \mu-y\right)  t\right\}
}{\mu-y}\phi\left(  y\right)  dy
\]
and recall%
\[
{\tilde{\psi}_{1,0}\left(  t,\mu;\mu^{\prime}\right)  ={\int_{\left[  -1,1\right]  }%
}\omega\left(  s\right)  \cos\left\{  ts \sigma^{-1} \left(  \mu-\mu^{\prime}\right)
\right\}  ds}
\]
and
$u_{m}=\min_{\tau\in\left\{  a,b\right\}  }\min_{\left\{  j:\mu_{j}%
\neq\tau\right\}  }\left\vert \mu_{j}-\tau\right\vert $ and $\left\Vert
\phi\right\Vert _{1,\infty}=\sup_{\mu\in \textcolor{black}{U}  }C_{\mu}\left(
\phi\right)  $.

Then we have:
\begin{lemma}
\label{SpeedOracleExt}
Define%
\[
\varphi_{m}\left(  t,\boldsymbol{\mu}\right)  =m^{-1}\sum_{i=1}^{m}\left[
\psi_{1}\left(  t,\mu_{i}\right)  -2^{-1}\left\{  \phi\left(  a\right)
{\tilde{\psi}}_{1,0}\left(  t,\mu;a\right)  +\phi\left(  b\right)  {\tilde{\psi}}_{1,0}\left(
t,\mu;b\right)  \right\}  \right]  .
\]
Then, for $t$ such that $tu_{m}\geq2$ and $t\left(  b-a\right)  \geq2$,%
\begin{equation}
\left\vert \check{\pi}_{0,m}^{-1}\varphi_{m}\left(  t,\boldsymbol{\mu}\right)
-1\right\vert \leq\frac{\textcolor{black}{20}\left\Vert \phi\right\Vert _{1,\infty}}{\pi
t\check{\pi}_{0,m}}+\frac{\textcolor{black}{20}\left\Vert \phi\right\Vert _{\infty}}{tu_{m}%
\check{\pi}_{0,m}}+\frac{\textcolor{black}{10}\left\Vert \phi\right\Vert _{\infty}}{t\left(
b-a\right)  \check{\pi}_{0,m}}+\frac{4\left(  \left\Vert \omega\right\Vert
_{\mathrm{TV}}+\left\Vert \omega\right\Vert _{\infty}\right)  \left\Vert
\phi\right\Vert _{\infty}}{tu_{m}\check{\pi}_{0,m} {\sigma^{-1}}}.\label{eqdx13xx}%
\end{equation}
\end{lemma}
\end{itemize}

\subsection{\textcolor{blue}{Concentrations of three empirical processes}}

\label{SecConcentrations}

Recall that for a nonnegative constant $\rho$ and $j\in\mathbb{N}_{+}$,%
\[
\mathcal{B}_{j,m}\left(  \rho\right)  =\left\{  \boldsymbol{\mu}\in
\mathbb{R}^{m}:\frac{1}{m}\sum\nolimits_{i=1}^{m}\left\vert \mu_{i}\right\vert
^{j}\leq\rho\right\}  .
\]

\begin{lemma}
\label{ConcentrationEmpiricalProcesses}Let $a_{0},b_{0},c_{0},c_{1},c_{2}$ and
$c_{3}$ be constants in $\mathbb{R}$ such that $c_{1}>0$ and $c_{2}\leq c_{3}%
$. Define the function%
\[
f\left(  x,v,y,a_{0},b_{0},c_{0}\right)  =\left(  a_{0}x+b_{0}\right)
\cos\left(  v\left(  x-y\right)  -c_{0}\right)
\]
for $\left(  x,y,v\right)  \in\mathbb{R}\times\left[  0,c_{1}\right]
\times\left[  c_{2},c_{3}\right]  $. Let $\left\{  z_{i}\right\}  _{i=1}^{m}$
be random variables whose CDFs are members of a location-shift family, such
that $z_{i}$ has CDF $F_{\mu_{i}}$. Assume $\mathbb{E}\left[  X_{\left(
0\right)  }^{2}\right]  <\infty$ where $X_{\left(  0\right)  }$ has CDF
$F_{0}$. Consider the empirical process%
\[
\mathbb{G}_{m}\left(  v,y\right)  =\frac{1}{m}\sum_{i=1}^{m}\left(  f\left(
z_{i},v,y,a_{0},c_{0}\right)  -\mathbb{E}\left[  f\left(  z_{i},v,y,a_{0}%
,c_{0}\right)  \right]  \right)  .
\]
There are positive constants $\Delta_{1},\Delta_{2},\kappa,\kappa_{1}$ and
$\kappa_{2}$ with $\kappa_{1}>\kappa_{2}$ and a nonnegative sequence $\left\{
\rho_{m}\right\}  _{m\geq1}$, such that the following hold:

\begin{enumerate}
\item Case 1: Assume $a_{0}=0,b_{0}=1$ and $c_{2}<c_{3}$. If $\boldsymbol{\mu
}\in\mathcal{B}_{1,m}\left(  \rho_{m}\right)  $ and%
\[
\left(  \max\left\{  \Delta_{1},\Delta_{2}\right\}  \right)  ^{-1}\kappa
_{2}-\left(  2c_{1}+R_{1}\left(  \rho_{m}\right)  \right)  >0,
\]
then%
\[
\Pr\left\{  \sup\nolimits_{\left(  v,y\right)  \in\left[  0,c_{1}\right]
\times\left[  c_{2},c_{3}\right]  }\left\vert \mathbb{G}_{m}\left(
v,y\right)  \right\vert \geq\kappa_{1}\right\}  \leq\hat{p}_{1,m}\left(
\Xi_{1}\right)  ,
\]
where $\Xi_{1}=\left(  \kappa_{1},\kappa_{2},c_{1},c_{2},c_{3},\Delta
_{1},\Delta_{2},F_{0},\rho_{m}\right)  ,R_{1}\left(  \rho_{m}\right)
=2\mathbb{E}\left[  \left\vert X_{(0)}\right\vert \right]  +2\max\left\{
\left\vert c_{2}\right\vert ,\left\vert c_{3}\right\vert \right\}  +2\rho_{m}$
and
\begin{align}
\hat{p}_{1,m}\left(  \Xi_{1}\right)   &  =2\left(  c_{1}\Delta_{1}%
^{-1}+1\right)  \left[  \left(  c_{3}-c_{2}\right)  \Delta_{2}^{-1}+1\right]
\exp\left(  -2^{-1}m\left(  \kappa_{1}-\kappa_{2}\right)  ^{2}\right)
\nonumber\\
&  \text{ \ \ }+m^{-1}\left[  \left(  \max\left\{  \Delta_{1},\Delta
_{2}\right\}  \right)  ^{-1}\kappa_{2}-\left(  2c_{1}+R_{1}\left(  \rho
_{m}\right)  \right)  \right]  ^{-2}\mathbb{V}\left[  \left\vert X_{\left(
0\right)  }\right\vert \right]  . \label{probExpressCase1}%
\end{align}

\item Case 2: Assume $a_{0}=0,b_{0}=1$ and $c_{2}=c_{3}$. If $\boldsymbol{\mu
}\in\mathcal{B}_{1,m}\left(  \rho_{m}\right)  $ and $\Delta_{1}^{-1}\kappa
_{2}-R_{1}\left(  \rho_{m}\right)  >0$, then
\[
\Pr\left\{  \sup\nolimits_{v\in\left[  0,c_{1}\right]  }\left\vert
\mathbb{G}_{m}\left(  v,c_{2}\right)  \right\vert \geq\kappa_{1}\right\}
\leq\hat{p}_{2,m}\left(  \Xi_{2}\right)  ,
\]
where $\Xi_{2}=\left(  \kappa_{1},\kappa_{2},c_{1},\Delta_{1},F_{0},\rho
_{m}\right)  $ and%
\begin{align}
\hat{p}_{2,m}\left(  \Xi_{2}\right)   &  =2\left(  c_{1}\Delta_{1}%
^{-1}+1\right)  \exp\left(  -2^{-1}m\left(  \kappa_{1}-\kappa_{2}\right)
^{2}\right) \nonumber\\
&  \text{ \ \ }+m^{-1}\left[  \Delta_{1}^{-1}\kappa_{2}-R_{1}\left(  \rho
_{m}\right)  \right]  ^{-2}\mathbb{V}\left[  \left\vert X_{\left(  0\right)
}\right\vert \right]  . \label{probExpressCase2}%
\end{align}

\item Case 3: Assume $a_{0}=1,b_{0}=0$ and $\left[  c_{2},c_{3}\right]
=\left\{  0\right\}  $. If $\boldsymbol{\mu}\in\mathcal{B}_{2,m}\left(
\rho_{m}\right)  $ and $2^{-1}\Delta_{1}^{-1}\kappa_{2}-2^{-1}R_{2}\left(
\rho_{m}\right)  >0$, then%
\[
\Pr\left\{  \sup\nolimits_{v\in\left[  0,c_{1}\right]  }\mathbb{G}_{m}\left(
v,0\right)  \geq\kappa_{1}\right\}  \leq\hat{p}_{3,m}\left(  \Xi_{3}\right)
,
\]
where $\Xi_{3}=\left(  \kappa_{1},\kappa_{2},c_{1},\Delta_{1},F_{0},\rho
_{m},\kappa\right)  $, $R_{2}\left(  \rho_{m}\right)  =4\mathbb{E}\left[
X_{\left(  0\right)  }^{2}\right]  +4\rho_{m}$,%
\begin{align}
\hat{p}_{3,m}\left(  \Xi_{3}\right)   &  =2\left(  c_{1}\Delta_{1}%
^{-1}+1\right)  \exp\left(  -2^{-1}m\kappa^{-2}\left(  \kappa_{1}-\kappa
_{2}\right)  ^{2}\right) \nonumber\\
&  \text{ \ \ \ }+m^{-1}\left[  2^{-1}\Delta_{1}^{-1}\kappa_{2}-2^{-1}%
R_{2}\left(  \rho_{m}\right)  \right]  ^{-2}\mathbb{V}\left[  X_{\left(
0\right)  }^{2}\right]  +1-C_{m,\boldsymbol{\mu},F_{0}}
\label{probExpressCase3}%
\end{align}
and $C_{m,\boldsymbol{\mu},F_{0}}=\left[  1-2\Pr\left\{  X_{\left(  0\right)
}>\kappa-\left\Vert \boldsymbol{\mu}\right\Vert _{\infty}\right\}  \right]
^{m}$.
\end{enumerate}
\end{lemma}

The proof of \autoref{ConcentrationEmpiricalProcesses} is given below.

\subsubsection{Proof of \autoref{ConcentrationEmpiricalProcesses}}

The whole proof has two parts: \textquotedblleft\textbf{Part I}:\ preparations
(including laying out a grid, writing down partial derivatives of various
related quantities, and deriving coarse upper bounds)\textquotedblright\ and
\textquotedblleft\textbf{Part II}: bounding the tail probabilities of maxima
related to three empirical processes\textquotedblright, for which \textbf{Part
II} has three parts, each for one empirical process induced by $\mathbb{G}%
_{m}\left(  v,y\right)  $. To reduce cumbersome notations, let us just write
$f_{i}\left(  v,y\right)  $ for $f\left(  z_{i},v,y,a_{0},c_{0}\right)  $, and
introduce
\[
\hat{G}_{m}\left(  v,y\right)  =m^{-1}\sum_{i=1}^{m}f_{i}\left(  v,y\right)
\text{ \ and \ \ }G_{m}\left(  v,y\right)  =\mathbb{E}\left[  \hat{G}%
_{m}\left(  v,y\right)  \right]  \text{.}%
\]
Then $\mathbb{G}_{m}\left(  v,y\right)  =\hat{G}_{m}\left(  v,y\right)
-G_{m}\left(  v,y\right)  $.

\textbf{Part I: preparations.} Set $D=\left[  0,c_{1}\right]  \times\left[
c_{2},c_{3}\right]  $. Let $\mathcal{K}_{1}=\left\{  v_{1},\ldots,v_{l_{\ast}%
}\right\}  $ for some $l_{\ast}\in\mathbb{N}_{+}$ be a partition of $\left[
0,c_{1}\right]  $ with norm $\Delta_{1}$ into $l_{\ast}-1$ subintervals of
equal lengths, and let $\mathcal{K}_{2}=\left\{  y_{1},\ldots,y_{l_{\#}%
}\right\}  $ for some $l_{\#}\in\mathbb{N}_{+}$ be a partition of $\left[
c_{2},c_{3}\right]  $ with norm $\Delta_{2}$ into $l_{\#}-1$ subintervals of
equal lengths. So, $\mathcal{K}_{1}\times\mathcal{K}_{2}$ forms a partition of
$D$. For each $\left(  v,y\right)  \in D$, pick $\left(  v_{i},y_{j}\right)  $
from $\mathcal{K}_{1}\times\mathcal{K}_{2}$ that is closest to $\left(
v,y\right)  $ in Euclidean norm. Let $\partial_{\cdot}$ denote the operator of
differentiation with respect to the subscript. Assume \textquotedblleft%
$a_{0}\neq0$ and $\max_{1\leq i\leq m}\mathbb{E}\left(  \left\vert
z_{i}\right\vert ^{2}\right)  <\infty$\textquotedblright\ or \textquotedblleft%
$a_{0}=0$ and $\max_{1\leq i\leq m}\mathbb{E}\left(  \left\vert z_{i}%
\right\vert \right)  <\infty$\textquotedblright. Then $\partial_{v}%
\mathbb{G}_{m}\left(  v,y\right)  $ and $\partial_{y}\mathbb{G}_{m}\left(
v,y\right)  $ are both defined, and by Lagrange mean value theorem, we obtain%
\begin{align}
\mathbb{G}_{m}\left(  v,y\right)   &  \leq\left\vert \mathbb{G}_{m}\left(
v_{i},y_{j}\right)  \right\vert +\left\vert \mathbb{G}_{m}\left(  v,y\right)
-\mathbb{G}_{m}\left(  v_{i},y_{j}\right)  \right\vert \nonumber\\
&  \leq\left\vert \mathbb{G}_{m}\left(  v_{i},y_{j}\right)  \right\vert
+\Delta_{1}\sup_{v\in\left[  0,c_{1}\right]  }\left\vert \partial
_{v}\mathbb{G}_{m}\left(  v,y\right)  \right\vert +\Delta_{2}\sup_{y\in\left[
c_{2},c_{3}\right]  }\left\vert \partial_{y}\mathbb{G}_{m}\left(  v,y\right)
\right\vert . \label{TotalBnd}%
\end{align}
Further, we can apply the dominated convergence theorem to compute
$\partial_{v}\mathbb{G}_{m}\left(  v,y\right)  $ and $\partial_{y}%
\mathbb{G}_{m}\left(  v,y\right)  $ by exchanging differentiation and
integration, and obtain%
\begin{equation}
\left\{
\begin{array}
[c]{c}%
\partial_{v}\mathbb{G}_{m}\left(  v,y\right)  =m^{-1}\sum_{i=1}^{m}\left(
\partial_{v}f_{i}\left(  v,y\right)  -\mathbb{E}\left[  \partial_{v}%
f_{i}\left(  v,y\right)  \right]  \right) \\
\partial_{y}\mathbb{G}_{m}\left(  v,y\right)  =m^{-1}\sum_{i=1}^{m}\left(
\partial_{y}f_{i}\left(  v,y\right)  -\mathbb{E}\left[  \partial_{y}%
f_{i}\left(  v,y\right)  \right]  \right)
\end{array}
\right.  \label{eqPartialsG}%
\end{equation}
based on%
\begin{equation}
\left\{
\begin{array}
[c]{l}%
\partial_{v}f\left(  x,v,y,a_{0},c_{0}\right)  =-\left(  a_{0}x+b_{0}\right)
\left(  x-y\right)  \sin\left(  v\left(  x-y\right)  -c_{0}\right) \\
\partial_{y}f\left(  x,v,y,a_{0},c_{0}\right)  =\left(  a_{0}x+b_{0}\right)
v\sin\left(  v\left(  x-y\right)  -c_{0}\right)
\end{array}
\right.  , \label{eqPartialsf}%
\end{equation}
where we recall%
\[
f\left(  x,v,y,a_{0},c_{0}\right)  =\left(  a_{0}x+b_{0}\right)  \cos\left(
v\left(  x-y\right)  -c_{0}\right)  .
\]
Combining (\ref{eqPartialsG}) and (\ref{eqPartialsf}) gives the following two
inequalities:
\begin{align*}
\left\vert \partial_{v}\mathbb{G}_{m}\left(  v,y\right)  \right\vert  &
\leq\frac{1}{m}\sum_{i=1}^{m}\left\{  \left(  \left\vert a_{0}z_{i}\right\vert
+\left\vert b_{0}\right\vert \right)  \left(  \left\vert z_{i}\right\vert
+\left\vert y\right\vert \right)  +\mathbb{E}\left[  \left(  \left\vert
a_{0}z_{i}\right\vert +\left\vert b_{0}\right\vert \right)  \left(  \left\vert
z_{i}\right\vert +\left\vert y\right\vert \right)  \right]  \right\} \\
&  =2\left\vert b_{0}y\right\vert +\frac{\left\vert a_{0}\right\vert }{m}%
\sum_{i=1}^{m}\left[  \left\vert z_{i}^{2}\right\vert +\mathbb{E}\left(
\left\vert z_{i}^{2}\right\vert \right)  \right]  +\frac{\left\vert
a_{0}y\right\vert +\left\vert b_{0}\right\vert }{m}\sum_{i=1}^{m}\left[
\left\vert z_{i}\right\vert +\mathbb{E}\left(  \left\vert z_{i}\right\vert
\right)  \right]
\end{align*}
and%
\[
\left\vert \partial_{y}\mathbb{G}_{m}\left(  v,y\right)  \right\vert \leq
\frac{\left\vert v\right\vert }{m}\sum_{i=1}^{m}\left[  \left(  \left\vert
a_{0}z_{i}\right\vert +\left\vert b_{0}\right\vert \right)  +\mathbb{E}\left(
\left\vert a_{0}z_{i}\right\vert +\left\vert b_{0}\right\vert \right)
\right]  =2\left\vert vb_{0}\right\vert +\frac{\left\vert va_{0}\right\vert
}{m}\sum_{i=1}^{m}\left[  \left\vert z_{i}\right\vert +\mathbb{E}\left(
\left\vert z_{i}\right\vert \right)  \right]  .
\]
Therefore,
\begin{equation}
\sup_{\left(  v,y\right)  \in D}\left\vert \partial_{v}\mathbb{G}_{m}\left(
v,y\right)  \right\vert \leq2\left\vert b_{0}\right\vert C_{c_{2},c_{3}}%
+\frac{\left\vert a_{0}\right\vert }{m}\sum_{i=1}^{m}\left[  \left\vert
z_{i}^{2}\right\vert +\mathbb{E}\left(  \left\vert z_{i}^{2}\right\vert
\right)  \right]  +\frac{\left\vert a_{0}\right\vert C_{c_{2},c_{3}%
}+\left\vert b_{0}\right\vert }{m}\sum_{i=1}^{m}\left[  \left\vert
z_{i}\right\vert +\mathbb{E}\left(  \left\vert z_{i}\right\vert \right)
\right]  \label{supPartialv}%
\end{equation}
where $C_{x,y}=\max\left\{  \left\vert x\right\vert ,\left\vert y\right\vert
\right\}  $ for all $x,y\in\mathbb{R}$, and%
\begin{equation}
\sup_{\left(  v,y\right)  \in D}\left\vert \partial_{y}\mathbb{G}_{m}\left(
v,y\right)  \right\vert \leq2c_{1}\left\vert b_{0}\right\vert +\frac
{c_{1}\left\vert a_{0}\right\vert }{m}\sum_{i=1}^{m}\left[  \left\vert
z_{i}\right\vert +\mathbb{E}\left(  \left\vert z_{i}\right\vert \right)
\right]  . \label{supPartialY}%
\end{equation}

\textbf{Part II:} bounding the tail probabilities of three empirical
processes. We will consider three cases:%
\[
\left\{
\begin{array}
[c]{l}%
\text{Case 1: }a_{0}=0,b_{0}=1,c_{2}<c_{3}\text{ and }\sup_{\left(
y,v\right)  \in D}\left\vert \mathbb{G}_{m}\left(  v,y\right)  \right\vert \\
\text{Case 2: }a_{0}=0,b_{0}=1,\left[  c_{2},c_{3}\right]  =\left\{
c_{2}\right\}  \text{ and }\sup_{v\in\left[  0,c_{1}\right]  }\left\vert
\mathbb{G}_{m}\left(  v,c_{2}\right)  \right\vert \\
\text{Case 3: }a_{0}=1,b_{0}=0,\left[  c_{2},c_{3}\right]  =\left\{
0\right\}  \text{ and }\sup_{v\in\left[  0,c_{1}\right]  }\left\vert
\mathbb{G}_{m}\left(  v,0\right)  \right\vert
\end{array}
\right.  .
\]
\textbf{Part II, Case 1:} $a_{0}=0,b_{0}=1,c_{2}<c_{3}$ and $\sup_{\left(
y,v\right)  \in D}\left\vert \mathbb{G}_{m}\left(  v,y\right)  \right\vert $.
In this case, $f$ becomes%
\[
f\left(  x,v,y,0,1,c_{0}\right)  =\cos\left(  v\left(  x-y\right)
-c_{0}\right)  ,
\]
and (\ref{supPartialv}) and (\ref{supPartialY}) respectively reduce to%
\begin{equation}
\left\{
\begin{array}
[c]{l}%
\sup_{\left(  v,y\right)  \in D}\left\vert \partial_{v}\mathbb{G}_{m}\left(
v,y\right)  \right\vert \leq2C_{c_{2},c_{3}}+\mathbb{H}_{1,m}\left(
\mathbf{z}\right) \\
\sup_{\left(  v,y\right)  \in D}\left\vert \partial_{y}\mathbb{G}_{m}\left(
v,y\right)  \right\vert \leq2c_{1}%
\end{array}
\right.  , \label{supPartialva}%
\end{equation}
where%
\begin{equation}
\mathbb{H}_{j,m}\left(  \mathbf{z}\right)  =\frac{1}{m}\sum_{i=1}^{m}\left[
\left\vert z_{i}\right\vert ^{j}+\mathbb{E}\left(  \left\vert z_{i}\right\vert
^{j}\right)  \right]  \text{ \ for }j\in\mathbb{N}_{+} \label{Hprocess}%
\end{equation}
whenever $\mathbb{H}_{j,m}\left(  \mathbf{z}\right)  $ is finite for a
$j\in\mathbb{N}_{+}$.

We can relax the upper bound for $\sup_{\left(  v,y\right)  \in D}\left\vert
\partial_{v}\mathbb{G}_{m}\left(  v,y\right)  \right\vert $ in
(\ref{supPartialva}) as follows. Since $\mathcal{F}$ is a location-shift
family, there are independent and identically distributed (i.i.d.) $\left\{
X_{i}\right\}  _{i=1}^{m}$ with common CDF $F_{0}$ such that $z_{i}=\mu
_{i}+X_{i}$ for $1\leq i\leq m$. So, $\mathbb{H}_{j,m}\left(  \mathbf{z}%
\right)  $ is finite for a $j\in\mathbb{N}_{+}$ whenever $\mathbb{E}\left(
\left\vert X_{\left(  0\right)  }\right\vert ^{j}\right)  <\infty$, where
$X_{\left(  0\right)  }$ has CDF $F_{0}$. Recall%
\[
\mathcal{B}_{j,m}\left(  \rho_{m}\right)  =\left\{  \boldsymbol{\mu}%
\in\mathbb{R}^{m}:\frac{1}{m}\sum\nolimits_{i=1}^{m}\left\vert \mu
_{i}\right\vert ^{j}\leq\rho_{m}\right\}  \text{ for }j\in\mathbb{N}_{+}.
\]
If $\boldsymbol{\mu}\in\mathcal{B}_{1,m}\left(  \rho_{m}\right)  $, then%
\begin{align}
2C_{c_{2},c_{3}}+\mathbb{H}_{1,m}\left(  \mathbf{z}\right)   &  \leq
2C_{c_{2},c_{3}}+\frac{2}{m}\sum_{i=1}^{m}\left\vert \mu_{i}\right\vert
+\frac{1}{m}\sum_{i=1}^{m}\left[  \left\vert X_{i}\right\vert +\mathbb{E}%
\left(  \left\vert X_{i}\right\vert \right)  \right] \nonumber\\
&  \leq2C_{c_{2},c_{3}}+2\rho_{m}+\frac{1}{m}\sum_{i=1}^{m}\left[  \left\vert
X_{i}\right\vert -\mathbb{E}\left(  \left\vert X_{i}\right\vert \right)
\right]  +2\mathbb{E}\left[  \left\vert X_{\left(  0\right)  }\right\vert
\right] \label{L2L1Connection}\\
&  =R_{1}\left(  \rho_{m}\right)  +T_{1,m}\left(  \mathbf{x},F_{0}\right)  ,
\label{BndHprocessRelax}%
\end{align}
where $\mathbf{x}=\left(  X_{1},\ldots,X_{m}\right)  $,%
\begin{equation}
R_{1}\left(  \rho_{m}\right)  =2\mathbb{E}\left[  \left\vert X_{(0)}%
\right\vert \right]  +2C_{c_{2},c_{3}}+2\rho_{m}=2\mathbb{E}\left[  \left\vert
X_{(0)}\right\vert \right]  +2\max\left\{  \left\vert c_{2}\right\vert
,\left\vert c_{3}\right\vert \right\}  +2\rho_{m}, \label{TheShiftRho}%
\end{equation}
and whenever $\mathbb{E}\left[  \left\vert X_{\left(  0\right)  }\right\vert
^{j}\right]  <\infty$,%
\begin{equation}
T_{j,m}\left(  \mathbf{x},F_{0}\right)  =\frac{1}{m}\sum_{i=1}^{m}\left[
\left\vert X_{i}\right\vert ^{j}-\mathbb{E}\left(  \left\vert X_{i}\right\vert
^{j}\right)  \right]  \text{ \ for \ }j\in\mathbb{N}_{+}. \label{Tprocess}%
\end{equation}

In summary, if $\boldsymbol{\mu}\in\mathcal{B}_{1,m}\left(  \rho_{m}\right)  $
and $\mathbb{E}\left[  \left\vert X_{\left(  0\right)  }\right\vert \right]
<\infty$, then
\[
\sup_{\left(  v,y\right)  \in D}\left\vert \partial_{v}\mathbb{G}_{m}\left(
v,y\right)  \right\vert \leq R_{1}\left(  \rho_{m}\right)  +T_{1,m}\left(
\mathbf{x},F_{0}\right)  ,
\]
and (\ref{TotalBnd}) is relaxed into%
\[
\mathbb{G}_{m}\left(  v,y\right)  \leq\left\vert \mathbb{G}_{m}\left(
v_{i},y_{j}\right)  \right\vert +\max\left\{  \Delta_{1},\Delta_{2}\right\}
\times\left[  2c_{1}+R_{1}\left(  \rho_{m}\right)  +T_{1,m}\left(
\mathbf{x},F_{0}\right)  \right]  ,
\]
and hence%
\begin{equation}
\sup_{\left(  v,y\right)  \in D}\left\vert \mathbb{G}_{m}\left(  v,y\right)
\right\vert \leq\sup_{\left(  v_{i},y_{j}\right)  \in\mathcal{K}_{1}%
\times\mathcal{K}_{2}}\left\vert \mathbb{G}_{m}\left(  v_{i},y_{j}\right)
\right\vert +C_{\Delta_{1},\Delta_{2}}\left[  2c_{1}+R_{1}\left(  \rho
_{m}\right)  +T_{1,m}\left(  \mathbf{x},F_{0}\right)  \right]
\label{TotalBnda}%
\end{equation}
where we recall $C_{\Delta_{1},\Delta_{2}}=\max\left\{  \Delta_{1},\Delta
_{2}\right\}  $.

We are ready to bound the tail probability of $\sup_{\left(  y,v\right)  \in
G}\left\vert \mathbb{G}_{m}\left(  v,y\right)  \right\vert $. Consider two
positive constants $\kappa_{1},\kappa_{2}$ such that $\kappa_{1}>\kappa_{2}$,
and define%
\begin{equation}
\left\{
\begin{array}
[c]{l}%
B_{1,m}=\Pr\left\{  \sup_{\left(  v_{i},y_{j}\right)  \in\mathcal{K}_{1}%
\times\mathcal{K}_{2}}\left\vert \mathbb{G}_{m}\left(  v_{i},y_{j}\right)
\right\vert \geq\kappa_{1}-\kappa_{2}\right\} \\
B_{2,m}=\Pr\left\{  T_{1,m}\left(  \mathbf{x},F_{0}\right)  \geq C_{\Delta
_{1},\Delta_{2}}^{-1}\kappa_{2}-\left(  2c_{1}+R_{1}\left(  \rho_{m}\right)
\right)  \right\}
\end{array}
\right.  . \label{eventsB}%
\end{equation}
Then (\ref{TotalBnda}) implies%
\begin{equation}
B_{0,m}=\Pr\left\{  \sup\nolimits_{\left(  v,y\right)  \in D}\left\vert
\mathbb{G}_{m}\left(  v,y\right)  \right\vert \geq\kappa_{1}\right\}  \leq
B_{1,m}+B_{2,m}, \label{eq10h1}%
\end{equation}
and it suffices to bound $B_{1,m}$ and $B_{2,m}$. For $B_{1,m}$, we employ
Hoeffding's inequality to obtain
\begin{equation}
\Pr\left\{  \left\vert \mathbb{G}_{m}\left(  v_{i},y_{j}\right)  \right\vert
\geq\lambda\right\}  \leq2\exp\left(  -2^{-1}m\lambda^{2}\right)  \text{ for
each }\lambda>0, \label{eq10g2}%
\end{equation}
and apply to $B_{1,m}$ the union bound on the grid $\mathcal{K}_{1}%
\times\mathcal{K}_{2}$ based on (\ref{eq10g2}) to obtain%
\begin{equation}
B_{1,m}\leq2l_{\ast}l_{\#}\exp\left(  -2^{-1}m\left(  \kappa_{1}-\kappa
_{2}\right)  ^{2}\right)  . \label{eq10g1}%
\end{equation}
For $B_{2,m}$, Chebyshev's inequality implies%
\begin{equation}
B_{2,m}\leq\frac{1}{m}\left[  C_{\Delta_{1},\Delta_{2}}^{-1}\kappa_{2}-\left(
2c_{1}+R_{1}\left(  \rho_{m}\right)  \right)  \right]  ^{-2}\mathbb{V}\left[
\left\vert X_{\left(  0\right)  }\right\vert \right]  . \label{bndProB2m}%
\end{equation}
Combining (\ref{eq10h1}), (\ref{eq10g1}) and (\ref{bndProB2m}), we see that,
if $\boldsymbol{\mu}\in\mathcal{B}_{m}\left(  \rho_{m}\right)  $ and
$\mathbb{E}\left[  \left\vert X_{\left(  0\right)  }\right\vert ^{2}\right]
<\infty$, then%
\begin{equation}
\Pr\left\{  \sup\nolimits_{\left(  v,y\right)  \in D}\left\vert \mathbb{G}%
_{m}\left(  v,y\right)  \right\vert \geq\kappa_{1}\right\}  \leq\hat{p}%
_{1,m}\left(  \kappa_{1},\kappa_{2},c_{1},c_{2},c_{3},\Delta_{1},\Delta
_{2},F_{0},\rho_{m}\right)  , \label{BndDoubleSup}%
\end{equation}
where we recall $\Xi_{1}=\left(  \kappa_{1},\kappa_{2},c_{1},c_{2}%
,c_{3},\Delta_{1},\Delta_{2},F_{0},\rho_{m}\right)  $ and
\begin{align*}
\hat{p}_{1,m}\left(  \Xi_{1}\right)   &  =2\left(  c_{1}\Delta_{1}%
^{-1}+1\right)  \left[  \left(  c_{3}-c_{2}\right)  \Delta_{2}^{-1}+1\right]
\exp\left(  -2^{-1}m\left(  \kappa_{1}-\kappa_{2}\right)  ^{2}\right) \\
&  \text{ \ \ }+\frac{m^{-1}\mathbb{V}\left[  \left\vert X_{\left(  0\right)
}\right\vert \right]  }{\left[  \left(  \max\left\{  \Delta_{1},\Delta
_{2}\right\}  \right)  ^{-1}\kappa_{2}-\left(  2c_{1}+R_{1}\left(  \rho
_{m}\right)  \right)  \right]  ^{2}}%
\end{align*}
and have used $l_{\ast}=c_{1}\Delta_{1}^{-1}+1$ and $l_{\#}=\left(
c_{3}-c_{2}\right)  \Delta_{2}^{-1}+1$ in (\ref{eq10g1}).

\textbf{Part II, Case 2}:\label{PartIIcase2} $a_{0}=0,b_{0}=1,c_{2}=c_{3}$ and
$\sup_{v\in\left[  0,c_{1}\right]  }\left\vert \mathbb{G}_{m}\left(
v,c_{2}\right)  \right\vert $. This is a special case of \textbf{Case
1}.\ that we have already considered. Specifically, in \textbf{Case 2}, $f$
becomes%
\[
f\left(  x,v,y,0,1,c_{0}\right)  =\cos\left(  v\left(  x-c_{2}\right)
-c_{0}\right)  ,
\]
$C_{c_{2},c_{3}}=\max\left\{  \left\vert c_{2}\right\vert ,\left\vert
c_{3}\right\vert \right\}  =\left\vert c_{2}\right\vert $, $\Delta_{2}=0$,
$l_{\#}=1,\mathbb{G}_{m}\left(  v,y\right)  \equiv\mathbb{G}_{m}\left(
v,c_{2}\right)  $ and%
\[
\sup_{\left(  v,y\right)  \in D}\left\vert \mathbb{G}_{m}\left(  v,y\right)
\right\vert =\sup_{v\in\left[  0,c_{1}\right]  }\left\vert \mathbb{G}%
_{m}\left(  v,c_{2}\right)  \right\vert .
\]
Further, we see the following:

\begin{itemize}
\item (\ref{TotalBnd}) reduces to%
\[
\mathbb{G}_{m}\left(  v,c_{2}\right)  \leq\left\vert \mathbb{G}_{m}\left(
v_{i},c_{2}\right)  \right\vert +\Delta_{1}\sup_{v\in\left[  0,c_{1}\right]
}\left\vert \partial_{v}\mathbb{G}_{m}\left(  v,c_{2}\right)  \right\vert ,
\]
and (\ref{supPartialv}) reduces to%
\[
\sup_{v\in\left[  0,c_{1}\right]  }\left\vert \partial_{v}\mathbb{G}%
_{m}\left(  v,c_{2}\right)  \right\vert \leq2\left\vert c_{2}\right\vert
+\frac{1}{m}\sum_{i=1}^{m}\left[  \left\vert z_{i}\right\vert +\mathbb{E}%
\left(  \left\vert z_{i}\right\vert \right)  \right]  =2\left\vert
c_{2}\right\vert +\mathbb{H}_{1,m}\left(  \mathbf{z}\right)  .
\]

\item (\ref{BndHprocessRelax}) and (\ref{TheShiftRho}) respectively reduce to%
\[
2\left\vert c_{2}\right\vert +\mathbb{H}_{1,m}\left(  \mathbf{z}\right)  \leq
R_{1}\left(  \rho_{m}\right)  +T_{1,m}\left(  \mathbf{x},F_{0}\right)  \text{
\ and }R_{1}\left(  \rho_{m}\right)  =2\mathbb{E}\left[  \left\vert
X_{(0)}\right\vert \right]  +2\left\vert c_{2}\right\vert +2\rho_{m}%
\]
when $\boldsymbol{\mu}\in\mathcal{B}_{1,m}\left(  \rho_{m}\right)  $ and
$\mathbb{E}\left[  \left\vert X_{\left(  0\right)  }\right\vert \right]
<\infty$.

\item (\ref{eventsB}) and (\ref{eq10g1}) reduce to%
\[
B_{1,m}=\Pr\left\{  \sup\nolimits_{1\leq i\leq l_{\ast}}\left\vert
\mathbb{G}_{m}\left(  v_{i},c_{2}\right)  \right\vert \geq\kappa_{1}%
-\kappa_{2}\right\}  \leq2l_{\ast}\exp\left(  -2^{-1}m\left(  \kappa
_{1}-\kappa_{2}\right)  ^{2}\right)  .
\]

\end{itemize}

With these in mind, we see that%
\[
\sup_{v\in\left[  0,c_{1}\right]  }\left\vert \mathbb{G}_{m}\left(
v,c_{2}\right)  \right\vert \leq\sup_{1\leq i\leq l_{\ast}}\left\vert
\mathbb{G}_{m}\left(  v_{i},c_{2}\right)  \right\vert +\Delta_{1}\left[
R_{1}\left(  \rho_{m}\right)  +T_{1,m}\left(  \mathbf{x},F_{0}\right)
\right]  .
\]
On the other hand, by Chebyshev's inequality,%
\[
\tilde{B}_{2,m}=\Pr\left\{  T_{1,m}\left(  \mathbf{x},F_{0}\right)  \geq
\Delta_{1}^{-1}\kappa_{2}-R_{1}\left(  \rho_{m}\right)  \right\}  \leq
\frac{m^{-1}\mathbb{V}\left[  \left\vert X_{\left(  0\right)  }\right\vert
\right]  }{\left(  \Delta_{1}^{-1}\kappa_{2}-R_{1}\left(  \rho_{m}\right)
\right)  ^{2}}.
\]
So,%
\[
B_{0,m}=\Pr\left\{  \sup\nolimits_{v\in\left[  0,c_{1}\right]  }\left\vert
\partial_{v}\mathbb{G}_{m}\left(  v,c_{2}\right)  \right\vert \geq\kappa
_{1}\right\}  \leq B_{1,m}+\tilde{B}_{2,m},
\]
and if $\boldsymbol{\mu}\in\mathcal{B}_{1,m}\left(  \rho_{m}\right)  $ and
$\mathbb{E}\left[  X_{\left(  0\right)  }^{2}\right]  <\infty$, then%
\[
\Pr\left\{  \sup\nolimits_{v\in\left[  0,c_{1}\right]  }\left\vert
\mathbb{G}_{m}\left(  v,c_{2}\right)  \right\vert \geq\kappa_{1}\right\}
\leq\hat{p}_{2,m}\left(  \kappa_{1},\kappa_{2},c_{1},c_{2},\Delta_{1}%
,F_{0},\rho_{m}\right)  ,
\]
where $\Xi_{2}=\left(  \kappa_{1},\kappa_{2},c_{1},c_{2},\Delta_{1},F_{0}%
,\rho_{m}\right)  $ and%
\begin{equation}
\hat{p}_{2,m}\left(  \Xi_{2}\right)  =2\left(  c_{1}\Delta_{1}^{-1}+1\right)
\exp\left(  -2^{-1}m\left(  \kappa_{1}-\kappa_{2}\right)  ^{2}\right)
+m^{-1}\mathbb{V}\left[  \left\vert X_{\left(  0\right)  }\right\vert \right]
\left(  \Delta_{1}^{-1}\kappa_{2}-R_{1}\left(  \rho_{m}\right)  \right)
^{-2}. \label{TotalBndCase2}%
\end{equation}
and we have used $l_{\ast}=c_{1}\Delta_{1}^{-1}+1$ in the upper bound for
$B_{1,m}$.

\textbf{Part II, Case 3}:\label{PartIIcase3} $a_{0}=1,b_{0}=0,c_{2}=c_{3}=0$
and $\sup_{v\in\left[  0,c_{1}\right]  }\left\vert \mathbb{G}_{m}\left(
v,0\right)  \right\vert $. In this case, $f$ becomes%
\[
f\left(  x,v,y,1,0,c_{0}\right)  =x\cos\left(  vx-c_{0}\right)  ,
\]
$C_{c_{2},c_{3}}=\max\left\{  \left\vert c_{2}\right\vert ,\left\vert
c_{3}\right\vert \right\}  =0$, $\Delta_{2}=0$ and $l_{\#}=1$. Note that
$f\left(  x,v,y,1,0,2^{-1}\pi\right)  =x\sin\left(  vx\right)  $. Further, we
see the following:

\begin{itemize}
\item (\ref{TotalBnd}) reduces to%
\begin{equation}
\mathbb{G}_{m}\left(  v,0\right)  \leq\left\vert \mathbb{G}_{m}\left(
v_{i},0\right)  \right\vert +\Delta_{1}\sup_{v\in\left[  0,c_{1}\right]
}\left\vert \partial_{v}\mathbb{G}_{m}\left(  v,0\right)  \right\vert .
\label{bndG1}%
\end{equation}

\item (\ref{supPartialv}) reduces to%
\begin{equation}
\sup_{v\in\left[  0,c_{1}\right]  }\left\vert \partial_{v}\mathbb{G}%
_{m}\left(  v,0\right)  \right\vert \leq\frac{1}{m}\sum_{i=1}^{m}\left[
\left\vert z_{i}^{2}\right\vert +\mathbb{E}\left(  \left\vert z_{i}%
^{2}\right\vert \right)  \right]  =\mathbb{H}_{2,m}\left(  \mathbf{z}\right)
, \label{bndG2}%
\end{equation}
where $\mathbb{H}_{2,m}\left(  \mathbf{z}\right)  $ is defined by
(\ref{Hprocess}).
\end{itemize}

Let us bound $\mathbb{H}_{2,m}\left(  \mathbf{z}\right)  $ from above, so as
to upper bound $\sup_{v\in\left[  0,c_{1}\right]  }\mathbb{G}_{m}\left(
v,0\right)  $. Again since $\mathcal{F}$ is a location-shift family, there are
independent and identically distributed (i.i.d.) $\left\{  X_{i}\right\}
_{i=1}^{m}$ with common CDF $F_{0}$ such that $z_{i}=\mu_{i}+X_{i}$ for $1\leq
i\leq m$. So, if $\boldsymbol{\mu}\in\mathcal{B}_{2,m}\left(  \rho_{m}\right)
$ and $\mathbb{E}\left[  \left\vert X_{\left(  0\right)  }\right\vert
^{2}\right]  <\infty$, then%
\begin{align}
\mathbb{H}_{2,m}\left(  \mathbf{z}\right)   &  =\frac{1}{m}\sum_{i=1}%
^{m}\left[  \left(  \mu_{i}+X_{i}\right)  ^{2}+\mathbb{E}\left(  \left(
\mu_{i}+X_{i}\right)  ^{2}\right)  \right]  \leq\frac{1}{m}\sum_{i=1}%
^{m}\left[  2\left(  \mu_{i}^{2}+X_{i}^{2}\right)  +2\mathbb{E}\left(  \mu
_{i}^{2}+X_{i}^{2}\right)  \right] \nonumber\\
&  \leq4\rho_{m}+2T_{2,m}\left(  \mathbf{x},F_{0}\right)  +4\mathbb{E}\left[
X_{\left(  0\right)  }^{2}\right] \nonumber\\
&  \leq2T_{2,m}\left(  \mathbf{x},F_{0}\right)  +R_{2}\left(  \rho_{m}\right)
, \label{bndG4}%
\end{align}
where we recall $R_{2}\left(  \rho_{m}\right)  =4\mathbb{E}\left[  X_{\left(
0\right)  }^{2}\right]  +4\rho_{m}$ and $T_{2,m}\left(  \mathbf{x}%
,F_{0}\right)  $ is defined by (\ref{Tprocess}). Combining (\ref{bndG1}),
(\ref{bndG2}) and (\ref{bndG4}) gives%
\begin{equation}
\sup_{v\in\left[  0,c_{1}\right]  }\mathbb{G}_{m}\left(  v,0\right)  \leq
\max_{1\leq i\leq l_{\ast}}\left\vert \mathbb{G}_{m}\left(  v_{i},0\right)
\right\vert +\Delta_{1}\left(  2T_{2,m}\left(  \mathbf{x},F_{0}\right)
+R_{2}\left(  \rho_{m}\right)  \right)  . \label{bndTotalCase3}%
\end{equation}

We are ready to obtain a tail bound on $\sup_{v\in\left[  0,c_{1}\right]
}\mathbb{G}_{m}\left(  v,0\right)  $. Since $f\left(  x,v,y,1,0,c_{0}\right)
=x\cos\left(  vx-c_{0}\right)  $ is not uniformly bounded by a constant and
$f\left(  x,v,y,1,0,c_{0}\right)  $ does not have uniformly bounded
coordinate-wise differences, neither Hoeffding's inequality nor McDiarmid's
inequality can be directly apply to bound the tail probability of
$\mathbb{G}_{m}\left(  v_{i},0\right)  $ for a fixed $i$. So, let us use the
truncation trick and define the event $\tilde{B}_{\kappa}=\Pr\left\{
\max\nolimits_{1\leq i\leq m}\left\vert z_{i}\right\vert <\kappa\right\}  $
for $\kappa>0$ and set%
\[
\tilde{B}_{1,m}=\Pr\left\{  \sup\nolimits_{1\leq i\leq l_{\ast}}\left\vert
\mathbb{G}_{m}\left(  v_{i},0\right)  \right\vert \geq\kappa_{1}-\kappa
_{2},1_{\tilde{B}_{\kappa}}\right\}  .
\]
Then
\begin{equation}
B_{1,m}=\Pr\left\{  \sup_{1\leq i\leq l_{\ast}}\left\vert \mathbb{G}%
_{m}\left(  v_{i},0\right)  \right\vert \geq\kappa_{1}-\kappa_{2}\right\}
\leq\tilde{B}_{1,m}+1-\Pr\left\{  \tilde{B}_{\kappa}\right\}  ,
\label{NewBndSupGrid}%
\end{equation}
and Hoeffding's inequality implies%
\[
\Pr\left\{  \left\vert \mathbb{G}_{m}\left(  v_{i},0\right)  \right\vert
\geq\kappa_{1}-\kappa_{2},1_{\tilde{B}_{\kappa}}\right\}  \leq2\exp\left(
-2^{-1}m\kappa^{-2}\left(  \kappa_{1}-\kappa_{2}\right)  ^{2}\right)  .
\]
So, a union bound (together with $l_{\ast}=c_{1}\Delta_{1}^{-1}+1$) implies%
\begin{equation}
\tilde{B}_{1,m}\leq2l_{\ast}\exp\left(  -2^{-1}m\kappa^{-2}\left(  \kappa
_{1}-\kappa_{2}\right)  ^{2}\right)  =2\left(  c_{1}\Delta_{1}^{-1}+1\right)
\exp\left(  -2^{-1}m\kappa^{-2}\left(  \kappa_{1}-\kappa_{2}\right)
^{2}\right)  . \label{bndG3}%
\end{equation}
Set%
\[
\hat{B}_{2,m}=\Pr\left\{  \left\vert T_{2,m}\left(  \mathbf{x},F_{0}\right)
\right\vert \geq2^{-1}\Delta_{1}^{-1}\kappa_{2}-2^{-1}R_{2}\left(  \rho
_{m}\right)  \right\}  .
\]
Then, when $\boldsymbol{\mu}\in\mathcal{B}_{2,m}\left(  \rho_{m}\right)  $ and
$\mathbb{E}\left[  X_{\left(  0\right)  }^{2}\right]  <\infty$, Chebyshev's
inequality implies%
\[
\hat{B}_{2,m}\leq\frac{1}{m}\left[  2^{-1}\Delta_{1}^{-1}\kappa_{2}%
-2^{-1}R_{2}\left(  \rho_{m}\right)  \right]  ^{-2}\mathbb{V}\left[
X_{\left(  0\right)  }^{2}\right]  .
\]

It is left to bound $1-\Pr\left\{  \tilde{B}_{\kappa}\right\}  =\Pr\left\{
\max_{1\leq i<m}\left\vert z_{i}\right\vert >\kappa\right\}  $. Recall
$\left\Vert \boldsymbol{\mu}\right\Vert _{\infty}=\max_{1\leq i\leq
m}\left\vert \mu_{i}\right\vert $. Since $z_{i}=\mu_{i}+X_{i}$ for $1\leq
i\leq m$ for i.i.d. $\left\{  X_{i}\right\}  _{i=1}^{m}$ with common CDF
$F_{0}$, then%
\[
\left\vert z_{i}\right\vert \leq\left\vert X_{i}\right\vert +\left\Vert
\boldsymbol{\mu}\right\Vert _{\infty}\text{ and }\Pr\left\{  \max_{1\leq i\leq
m}\left\vert z_{i}\right\vert \leq\kappa\right\}  \geq\Pr\left\{  \max_{1\leq
i\leq m}\left\vert X_{i}\right\vert +\left\Vert \boldsymbol{\mu}\right\Vert
_{\infty}\leq\kappa\right\}  .
\]
Further, $F_{0}$ is an even function. So, $\Pr\left\{  X_{\left(  0\right)
}<-t\right\}  =\Pr\left(  X_{\left(  0\right)  }>t\right)  $ and%
\begin{align*}
1-\Pr\left\{  \tilde{B}_{\kappa}\right\}   &  \leq1-\Pr\left\{  \max_{1\leq
i\leq m}\left\vert X_{i}\right\vert \leq\kappa-\left\Vert \boldsymbol{\mu
}\right\Vert _{\infty}\right\}  =1-\left(  \Pr\left\{  \left\vert X_{\left(
0\right)  }\right\vert \leq\kappa-\left\Vert \boldsymbol{\mu}\right\Vert
_{\infty}\right\}  \right)  ^{m}\\
&  =1-\left[  1-2\Pr\left\{  X_{\left(  0\right)  }>\kappa-\left\Vert
\boldsymbol{\mu}\right\Vert _{\infty}\right\}  \right]  ^{m}%
=1-C_{m,\boldsymbol{\mu},F_{0}},
\end{align*}
where $C_{m,\boldsymbol{\mu},F_{0}}=\left[  1-2\Pr\left\{  X_{\left(
0\right)  }>\kappa-\left\Vert \boldsymbol{\mu}\right\Vert _{\infty}\right\}
\right]  ^{m}$. Thus, when $\boldsymbol{\mu}\in\mathcal{B}_{2,m}\left(
\rho_{m}\right)  $ and $\mathbb{E}\left[  X_{\left(  0\right)  }^{2}\right]
<\infty$, we have%
\[
B_{0,m}=\Pr\left\{  \sup\nolimits_{v\in\left[  0,c_{1}\right]  }\mathbb{G}%
_{m}\left(  v,0\right)  \geq\kappa_{1}\right\}  \leq B_{1,m}+\hat{B}%
_{2,m}+1-\Pr\left\{  \tilde{B}_{\kappa}\right\}  \leq\hat{p}_{3,m}\left(
\Xi_{3}\right)  ,
\]
where $\Xi_{3}=\left(  \kappa_{1},\kappa_{2},c_{1},\Delta_{1},F_{0},\rho
_{m},\kappa\right)  $ and%
\begin{align}
\hat{p}_{3,m}\left(  \Xi_{3}\right)   &  =2\left(  c_{1}\Delta_{1}%
^{-1}+1\right)  \exp\left(  -2^{-1}m\kappa^{-2}\left(  \kappa_{1}-\kappa
_{2}\right)  ^{2}\right) \nonumber\\
&  \text{ \ \ \ }+m^{-1}\left[  2^{-1}\Delta_{1}^{-1}\kappa_{2}-2^{-1}%
R_{2}\left(  \rho_{m}\right)  \right]  ^{-2}\mathbb{V}\left[  X_{\left(
0\right)  }^{2}\right]  +1-C_{m,\boldsymbol{\mu},F_{0}}. \label{probBnd3}%
\end{align}
\qed

\section{Proofs Related to Construction I}

\label{ProofsLocA}

\subsection{Proof of \autoref{ThmIV}}

Pick any $\mu^{\prime}\in U$. Let $\tilde{g}\left(  t;\mu^{\prime}\right)
=\int_{a}^{b}\exp\left\{  -\iota\left(  y-\mu^{\prime}\right)  t\right\}  dy$
and%
\[
Q_{1}\left(  t,x;\mu^{\prime}\right)  =\frac{t}{2\pi}\int_{\left[
-1,1\right]  }\frac{1}{\hat{F}_{\mu^{\prime}}\left(  ts\right)  }\tilde
{g}\left(  ts;\mu^{\prime}\right)  \exp\left(  \iota txs\right)  ds.
\]
Then,%
\begin{align*}
&  \int Q_{1}\left(  t,x;\mu^{\prime}\right)  dF_{\mu}\left(  x\right)
=\frac{t}{2\pi}\int_{-1}^{1}\frac{\hat{F}_{\mu}\left(  ts\right)  }{\hat
{F}_{\mu^{\prime}}\left(  ts\right)  }\tilde{g}\left(  ts;\mu^{\prime}\right)
ds\\
&  =\frac{t}{2\pi}\int_{-1}^{1}\exp\left\{  \iota ts\left(  \mu-\mu^{\prime
}\right)  \right\}  \tilde{g}\left(  ts;\mu^{\prime}\right)  ds\\
&  =\frac{t}{2\pi}\int_{-1}^{1}ds\int_{a}^{b}\exp\left\{  -\iota\left(
y-\mu^{\prime}\right)  ts\right\}  \exp\left\{  \iota ts\left(  \mu
-\mu^{\prime}\right)  \right\}  dy\\
&  =\frac{1}{2\pi}\int_{-1}^{1}ds\int_{a}^{b}t\exp\left\{  \iota ts\left(
\mu-y\right)  \right\}  dy.
\end{align*}
So, from \autoref{LmDirichlet} we have%
\[
\int Q_{1}\left(  t,x;\mu^{\prime}\right)  dF_{\mu}\left(  x\right)  =\frac
{1}{\pi}\int_{\left(  \mu-b\right)  t}^{\left(  \mu-a\right)  t}\frac{\sin
v}{v}dv.
\]
Now set $\mu^{\prime}=0$. Since $\mathcal{F}$ is a Type I location-shift
family, $\hat{F}_{0}\equiv r_{0}$ holds,
\textcolor{blue}{$r_0(\cdot)$ is an even function}, and
\begin{align*}
Q_{1}\left(  t,x;0\right)   &  =\frac{t}{2\pi}\int_{\left[  -1,1\right]
}\frac{1}{r_{0}\left(  ts\right)  }\exp\left(  \iota txs\right)  ds\int%
_{a}^{b}\exp\left(  -\iota yts\right)  dy\\
&  =\frac{t}{2\pi}\int_{a}^{b}dy\int_{\left[  -1,1\right]  }\frac{\cos\left\{
ts\left(  x-y\right)  \right\}  }{r_{0}\left(  ts\right)  }ds.
\end{align*}
Namely, $Q_{1}\left(  t,x;0\right)  $ is exactly $K_{1}\left(  t,x\right)  $.
Finally, we only need to capture the contributions of the end points $a$ and
$b$ to estimating $\pi_{1,m}$. By \autoref{ThmPoinNull}, we only need to set
$\left(  K,\psi\right)  $ as given by (\ref{IV-a}).\qed

\subsection{Proof of \autoref{ThmIVa}}

The proof is split into three parts in three subsections: bounding the
variance of the error $e_{m}\left(  t\right)  =\hat{\varphi}_{m}\left(
t,\mathbf{z}\right)  -\varphi_{m}\left(  t,\boldsymbol{\mu}\right)  $,
deriving a uniform bound on $\left\vert e_{m}\left(  t\right)  \right\vert $
and then the consistency of $\hat{\varphi}_{m}\left(  t,\mathbf{z}\right)  $,
and specializing these results to Gaussian family.

\subsubsection{Bounding the variance of the error}

Define $\hat{\varphi}_{1,m}\left(  t,\mathbf{z}\right)  =m^{-1}\sum_{i=1}%
^{m}K_{1}\left(  t,z_{i}\right)  $\ and $\varphi_{1,m}\left(
t,\boldsymbol{\mu}\right)  =\mathbb{E}\left\{  \hat{\varphi}_{1,m}\left(
t,\mathbf{z}\right)  \right\}  $. First, we study the
\textcolor{blue}{variance} of $e_{1,m}\left(  t\right)  =\hat{\varphi}%
_{1,m}\left(  t,\mathbf{z}\right)  -\varphi_{1,m}\left(  t,\boldsymbol{\mu
}\right)  $. Recall
\[
K_{1}\left(  t,x\right)  =\frac{t}{2\pi}\int_{a}^{b}dy\int_{\left[
-1,1\right]  }\frac{\cos\left\{  ts\left(  x-y\right)  \right\}  }%
{r_{0}\left(  ts\right)  }ds
\]
and $\psi_{1}\left(  t,\mu\right)  =\int K_{1}\left(  t,x\right)  dF_{\mu
}\left(  x\right)  $. Set $w_{1,i}\left(  v,y\right)  =\cos\left\{  v\left(
z_{i}-y\right)  \right\}  $ for each $i$ and $v\in\mathbb{R}$ and $y\in\left[
a,b\right]  $. Define%
\begin{equation}
S_{1,m}\left(  v,y\right)  =\frac{1}{m}\sum_{i=1}^{m}\left[  w_{1,i}\left(
v,y\right)  -\mathbb{E}\left\{  w_{1,i}\left(  v,y\right)  \right\}  \right]
. \label{eq2e}%
\end{equation}
\textcolor{blue}{Since $S_{1,m}\left(  \cdot,\cdot\right)  $ is even in its first argument and
$r_{0}\left(  \cdot\right)  $ is even,} then%
\begin{equation}
e_{1,m}\left(  t\right)  =\frac{t}{2\pi}\int_{a}^{b}dy\int_{\left[
-1,1\right]  }\frac{S_{1,m}\left(  ts,y\right)  }{r_{0}\left(  ts\right)
}%
ds=\textcolor{blue}{\frac{2t}{2\pi}\int_{a}^{b}dy\int_{\left[ 0,1\right] }\frac {S_{1,m}\left( ts,y\right) }{r_{0}\left( ts\right) }ds}.
\label{phiA}%
\end{equation}
Since $\left\vert w_{1,i}\left(  ts,y\right)  \right\vert \leq1$ uniformly in
$\left(  t,s,y,z_{i},i\right)  $ and $\left\{  z_{i}\right\}  _{i=1}^{m}$ are
independent, we have
\begin{align}
\mathbb{V}\left\{  \hat{\varphi}_{1,m}\left(  t,\mathbf{z}\right)  \right\}
&  =\frac{1}{m^{2}}\sum_{i=1}^{m}\mathbb{V}\left(  \frac{t}{2\pi}\int_{a}%
^{b}dy\int_{\left[  -1,1\right]  }\frac{w_{1,i}\left(  ts,y\right)
-\mathbb{E}\left\{  w_{1,i}\left(  ts,y\right)  \right\}  }{r_{0}\left(
ts\right)  }ds\right) \nonumber\\
&  \leq\frac{1}{m^{2}}\sum_{i=1}^{m}\int\left(  \frac{t}{2\pi}\int_{a}%
^{b}dy\int_{\left[  -1,1\right]  }\frac{2}{r_{0}\left(  ts\right)  }ds\right)
^{2}dF_{\mu_{i}}\left(  z\right) \nonumber\\
&  \leq\frac{t^{2}\left(  b-a\right)  ^{2}}{\pi^{2}m}g^{2}\left(  t,0\right)
\label{eq3f}%
\end{align}
where we recall $g\left(  t,\mu\right)  =\int_{\left[  -1,1\right]  }1/r_{\mu
}\left(  ts\right)  ds$.

Recall%
\[
\psi\left(  t,\mu\right)  =\psi_{1}\left(  t,\mu\right)  -2^{-1}\left\{
\psi_{1,0}\left(  t,\mu;a\right)  +\psi_{1,0}\left(  t,\mu;b\right)  \right\}
\]
and the functions \textcolor{blue}{
\[
K_{1,0}\left(  t,x;\mu^{\prime}\right)  =\int_{\left[  -1,1\right]  }\frac
{\omega\left(  s\right)  \cos\left(  ts\left(  x-\mu^{\prime}\right)  \right)
}{r_{\mu^{\prime}}\left(  ts\right)  }ds.
\]
}and $\psi_{1,0}$ from \autoref{ThmPoinNull}. Define for $\tau\in
\textcolor{blue}{U}$%
\begin{equation}
\hat{\varphi}_{1,0,m}\left(  t,\mathbf{z};\tau\right)  =m^{-1}\sum_{i=1}%
^{m}K_{1,0}\left(  t,z_{i};\tau\right)  \text{ and }\varphi_{1,0,m}\left(
t,\boldsymbol{\mu};\tau\right)  =\mathbb{E}\left\{  \hat{\varphi}%
_{1,0,m}\left(  t,\mathbf{z};\tau\right)  \right\}  . \label{eq1d}%
\end{equation}
\textcolor{blue}{
Since $r_{\mu}\equiv r_{0}$ for all $\mu\in U$, both $r_{0}\left(
\cdot\right)  $ and $\omega\left(  \cdot\right)  $ are even, and
$S_{1,m}\left(  \cdot,\cdot\right)  $ is even in its first argument, we obtain\[
K_{1,0}\left(  t,x;\mu^{\prime}\right)  =2\int_{\left[  0,1\right]  }\frac
{\omega\left(  s\right)  \cos\left(  ts\left(  x-\mu^{\prime}\right)  \right)
}{r_{0}\left(  ts\right)  }ds
\]
and with $e_{1,0,m}\left(  t,\mu^{\prime}\right)  =\hat{\varphi}_{1,0,m}\left(
t,\mathbf{z};\mu^{\prime}\right)  -\varphi_{1,0,m}\left(  t,\boldsymbol{\mu};\mu^{\prime}\right)  $, we have\begin{equation}
e_{1,0,m}\left(  t,\tau\right)  =\int_{\left[  -1,1\right]  }\frac{\omega\left(
s\right)  }{r_{\tau}\left(  ts\right)  }S_{1,m}\left(  ts,\tau\right)
ds=2\int_{\left[  0,1\right]  }\frac{\omega\left(  s\right)  }{r_{0}\left(
ts\right)  }S_{1,m}\left(  ts,\tau\right)  ds.\label{phiB}\end{equation}
}By Theorem 2 of \cite{Chen:2018a},%
\begin{equation}
\mathbb{V}\left\{  \hat{\varphi}_{1,0,m}\left(  t,\mathbf{z};\tau\right)
\right\}  \leq\left\Vert \omega\right\Vert _{\infty}^{2}m^{-1}g^{2}\left(
t,\tau\right)  \label{eq3e}%
\end{equation}
Combining (\ref{eq3f}) and (\ref{eq3e}) and applying twice the inequality
$(a_{\ast}+b_{\ast})^{2} \le2 a_{\ast}^{2}+ 2 b_{\ast}^{2}$ for $a_{\ast
},b_{\ast} \in\mathbb{R}$, we have%
\[
\mathbb{V}\left\{  \hat{\varphi}_{m}\left(  t,\mathbf{z}\right)  \right\}
\leq\frac{2\left\Vert \omega\right\Vert _{\infty}^{2}\max_{\tau\in\left\{
a,b\right\}  }g^{2}\left(  t,\tau\right)  }{m}+\frac{2t^{2}\left(  b-a\right)
^{2}}{\pi^{2}m}g^{2}\left(  t,0\right)  .
\]
However, $\mathcal{F}$ is a location-shift family. So, $g\left(  t,\mu\right)
$ is independent of $\mu$, i.e., $g\left(  t,0\right)  =g\left(  t,\mu\right)
$ for all $\mu\in U$. So,%
\begin{equation}
\mathbb{V}\left\{  \hat{\varphi}_{m}\left(  t,\mathbf{z}\right)  \right\}
\leq\left(  \left\Vert \omega\right\Vert _{\infty}^{2}+ \pi^{-2}\left(
b-a\right)  ^{2}t^{2}\right)  \frac{2 g^{2}\left(  t,0\right)  }{m}.
\label{eqLsosdx10}%
\end{equation}

\subsubsection{\textcolor{blue}{Uniform bound on error and consistency}}

Second, let us show the uniform bound on $e_{m}\left(  t\right)  =$
$\hat{\varphi}_{m}\left(  t,\mathbf{z}\right)  -\varphi_{m}\left(
t,\boldsymbol{\mu}\right)  $, for which we will use
\autoref{ConcentrationEmpiricalProcesses}. Let $\left\{  \tau_{m}\right\}
_{m\geq1}$, $\left\{  \gamma_{m}\right\}  _{m\geq1}$ and $\left\{  \rho
_{m}\right\}  _{m\geq1}$ be three positive sequences to be determined later.
In the rest of the proof, $q,\gamma$ and $\vartheta$ are positive constants
and $\vartheta^{\prime}$ is a non-negative constant whose ranges will be
determined later. Let us bound the tail probabilities of $\sup_{v\in\left[
0,\tau_{m}\right]  }S_{1,m}\left(  v,\tau\right)  $ for a fixed $\tau
\in\left[  a,b\right]  $ and $\sup_{\left(  v,y\right)  \in\left[  0,\tau
_{m}\right]  \times\left[  a,b\right]  }S_{1,m}\left(  v,y\right)  $, as follows:

\begin{itemize}
\item \textbf{Case 2a} \textquotedblleft Bound the tail probability of
$\sup_{v\in\left[  0,\tau_{m}\right]  }S_{1,m}\left(  v,\tau\right)  $ when
$\boldsymbol{\mu}\in\mathcal{B}_{1,m}\left(  \rho_{m}\right)  $ and
$\mathbb{E}\left[  X_{\left(  0\right)  }^{2}\right]  <\infty$%
\textquotedblright: we can use Case 2 of
\autoref{ConcentrationEmpiricalProcesses} by setting%
\begin{equation}
\left\{
\begin{array}
[c]{l}%
c_{0}=0,c_{1}=\tau_{m}\rightarrow\infty;c_{2}=c_{3}=\tau;\Delta_{1}%
=m^{-\vartheta}\\
\kappa_{1}=\gamma_{m}m^{-1/2},\kappa_{2}=\gamma_{m}^{-1}m^{-1/2}\ln
c_{1};\gamma_{m}^{2}>\ln\tau_{m};\Delta_{1}^{-1}\kappa_{2}/R_{1}\left(
\rho_{m}\right)  \rightarrow\infty
\end{array}
\right.  , \label{ParConditionsCase2}%
\end{equation}
where $R_{1}\left(  \rho_{m}\right)  $ is defined by (\ref{TheShiftRho}) and
becomes%
\begin{equation}
R_{1}\left(  \rho_{m}\right)  =2\mathbb{E}\left[  \left\vert X_{(0)}%
\right\vert \right]  +2\left\vert \tau\right\vert +2\rho_{m}\text{ for a fixed
}\tau\in\left\{  a,b\right\}  . \label{Case2aRho1}%
\end{equation}
Then $\kappa_{1}>\kappa_{2}$ for all $m\geq1$ due to $\gamma_{m}^{2}>\ln
\tau_{m}$, $c_{1}\Delta_{1}^{-1}\rightarrow\infty$ and%
\begin{align}
\exp\left(  -2^{-1}m\left(  \kappa_{1}-\kappa_{2}\right)  ^{2}\right)   &
=\exp\left(  -2^{-1}\gamma_{m}^{2}-2^{-1}\gamma_{m}^{-2}\ln^{2}c_{1}+\ln
c_{1}\right) \nonumber\\
&  \leq\exp\left(  -2^{-1}\gamma_{m}^{2}+\ln c_{1}\right)  =c_{1}\exp\left(
-2^{-1}\gamma_{m}^{2}\right)  . \label{BndExpProb}%
\end{align}
Let us use $\left\vert o\left(  1\right)  \right\vert $ to emphasize that the
sequence $o\left(  1\right)  $ is non-negative. Then $\hat{p}_{2,m}\left(
\Xi_{2}\right)  $ defined in (\ref{probExpressCase2}) as%
\[
\hat{p}_{2,m}\left(  \Xi_{2}\right)  =2\left(  c_{1}\Delta_{1}^{-1}+1\right)
\exp\left(  -2^{-1}m\left(  \kappa_{1}-\kappa_{2}\right)  ^{2}\right)
+m^{-1}\left[  \Delta_{1}^{-1}\kappa_{2}-R_{1}\left(  \rho_{m}\right)
\right]  ^{-2}\mathbb{V}\left[  \left\vert X_{\left(  0\right)  }\right\vert
\right]
\]
\ satisfies, for all $m$ large enough,%
\begin{align}
\hat{p}_{2,m}\left(  \Xi_{2}\right)   &  \leq2\left(  c_{1}\Delta_{1}%
^{-1}+1\right)  c_{1}\exp\left(  -2^{-1}\gamma_{m}^{2}\right)  +m^{-1}%
\Delta_{1}^{2}\kappa_{2}^{-2}\left(  1-\left\vert o\left(  1\right)
\right\vert \right)  ^{-2}\mathbb{V}\left[  \left\vert X_{\left(  0\right)
}\right\vert \right] \label{BndProbhatp2}\\
&  \leq\hat{p}_{2,m}^{\ast}\left(  c_{1},\Delta_{1},\gamma_{m}\right)
:=C\left[  c_{1}^{2}\Delta_{1}^{-1}\exp\left(  -2^{-1}\gamma_{m}^{2}\right)
+m^{-1}\Delta_{1}^{2}\kappa_{2}^{-2}\right]  \label{probCase2}%
\end{align}
where we have used (\ref{BndExpProb}) and $\Delta_{1}^{-1}\kappa_{2}%
/R_{1}\left(  \rho_{m}\right)  \rightarrow\infty$ in deriving
(\ref{BndProbhatp2}). Namely,%
\begin{equation}
\Pr\left\{  \sup\nolimits_{v\in\left[  0,\tau_{m}\right]  }S_{1,m}\left(
v,\tau\right)  \geq\kappa_{1}\right\}  \leq\hat{p}_{2,m}^{\ast}\left(
c_{1},\Delta_{1},\gamma_{m}\right)  \label{supProbCase2}%
\end{equation}

\item \textbf{Case 1a} \textquotedblleft Bound the tail probability of
$\sup_{\left(  v,y\right)  \in\left[  0,\tau_{m}\right]  \times\left[
a,b\right]  }S_{1,m}\left(  v,y\right)  $ when $\boldsymbol{\mu}\in
\mathcal{B}_{1,m}\left(  \rho_{m}\right)  $ and $\mathbb{E}\left[  X_{\left(
0\right)  }^{2}\right]  <\infty$\textquotedblright: we can use Case 1 of
\autoref{ConcentrationEmpiricalProcesses} by setting%
\begin{equation}
\left\{
\begin{array}
[c]{l}%
c_{0}=0,c_{1}=\tau_{m}\rightarrow\infty,\left[  c_{2},c_{3}\right]  =\left[
a,b\right]  ,\Delta_{1}=\Delta_{2}=m^{-\vartheta}\\
\kappa_{1}=\gamma_{m}m^{-1/2},\kappa_{2}=\gamma_{m}^{-1}m^{-1/2}\ln
c_{1};\gamma_{m}^{2}>\ln\tau_{m};\\
\Delta_{1}^{-1}\kappa_{2}/\left(  2c_{1}+R_{1}\left(  \rho_{m}\right)
\right)  \rightarrow\infty
\end{array}
\right.  , \label{ParConditionsCase1}%
\end{equation}
where $R_{1}\left(  \rho_{m}\right)  $ is defined by (\ref{TheShiftRho}) and
becomes%
\begin{equation}
R_{1}\left(  \rho_{m}\right)  =2\mathbb{E}\left[  \left\vert X_{(0)}%
\right\vert \right]  +2\max\left\{  \left\vert a\right\vert ,\left\vert
b\right\vert \right\}  +2\rho_{m}. \label{Case1aRho1}%
\end{equation}
Then $\hat{p}_{1,m}\left(  \Xi_{1}\right)  $ defined in
(\ref{probExpressCase1}) as%
\begin{align*}
\hat{p}_{1,m}\left(  \Xi_{1}\right)   &  =2\left(  c_{1}\Delta_{1}%
^{-1}+1\right)  \left[  \left(  c_{3}-c_{2}\right)  \Delta_{2}^{-1}+1\right]
\exp\left(  -2^{-1}m\left(  \kappa_{1}-\kappa_{2}\right)  ^{2}\right) \\
&  \text{ \ \ }+m^{-1}\left[  \left(  \max\left\{  \Delta_{1},\Delta
_{2}\right\}  \right)  ^{-1}\kappa_{2}-\left(  2c_{1}+R_{1}\left(  \rho
_{m}\right)  \right)  \right]  ^{-2}\mathbb{V}\left[  \left\vert X_{\left(
0\right)  }\right\vert \right]
\end{align*}
satisfies, for all $m$ large enough,%
\begin{align}
\hat{p}_{1,m}\left(  \Xi_{1}\right)   &  \leq Cc_{1}\Delta_{1}^{-2}c_{1}%
\exp\left(  -2^{-1}\gamma_{m}^{2}\right)  +m^{-1}\Delta_{1}^{2}\kappa_{2}%
^{-2}\left(  1-\left\vert o\left(  1\right)  \right\vert \right)
^{-2}\mathbb{V}\left[  \left\vert X_{\left(  0\right)  }\right\vert \right]
\label{BndProbhatp1}\\
&  \leq\hat{p}_{1,m}^{\ast}\left(  c_{1},\Delta_{1},\gamma_{m}\right)
:=C\left[  c_{1}^{2}\Delta_{1}^{-2}\exp\left(  -2^{-1}\gamma_{m}^{2}\right)
+m^{-1}\Delta_{1}^{2}\kappa_{2}^{-2}\right]  \label{probCase1}%
\end{align}
where we have used (\ref{BndExpProb}) and $\Delta_{1}^{-1}\kappa_{2}/\left(
2c_{1}+R_{1}\left(  \rho_{m}\right)  \right)  \rightarrow\infty$ in deriving
(\ref{BndProbhatp1}). Namely,
\begin{equation}
\Pr\left\{  \sup\nolimits_{\left(  v,y\right)  \in\left[  0,\tau_{m}\right]
\times\left[  a,b\right]  }S_{1,m}\left(  v,y\right)  \geq\kappa_{1}\right\}
\leq\hat{p}_{1,m}^{\ast}\left(  c_{1},\Delta_{1},\gamma_{m}\right)  .
\label{supProbCase1}%
\end{equation}

\end{itemize}

When no confusion should arise, we will just write $\hat{p}_{1,m}^{\ast
}\left(  c_{1},\Delta_{1},\gamma_{m}\right)  $ as $\hat{p}_{1,m}^{\ast}$ and
$\hat{p}_{2,m}^{\ast}\left(  c_{1},\Delta_{1},\gamma_{m}\right)  $ as $\hat
{p}_{2,m}^{\ast}$. Since $s\in\left[  0,1\right]  $ and $S_{1,m}\left(
\cdot,\cdot\right)  $ is continuous, we must have%
\begin{equation}
\left\{
\begin{array}
[c]{l}%
\sup_{\left(  t,s\right)  \in\left[  0,\tau_{m}\right]  \times\left[
0,1\right]  }S_{1,m}\left(  ts,y\right)  =\sup_{v\in\left[  0,\tau_{m}\right]
}S_{1,m}\left(  v,y\right)  \text{ for each fixed }y\in\mathbb{R}\\
\sup_{\left(  t,s,y\right)  \in\left[  0,\tau_{m}\right]  \times\left[
0,1\right]  \times\left[  a,b\right]  }S_{1,m}\left(  ts,y\right)
=\sup_{\left(  v,y\right)  \in\left[  0,\tau_{m}\right]  \times\left[
a,b\right]  }S_{1,m}\left(  v,y\right)
\end{array}
\right.  . \label{SupEquiv}%
\end{equation}
Recall the settings (\ref{ParConditionsCase2}) and (\ref{ParConditionsCase1})
for the parameters $\kappa_{1}=\gamma_{m}m^{-1/2},\kappa_{2}=\gamma_{m}%
^{-1}m^{-1/2}\ln\tau_{m}$ and $\Delta_{1}=m^{-\vartheta}$. Then, (\ref{phiA}),
(\ref{supProbCase1}) and (\ref{SupEquiv})\ together imply that, with
probability at least $1-\hat{p}_{1,m}^{\ast}$,%
\begin{align}
&  \sup_{\boldsymbol{\mu}\in\mathcal{B}_{1,m}\left(  \rho_{m}\right)  }%
\sup_{t\in\left[  0,\tau_{m}\right]  }\left\vert e_{1,m}\left(  t\right)
\right\vert \nonumber\\
&  \leq\frac{2\left(  b-a\right)  \tau_{m}\gamma_{m}}{2\pi\sqrt{m}}\sup
_{t\in\left[  0,\tau_{m}\right]  }\int_{\left[  0,1\right]  }\frac{ds}%
{r_{0}\left(  ts\right)  }=\left(  2\pi\right)  ^{-1}\left(  b-a\right)
{\tau_{m}\Upsilon\left(  \gamma_{m}m^{-1/2},\tau_{m}\right)  },
\label{BndNullPhi1Sup}%
\end{align}
where%
\begin{equation}
\textcolor{blue}{\Upsilon\left( \lambda,\tilde{\tau}\right)}=2\lambda
\sup_{t\in\left[  0,\tilde{\tau}\right]  }\int_{\left[  0,1\right]  }\frac
{ds}{r_{0}\left(  ts\right)  }\text{ for }\lambda,\tilde{\tau}\geq0\text{.}
\label{UpsilonA}%
\end{equation}
Further, (\ref{phiB}), (\ref{supProbCase2}) and (\ref{SupEquiv})\ together
imply that, with probability at least $1-\hat{p}_{2,m}^{\ast}$,
\begin{align}
&  \sup_{\boldsymbol{\mu}\in\mathcal{B}_{1,m}\left(  \rho_{m}\right)  }%
\sup_{t\in\left[  0,\tau_{m}\right]  }\left\vert e_{1,0,m}\left(
t,\tau\right)  \right\vert \nonumber\\
&  \leq\frac{2\left\Vert \omega\right\Vert _{\infty}\gamma_{m}}{\sqrt{m}}%
\sup_{t\in\left[  0,\tau_{m}\right]  }\int_{\left[  0,1\right]  }\frac
{ds}{r_{0}\left(  ts\right)  }=\left\Vert \omega\right\Vert _{\infty
}\textcolor{blue}{\Upsilon\left(  \gamma_{m}m^{-1/2},\tau_{m}\right)  }
\label{PointNullPhi10Sup}%
\end{align}
for $\tau\in\left\{  a,b\right\}  $. Therefore, in view of%
\[
\hat{\varphi}_{m}\left(  t,\mathbf{z}\right)  =\hat{\varphi}_{1,m}\left(
t,\mathbf{z}\right)  -2^{-1}\sum_{\tau\in\left\{  a,b\right\}  }\hat{\varphi
}_{1,0,m}\left(  t,\mathbf{z};\tau\right)  ,
\]
(\ref{BndNullPhi1Sup}) and (\ref{PointNullPhi10Sup}), a union bound for
probability implies that%
\begin{equation}
\sup_{\boldsymbol{\mu}\in\mathcal{B}_{1,m}\left(  \rho_{m}\right)  }\sup
_{t\in\left[  0,\tau_{m}\right]  }\left\vert e_{m}\left(  t\right)
\right\vert \leq\left\{  \left(  b-a\right)  \frac{\tau_{m}}{2\pi}+\left\Vert
\omega\right\Vert _{\infty}\right\}  \textcolor{blue}{\Upsilon\left(  \gamma
_{m}m^{-1/2},\tau_{m}\right)  } \label{BndNullErrorSup}%
\end{equation}
with probability at least $1-2\hat{p}_{2,m}^{\ast}-\hat{p}_{1,m}^{\ast}$.

\textcolor{blue}{ Third, we derive conditions that ensure $\hat{\varphi}_{m}\left(
t,\mathbf{z}\right)  $ is consistent and uniformly consistent. The settings in
(\ref{ParConditionsCase2}) and (\ref{ParConditionsCase1}) and definitions of
$\hat{p}_{2,m}^{\ast}$ and $\hat{p}_{1,m}^{\ast}$ respectively in
(\ref{probCase2}) and (\ref{probCase1}) imply}%
\[
c_{1}^{2}\Delta_{1}^{-2}\exp\left(  -2^{-1}\gamma_{m}^{2}\right)
+m^{-1}\Delta_{1}^{2}\kappa_{2}^{-2}=\tau_{m}^{2}m^{2\vartheta}\exp\left(
-2^{-1}\gamma_{m}^{2}\right)  +m^{-2\vartheta}\gamma_{m}^{2}\ln^{-2}\tau_{m}.
\]
Recall $\Delta_{1}^{-1}\kappa_{2}/\left(  2c_{1}+R_{1}\left(  \rho_{m}\right)
\right)  \rightarrow\infty$ in (\ref{ParConditionsCase1}) with $R_{1}\left(
\rho_{m}\right)  $ as (\ref{Case1aRho1}) and $\Delta_{1}^{-1}\kappa_{2}%
/R_{1}\left(  \rho_{m}\right)  \rightarrow\infty$ in (\ref{ParConditionsCase2}%
) with $R_{1}\left(  \rho_{m}\right)  $ as (\ref{Case2aRho1}). We see that the
difference between the two $R_{1}\left(  \rho_{m}\right)  $'s is a constant,
that the two $R_{1}\left(  \rho_{m}\right)  $'s are of the same order as that
of $\rho_{m}$, and that for the rest of the proof, it suffices to take either
$R_{1}\left(  \rho_{m}\right)  $ and work with it. However, $c_{1}=\tau
_{m}\geq0$. So, $\Delta_{1}^{-1}\kappa_{2}/\left(  2c_{1}+R_{1}\left(
\rho_{m}\right)  \right)  \rightarrow\infty$ implies $\Delta_{1}^{-1}%
\kappa_{2}/R_{1}\left(  \rho_{m}\right)  \rightarrow\infty$, and {%
\begin{equation}
\left\{
\begin{array}
[c]{l}%
\left(  2\tau_{m}+R_{1}\left(  \rho_{m}\right)  \right)  ^{-1}m^{\vartheta
-1/2}\gamma_{m}^{-1}\ln\tau_{m}\rightarrow\infty\\
\multicolumn{1}{c}{\tau_{m}^{2}m^{2\vartheta}\exp\left(  -2^{-1}\gamma_{m}%
^{2}\right)  +m^{-2\vartheta}\gamma_{m}^{2}\ln^{-2}\tau_{m}=o\left(  1\right)
}%
\end{array}
\right.  \Longrightarrow\hat{p}_{2,m}^{\ast}=o\left(  1\right)  =\hat{p}%
_{1,m}^{\ast}. \label{BndNullProbOnErrorBnd}%
\end{equation}
}Consider the decomposition%
\begin{align}
\hat{\varphi}_{m}\left(  t,\mathbf{z}\right)  =  &  -\left\{  \hat{\varphi
}_{1,m}\left(  t,\mathbf{z}\right)  -\varphi_{1,m}\left(  t,\boldsymbol{\mu
}\right)  \right\}  +\frac{1}{2}\left\{  \hat{\varphi}_{1,0,m}\left(
t,\mathbf{z};\tau\right)  -\varphi_{1,0,m}\left(  t,\boldsymbol{\mu}%
;\tau\right)  \right\} \nonumber\\
&  +\frac{1}{2}\left\{  \hat{\varphi}_{1,0,m}\left(  t,\mathbf{z};b\right)
-\varphi_{1,0,m}\left(  t,\boldsymbol{\mu};b\right)  \right\}  +\varphi
_{m}\left(  t,\boldsymbol{\mu}\right)  , \label{BndNullDecomposition}%
\end{align}
where%
\[
\varphi_{m}\textcolor{blue}{\left(t,\boldsymbol{\mu}\right)}=1-\varphi
_{1,m}\left(  \textcolor{blue}{t},\boldsymbol{\mu}\right)  +2^{-1}%
\varphi_{1,0,m}\left(  \textcolor{blue}{t},\boldsymbol{\mu};a\right)
+2^{-1}\varphi_{1,0,m}\left(  \textcolor{blue}{t},\boldsymbol{\mu};b\right)
.
\]
\textcolor{blue}{\autoref{SpeedOracleBoundedNull}, whose $\tilde{\psi}_{1,0}$ has parameter $\sigma>0$ and
whose $\tilde{\psi}_{1,0}$ with $\sigma=1$ reduces to $\psi_{1,0}$}, asserts that, for all positive $t$ such that $t\left(  b-a\right)
\geq2$ and $tu_{m}\geq2$,
\begin{equation}
\left\vert \pi_{1,m}^{-1}\varphi_{m}\left(  t,\boldsymbol{\mu}\right)
-1\right\vert \leq C\left(  \frac{1}{\pi_{1,m}u_{m}t}+\frac{1}{t\pi_{1,m}%
}\right)  , \label{eq10ca}%
\end{equation}
where $u_{m}=\min_{\tau\in\left\{  a,b\right\}  }\min_{\left\{  j:\mu_{j}%
\neq\tau\right\}  }\left\vert \mu_{j}-\tau\right\vert $. So, $\pi_{1,m}%
^{-1}\varphi_{m}\left(  \tau_{m},\boldsymbol{\mu}\right)  =1+o\left(
1\right)  $ when $\tau_{m}^{-1}\left(  1+u_{m}^{-1}\right)  =o\left(
\pi_{1,m}\right)  $, and it suffices to show $\pi_{1,m}^{-1}\left\vert
\hat{\varphi}_{m}\left(  t,\mathbf{z}\right)  -\varphi_{m}\left(
t,\boldsymbol{\mu}\right)  \right\vert \rightsquigarrow0$ for $t$ in a
suitable range. From (\ref{BndNullErrorSup}) and $\left\Vert \omega\right\Vert
_{\infty}<\infty$, we obtain%
\begin{equation}
\sup_{t\in\left[  0,\tau_{m}\right]  }\pi_{1,m}^{-1}\left\vert e_{m}\left(
t\right)  \right\vert \leq C\pi_{1,m}^{-1}\tau_{m}{\Upsilon\left(  \gamma
_{m}m^{-1/2},\tau_{m}\right)  }. \label{BndNullErrorRate}%
\end{equation}
\textcolor{blue}{ Therefore, from (\ref{BndNullProbOnErrorBnd}), (\ref{eq10ca}),
(\ref{BndNullErrorRate}) and the decomposition}%
\[
\left\vert \pi_{1,m}^{-1}\hat{\varphi}_{m}\left(  t,\mathbf{z}\right)
-1\right\vert \leq\pi_{1,m}^{-1}\left\vert \hat{\varphi}_{m}\left(
t,\mathbf{z}\right)  -\varphi_{m}\left(  t,\boldsymbol{\mu}\right)
\right\vert +\left\vert \pi_{1,m}^{-1}\varphi_{m}\left(  t_{m},\boldsymbol{\mu
}\right)  -1\right\vert ,
\]
{we conclude:%
\begin{equation}
\left\{
\begin{array}
[c]{l}%
\left(  2\tau_{m}+R_{1}\left(  \rho_{m}\right)  \right)  ^{-1}m^{\vartheta
-1/2}\gamma_{m}^{-1}\ln\tau_{m}\rightarrow\infty\\
\tau_{m}^{2}m^{2\vartheta}\exp\left(  -2^{-1}\gamma_{m}^{2}\right)
+m^{-2\vartheta}\gamma_{m}^{2}\ln^{-2}\tau_{m}=o\left(  1\right) \\
\tau_{m}^{-1}\left(  1+u_{m}^{-1}\right)  =o\left(  \pi_{1,m}\right)
,\pi_{1,m}^{-1}\tau_{m}{\Upsilon\left(  \gamma_{m}m^{-1/2},\tau_{m}\right)
}=o\left(  1\right)
\end{array}
\right.  \Longrightarrow\left\vert \pi_{1,m}^{-1}\hat{\varphi}_{m}\left(
\tau_{m},\mathbf{z}\right)  -1\right\vert \rightsquigarrow0.
\label{BndNullConsistency}%
\end{equation}
In order to obtain a uniform consistency class, we just need to ensure that
the conditions in (\ref{BndNullConsistency}) is uniform to some extent. To
this tend, for each $m\geq1$ let $\Omega_{m}\left(  \rho_{m}\right)  $ be a
subset of $\mathcal{B}_{1,m}\left(  \rho_{m}\right)  $. Notice that we already
have (\ref{BndNullErrorSup}) and (\ref{eq10ca}). So, a uniform consistency
class for $\hat{\varphi}_{m}\left(  t,\mathbf{z}\right)  $ is%
\begin{equation}
\mathcal{Q}\left(  \mathcal{F}\right)  =\left\{
\begin{array}
[c]{l}%
\left(  2\tau_{m}+R_{1}\left(  \rho_{m}\right)  \right)  ^{-1}m^{\vartheta
-1/2}\gamma_{m}^{-1}\ln\tau_{m}\rightarrow\infty\\
\tau_{m}^{2}m^{2\vartheta}\exp\left(  -2^{-1}\gamma_{m}^{2}\right)
+m^{-2\vartheta}\gamma_{m}^{2}\ln^{-2}\tau_{m}=o\left(  1\right) \\
\tau_{m}^{-1}\sup_{\boldsymbol{\mu}\in\Omega_{m}\left(  \rho_{m}\right)  }%
\pi_{1,m}^{-1}\left(  1+u_{m}^{-1}\right)  =o\left(  1\right) \\
\tau_{m}\Upsilon\left(  \gamma_{m}m^{-1/2},\tau_{m}\right)  \sup
_{\boldsymbol{\mu}\in\Omega_{m}\left(  \rho_{m}\right)  }\pi_{1,m}%
^{-1}=o\left(  1\right)
\end{array}
\right\}  . \label{BndNullUniformConsistency}%
\end{equation}
}

\subsubsection{\textcolor{blue}{Specialization to Gaussian family}}

Finally, we specialize the previous results to Gaussian family. Consider the
Gaussian family $\mathcal{F}$ with scale parameter $\sigma>0$. In this
setting, $r_{0}\left(  t\right)  =\exp\left(  -2^{-1}t^{2}\sigma^{2}\right)
.$ Consider $m\geq2$. Let us specify when the conditions in
(\ref{BndNullConsistency}) hold. Set $\tau_{m}=\sqrt{2\gamma\ln m}$ and
$\gamma_{m}=\sqrt{2q\ln m}$ for some positive constants $q$ and $\gamma$, such
that $q>2\vartheta>0$. Then,%
\begin{equation}
\tau_{m}^{2}m^{2\vartheta}\exp\left(  -2^{-1}\gamma_{m}^{2}\right)
+m^{-2\vartheta}\gamma_{m}^{2}\ln^{-2}\tau_{m}\leq C\left(  m^{2\vartheta
-q}\ln m+m^{-2\vartheta}\ln m\times o\left(  1\right)  \right)  =o\left(
1\right)  , \label{BndNullProbToZeroGaussian}%
\end{equation}
which implies $\hat{p}_{2,m}^{\ast}=o\left(  1\right)  =\hat{p}_{1,m}^{\ast}$.
In addition, $\vartheta>2^{-1}$ and $R_{1}\left(  \rho_{m}\right)  =O\left(
m^{\vartheta^{\prime}}\right)  $ with $0\leq\vartheta^{\prime}<\vartheta
-2^{-1}$ imply%
\begin{equation}
\frac{\Delta_{1}^{-1}\kappa_{2}}{2c_{1}+R_{1}\left(  \rho_{m}\right)  }%
=\frac{m^{\vartheta}\ln\left(  \sqrt{2\gamma\ln m}\right)  }{\sqrt{2q\ln
m}\left[  2\sqrt{2\gamma\ln m}+R_{1}\left(  \rho_{m}\right)  \right]
}\rightarrow\infty\text{.} \label{TheOrderForShiftRho}%
\end{equation}

On the other hand, recall $g\left(  t,0\right)  =2\int_{\left[  0,1\right]
}\left(  1/r_{0}\left(  ts\right)  \right)  ds$. Then $\Upsilon\left(
\lambda,\tilde{\tau}\right)  =\lambda\sup_{t\in\left[  0,\tilde{\tau}\right]
}g\left(  \tilde{\tau},0\right)  $. Since $1/r_{0}\left(  t\right)
=\exp\left(  2^{-1}t^{2}\sigma^{2}\right)  $ and $s^{2}\leq s$ for
$s\in\left[  0,1\right]  $,%
\begin{align}
\sup_{t\in\left[  0,\tilde{\tau}\right]  }g\left(  t,0\right)   &  \leq
2\int_{0}^{1}\sup_{t\in\left[  0,\tilde{\tau}\right]  }\exp\left(  2^{-1}%
t^{2}s\sigma^{2}\right)  ds=2\int_{0}^{1}\exp\left(  2^{-1}s\tilde{\tau}%
^{2}\sigma^{2}\right)  ds\nonumber\\
&  =\frac{4\left(  \exp\left(  2^{-1}\tilde{\tau}^{2}\right)  -1\right)
}{\tilde{\tau}^{2}\sigma^{2}}\leq\frac{4\exp\left(  2^{-1}\tilde{\tau}%
^{2}\sigma^{2}\right)  }{\tilde{\tau}^{2}\sigma^{2}}. \label{supofg(t,0)}%
\end{align}
This implies
\begin{equation}
\Upsilon\left(  \lambda,\tilde{\tau}\right)  \leq\lambda\frac{4\exp\left(
2^{-1}\tilde{\tau}^{2}\sigma^{2}\right)  }{\tilde{\tau}^{2}\sigma^{2}}\text{
\ and \ }{\Upsilon\left(  \gamma_{m}m^{-1/2},\tau_{m}\right)  }\leq
\frac{4\sqrt{2q\ln m}}{2\gamma\ln m\times\sigma^{2}}m^{\sigma^{2}\gamma-1/2},
\label{boundOnUpsilon}%
\end{equation}
and $\tau_{m}\pi_{1,m}^{-1}{\Upsilon\left(  \gamma_{m}m^{-1/2},\tau
_{m}\right)  }=o\left(  1\right)  $ if\ $\pi_{1,m}^{-1}m^{\sigma^{2}%
\gamma-1/2}=o\left(  1\right)  $. Namely, (\ref{BndNullErrorSup}) becomes%
\begin{equation}
\sup_{\boldsymbol{\mu}\in\mathcal{B}_{1,m}\left(  \rho_{m}\right)  }\sup
_{t\in\left[  0,\sqrt{2\gamma\ln m}\right]  }\left\vert e_{m}\left(  t\right)
\right\vert \leq C\frac{4\sqrt{q}m^{\sigma^{2}\gamma-1/2}}{\sigma^{2}%
\sqrt{\gamma}}\leq Cm^{\sigma^{2}\gamma-1/2} \label{ErrorBoundGaussianBndNull}%
\end{equation}
with probability at least $1-o\left(  1\right)  $, and the right most upper
bound in (\ref{ErrorBoundGaussianBndNull})\ is $o\left(  1\right)  $ when
$0<\gamma<2^{-1}\sigma^{-2}$.

Note that $\pi_{1,m}^{-1}m^{\sigma^{2}\gamma-1/2}=o\left(  1\right)  $ for
$\pi_{1,m}^{-1}=O\left(  \sqrt{\ln m}\right)  $ and $0<\gamma<2^{-1}%
\sigma^{-2}$. Further, the condition $\tau_{m}^{-1}\left(  1+u_{m}%
^{-1}\right)  =o\left(  \pi_{1,m}\right)  $ implies $\tau_{m}u_{m}%
\rightarrow\infty$, which implies $\pi_{1,m}^{-1}\leq C\sqrt{\ln m}$ for a
constant $C>0$ for all sufficiently large $m$. For example, we can set
$u_{m}\geq\left(  \tau_{m}\right)  ^{-1}\ln\ln m$ to ensure $\tau_{m}%
u_{m}\rightarrow\infty$. With these preparations, we obtain a uniform
consistency class as%
\[
\mathcal{Q}\left(  \mathcal{F}\right)  =\left\{
\begin{array}
[c]{l}%
0<\gamma<2^{-1}\sigma^{-2};\tau_{m}=\sqrt{2\gamma\ln m}\\
q>2\vartheta;0\leq\vartheta^{\prime}<\vartheta-2^{-1};R_{1}\left(  \rho
_{m}\right)  =O\left(  m^{\vartheta^{\prime}}\right) \\
\pi_{1,m}^{-1}=O\left(  \sqrt{\ln m}\right)  ;u_{m}\geq\left(  \tau
_{m}\right)  ^{-1}\ln\ln m
\end{array}
\right\}  .
\]
The set $\mathcal{Q}\left(  \mathcal{F}\right)  $ can be interpreted as
follows: let $\tilde{C}>0$ be a constant, and define%
\[
\mathcal{U}_{m}=\left\{  \boldsymbol{\mu}\in\mathcal{B}_{1,m}\left(  \rho
_{m}\right)  :%
\begin{array}
[c]{l}%
\vartheta=3/4;q=7/4;\tau_{m}=\sqrt{2^{-1}\sigma^{-2}\ln m}\\
0\leq\vartheta^{\prime}<1/4;R_{1}\left(  \rho_{m}\right)  \leq\tilde
{C}m^{\vartheta^{\prime}}\\
\pi_{1,m}^{-1}\leq\tilde{C}\sqrt{\ln m};u_{m}\geq\left(  \tau_{m}\right)
^{-1}\ln\ln m
\end{array}
\right\}  ,
\]
then%
\[
\sup\nolimits_{\boldsymbol{\mu}\in\mathcal{U}_{m}}\left(  \left\vert \pi
_{1,m}^{-1}\sup\nolimits_{t\in\left\{  \tau_{m}\right\}  }\hat{\varphi}%
_{m}\left(  t,\mathbf{z}\right)  -1\right\vert \right)  \rightsquigarrow
0\text{ \ as \ }m\rightarrow\infty.
\]
\qed

\section{Proofs Related to Construction II}

\label{AppProofsLocB}

\subsection{Proof of \autoref{ThmTypeIOneSided}}

Since $F_{0}\left(  x\right)  $ is differentiable in $x$ and $\int\left\vert
x\right\vert dF_{\mu}\left(  x\right)  <\infty$ for all $\mu\in U$,
$r_{0}\left(  t\right)  $ is differentiable in $t\in\mathbb{R}$. Assume
(\ref{CondVc1}), i.e.,%
\[
\int_{0}^{t}\frac{1}{y}dy\int_{-1}^{1}\left\vert \frac{d}{ds}\frac{1}%
{r_{0}\left(  ys\right)  }\right\vert ds<\infty\text{ \ for each }t>0.
\]
Then%
\begin{equation}
\frac{1}{2\pi}\int_{0}^{t}dy\int_{-1}^{1}\frac{1}{\iota y}\frac{d}{ds}%
\frac{\exp\left(  \iota ysx\right)  }{r_{0}\left(  ys\right)  }ds=\frac
{1}{2\pi}\int_{0}^{1}dy\int_{-1}^{1}\frac{1}{\iota y}\frac{d}{ds}\frac
{\exp\left(  \iota tysx\right)  }{r_{0}\left(  tys\right)  }ds, \label{eq10e}%
\end{equation}
and%
\[
K_{1}^{\dagger}\left(  t,x\right)  =\frac{1}{2\pi}\int_{0}^{1}dy\int_{-1}%
^{1}\frac{1}{\iota y}\frac{d}{ds}\frac{\exp\left(  \iota tysx\right)  }%
{r_{0}\left(  tys\right)  }ds
\]
is well-defined and equal to the left-hand side (LHS) of (\ref{eq10e}).
Further,%
\begin{align*}
\int K_{1}^{\dagger}\left(  t,x\right)  dF_{\mu}\left(  x\right)  =  &
\frac{1}{2\pi}\int dF_{\mu}\left(  x\right)  \int_{0}^{t}dy\int_{-1}^{1}%
\frac{1}{\iota y}\left\{  \frac{d}{ds} \frac{\exp\left(  \iota ysx\right)
}{r_{0}\left(  ys\right)  } \right\}  ds\\
=  &  \frac{1}{2\pi}\int_{0}^{t}dy\int_{-1}^{1}ds\int\frac{1}{\iota y}\left\{
\frac{d}{ds} \frac{\exp\left(  \iota ysx\right)  }{r_{0}\left(  ys\right)  }
\right\}  dF_{\mu}\left(  x\right) \\
=  &  \frac{1}{2\pi}\int_{0}^{t}dy\int_{-1}^{1}ds\frac{1}{\iota y}\left\{
\frac{d}{ds}\frac{1}{r_{0}\left(  ys\right)  }\right\}  \int\exp\left(  \iota
ysx\right)  dF_{\mu}\left(  x\right) \\
&  +\frac{1}{2\pi}\int_{0}^{t}dy\int_{-1}^{1}ds\frac{1}{\iota y}\frac{1}%
{r_{0}\left(  ys\right)  }\int\left\{  \frac{d}{ds}\exp\left(  \iota
ysx\right)  \right\}  dF_{\mu}\left(  x\right) \\
=  &  \frac{1}{2\pi}\int_{0}^{t}dy\int_{-1}^{1}ds\frac{1}{\iota y}\left\{
\frac{d}{ds}\frac{1}{r_{0}\left(  ys\right)  }\right\}  \int\exp\left(  \iota
ysx\right)  dF_{\mu}\left(  x\right) \\
&  +\frac{1}{2\pi}\int_{0}^{t}dy\int_{-1}^{1}ds\frac{1}{\iota y}\frac{1}%
{r_{0}\left(  ys\right)  }\left\{  \frac{d}{ds}\int\exp\left(  \iota
ysx\right)  dF_{\mu}\left(  x\right)  \right\}  ,
\end{align*}
where we have invoked Fubini's theorem due to (\ref{CondVc1}) and the identity%
\[
\frac{d}{ds}\frac{\exp\left(  \iota ysx\right)  }{r_{0}\left(  ys\right)
}=\left\{  \frac{d}{ds}\frac{1}{r_{0}\left(  ys\right)  }\right\}  e^{\iota
ysx}+\frac{1}{r_{0}\left(  ys\right)  }\left\{  \frac{d}{ds}\exp\left(  \iota
ysx\right)  \right\}
\]
to obtain the second and third equalities, and the condition $\int\left\vert
x\right\vert dF_{\mu}\left(  x\right)  <\infty$ for all $\mu\in U$ and
$\left\vert s\right\vert \leq1$ to assert%
\[
\int\frac{d}{ds}\left\{  \frac{1}{y}\exp\left(  \iota ysx\right)  \right\}
dF_{\mu}\left(  x\right)  =\frac{d}{ds}\left\{  \int\frac{1}{y}\exp\left(
\iota ysx\right)  dF_{\mu}\left(  x\right)  \right\}
\]
to obtain the fourth equality. In other words, we have shown%
\begin{equation}
\int K_{1}^{\dagger}\left(  t,x\right)  dF_{\mu}\left(  x\right)  =\frac
{1}{2\pi}\int_{0}^{t}dy\int_{-1}^{1}ds\left[  \frac{1}{\iota y}\frac{d}%
{ds}\left\{  \frac{1}{r_{0}\left(  ys\right)  }\int\exp\left(  \iota
ysx\right)  dF_{\mu}\left(  x\right)  \right\}  \right]  . \label{eq2a}%
\end{equation}
However, since $\mathcal{F}$ is a Type I location-shift family, we must have
$\hat{F}_{\mu}\left(  t\right)  =\hat{F}_{0}\left(  t\right)  \exp\left(
\iota t\mu\right)  $ and $\hat{F}_{0}=r_{0}$. Therefore, the RHS (i.e.,
right-hand side) of (\ref{eq2a}) is equal to%
\[
\frac{1}{2\pi}\int_{0}^{t}dy\int_{-1}^{1}ds\left\{  \frac{1}{\iota y}\frac
{d}{ds}\exp\left(  \iota ys\mu\right)  \right\}  =\frac{1}{2\pi}\int_{0}%
^{t}dy\int_{-1}^{1}\mu\exp\left(  \iota ys\mu\right)  ds.
\]
Namely,%
\[
\int K_{1}^{\dagger}\left(  t,x\right)  dF_{\mu}\left(  x\right)  =\psi
_{1}\left(  t,\mu\right)  =\frac{1}{2\pi}\int_{0}^{t}dy\int_{-1}^{1}\mu
\exp\left(  \iota ys\mu\right)  ds.
\]
Since $\psi_{1}$ is real, $\psi_{1}\left(  t,\mu\right)  =\int K_{1}\left(
t,x\right)  dF_{\mu}\left(  x\right)  $ has to hold, where $K_{1}\left(
t,x\right)  =\Re\left\{  K_{1}^{\dagger}\left(  t,x\right)  \right\}  $. By
\autoref{ThmPoinNull} and \autoref{lm:Dirichlet}, the pair (\ref{V-a}) is as desired.

\textcolor{blue}{Let $\tilde{r}_{0}^{\prime}\left(  t\right) = \frac{d}{dt} \left(1/r_0\left(t\right)\right)$. Then $\tilde{r}_{0}^{\prime}\left(  \cdot\right)$ is odd since $r_0\left(\cdot\right)$ is even.}
Further,%
\begin{align*}
K_{1}\left(  t,x\right)   &  =\Re\left\{  \frac{1}{2\pi}\int_{0}^{t}%
dy\int_{-1}^{1}\frac{1}{\iota y}\frac{d}{ds}\frac{\exp\left(  \iota
ysx\right)  }{r_{0}\left(  ys\right)  }ds\right\} \\
&  =\Re\left[  \frac{1}{2\pi}\int_{0}^{t}dy\int_{-1}^{1}\frac{\exp\left(
\iota ysx\right)  }{\iota y}\tilde{r}_{0}^{\prime}\left(  ys\right)
ds\right]  +\Re\left\{  \frac{1}{2\pi}\int_{0}^{t}dy\int_{-1}^{1}\frac
{x\exp\left(  \iota ysx\right)  }{r_{0}\left(  ys\right)  }ds\right\} \\
&  =\Re\left[  \frac{1}{2\pi}\int_{0}^{t}dy\int_{-1}^{1}\frac{\exp\left(
\iota ysx\right)  }{\iota y}\tilde{r}_{0}^{\prime}\left(  ys\right)
ds\right]  +\frac{1}{2\pi}\int_{0}^{t}dy\int_{-1}^{1}\frac{x\cos\left(
ysx\right)  }{r_{0}\left(  ys\right)  }ds\\
&  =\frac{1}{2\pi}\int_{0}^{t}dy\int_{-1}^{1}\frac{\sin\left(  ysx\right)
}{y}\tilde{r}_{0}^{\prime}\left(  ys\right)  ds+\frac{1}{2\pi}\int_{0}%
^{t}dy\int_{-1}^{1}\frac{x\cos\left(  ysx\right)  }{r_{0}\left(  ys\right)
}ds\\
&  =\frac{1}{2\pi}\int_{0}^{1}dy\int_{-1}^{1}\left[  \frac{\sin\left(
ytsx\right)  }{y}\left\{  \frac{d}{ds}\frac{1}{r_{0}\left(  tys\right)
}\right\}  +\frac{tx\cos\left(  tysx\right)  }{r_{0}\left(  tys\right)
}\right]  ds.
\end{align*}
\qed

\subsection{Proof of \autoref{ThmTypeI-V}}

The proof is split into three parts in three subsections: bounding the
variance of the error $e_{m}\left(  t\right)  =\hat{\varphi}_{m}\left(
t,\mathbf{z}\right)  -\varphi_{m}\left(  t,\boldsymbol{\mu}\right)  $,
deriving a uniform bound on $\left\vert e_{m}\left(  t\right)  \right\vert $
and then the consistency of $\hat{\varphi}_{m}\left(  t,\mathbf{z}\right)  $,
and specializing these results to Gaussian family.

\subsubsection{Upper bound on the variance of error}

First, we derive an upper bound for $\mathbb{V}\left\{  \hat{\varphi}%
_{m}\left(  t,\mathbf{z}\right)  \right\}  $, and we start with $\mathbb{V}%
\left\{  \hat{\varphi}_{1,m}\left(  t,\mathbf{z}\right)  \right\}  $, where
$\hat{\varphi}_{1,m}\left(  t,\mathbf{z}\right)  =m^{-1}\sum_{i=1}^{m}%
K_{1}\left(  t,z_{i}\right)  $ and\ $\varphi_{1,m}\left(  t,\boldsymbol{\mu
}\right)  =\mathbb{E}\left\{  \hat{\varphi}_{1,m}\left(  t,\mathbf{z}\right)
\right\}  $. Set $\tilde{r}_{0}\left(  tys\right)  =y^{-1}\partial_{s}\left\{
1/r_{0}\left(  tys\right)  \right\}  $. Then
\begin{align*}
K_{1}\left(  t,x\right)   &  =\frac{1}{2\pi}\int_{0}^{1}dy\int_{-1}%
^{1}\left\{  \sin\left(  ytsx\right)  \tilde{r}_{0}\left(  tys\right)
+\frac{tx\cos\left(  tysx\right)  }{r_{0}\left(  tys\right)  }\right\}  ds\\
&
=\textcolor{blue}{\frac{2}{2\pi}\int_{0}^{1}dy\int_{0}^{1}\left\{ \sin\left( ytsx\right) \tilde{r}_{0}\left( tys\right) +\frac{tx\cos\left( tysx\right) }{r_{0}\left( tys\right) }\right\} ds},
\end{align*}
\textcolor{blue}{
since both $\sin\left(  \cdot\right)  \tilde{r}_{0}\left(  \cdot\right)  $ and
$\cos\left(  \cdot\right)  /r_{0}\left(  \cdot\right)  $ are even functions.
}Define%
\[
\left\{
\begin{array}
[c]{l}%
\tilde{S}_{1,m,0}\left(  t,y\right)  =m^{-1}\sum_{i=1}^{m}\sin\left(
ytz_{i}\right)  \text{ \ and }\ \tilde{S}_{1,m,1}\left(  t,y\right)
=m^{-1}\sum_{i=1}^{m}z_{i}\cos\left(  tyz_{i}\right) \\
\textcolor{blue}{\tilde{\Delta}_{1,m,0}\left( x,y\right) =\tilde{S}_{1,m,0}\left( x,y\right) -\mathbb{E}\left\{ \tilde{S}_{1,m,0}\left( x,y\right) \right\}}\\
\textcolor{blue}{\tilde{\Delta}_{1,m,1}\left( x,y\right) =\tilde{S}_{1,m,1}\left( x,y\right) -\mathbb{E}\left\{ \tilde{S}_{1,m,1}\left( x,y\right) \right\}}
\end{array}
\right.  .
\]
\textcolor{blue}{
Then with $e_{1,m}\left(  t\right)  =\hat{\varphi}_{1,m}\left(  t,\mathbf{z}\right)  -\varphi_{1,m}\left(  t,\boldsymbol{\mu}\right)  $, we have\begin{equation}
e_{1,m}\left(  t\right)    =\frac{2}{2\pi}\int_{0}^{1}dy\int_{0}^{1}\left[
\tilde{r}_{0}\left(  tys\right)  \tilde{\Delta}_{1,m,0}\left(  ts,y\right)
+\frac{t}{r_{0}\left(  tys\right)  }\tilde{\Delta}_{1,m,1}\left(  ts,y\right)
\right]  ds.\label{phi1OneSide}\end{equation}
}

By the inequality $\mathbb{V}\left(  X\right)  \leq\mathbb{E}\left(
X^{2}\right)  $ for any random variable $X$, $\mathbb{V}\left\{  \tilde
{S}_{1,m,0}\left(  ts,y\right)  \right\}  \leq m^{-1}$ and%
\begin{equation}
\mathbb{V}\left\{  \tilde{S}_{1,m,1}\left(  ts,y\right)  \right\}  \leq
\frac{1}{m^{2}}\sum_{i=1}^{m}\mathbb{E}\left(  z_{i}^{2}\right)  =\frac
{1}{m^{2}}\sum_{i=1}^{m}\left(  \sigma_{i}^{2}+\mu_{i}^{2}\right)  ,
\label{eq14e}%
\end{equation}
where $\sigma_{i}^{2}$ is the variance of $z_{i}$ and $\mu_{i}$ the mean of
$z_{i}$. Further, $\mathbb{V}\left\{  \hat{\varphi}_{1,m}\left(
t,\mathbf{z}\right)  \right\}  \leq2\tilde{I}_{1,m,0}+2\tilde{I}_{1,m,1}$,
where%
\[
\tilde{I}_{1,m,0}=\mathbb{E}\left(  \left\{  \frac{1}{2\pi}\int_{0}^{1}%
dy\int_{-1}^{1}\tilde{r}_{0}\left(  tys\right)
\textcolor{blue}{\tilde{\Delta}_{1,m,0}\left(  ts,y\right)} ds\right\}
^{2}\right)
\]
and%
\[
\tilde{I}_{1,m,1}=\mathbb{E}\left(  \left\{  \frac{1}{2\pi}\int_{0}^{1}%
dy\int_{-1}^{1}\frac{t}{r_{0}\left(  tys\right)  }
\textcolor{blue}{\tilde{\Delta}_{1,m,1}\left(  ts,y\right)} ds\right\}
^{2}\right)  .
\]
Set
\begin{equation}
\bar{r}_{0}\left(  t\right)  =\sup_{\left(  y,s\right)  \in\left[  0,1\right]
\times\left[  -1,1\right]  }%
\textcolor{blue}{\left\vert\tilde{r}_{0}\left( tys\right)\right\vert} \text{
\ and \ }\check{r}_{0}\left(  t\right)  =\sup_{\left(  y,s\right)  \in\left[
0,1\right]  \times\left[  -1,1\right]  }\frac{1}{r_{0}\left(  tys\right)  }.
\label{eq13d}%
\end{equation}
Then, applying Holder's inequality to $\int_{-1}^{1}\left\vert \tilde{\Delta
}_{1,m,0}\left(  ts,y\right)  \right\vert ds$ and $\int_{0}^{1}dy\left(
\int_{-1}^{1}\left\vert \tilde{\Delta}_{1,m,0}\left(  ts,y\right)  \right\vert
^{2}ds\right)  ^{1/2}$ yields
\begin{align*}
\tilde{I}_{1,m,0}  &  \leq\frac{\bar{r}_{0}^{2}\left(  t\right)  }{4\pi^{2}%
}\mathbb{E}\left(  \left\{  \int_{0}^{1}dy\int_{-1}^{1}\left\vert
\tilde{\Delta}_{1,m,0}\left(  ts,y\right)  \right\vert ds\right\}  ^{2}\right)
\\
&  \leq\frac{\bar{r}_{0}^{2}\left(  t\right)  }{4\pi^{2}}\mathbb{E}\left[
\left(  2^{1/2}\int_{0}^{1}dy\left(  \int_{-1}^{1}\left\vert \tilde{\Delta
}_{1,m,0}\left(  ts,y\right)  \right\vert ^{2}ds\right)  ^{1/2}\right)
^{2}\right] \\
&  \leq\frac{\bar{r}_{0}^{2}\left(  t\right)  }{2\pi^{2}}\mathbb{E}\left[
\int_{0}^{1}dy\left(  \int_{-1}^{1}\left\vert \tilde{\Delta}_{1,m,0}\left(
ts,y\right)  \right\vert ^{2}ds\right)  \right] \\
&  =\frac{\bar{r}_{0}^{2}\left(  t\right)  }{2\pi^{2}}\int_{0}^{1}dy\int%
_{-1}^{1}\mathbb{V}\left\{  \tilde{S}_{1,m,0}\left(  ts,y\right)  \right\}
ds\leq\frac{\bar{r}_{0}^{2}\left(  t\right)  }{2\pi^{2}}\frac{2}{m}=\frac
{\bar{r}_{0}^{2}\left(  t\right)  }{\pi^{2}m}.
\end{align*}
By similar arguments given to derive $\tilde{I}_{1,m,0}\leq\frac{\bar{r}%
_{0}^{2}\left(  t\right)  }{\pi^{2}m}$, we have%
\begin{align*}
\tilde{I}_{1,m,1}  &  \leq\frac{t^{2}\check{r}_{0}^{2}\left(  t\right)  }%
{4\pi^{2}}\mathbb{E}\left(  \left\{  \int_{0}^{1}dy\int_{-1}^{1}
\textcolor{blue}{\tilde{\Delta}_{1,m,1}\left( ts,y\right)} ds\right\}
^{2}\right) \\
&  \leq\frac{t^{2}\check{r}_{0}^{2}\left(  t\right)  }{2\pi^{2}}\int_{0}%
^{1}dy\int_{-1}^{1}\mathbb{V}\left(  \tilde{S}_{1,m,1}\left(  ts,y\right)
\right)  ds\leq\frac{t^{2}\check{r}_{0}^{2}\left(  t\right)  }{\pi^{2}}%
\frac{1}{m^{2}}\sum_{i=1}^{m}\left(  \sigma_{i}^{2}+\mu_{i}^{2}\right)  .
\end{align*}
Therefore,%
\begin{equation}
\mathbb{V}\left\{  \hat{\varphi}_{1,m}\left(  t,\mathbf{z}\right)  \right\}
\leq\frac{2}{\pi^{2}m}\left[  \bar{r}_{0}^{2}\left(  t\right)  +t^{2}\check
{r}_{0}^{2}\left(  t\right)  \tilde{D}_{m}\right]  \text{ \ with \ }\tilde
{D}_{m}=m^{-1}\sum_{i=1}^{m}\left(  \sigma_{i}^{2}+\mu_{i}^{2}\right)  .
\label{eq13f}%
\end{equation}
When $\sigma_{i}^{2}=\sigma^{2}$ for all $i\in\left\{  1,\ldots,m\right\}  $,
$\tilde{D}_{m}$ becomes $D_{m}=\sigma^{2}+m^{-1}\sum_{i=1}^{m}\mu_{i}^{2}$.

Recall $K\left(  t,x\right)  =2^{-1}-K_{1}\left(  t,x\right)  -2^{-1}%
K_{1,0}\left(  t,x;0\right)  $ and for $\tau\in U$%
\[
\hat{\varphi}_{1,0,m}\left(  t,\mathbf{z};\tau\right)  =m^{-1}\sum_{i=1}%
^{m}K_{1,0}\left(  t,z_{i};\tau\right)  \text{ and \ }\varphi_{1,0,m}\left(
t,\boldsymbol{\mu};\tau\right)  =\mathbb{E}\left\{  \hat{\varphi}%
_{1,0,m}\left(  t,\mathbf{z};\tau\right)  \right\}  .
\]
By Theorem 2 of \cite{Chen:2018a} \textcolor{blue}{or (\ref{eq3e})},%
\begin{equation}
\mathbb{V}\left\{  \hat{\varphi}_{1,0,m}\left(  t,\mathbf{z};0\right)
\right\}  \leq m^{-1}\left\Vert \omega\right\Vert _{\infty}^{2}g^{2}\left(
t,0\right)  . \label{eq14d}%
\end{equation}
where $g\left(  t,\mu\right)  =\int_{\left[  -1,1\right]  }r_{\mu}^{-1}\left(
ts\right)  ds$. So, combining (\ref{eq13f}) with (\ref{eq14d}) gives
\begin{align}
\mathbb{V}\left\{  \hat{\varphi}_{m}\left(  t,\mathbf{z}\right)  \right\}   &
\leq2\mathbb{V}\left\{  \hat{\varphi}_{1,m}\left(  t,\mathbf{z}\right)
\right\}  +2\mathbb{V}\left\{  2^{-1}\hat{\varphi}_{1,0,m}\left(
t,\mathbf{z};0\right)  \right\} \nonumber\\
&  \leq\frac{4}{\pi^{2}m}\left[  \bar{r}_{0}^{2}\left(  t\right)  +t^{2}%
\check{r}_{0}^{2}\left(  t\right)  \tilde{D}_{m}\right]  +\frac{\left\Vert
\omega\right\Vert _{\infty}^{2}}{2m}g^{2}\left(  t,0\right)  . \label{eq15e}%
\end{align}

\subsubsection{\textcolor{blue}{Uniform bound on error and consistency}}

Second, we derive a uniform bound for $e_{m}\left(  t\right)  =\hat{\varphi
}_{m}\left(  t,\mathbf{z}\right)  -\varphi_{m}\left(  t,\boldsymbol{\mu
}\right)  $, for which we will use \autoref{ConcentrationEmpiricalProcesses}.
Let $\left\{  \tau_{m}\right\}  _{m\geq1},\left\{  \gamma_{m}\right\}
_{m\geq1},\left\{  \rho_{m}\right\}  _{m\geq1}$ and $\left\{  \tilde{\rho}%
_{m}\right\}  _{m\geq1}$ be four positive sequences to be determined later. In
the rest of the proof, $q,\gamma,\vartheta,\vartheta_{1}$ and $\vartheta_{2}$
are positive constants and $\vartheta^{\prime}$ and $\vartheta^{^{\prime
\prime}}$ are non-negative constants whose ranges will be determined later.
Set $D=\left[  0,\tau_{m}\right]  \times\left[  0,1\right]  \times\left[
0,1\right]  $. The continuity of $\tilde{\Delta}_{1,m,0}\left(  \cdot
,\cdot\right)  $ and that of $\tilde{\Delta}_{1,m,1}\left(  \cdot
,\cdot\right)  $ respectively imply%
\begin{equation}
\left\{
\begin{array}
[c]{c}%
\sup_{\left(  t,s,y\right)  \in D}\tilde{\Delta}_{1,m,0}\left(  ts,y\right)
=\sup_{t\in\left[  0,\tau_{m}\right]  }\tilde{\Delta}_{1,m,0}\left(
t,1\right) \\
\sup_{\left(  t,s,y\right)  \in D}\tilde{\Delta}_{1,m,1}\left(  ts,y\right)
=\sup_{t\in\left[  0,\tau_{m}\right]  }\tilde{\Delta}_{1,m,1}\left(
t,1\right)
\end{array}
\right.  . \label{TildeDeltaSupEquiv}%
\end{equation}
So, it suffices to bound the tail probabilities of $\sup_{t\in\left[
0,\tau_{m}\right]  }\tilde{\Delta}_{1,m,0}\left(  t,1\right)  $ and
$\sup_{t\in\left[  0,\tau_{m}\right]  }\tilde{\Delta}_{1,m,1}\left(
t,1\right)  $ respectively, as done below:

\begin{itemize}
\item \textbf{Case 3a} \textquotedblleft Bound the tail probability of
$\sup\limits_{t\in\left[  0,\tau_{m}\right]  }\tilde{\Delta}_{1,m,1}\left(
t,1\right)  $ when $\boldsymbol{\mu}\in\mathcal{B}_{2,m}\left(  \tilde{\rho
}_{m}\right)  $ and $\mathbb{E}\left[  \left\vert X_{\left(  0\right)
}\right\vert ^{2}\right]  <\infty$\textquotedblright: we can use Case 3 of
\autoref{ConcentrationEmpiricalProcesses} (where $f\left(  x,v,y,1,0,c_{0}%
\right)  =x\cos\left(  vx-c_{0}\right)  $) by setting%
\begin{equation}
\left\{
\begin{array}
[c]{l}%
c_{0}=0;c_{1}=\tau_{m}\rightarrow\infty;c_{2}=c_{3}=0;\Delta_{1}%
=m^{-\vartheta}\\
\kappa_{1}=\gamma_{m}m^{-\vartheta_{1}};\kappa_{2}=\gamma_{m}^{-1}%
m^{-\vartheta_{1}}\ln c_{1};\Delta_{1}^{-1}\kappa_{2}/R_{2}\left(  \tilde
{\rho}_{m}\right)  \rightarrow\infty
\end{array}
\right.  , \label{ParSetting4}%
\end{equation}
where $R_{2}\left(  \tilde{\rho}_{m}\right)  =4\mathbb{E}\left[  X_{\left(
0\right)  }^{2}\right]  +4\tilde{\rho}_{m}$. Recall $C_{m,\boldsymbol{\mu
},F_{0}}=\left[  1-2\Pr\left\{  X_{\left(  0\right)  }>\kappa-\left\Vert
\boldsymbol{\mu}\right\Vert _{\infty}\right\}  \right]  ^{m}$ and
(\ref{probExpressCase3}) as%
\begin{align*}
\hat{p}_{3,m}\left(  \Xi_{3}\right)   &  =2\left(  c_{1}\Delta_{1}%
^{-1}+1\right)  \exp\left(  -2^{-1}m\kappa^{-2}\left(  \kappa_{1}-\kappa
_{2}\right)  ^{2}\right) \\
&  \text{ \ \ \ }+m^{-1}\left[  2^{-1}\Delta_{1}^{-1}\kappa_{2}-2^{-1}%
R_{2}\left(  \tilde{\rho}_{m}\right)  \right]  ^{-2}\mathbb{V}\left[
X_{\left(  0\right)  }^{2}\right]  +1-C_{m,\boldsymbol{\mu},F_{0}},
\end{align*}
where $X_{\left(  0\right)  }$ has CDF $F_{0}$ and $\left\Vert \boldsymbol{\mu
}\right\Vert _{\infty}=\max_{1\leq i\leq m}\left\vert \mu_{i}\right\vert $.
Then, for all $m$ large enough, $\hat{p}_{3,m}\left(  \Xi_{3}\right)  \leq
\hat{p}_{3,m}^{\ast}\left(  c_{1},\Delta_{1},\kappa,\vartheta_{1},\gamma
_{m}\right)  $, where
\begin{align}
\hat{p}_{3,m}^{\ast}\left(  c_{1},\Delta_{1},\kappa,\vartheta_{1},\gamma
_{m}\right)   &  =C\left[  c_{1}\Delta_{1}^{-1}\exp\left(  -2^{-1}\kappa
^{-2}m^{1-2\vartheta_{1}}\left(  \gamma_{m}^{2}-2\ln\tau_{m}\right)  \right)
\right. \nonumber\\
&  \text{ \ \ \ \ }\left.  +m^{-1}\Delta_{1}^{2}\kappa_{2}^{-2}%
+1-C_{m,\boldsymbol{\mu},F_{0}}\right]  . \label{ProbTailCase3AX}%
\end{align}
Namely, for all $m$ large enough,%
\begin{equation}
\Pr\left\{  \sup\nolimits_{t\in\left[  0,\tau_{m}\right]  }\tilde{\Delta
}_{1,m,1}\left(  t,1\right)  \geq\frac{\gamma_{m}}{m^{\vartheta_{1}}}\right\}
\leq\hat{p}_{3,m}^{\ast}\left(  c_{1},\Delta_{1},\kappa,\vartheta_{1}%
,\gamma_{m}\right)  . \label{supBndCase3AX}%
\end{equation}
However, when

\item \textbf{Case 2b} \textquotedblleft Bound the tail probability of
$\sup\limits_{t\in\left[  0,\tau_{m}\right]  }\tilde{\Delta}_{1,m,0}\left(
t,1\right)  $ when $\boldsymbol{\mu}\in\mathcal{B}_{1,m}\left(  \rho
_{m}\right)  $ and $\mathbb{E}\left[  \left\vert X_{\left(  0\right)
}\right\vert ^{2}\right]  <\infty$\textquotedblright: we can use Case 2 of
\autoref{ConcentrationEmpiricalProcesses} (where $f\left(  x,v,y,0,1,c_{0}%
\right)  =\cos\left(  v\left(  x-c_{2}\right)  -c_{0}\right)  $) by setting%
\begin{equation}
\left\{
\begin{array}
[c]{l}%
c_{0}=2^{-1}\pi,c_{1}=\tau_{m}\rightarrow\infty,c_{2}=c_{3}=0;\Delta
_{1}=m^{-\vartheta}\\
\kappa_{1}=\gamma_{m}m^{-1/2},\kappa_{2}=\gamma_{m}^{-1}m^{-1/2}\ln
c_{1};\gamma_{m}^{2}>\ln\tau_{m};\Delta_{1}^{-1}\kappa_{2}/R_{1}\left(
\rho_{m}\right)  \rightarrow\infty
\end{array}
\right.  . \label{SettingCase2b}%
\end{equation}
where $R_{1}\left(  \rho_{m}\right)  $ in (\ref{TheShiftRho}) becomes
$R_{1}\left(  \rho_{m}\right)  =2\mathbb{E}\left[  \left\vert X_{(0)}%
\right\vert \right]  +2\rho_{m}$. Then $f\left(  x,v,y,0,1,2^{-1}\pi\right)
=\sin\left(  vx\right)  $. We now have a special case of the analysis on
bounding the tail probability of $\sup_{v\in\left[  0,\tau_{m}\right]
}S_{1,m}\left(  v,\tau\right)  $ that has been marked as Case 2a and provided
in the proof of \autoref{ThmIVa}. So, arguments there imply that $\hat
{p}_{2,m}\left(  \Xi_{2}\right)  $ in (\ref{probExpressCase2})\ satisfies, for
all sufficiently large $m$,%
\begin{equation}
\hat{p}_{2,m}\left(  \Xi_{2}\right)  \leq\hat{p}_{2,m}^{\ast}\left(
c_{1},\Delta_{1},\gamma_{m}\right)  =C\left[  c_{1}^{2}\Delta_{1}^{-1}%
\exp\left(  -2^{-1}\gamma_{m}^{2}\right)  +m^{-1}\Delta_{1}^{2}\kappa_{2}%
^{-2}\right]  \label{probCase2Copy}%
\end{equation}
where $\hat{p}_{2,m}^{\ast}\left(  c_{1},\Delta_{1}\right)  $ has been defined
by (\ref{probCase2}) in the proof of \autoref{ThmIVa}. Namely, when $m$ is
large,%
\begin{equation}
\Pr\left\{  \sup\nolimits_{t\in\left[  0,c_{1}\right]  }\tilde{\Delta}%
_{1,m,0}\left(  t,1\right)  \geq\frac{\gamma_{m}}{\sqrt{m}}\right\}  \leq
\hat{p}_{2,m}^{\ast}\left(  c_{1},\Delta_{1},\gamma_{m}\right)  .
\label{supBndSinAX}%
\end{equation}

\end{itemize}

Before proceeding futher, we need to pause to make an important remark as
follows. Let $\overline{\mu_{j,m}}=m^{-1}\sum_{i=1}^{m}\left\vert \mu
_{i}\right\vert ^{j}$. Then the Cauchy-Schwarz inequality implies
$\overline{\mu_{1,m}}^{2}\leq\overline{\mu_{2,m}}$, and so $\overline
{\mu_{2,m}}\leq\tilde{\rho}_{m}$ implies $\overline{\mu_{1,m}}\leq\sqrt
{\tilde{\rho}_{m}}$. Namely, we can set as $\sqrt{\tilde{\rho}_{m}}$ the
$\rho_{m}$ in upper bound in (\ref{L2L1Connection}) for Case 2 in the proof of
\autoref{ConcentrationEmpiricalProcesses}, and set as $R_{1}\left(
\sqrt{\tilde{\rho}_{m}}\right)  $ the $R_{1}\left(  \rho_{m}\right)  $ in
(\ref{BndHprocessRelax}) and (\ref{TheShiftRho}) without affecting any
conclusions for Case 2 of \autoref{ConcentrationEmpiricalProcesses} or of Case
2a in the proof of \autoref{ThmIVa}, and hence without affecting any arguments
presented thus far in this proof when we set $R_{1}\left(  \rho_{m}\right)  $
in (\ref{SettingCase2b}) as $R_{1}\left(  \sqrt{\tilde{\rho}_{m}}\right)  $.
This remark will be used in the rest of the proof, where we take $\rho
_{m}=\sqrt{\tilde{\rho}_{m}}$.

Again, we will simply write $\hat{p}_{3,m}^{\ast}\left(  c_{1},\Delta
_{1},\kappa,\vartheta_{1},\gamma_{m}\right)  $ as $\hat{p}_{3,m}^{\ast}$ when
no confusion arises. Let $\mathcal{\tilde{B}}_{m}\left(  \tilde{\rho}%
_{m}\right)  =\mathcal{B}_{1,m}\left(  \sqrt{\tilde{\rho}_{m}}\right)
\cap\mathcal{B}_{2,m}\left(  \tilde{\rho}_{m}\right)  $. Then,
(\ref{phi1OneSide}), (\ref{TildeDeltaSupEquiv}), (\ref{supBndSinAX}) and
(\ref{supBndCase3AX}) implies that, for all $m$ large enough, with probability
at least $1-\hat{p}_{2,m}^{\ast}-\hat{p}_{3,m}^{\ast}$,%
\begin{align*}
\sup_{\boldsymbol{\mu}\in\mathcal{\tilde{B}}_{m}\left(  \tilde{\rho}%
_{m}\right)  }\sup_{t\in\left[  0,\tau_{m}\right]  }\left\vert e_{1,m}\left(
t\right)  \right\vert  &  \leq\frac{2}{2\pi}\sup_{\boldsymbol{\mu}%
\in\mathcal{\tilde{B}}_{m}\left(  \tilde{\rho}_{m}\right)  }\sup_{t\in\left[
0,\tau_{m}\right]  }\int_{0}^{1}dy\int_{0}^{1}\left\vert \tilde{r}_{0}\left(
tys\right)  \right\vert \sup_{\left(  t,s,y\right)  \in D}\left\vert
\tilde{\Delta}_{1,m,0}\left(  ts,y\right)  \right\vert ds\\
&  \text{ \ \ }+\frac{2}{2\pi}\sup_{\boldsymbol{\mu}\in\mathcal{\tilde{B}}%
_{m}\left(  \tilde{\rho}_{m}\right)  }\sup_{t\in\left[  0,\tau_{m}\right]
}\int_{0}^{1}dy\int_{0}^{1}\frac{t}{r_{0}\left(  tys\right)  }\sup_{\left(
t,s,y\right)  \in D}\left\vert \tilde{\Delta}_{1,m,1}\left(  ts,y\right)
\right\vert ds\\
&  \leq \Upsilon_{1}\left(  \tau_{m}m^{-1/2},\tau_{m}\right)  +\tau_{m}%
\Upsilon_{2}\left(  \tau_{m}m^{-\vartheta_{1}},\tau_{m}\right)
\end{align*}
where we recall $\tilde{r}_{0}\left(  tys\right)  =y^{-1}\partial_{s}\left\{
1/r_{0}\left(  tys\right)  \right\}  $ and that $r_{0}\left(  \cdot\right)  $
is a positive function, and have set%
\begin{equation}
\left\{
\begin{array}
[c]{l}%
\Upsilon_{1}\left(  \lambda,\tilde{\tau}\right)  =\dfrac{\lambda}{\pi}%
\sup_{t\in\left[  0,\tilde{\tau}\right]  }%
{\displaystyle\int_{0}^{1}}
dy%
{\displaystyle\int_{0}^{1}}
\left\vert \tilde{r}_{0}\left(  tys\right)  \right\vert ds\smallskip\\
\Upsilon_{2}\left(  \lambda,\tilde{\tau}\right)  =\dfrac{\lambda}{\pi}%
\sup_{t\in\left[  0,\tilde{\tau}\right]  }%
{\displaystyle\int_{0}^{1}}
dy%
{\displaystyle\int_{0}^{1}}
\dfrac{ds}{r_{0}\left(  tys\right)  }%
\end{array}
\right.  \text{ for }\lambda,\tilde{\tau}\geq0. \label{upsilon1and2X}%
\end{equation}
Namely, for all $m$ sufficiently large, with probability at least $1-\hat
{p}_{2,m}^{\ast}-\hat{p}_{3,m}^{\ast}$,%
\begin{equation}
\sup_{\boldsymbol{\mu}\in\mathcal{\tilde{B}}_{m}\left(  \tilde{\rho}%
_{m}\right)  }\sup_{t\in\left[  0,\tau_{m}\right]  }\left\vert e_{1,m}\left(
t\right)  \right\vert \leq \Upsilon_{1}\left(  \gamma_{m}m^{-1/2},\tau
_{m}\right)  +\tau_{m}\Upsilon_{2}\left(  \gamma_{m}m^{-\vartheta_{1}}%
,\tau_{m}\right)  . \label{SupBoundPhi1mX}%
\end{equation}
Recall $e_{1,0,m}\left(  t,\mu^{\prime}\right)  =\hat{\varphi}_{1,0,m}\left(
t,\mathbf{z};\mu^{\prime}\right)  -\varphi_{1,0,m}\left(  t,\boldsymbol{\mu
};\mu^{\prime}\right)  $. However, analogous to (\ref{PointNullPhi10Sup}), the
proof of \autoref{ThmIVa} already provides, for all large $m$,%
\begin{equation}
\sup_{\boldsymbol{\mu}\in\mathcal{\tilde{B}}_{m}\left(  \tilde{\rho}%
_{m}\right)  }\sup_{t\in\left[  0,\tau_{m}\right]  }\left\vert e_{1,0,m}%
\left(  t,0\right)  \right\vert \leq\left\Vert \omega\right\Vert _{\infty
}\Upsilon\left(  \gamma_{m}m^{-1/2},\tau_{m}\right)
\label{BndErrorPhi10TauAsZero}%
\end{equation}
with probability at least $1-\hat{p}_{2,m}^{\ast}$, where $\hat{p}_{2,m}%
^{\ast}$ is defined by (\ref{probCase2}) and $\Upsilon\left(  \lambda
,\tilde{\tau}\right)  $ is defined by (\ref{UpsilonA}). Since%
\[
\hat{\varphi}_{m}\left(  t,\mathbf{z}\right)  -\varphi_{m}\left(
t,\boldsymbol{\mu}\right)  =\hat{\varphi}_{1,m}\left(  t,\mathbf{z}\right)
-\varphi_{1,m}\left(  t,\boldsymbol{\mu}\right)  +2^{-1}\hat{\varphi}%
_{1,0,m}\left(  t,\mathbf{z};0\right)  -2^{-1}\mathbb{E}\left\{
\varphi_{1,0,m}\left(  t,\boldsymbol{\mu};0\right)  \right\}  ,
\]
(\ref{BndErrorPhi10TauAsZero}) and (\ref{SupBoundPhi1mX}) together imply that,
for all $m$ large enough,
\begin{equation}
\sup_{\boldsymbol{\mu}\in\mathcal{\tilde{B}}_{m}\left(  \tilde{\rho}%
_{m}\right)  }\sup_{t\in\left[  0,\tau_{m}\right]  }\left\vert e_{m}\left(
t\right)  \right\vert \leq\hat{\Upsilon}\left(  \tau_{m},m,\vartheta
_{1}\right)  \label{BoundSupWholeErrorX}%
\end{equation}
with probability at least $1-p_{m}^{\ast}\left(  c_{1},\Delta_{1}%
,\kappa,\vartheta_{1},\gamma_{m}\right)  $, where%
\begin{equation}
\left\{
\begin{array}
[c]{l}%
\hat{\Upsilon}\left(  \tau_{m},\gamma_{m},\vartheta_{1}\right)  =2^{-1}%
\left\Vert \omega\right\Vert _{\infty}\Upsilon\left(  \gamma_{m}m^{-1/2}%
,\tau_{m}\right)  +\Upsilon_{1}\left(  \gamma_{m}m^{-1/2},\tau_{m}\right)
+\tau_{m}\Upsilon_{2}\left(  \gamma_{m}m^{-\vartheta_{1}},\tau_{m}\right) \\
p_{m}^{\ast}\left(  c_{1},\Delta_{1},\kappa,\vartheta_{1},\gamma_{m}\right)
=2\hat{p}_{2,m}^{\ast}\left(  c_{1},\Delta_{1},\gamma_{m}\right)  +\hat
{p}_{3,m}^{\ast}\left(  c_{1},\Delta_{1},\kappa,\vartheta_{1},\gamma
_{m}\right)  .
\end{array}
\right.  \label{UpperBndTotalErrorAndProbOneSide}%
\end{equation}
We will write $p_{m}^{\ast}\left(  c_{1},\Delta_{1},\kappa,\vartheta
_{1},\gamma_{m}\right)  $ as $p_{m}^{\ast}\ $when no confusion arises.

Third, we show consistency of $\hat{\varphi}_{m}\left(  t,\mathbf{z}\right)
$. In the decomposition%
\[
\left\vert \pi_{1,m}^{-1}\hat{\varphi}_{m}\left(  t,\mathbf{z}\right)
-1\right\vert \leq\pi_{1,m}^{-1}\left\vert \hat{\varphi}_{m}\left(
t,\mathbf{z}\right)  -\varphi_{m}\left(  t,\boldsymbol{\mu}\right)
\right\vert +\left\vert \pi_{1,m}^{-1}\varphi_{m}\left(  t_{m},\boldsymbol{\mu
}\right)  -1\right\vert ,
\]
\autoref{SpeedOneSidedNull}, whose $\tilde{\psi}_{1,0}$ has parameter $\sigma>0$ and
whose $\tilde{\psi}_{1,0}$ with $\sigma=1$ reduces to $\psi_{1,0}$,
asserts that, when $t\tilde{u}_{m}\geq2$,
\[
\left\vert \pi_{1,m}^{-1}\varphi_{m}\left(  t,\boldsymbol{\mu}\right)
-1\right\vert \leq\frac{4+2\left(  \left\Vert \omega\right\Vert _{\infty
}+\left\Vert \omega\right\Vert _{\mathrm{TV}}\right)  }{t\tilde{u}_{m}%
\pi_{1,m}}\leq\frac{C}{t\tilde{u}_{m}\pi_{1,m}},
\]
where $\tilde{u}_{m}=\min_{\left\{  j:\mu_{j}\neq0\right\}  }\left\vert
\mu_{j}\right\vert $. Therefore,%
\begin{equation}
\tau_{m}^{-1}\left(  1+\tilde{u}_{m}^{-1}\right)  =o\left(  \pi_{1,m}\right)
\Longrightarrow\pi_{1,m}^{-1}\varphi_{m}\left(  \tau_{m},\boldsymbol{\mu
}\right)  \rightarrow1. \label{eqdx4aA1}%
\end{equation}
In addition, (\ref{BoundSupWholeErrorX}) and
(\ref{UpperBndTotalErrorAndProbOneSide}) imply that
\begin{equation}
\sup_{\boldsymbol{\mu}\in\mathcal{\tilde{B}}_{m}\left(  \tilde{\rho}%
_{m}\right)  }\sup_{t\in\left[  0,\tau_{m}\right]  }\left\vert e_{m}\left(
t\right)  \right\vert \rightsquigarrow0\text{ if }\Upsilon\left(  \tau
_{m},m,\vartheta_{1}\right)  =o\left(  1\right)
\label{OneSideConsistencyCondB}%
\end{equation}
Recall $p_{m}^{\ast}$ defined by (\ref{UpperBndTotalErrorAndProbOneSide}). We
need one more ingradient, i.e., $p_{m}^{\ast}=o\left(  1\right)  $. Recall the
definitions of $\hat{p}_{2,m}^{\ast}$ in (\ref{probCase2Copy}) or
(\ref{probCase2}) and $\hat{p}_{3,m}^{\ast}$ in (\ref{ProbTailCase3AX}), the
restrictions $\Delta_{1}^{-1}\kappa_{2}/R_{1}\left(  \rho_{m}\right)
\rightarrow\infty$ in (\ref{SettingCase2b})\ and $\Delta_{1}^{-1}\kappa
_{2}/R_{2}\left(  \tilde{\rho}_{m}\right)  \rightarrow\infty$ in
(\ref{ParSetting4}), and%
\[
C_{m,\boldsymbol{\mu},F_{0}}=\left[  1-2\Pr\left\{  X_{\left(  0\right)
}>\kappa-\left\Vert \boldsymbol{\mu}\right\Vert _{\infty}\right\}  \right]
^{m}.
\]
By the remark right after (\ref{supBndSinAX}), we see $R_{1}\left(  \rho
_{m}\right)  $ in (\ref{SettingCase2b}) has been set to be $R_{1}\left(
\sqrt{\tilde{\rho}_{m}}\right)  $, and it suffices to ensure $\Delta_{1}%
^{-1}\kappa_{2}/R_{2}\left(  \tilde{\rho}_{m}\right)  \rightarrow\infty$ for
the aforementioned restrictions to be met. Then, with $m^{-1}\Delta_{1}%
^{2}\kappa_{2}^{-2}=m^{-1-2\vartheta+2\vartheta_{1}}\gamma_{m}^{2}\ln^{-2}%
\tau_{m}$ via setting (\ref{ParSetting4}), we obtain
\begin{equation}
\left\{
\begin{array}
[c]{l}%
m^{-2\vartheta}\gamma_{m}^{2}\ln^{-2}\tau_{m}=o\left(  1\right)  ;\tau_{m}%
^{2}m^{\vartheta}\exp\left(  -2^{-1}\gamma_{m}^{2}\right)  =o\left(  1\right)
\\
\left(  \gamma_{m}R_{1}\left(  \sqrt{\tilde{\rho}_{m}}\right)  \right)
^{-1}m^{\vartheta-1/2}\ln\tau_{m}\rightarrow\infty;\left(  \gamma_{m}%
R_{2}\left(  \tilde{\rho}_{m}\right)  \right)  ^{-1}m^{\vartheta-\vartheta
_{1}}\ln\tau_{m}\rightarrow\infty\\
m^{-1-2\vartheta+2\vartheta_{1}}\gamma_{m}^{2}\ln^{-2}\tau_{m}=o\left(
1\right) \\
\tau_{m}m^{\vartheta}\exp\left(  -2^{-1}\kappa^{-2}m^{1-2\vartheta_{1}}\left(
\gamma_{m}^{2}-2\ln\tau_{m}\right)  \right)  =o\left(  1\right) \\
\left[  1-2\Pr\left\{  X_{\left(  0\right)  }>\kappa-\left\Vert
\boldsymbol{\mu}\right\Vert _{\infty}\right\}  \right]  ^{m}=1+o\left(
1\right)
\end{array}
\right.  \Rightarrow p_{m}^{\ast}=o\left(  1\right)  .
\label{OneSideTotalErrorProb}%
\end{equation}
So, (\ref{eqdx4aA1}), (\ref{OneSideConsistencyCondB}) and
(\ref{OneSideTotalErrorProb}) together are equivalent to $\left\vert \pi
_{1,m}^{-1}\hat{\varphi}_{m}\left(  \tau_{m},\mathbf{z}\right)  -1\right\vert
\rightsquigarrow0$, proving consistency.

To achieve uniform consistency, since we already have (\ref{eqdx4aA1}),
(\ref{UpperBndTotalErrorAndProbOneSide}) and (\ref{OneSideTotalErrorProb}), we
can define $\Omega_{m}\left(  \tilde{\rho}_{m}\right)  $ be a subset of
$\mathcal{\tilde{B}}_{m}\left(  \tilde{\rho}_{m}\right)  =\mathcal{B}%
_{1,m}\left(  \sqrt{\tilde{\rho}_{m}}\right)  \cap\mathcal{B}_{2,m}\left(
\tilde{\rho}_{m}\right)  $ for each $m\geq1$. Then a uniform consistency class
for $\hat{\varphi}_{m}\left(  t,\mathbf{z}\right)  $ is%
\[
\mathcal{Q}\left(  \mathcal{F}\right)  =\left\{
\begin{array}
[c]{l}%
m^{-2\vartheta}\gamma_{m}^{2}\ln^{-2}\tau_{m}=o\left(  1\right)  ;\tau_{m}%
^{2}m^{\vartheta}\exp\left(  -2^{-1}\gamma_{m}^{2}\right)  =o\left(  1\right)
\\
\left(  \gamma_{m}R_{1}\left(  \sqrt{\tilde{\rho}_{m}}\right)  \right)
^{-1}m^{\vartheta-1/2}\ln\tau_{m}\rightarrow\infty;\left(  \gamma_{m}%
R_{2}\left(  \tilde{\rho}_{m}\right)  \right)  ^{-1}m^{\vartheta-\vartheta
_{1}}\ln\tau_{m}\rightarrow\infty\\
m^{-1-2\vartheta+2\vartheta_{1}}\gamma_{m}^{2}\ln^{-2}\tau_{m}=o\left(
1\right) \\
\tau_{m}m^{\vartheta}\exp\left(  -2^{-1}\kappa^{-2}m^{1-2\vartheta_{1}}\left(
\gamma_{m}^{2}-2\ln\tau_{m}\right)  \right)  =o\left(  1\right) \\
\left[  1-2\Pr\left\{  X_{\left(  0\right)  }>\kappa-\left\Vert
\boldsymbol{\mu}\right\Vert _{\infty}\right\}  \right]  ^{m}=1+o\left(
1\right) \\
\tau_{m}^{-1}\sup_{\boldsymbol{\mu}\in\Omega_{m}\left(  \tilde{\rho}%
_{m}\right)  }\left(  1+\tilde{u}_{m}^{-1}\right)  \pi_{1,m}^{-1}=o\left(
\pi_{1,m}\right) \\
\left[  \Upsilon\left(  \gamma_{m}m^{-1/2},\tau_{m}\right)  +\Upsilon
_{1}\left(  \gamma_{m}m^{-1/2},\tau_{m}\right)  \right]  \sup_{\boldsymbol{\mu
}\in\Omega_{m}\left(  \tilde{\rho}_{m}\right)  }\pi_{1,m}^{-1}=o\left(
1\right) \\
\tau_{m}\Upsilon_{2}\left(  \gamma_{m}m^{-\vartheta_{1}},\tau_{m}\right)
\sup_{\boldsymbol{\mu}\in\Omega_{m}\left(  \tilde{\rho}_{m}\right)  }\pi
_{1,m}^{-1}=o\left(  1\right)
\end{array}
\right\}  .
\]

\subsubsection{\textcolor{blue}{Specialization to Gaussian family}}

Finally, we specialize the above results to Gaussian family with scale
parameter $\sigma>0$. There are three key steps to achieve this. First, we take
look at the variance. Recall $r_{\mu}^{-1}\left(  tys\right)  =\exp\left(
2^{-1}t^{2}y^{2}s^{2}\sigma^{2}\right)  $ and%
\begin{equation}
\tilde{r}_{0}\left(  tys\right)  =\frac{1}{y}\frac{d}{ds}\frac{1}{r_{0}\left(
tys\right)  }=syt^{2}\sigma^{2}\exp\left(  2^{-1}y^{2}t^{2}s^{2}\sigma
^{2}\right)  . \label{bndPartialRecipCF}%
\end{equation}
We see from definition (\ref{eq13d}) that $\bar{r}_{0}\left(  t\right)  \leq
t^{2}\sigma^{2}\exp\left(  2^{-1}t^{2}\sigma^{2}\right)  $ and $\check{r}%
_{0}\left(  t\right)  \leq\exp\left(  2^{-1}t^{2}\sigma^{2}\right)  $.
Therefore, from (\ref{eq15e}) we obtain%
\begin{align*}
\mathbb{V}\left\{  \hat{\varphi}_{m}\left(  t,\mathbf{z}\right)  \right\}   &
\leq\frac{4}{\pi^{2}m}\left[  t^{4}\sigma^{4}\exp\left(  t^{2}\sigma
^{2}\right)  +t^{2}\exp\left(  t^{2}\sigma^{2}\right)  \tilde{D}_{m}\right]
+\frac{\left\Vert \omega\right\Vert _{\infty}^{2}}{2m}g^{2}\left(  t,0\right)
\\
&  =\frac{4t^{2}\exp\left(  t^{2}\sigma^{2}\right)  }{\pi^{2}m}\left(
t^{2}\sigma^{4}+\tilde{D}_{m}\right)  +\frac{\left\Vert \omega\right\Vert
_{\infty}^{2}}{2m}g^{2}\left(  t,0\right)  .
\end{align*}
In addition, when $\sigma_{i}^{2}=\sigma^{2}$ for all $i\in\left\{
1,\ldots,m\right\}  $, $\tilde{D}_{m}$ becomes $D_{m}=\sigma^{2}+m^{-1}%
\sum_{i=1}^{m}\mu_{i}^{2}$.

Second, let us upper bound $\Upsilon_{1}\left(  \lambda,\tilde{\tau}\right)  $
and $\Upsilon_{2}\left(  \lambda,\tilde{\tau}\right)  $ defined by
(\ref{upsilon1and2X}). Clearly, (\ref{bndPartialRecipCF}) implies%
\begin{equation}
\int_{0}^{1}dy\int_{0}^{1}\left\vert \tilde{r}_{0}\left(  tys\right)
\right\vert ds\leq\int_{0}^{1}st^{2}\sigma^{2}\exp\left(  2^{-1}t^{2}%
s^{2}\sigma^{2}\right)  ds=\exp\left(  2^{-1}t^{2}\sigma^{2}\right)  ,
\label{BndTildeCFInduce}%
\end{equation}
and due to $s^{2}\leq s$ for $s\in\left[  0,1\right]  $,%
\begin{equation}
\int_{0}^{1}dy\int_{0}^{1}\frac{1}{r_{0}\left(  tys\right)  }ds\leq\int%
_{0}^{1}\exp\left(  2^{-1}t^{2}s\sigma^{2}\right)  ds=\frac{2\left\{
\exp\left(  2^{-1}t^{2}\sigma^{2}\right)  -1\right\}  }{t^{2}\sigma^{2}}
\label{BndRecipCF}%
\end{equation}
and%
\[
g\left(  t,0\right)  =2\int_{0}^{1}\exp\left(  2^{-1}t^{2}s^{2}\sigma
^{2}\right)  ds\leq2\int_{0}^{1}\exp\left(  2^{-1}t^{2}s\sigma^{2}\right)
ds=\frac{4\left\{  \exp\left(  2^{-1}t^{2}\sigma^{2}\right)  -1\right\}
}{t^{2}\sigma^{2}}.
\]
Therefore, (\ref{BndTildeCFInduce}) and (\ref{BndRecipCF}) imply respectively%
\begin{equation}
\left\{
\begin{array}
[c]{l}%
\Upsilon_{1}\left(  \lambda,\tilde{\tau}\right)  =\dfrac{\lambda}{\pi}%
\sup_{t\in\left[  0,\tilde{\tau}\right]  }%
{\displaystyle\int_{0}^{1}}
dy%
{\displaystyle\int_{0}^{1}}
\left\vert \tilde{r}_{0}\left(  tys\right)  \right\vert ds\leq\dfrac{\lambda
}{\pi}\exp\left(  2^{-1}\tilde{\tau}^{2}\sigma^{2}\right) \\
\Upsilon_{2}\left(  \lambda,\tilde{\tau}\right)  =\dfrac{\lambda}{\pi}%
\sup_{t\in\left[  0,\tilde{\tau}\right]  }%
{\displaystyle\int_{0}^{1}}
dy%
{\displaystyle\int_{0}^{1}}
\dfrac{ds}{r_{0}\left(  tys\right)  }\leq\dfrac{\lambda}{\pi}\dfrac
{2\exp\left(  2^{-1}\tilde{\tau}^{2}\sigma^{2}\right)  }{\tilde{\tau}%
^{2}\sigma^{2}}%
\end{array}
\right.  . \label{BndTwoUpsilons}%
\end{equation}
Now set $\tau_{m}=\sqrt{2\gamma\ln m}$ and $\gamma_{m}=\sqrt{2q\ln m}$ for
some positive constants $q$ and $\gamma$, such that $q>2\vartheta>0$. Then
(\ref{BndTwoUpsilons}) implies%
\begin{equation}
\left\{
\begin{array}
[c]{l}%
\Upsilon_{1}\left(  \gamma_{m}m^{-1/2},\tau_{m}\right)  \leq\dfrac{\gamma
_{m}m^{-1/2}}{\pi}\exp\left(  2^{-1}\tau_{m}^{2}\sigma^{2}\right)
=\dfrac{\sqrt{2q\ln m}}{\pi}m^{\gamma\sigma^{2}-1/2}\\
\Upsilon_{2}\left(  \tau_{m}m^{-\vartheta_{1}},\tau_{m}\right)  \leq
\dfrac{\gamma_{m}m^{-\vartheta_{1}}}{\pi}\exp\left(  2^{-1}\tau_{m}^{2}%
\sigma^{2}\right)  =\dfrac{\sqrt{2q\ln m}}{\pi}m^{\gamma\sigma^{2}%
-\vartheta_{1}}%
\end{array}
\right.  . \label{TwoUpsilonBndRealized}%
\end{equation}
This, together with the inequality (\ref{boundOnUpsilon}), i.e.,%
\[
\Upsilon\left(  \lambda,\tilde{\tau}\right)  \leq\lambda\frac{4\exp\left(
2^{-1}\tilde{\tau}^{2}\sigma^{2}\right)  }{\tilde{\tau}^{2}\sigma^{2}}\text{
\ and }{\Upsilon\left(  \gamma_{m}m^{-1/2},\tau_{m}\right)  }\leq\frac
{4\sqrt{2q\ln m}}{2\gamma\ln m\times\sigma^{2}}m^{\sigma^{2}\gamma-1/2},
\]
imply that, the total error $\hat{\Upsilon}\left(  \tau_{m},\gamma
_{m},\vartheta_{1}\right)  $ defined by
(\ref{UpperBndTotalErrorAndProbOneSide}) satisfies $\hat{\Upsilon}\left(
\tau_{m},\gamma_{m},\vartheta_{1}\right)  \leq\bar{\Upsilon}\left(
q,\gamma,\vartheta_{1}\right)  $, where%
\[
\bar{\Upsilon}\left(  q,\gamma,\vartheta_{1}\right)  =\frac{\left\Vert
\omega\right\Vert _{\infty}\sqrt{2q\ln m}}{\gamma\ln m\times\sigma^{2}%
}m^{\sigma^{2}\gamma-1/2}+\dfrac{\sqrt{2q\ln m}}{\pi}m^{\gamma\sigma^{2}%
-1/2}+\dfrac{2\sqrt{\gamma q}\ln m}{\pi}m^{\gamma\sigma^{2}-\vartheta_{1}}.
\]
So, the uniform bound (\ref{BoundSupWholeErrorX}) becomes%
\begin{equation}
\Pr\left\{  \sup\nolimits_{\boldsymbol{\mu}\in\mathcal{\tilde{B}}_{m}\left(
\tilde{\rho}_{m}\right)  }\sup\nolimits_{t\in\left[  0,\tau_{m}\right]
}\left\vert e_{m}\left(  t\right)  \right\vert \leq\bar{\Upsilon}\left(
\gamma,\vartheta,\vartheta_{1}\right)  \right\}  \geq1-p_{m}^{\ast}\left(
c_{1},\Delta_{1},\kappa,\vartheta_{1},\gamma_{m}\right)
\label{GaussianTotalErrorUpperBnd}%
\end{equation}
for all large $m$, where $p_{m}^{\ast}\left(  c_{1},\Delta_{1},\kappa
,\vartheta_{1},\gamma_{m}\right)  $ is defined by
(\ref{UpperBndTotalErrorAndProbOneSide}) as%
\[
p_{m}^{\ast}\left(  c_{1},\Delta_{1},\kappa,\vartheta_{1},\gamma_{m}\right)
=2\hat{p}_{2,m}^{\ast}\left(  c_{1},\Delta_{1},\gamma_{m}\right)  +\hat
{p}_{3,m}^{\ast}\left(  c_{1},\Delta_{1},\kappa,\vartheta_{1},\gamma
_{m}\right)  .
\]
When $0<\gamma<\min\left\{  2^{-1},\vartheta_{1}\right\}  \times\sigma^{-2}$
and $0<\vartheta_{1}<2^{-1}$, i.e., when $0<\gamma<\vartheta_{1}\sigma^{-2}$
and $0<\vartheta_{1}<2^{-1}$,
\begin{equation}
\bar{\Upsilon}\left(  q,\gamma,\vartheta_{1}\right)  \leq Cm^{\sigma^{2}%
\gamma-\vartheta_{1}}\ln m=o\left(  1\right)
\label{GaussianErroBndSimplified}%
\end{equation}
and%
\begin{equation}
\pi_{1,m}^{-1}\sup\nolimits_{t\in\left[  0,\tau_{m}\right]  }\left\vert
e_{m}\left(  t\right)  \right\vert \leq\pi_{1,m}^{-1}\bar{\Upsilon}\left(
q,\gamma,\vartheta_{1}\right)  \leq\pi_{1,m}^{-1}Cm^{\sigma^{2}\gamma
-\vartheta_{1}}\ln m. \label{ErrorviaPi1Recip}%
\end{equation}

Third, let us show $p_{m}^{\ast}=o\left(  1\right)  $ by verifying the
conditions (\ref{OneSideTotalErrorProb}). Let us look at%
\[
C_{m,\boldsymbol{\mu},F_{0}}=\left[  1-2\Pr\left\{  X_{\left(  0\right)
}>\kappa-\left\Vert \boldsymbol{\mu}\right\Vert _{\infty}\right\}  \right]
^{m},
\]
where $X_{\left(  0\right)  }$ follows the Gaussian distribution with mean $0$
and variance $\sigma^{2}$, i.e., $X_{\left(  0\right)  }\sim\mathcal{N}%
_{1}\left(  0,\sigma^{2}\right)  $. By Mill's ratio (see, e.g., inequality (9)
of \cite{Gasull2015}), when $X\sim\mathcal{N}_{1}\left(  0,1\right)  $,%
\[
\frac{t}{1+t^{2}}\phi\left(  t\right)  <\Pr\left(  X>t\right)  <\frac
{\phi\left(  t\right)  }{t}\text{ for }t>0,
\]
where $\phi\left(  x\right)  =\left(  \sqrt{2\pi}\right)  ^{-1}\exp\left(
-2^{-1}x^{2}\right)  $ is the density of $X$, we see from the transform
$X\mapsto Y=\sigma X$ that%
\begin{equation}
\frac{t/\sigma}{1+\left(  t/\sigma\right)  ^{2}}\phi\left(  t/\sigma\right)
<\Pr\left(  Y>t\right)  <\frac{\phi\left(  t/\sigma\right)  }{t/\sigma}\text{
for \ }t>0. \label{MillRatioSigmaSquare}%
\end{equation}
Let $\tilde{\kappa}=\kappa-\left\Vert \boldsymbol{\mu}\right\Vert _{\infty}$
and $p_{0}\left(  \tilde{\kappa}\right)  :=\Pr\left\{  X_{\left(  0\right)
}>\tilde{\kappa}\right\}  $. Then, (\ref{MillRatioSigmaSquare}) implies that,
as $\tilde{\kappa}\rightarrow\infty$,
\[
p_{0}\left(  \tilde{\kappa}\right)  =\frac{\sigma+o\left(  1\right)  }%
{\tilde{\kappa}}\phi\left(  \frac{\tilde{\kappa}}{\sigma}\right)
=\frac{\sigma+o\left(  1\right)  }{\tilde{\kappa}}\exp\left(  -\left(
2\sigma^{2}\right)  ^{-1}\tilde{\kappa}^{2}\right)  .
\]
Note that we want $C_{m,\boldsymbol{\mu},F_{0}}=1+o\left(  1\right)  $, i.e.,
$m\ln\left(  1-2p_{0}\left(  \tilde{\kappa}\right)  \right)  \rightarrow0$.
Let $\tilde{\kappa}=\sqrt{2\sigma^{2}\vartheta_{2}\ln m}$ for $\vartheta
_{2}\geq1$. Then $p_{0}\left(  \tilde{\kappa}\right)  \rightarrow0$,
$\tilde{\kappa}\rightarrow\infty$ and%
\[
m\ln\left(  1-2p_{0}\left(  \tilde{\kappa}\right)  \right)  =m\left(
-2p_{0}\left(  \tilde{\kappa}\right)  +o\left(  p_{0}\left(  \tilde{\kappa
}\right)  \right)  \right)  =-2\times\frac{\sigma+o\left(  1\right)  }%
{\sqrt{2\sigma^{2}\vartheta_{2}\ln m}}m^{1-\vartheta_{2}}+o\left(  1\right)
=o\left(  1\right)  .
\]
Namely, $\kappa-\left\Vert \boldsymbol{\mu}\right\Vert _{\infty}=\sqrt
{2\sigma^{2}\vartheta_{2}\ln m}$ for $\vartheta_{2}\geq1$ forces
$C_{m,\boldsymbol{\mu},F_{0}}=1+o\left(  1\right)  $.

In addition, by the remark right after (\ref{supBndSinAX}), $R_{2}\left(
\tilde{\rho}_{m}\right)  =O\left(  m^{\vartheta^{\prime\prime}}\right)  $
implies $R_{1}\left(  \sqrt{\tilde{\rho}_{m}}\right)  =O\left(  m^{\vartheta
^{\prime\prime}/2}\right)  $. So,%
\[
\left\{
\begin{array}
[c]{l}%
\left(  \gamma_{m}R_{2}\left(  \tilde{\rho}_{m}\right)  \right)
^{-1}m^{\vartheta-\vartheta_{1}}\ln\tau_{m}\rightarrow\infty\ \text{when
}R_{2}\left(  \tilde{\rho}_{m}\right)  =O\left(  m^{\vartheta^{\prime\prime}%
}\right)  \text{ and }0\leq\vartheta^{\prime\prime}<\vartheta-\vartheta_{1},\\
\left(  \gamma_{m}R_{1}\left(  \sqrt{\tilde{\rho}_{m}}\right)  \right)
^{-1}m^{\vartheta-1/2}\ln\tau_{m}\rightarrow\infty\text{ when }0\leq
2^{-1}\vartheta^{\prime\prime}<\vartheta-2^{-1}%
\end{array}
\right.
\]
Further, $m^{-1-2\vartheta+2\vartheta_{1}}\gamma_{m}^{2}\ln^{-2}\tau
_{m}=o\left(  1\right)  $ when $\vartheta_{1}-\vartheta<2^{-1}$ (which holds
when $0<\vartheta_{1}<2^{-1}$ and $\vartheta>2^{-1}$), and $\tau_{m}%
^{2}m^{\vartheta}\exp\left(  -2^{-1}\gamma_{m}^{2}\right)  =o\left(  1\right)
$ when $q>\vartheta$. Moreover, when $0<\vartheta_{1}<2^{-1}$, then%
\begin{align*}
&  \ln\left[  \tau_{m}m^{\vartheta}\exp\left(  -2^{-1}\kappa^{-2}%
m^{1-2\vartheta_{1}}\left(  \gamma_{m}^{2}-2\ln\tau_{m}\right)  \right)
\right] \\
&  =C+2^{-1}\ln\ln m+\vartheta\ln m-m^{1-2\vartheta_{1}}\frac{\left(
2q\gamma\ln m\right)  \left(  1-o\left(  1\right)  \right)  }{2\times\left(
4\sigma^{4}\vartheta_{2}^{2}\ln^{2}m\right)  }\rightarrow-\infty,
\end{align*}
i.e., $\tau_{m}m^{\vartheta}\exp\left(  -2^{-1}\kappa^{-2}m^{1-2\vartheta_{1}%
}\left(  \tau_{m}^{2}-2\ln\tau_{m}\right)  \right)  =o\left(  1\right)  $. So,
as desired, we have%
\begin{equation}
\left\{
\begin{array}
[c]{l}%
0<\vartheta_{1}<2^{-1};\vartheta>2^{-1}\\
q>\vartheta;\tau_{m}=\sqrt{2\gamma\ln m}\\
\vartheta_{2}\geq1;\kappa-\left\Vert \boldsymbol{\mu}\right\Vert _{\infty
}=\sqrt{2\sigma^{2}\vartheta_{2}\ln m}\\
R_{2}\left(  \tilde{\rho}_{m}\right)  =O\left(  m^{\vartheta^{\prime\prime}%
}\right)  ;0\leq\vartheta^{\prime\prime}<\vartheta-2^{-1}%
\end{array}
\right.  \Longrightarrow p_{m}^{\ast}=o\left(  1\right)  .
\label{GaussianOneSideProbTozero}%
\end{equation}

Summarizing (\ref{GaussianTotalErrorUpperBnd}), (\ref{ErrorviaPi1Recip}),
(\ref{GaussianErroBndSimplified}) and (\ref{GaussianOneSideProbTozero}) and
(\ref{eqdx4aA1}), we have the desired consistency as%
\begin{equation}
\left\{
\begin{array}
[c]{l}%
0<\vartheta_{1}<2^{-1};\vartheta_{2}\geq1;q>\vartheta;0<\gamma<\vartheta
_{1}\sigma^{-2}\\
0\leq\vartheta^{\prime\prime}<\vartheta-2^{-1};R_{2}\left(  \tilde{\rho}%
_{m}\right)  =O\left(  m^{\vartheta^{\prime\prime}}\right) \\
\tau_{m}=\sqrt{2\gamma\ln m},\kappa-\left\Vert \boldsymbol{\mu}\right\Vert
_{\infty}=\sqrt{2\sigma^{2}\vartheta_{2}\ln m}\\
\tau_{m}^{-1}\left(  1+\tilde{u}_{m}^{-1}\right)  =o\left(  \pi_{1,m}\right)
;\pi_{1,m}^{-1}m^{\sigma^{2}\gamma-\vartheta_{1}}\ln m=o\left(  1\right)
\end{array}
\right.  \Rightarrow\left\vert \pi_{1,m}^{-1}\hat{\varphi}_{m}\left(  \tau
_{m},\mathbf{z}\right)  -1\right\vert \rightsquigarrow0.
\label{GaussianOneSideConsistency}%
\end{equation}
To achieve uniform consistency, we can define $\Omega_{m}\left(  \tilde{\rho
}_{m}\right)  $ be a subset of $\mathcal{\tilde{B}}_{m}\left(  \tilde{\rho
}_{m}\right)  $ for each $m\geq1$ and upgrade
(\ref{GaussianOneSideConsistency}) into%
\[
\mathcal{Q}\left(  \mathcal{F}\right)  =\left\{
\begin{array}
[c]{l}%
0<\vartheta_{1}<2^{-1};q>\vartheta\\
0\leq\vartheta^{\prime\prime}<\vartheta-2^{-1};R_{2}\left(  \tilde{\rho}%
_{m}\right)  =O\left(  m^{\vartheta^{\prime\prime}}\right) \\
0<\gamma<\vartheta_{1}\sigma^{-2};\tau_{m}=\sqrt{2\gamma\ln m}\\
\vartheta_{2}\geq1;\kappa-\left\Vert \boldsymbol{\mu}\right\Vert _{\infty
}=\sqrt{2\sigma^{2}\vartheta_{2}\ln m}\\
\tau_{m}^{-1}\sup_{\boldsymbol{\mu}\in\Omega_{m}\left(  \tilde{\rho}%
_{m}\right)  }\pi_{1,m}^{-1}\left(  1+\tilde{u}_{m}^{-1}\right)  =o\left(
1\right) \\
m^{\sigma^{2}\gamma-\vartheta_{1}}\ln m\sup_{\boldsymbol{\mu}\in\Omega
_{m}\left(  \tilde{\rho}_{m}\right)  }\pi_{1,m}^{-1}=o\left(  1\right)
\end{array}
\right\}  .
\]
Note again that $\tau_{m}\pi_{1,m}\rightarrow\infty$ holds when $\tau_{m}%
^{-1}\left(  1+\tilde{u}_{m}^{-1}\right)  =o\left(  \pi_{1,m}\right)  $, which
implies that $\pi_{1,m}^{-1}\leq C\sqrt{\ln m}$ for some constant $C>0$ for
all sufficiently large $m$. The set $\mathcal{Q}\left(  \mathcal{F}\right)  $
can be interpreted as follows: let $\tilde{C}>0$ be a constant, set
$\vartheta_{2}=1$, $\vartheta_{1}=3/8$, $q=2$, $\vartheta=5/8$ and
$\gamma=\left(  1/8\right)  \sigma^{-2}$, and define%
\[
\mathcal{U}_{m}=\left\{  \boldsymbol{\mu}\in\mathcal{\tilde{B}}_{m}\left(
\tilde{\rho}_{m}\right)  :%
\begin{array}
[c]{l}%
0\leq\vartheta^{\prime\prime}<1/8;R_{2}\left(  \tilde{\rho}_{m}\right)
\leq\tilde{C}m^{\vartheta^{\prime\prime}}\\
\tau_{m}=\sqrt{\left(  1/4\right)  \sigma^{-2}\ln m};\kappa-\left\Vert
\boldsymbol{\mu}\right\Vert _{\infty}=\sqrt{2\sigma^{2}\ln m}\\
\pi_{1,m}^{-1}\leq\tilde{C}\sqrt{\ln m};\tilde{u}_{m}\geq\left(  \tau
_{m}\right)  ^{-1}\ln\ln m
\end{array}
\right\}  ,
\]
then $\sup\nolimits_{\boldsymbol{\mu}\in\mathcal{U}_{m}}\left(  \left\vert
\pi_{1,m}^{-1}\sup\nolimits_{t\in\left\{  \tau_{m}\right\}  }\hat{\varphi}%
_{m}\left(  t,\mathbf{z}\right)  -1\right\vert \right)  \rightsquigarrow0$
as\ $m\rightarrow\infty.$\qed

\section{Proofs Related to the Extension}

\label{AppProofsExt}

\subsection{Proof of \autoref{ThmExtA}}

We show the first claim. Recall $\hat{\phi}\left(  s\right)  =\int_{a}^{b}%
\phi\left(  y\right)  \exp\left(  -\iota ys\right)  dy$ and%
\[
\mathcal{D}_{\phi}\left(  t,\mu;a,b\right)  =\frac{t}{2\pi}\int_{-1}^{1}%
\hat{\phi}\left(  ts\right)  \exp\left(  \iota\mu ts\right)  ds.
\]
Set%
\[
K_{1}\left(  t,x\right)  =\frac{t}{2\pi}\int_{\left[  -1,1\right]  }\frac
{1}{\hat{F}_{0}\left(  ts\right)  }\hat{\phi}\left(  ts\right)  \exp\left(
\iota txs\right)  ds.
\]
Then%
\begin{align*}
\int K_{1}\left(  t,x\right)  dF_{\mu}\left(  x\right)   &  =\frac{t}{2\pi
}\int_{\left[  -1,1\right]  }\frac{\hat{\phi}\left(  ts\right)  }{\hat{F}%
_{0}\left(  ts\right)  }\left\{  \int\exp\left(  \iota txs\right)  dF_{\mu
}\left(  x\right)  \right\}  ds\\
&  =\frac{t}{2\pi}\int_{-1}^{1}\frac{\hat{F}_{\mu}\left(  ts\right)  }{\hat
{F}_{0}\left(  ts\right)  }\hat{\phi}\left(  ts\right)  ds=\frac{t}{2\pi}%
\int_{-1}^{1}\exp\left\{  \iota ts\mu\right\}  \hat{\phi}\left(  ts\right)
ds\\
&  =\frac{t}{2\pi}\int_{a}^{b}\phi\left(  y\right)  dy\int_{-1}^{1}
\exp\left(  -\iota yts\right)  \exp\left\{  \iota ts\mu\right\}  ds.
\end{align*}
Namely, $\int K_{1}\left(  t,x\right)  dF_{\mu}\left(  x\right)
=\mathcal{D}_{\phi}\left(  t,\mu;a,b\right)  $ as desired. Since $\mathcal{F}$
is a Type I location-shift family, $\hat{F}_{0}\equiv r_{0}$ holds, $r_{0}$ is
an even function, and
\begin{align*}
K_{1}\left(  t,x\right)   &  =\frac{t}{2\pi}\int_{\left[  -1,1\right]  }%
\frac{1}{r_{0}\left(  ts\right)  }\hat{\phi}\left(  ts\right)  \exp\left(
\iota txs\right)  ds\\
&  =\frac{t}{2\pi}\int_{\left[  -1,1\right]  }\frac{1}{r_{0}\left(  ts\right)
}\exp\left(  \iota txs\right)  ds\int_{a}^{b}\phi\left(  y\right)  \exp\left(
-\iota yts\right)  dy\\
&  =\frac{t}{2\pi}\int_{a}^{b}\phi\left(  y\right)  dy\int_{\left[
-1,1\right]  }\frac{\exp\left\{  \iota ts\left(  x-y\right)  \right\}  }%
{r_{0}\left(  ts\right)  }ds\\
&  = \frac{t}{2\pi}\int_{a}^{b}\phi\left(  y\right)  dy\int_{\left[
-1,1\right]  }\frac{\cos\left\{  ts\left(  x-y\right)  \right\}  }%
{r_{0}\left(  ts\right)  }ds.
\end{align*}
Finally, we only need to capture the contributions of the end points $a$ and
$b$ to estimating $\check{\pi}_{0,m}$. By \autoref{ThmPoinNull}, we only need
to set $\left(  K,\psi\right)  $ as given by (\ref{eq2b}).\qed

\subsection{Proof of \autoref{ThmFinal}}

The proof is split into three parts in three subsections: bounding the
variance of the error $e_{m}\left(  t\right)  =\hat{\varphi}_{m}\left(
t,\mathbf{z}\right)  -\varphi_{m}\left(  t,\boldsymbol{\mu}\right)  $,
deriving a uniform bound on $\left\vert e_{m}\left(  t\right)  \right\vert $
and then the consistency of $\hat{\varphi}_{m}\left(  t,\mathbf{z}\right)  $,
and specializing these results to Gaussian family. The proof uses almost
identical arguments as those for the proof of \autoref{ThmIVa}, because the
only difference between the construction in \autoref{ThmIVa} and the
construction in \autoref{ThmFinal} (this Theorem) is that the integrand of the
latter construction has a factor $\phi\left(  y\right)  $ compared to the
former, so that the variance bounds (for $e_{m}\left(  t\right)  $) for the
latter will differ by at most a factor of $\left\Vert \phi\right\Vert
_{\infty}^{2}$ from those for the former and the error bounds (for
$e_{m}\left(  t\right)  $) for the latter will differ by at most a factor of
$\left\Vert \phi\right\Vert _{\infty}^{2}$ from those for the former.

\subsubsection{Bound on the variance of error}

First, let us bound the variance of $e_{m}\left(  t\right)  $. Take $t>0$.
Recall%
\[
K_{1}\left(  t,x\right)  =\frac{t}{2\pi}\int_{a}^{b}\phi\left(  y\right)
dy\int_{\left[  -1,1\right]  }\frac{\cos\left\{  ts\left(  x-y\right)
\right\}  }{r_{0}\left(  ts\right)  }ds=\frac{2t}{2\pi}\int_{a}^{b}\phi\left(
y\right)  dy\int_{0}^{1}\frac{\cos\left\{  ts\left(  x-y\right)  \right\}
}{r_{0}\left(  ts\right)  }ds
\]
(because both $\cos\left(  \cdot\right)  $ and $r_{0}\left(  \cdot\right)  $
are even) and $\psi_{1}\left(  t,\mu\right)  =\int K_{1}\left(  t,x\right)
dF_{\mu}\left(  x\right)  =\mathcal{D}_{\phi}\left(  t,\mu;a,b\right)  $.
Recall $w_{1,i}\left(  v,y\right)  =\cos\left\{  v\left(  z_{i}-y\right)
\right\}  $,%
\[
S_{1,m}\left(  v,y\right)  =\frac{1}{m}\sum_{i=1}^{m}\left[  w_{1,i}\left(
v,y\right)  -\mathbb{E}\left\{  w_{1,i}\left(  v,y\right)  \right\}  \right]
\]
and $g\left(  t,\mu\right)  =\int_{\left[  -1,1\right]  }1/r_{\mu}\left(
ts\right)  ds$ from the proof of \autoref{ThmIVa}. Let $e_{1,m}\left(
t\right)  =\hat{\varphi}_{1,m}\left(  t,\mathbf{z}\right)  -\varphi
_{1,m}\left(  t,\boldsymbol{\mu}\right)  $ with $\varphi_{1,m}\left(
t,\boldsymbol{\mu}\right)  =\mathbb{E}\left\{  \hat{\varphi}_{1,m}\left(
t,\mathbf{z}\right)  \right\}  $. Since $\hat{\varphi}_{1,m}\left(
t,\mathbf{z}\right)  =m^{-1}\sum_{i=1}^{m}K_{1}\left(  t,z_{i}\right)  $, we
immediately see
\begin{equation}
e_{1,m}\left(  t\right)  =\frac{t}{2\pi}\int_{a}^{b}\phi\left(  y\right)
dy\int_{\left[  -1,1\right]  }\frac{S_{1,m}\left(  ts,y\right)  }{r_{0}\left(
ts\right)  }ds=\frac{2t}{2\pi}\int_{a}^{b}\phi\left(  y\right)  dy\int_{0}%
^{1}\frac{S_{1,m}\left(  ts,y\right)  }{r_{0}\left(  ts\right)  }ds
\label{ExtensionPhiRepresentation}%
\end{equation}
and
\begin{align}
\mathbb{V}\left\{  \left\vert \hat{\varphi}_{1,m}\left(  t,\mathbf{z}\right)
\right\vert \right\}   &  =\mathbb{V}\left(  \frac{t}{2\pi}\int_{a}^{b}%
\phi\left(  y\right)  dy\int_{\left[  -1,1\right]  }\frac{S_{1,m}\left(
ts,y\right)  }{r_{0}\left(  ts\right)  }ds\right) \nonumber\\
&  \leq\frac{1}{m}\left(  \frac{t}{2\pi}\int_{a}^{b}\phi\left(  y\right)
dy\int_{\left[  -1,1\right]  }\frac{2}{r_{0}\left(  ts\right)  }ds\right)
^{2}\nonumber\\
&  \leq\frac{t^{2}\left(  b-a\right)  ^{2}\left\Vert \phi\right\Vert _{\infty
}^{2}}{\pi^{2}m}g^{2}\left(  t,0\right)  , \label{eq16b}%
\end{align}
where $g\left(  t,0\right)  =2\int_{\left[  0,1\right]  }\left(
1/r_{0}\left(  ts\right)  \right)  ds$.

Recall from the proof of \autoref{ThmIVa} the following quantities:
$\hat{\varphi}_{1,0,m}\left(  t,\mathbf{z};\tau\right)  =m^{-1}\sum_{i=1}%
^{m}K_{1,0}\left(  t,z_{i};\tau\right)  $ and $\varphi_{1,0,m}\left(
t,\boldsymbol{\mu};\tau\right)  =\mathbb{E}\left\{  \hat{\varphi}%
_{1,0,m}\left(  t,\mathbf{z};\tau\right)  \right\}  $ for $\tau\in\left\{
a,b\right\}  $, where%
\[
K_{1,0}\left(  t,x;\tau\right)  =2\int_{\left[  0,1\right]  }\frac
{\omega\left(  s\right)  \cos\left(  ts\left(  x-\tau\right)  \right)  }%
{r_{0}\left(  ts\right)  }ds.
\]
and%
\begin{equation}
e_{1,0,m}\left(  t,\tau\right)  =\hat{\varphi}_{1,0,m}\left(  t,\mathbf{z}%
;\tau\right)  -\varphi_{1,0,m}\left(  t,\boldsymbol{\mu};\tau\right)
=2\int_{\left[  0,1\right]  }\frac{\omega\left(  s\right)  }{r_{0}\left(
ts\right)  }S_{1,m}\left(  ts,\tau\right)  ds. \label{RepresentationPointNull}%
\end{equation}
Then%
\[
\mathbb{V}\left\{  \hat{\varphi}_{1,0,m}\left(  t,\mathbf{z};\phi\left(
\tau\right)  \right)  \right\}  \leq\left\Vert \omega\right\Vert _{\infty}%
^{2}m^{-1}g^{2}\left(  t, \textcolor{blue}{  \tau}  \right)  =\left\Vert
\omega\right\Vert _{\infty}^{2}m^{-1}g^{2}\left(  t,0\right)  ,
\]

Since%
\begin{equation}
\left\{
\begin{array}
[c]{l}%
K\left(  t,x\right)  =K_{1}\left(  t,x\right)  -2^{-1}\left\{  \phi\left(
a\right)  K_{1,0}\left(  t,x;a\right)  +\phi\left(  b\right)  K_{1,0}\left(
t,x;b\right)  \right\} \\
\psi\left(  t,\mu\right)  =\psi_{1}\left(  t,\mu\right)  -2^{-1}\left\{
\phi\left(  a\right)  \psi_{1,0}\left(  t,\mu;a\right)  +\phi\left(  b\right)
\psi_{1,0}\left(  t,\mu;b\right)  \right\} \\
\mathbb{E}\left\{  K\left(  t,Z\right)  \right\}  =\psi\left(  t,\mu\right)
\text{ for }Z\text{ with CDF }F_{\mu}\\
\hat{\varphi}_{m}\left(  t,\mathbf{z}\right)  =1-\sum_{i=1}^{m}K\left(
t,z_{i}\right)  \text{ and }\varphi_{m}\left(  t,\boldsymbol{\mu}\right)
=\mathbb{E}\left\{  \hat{\varphi}_{m}\left(  t,\mathbf{z}\right)  \right\}
\end{array}
\right.  , \label{ExtensionConstruction}%
\end{equation}
we have%
\begin{equation}
\hat{\varphi}_{m}\left(  t,\mathbf{z}\right)  =\hat{\varphi}_{1,m}\left(
t,\mathbf{z}\right)  -2^{-1}\left[  \phi\left(  a\right)  \hat{\varphi
}_{1,0,m}\left(  t,\mathbf{z};a\right)  +\phi\left(  b\right)  \hat{\varphi
}_{1,0,m}\left(  t,\mathbf{z};b\right)  \right]
\label{ExtensionDecomposition}%
\end{equation}
and hence
\begin{align}
\mathbb{V}\left\{  \hat{\varphi}_{m}\left(  t,\mathbf{z}\right)  \right\}   &
\leq2\mathbb{V}\left(  \hat{\varphi}_{1,m}\left(  t,\mathbf{z}\right)
\right)  +\mathbb{V}\left\{  \hat{\varphi}_{1,0,m}\left(  t,\mathbf{z}%
;a\right)  \right\}  +\mathbb{V}\left\{  \hat{\varphi}_{1,0,m}\left(
t,\mathbf{z};b\right)  \right\} \nonumber\\
&  \leq\frac{2\left\Vert \phi\right\Vert _{\infty}^{2}g^{2}\left(  t,0\right)
}{m}\left\{  \frac{t^{2}}{\pi^{2}}\left(  b-a\right)  ^{2}+\left\Vert
\omega\right\Vert _{\infty}^{2}\right\}  . \label{eqdx10}%
\end{align}
Note that (\ref{eqdx10}) is analogous to (\ref{eqLsosdx10}).

\subsubsection{\textcolor{blue}{Uniform bound on error and consistency}}

Second, let us derive a uniform bound on $\left\vert e_{m}\left(  t\right)
\right\vert $. Set $D=\left[  0,\tau_{m}\right]  \times\left[  0,1\right]
\times\left[  a,b\right]  $, Due to the representations
(\ref{ExtensionPhiRepresentation}) and (\ref{RepresentationPointNull}) and the
construction (\ref{ExtensionConstruction}), the settings
(\ref{ParConditionsCase2}) and (\ref{ParConditionsCase1}) and the bounds on
the tail probabilities of%
\[
\sup\nolimits_{\left(  t,s\right)  \in\left[  0,\tau_{m}\right]  \times\left[
0,1\right]  }S_{1,m}\left(  ts,\tau\right)  \text{ and }\sup\nolimits_{\left(
t,s,y\right)  \in D}S_{1,m}\left(  ts,y\right)
\]
in the proof of \autoref{ThmIVa} can be directly used here. Specifically, for
all $m$ large enough, the analog of (\ref{BndNullPhi1Sup}) is, with
probability at least $1-\hat{p}_{1,m}^{\ast}$,
\begin{align}
\sup_{\boldsymbol{\mu}\in\mathcal{B}_{1,m}\left(  \rho_{m}\right)  }\sup
_{t\in\left[  0,\tau_{m}\right]  }\left\vert e_{1,m}\left(  t\right)
\right\vert  &  \leq\frac{\gamma_{m}}{\sqrt{m}}\sup_{\boldsymbol{\mu}%
\in\mathcal{B}_{1,m}\left(  \rho_{m}\right)  }\left(  \sup_{t\in\left[
0,\tau_{m}\right]  }\frac{t}{2\pi}\int_{a}^{b}\left\vert \phi\left(  y\right)
\right\vert dy\int_{\left[  -1,1\right]  }\frac{1}{r_{0}\left(  ts\right)
}ds\right) \nonumber\\
&  \leq\frac{\left(  b-a\right)  \tau_{m}\left\Vert \phi\right\Vert _{\infty}%
}{2\pi}\Upsilon\left(  \tau_{m}m^{-1/2},\tau_{m}\right)  , \label{eq10dx3}%
\end{align}
and with probability at least $1-\hat{p}_{2,m}^{\ast}$,
(\ref{PointNullPhi10Sup}) holds, i.e.,%
\begin{equation}
\sup_{\boldsymbol{\mu}\in\mathcal{B}_{1,m}\left(  \rho_{m}\right)  }\sup
_{t\in\left[  0,\tau_{m}\right]  }\left\vert e_{1,0,m}\left(  t,\tau\right)
\right\vert \leq\frac{2\left\Vert \omega\right\Vert _{\infty}\gamma_{m}}%
{\sqrt{m}}\sup_{t\in\left[  0,\tau_{m}\right]  }\int_{\left[  0,1\right]
}\frac{ds}{r_{0}\left(  ts\right)  }=\left\Vert \omega\right\Vert _{\infty
}\Upsilon\left(  \gamma_{m}m^{-1/2},\tau_{m}\right)
\label{supBndPointNullExtension}%
\end{equation}
for $\tau\in\left\{  a,b\right\}  $. So, in view of
(\ref{ExtensionDecomposition}), (\ref{eq10dx3}) and
(\ref{supBndPointNullExtension}), the analog of (\ref{BndNullErrorSup}) is%
\begin{align}
&  \sup_{\boldsymbol{\mu}\in\mathcal{B}_{1,m}\left(  \rho_{m}\right)  }%
\sup_{t\in\left[  0,\tau_{m}\right]  }\left\vert e_{m}\left(  t\right)
\right\vert \nonumber\\
&  \leq\sup_{\boldsymbol{\mu}\in\mathcal{B}_{1,m}\left(  \rho_{m}\right)
}\sup_{t\in\left[  0,\tau_{m}\right]  }\left\vert e_{1,m}\left(  t\right)
\right\vert +2^{-1}\left\Vert \phi\right\Vert _{\infty}\sum_{\tau\in\left\{
a,b\right\}  }\sup_{\boldsymbol{\mu}\in\mathcal{B}_{m}\left(  \rho_{m}\right)
}\sup_{t\in\left[  0,\tau_{m}\right]  }\left\vert e_{1,0,m}\left(
t,\tau\right)  \right\vert \nonumber\\
&  \leq\left\{  \frac{\left(  b-a\right)  }{2\pi}\tau_{m}+\left\Vert
\omega\right\Vert _{\infty}\right\}  \left\Vert \phi\right\Vert _{\infty
}\Upsilon\left(  \gamma_{m}m^{-1/2},\tau_{m}\right)  , \label{eq11exh}%
\end{align}
with probability at least $1-2\hat{p}_{2,m}^{\ast}-\hat{p}_{1,m}^{\ast}$ for
all $m$ large enough, where $\hat{p}_{2,m}^{\ast}$ and $\hat{p}_{1,m}^{\ast}$
are exactly as in the proof of \autoref{ThmIVa}, and $\Upsilon\left(
\lambda,\tilde{\tau}\right)  $ is defined by (\ref{UpsilonA}) as%
\[
\Upsilon\left(  \lambda,\tilde{\tau}\right)  =2\lambda\sup_{t\in\left[
0,\tilde{\tau}\right]  }\int_{\left[  0,1\right]  }\frac{ds}{r_{0}\left(
ts\right)  }\text{ for }\lambda,\tilde{\tau}>0\text{.}%
\]

Third, we show the consistency of $\hat{\varphi}_{m}\left(  t,\mathbf{z}%
\right)  $. Again the results of in the proof of \autoref{ThmIVa} can be
directly used here. In particular, since (\ref{BndNullProbOnErrorBnd}) there
holds here, i.e.,%
\[
\left\{
\begin{array}
[c]{l}%
\left(  2\tau_{m}+R_{1}\left(  \rho_{m}\right)  \right)  ^{-1}m^{\vartheta
-1/2}\gamma_{m}^{-1}\ln\tau_{m}\rightarrow\infty\\
\multicolumn{1}{c}{\tau_{m}^{2}m^{2\vartheta}\exp\left(  -2^{-1}\gamma_{m}%
^{2}\right)  +m^{-2\vartheta}\gamma_{m}^{2}\ln^{-2}\tau_{m}=o\left(  1\right)
}%
\end{array}
\right.  \Longrightarrow\hat{p}_{2,m}^{\ast}=o\left(  1\right)  =\hat{p}%
_{1,m}^{\ast},
\]
and (\ref{eq11exh}) holds, we assert that%
\begin{equation}
\left\{
\begin{array}
[c]{l}%
\left(  2\tau_{m}+R_{1}\left(  \rho_{m}\right)  \right)  ^{-1}m^{\vartheta
-1/2}\gamma_{m}^{-1}\ln\tau_{m}\rightarrow\infty\\
\multicolumn{1}{c}{\tau_{m}^{2}m^{2\vartheta}\exp\left(  -2^{-1}\gamma_{m}%
^{2}\right)  +m^{-2\vartheta}\gamma_{m}^{2}\ln^{-2}\tau_{m}=o\left(  1\right)
}\\
\tau_{m}\pi_{1,m}^{-1}\Upsilon\left(  \gamma_{m}m^{-1/2},\tau_{m}\right)
=o\left(  1\right)
\end{array}
\right.  \Longrightarrow\sup_{\boldsymbol{\mu}\in\mathcal{B}_{1,m}\left(
\rho_{m}\right)  }\sup_{t\in\left[  0,\tau_{m}\right]  }\left\vert
e_{m}\left(  t\right)  \right\vert \rightsquigarrow0\text{.}
\label{ExtentionErrorConvergingInProb}%
\end{equation}
On the other hand, Lemma .. in, whose $\tilde{\psi}_{1,0}$ has parameter
$\sigma>0$ and whose $\tilde{\psi}_{1,0}$ with $\sigma=1$ reduces to
$\psi_{1,0}$ in this paper, asserts that, for all positive $t$ such that
$tu_{m}\geq2$ and $t\left(  b-a\right)  \geq2$,%
\begin{align}
\left\vert \check{\pi}_{0,m}^{-1}\varphi_{m}\left(  t,\boldsymbol{\mu}\right)
-1\right\vert  &  \leq\frac{20\left\Vert \phi\right\Vert _{1,\infty}}{\pi
t\check{\pi}_{0,m}}+\frac{20\left\Vert \phi\right\Vert _{\infty}}{tu_{m}%
\check{\pi}_{0,m}}+\frac{10\left\Vert \phi\right\Vert _{\infty}}{t\left(
b-a\right)  \check{\pi}_{0,m}}+\frac{4\left(  \left\Vert \omega\right\Vert
_{\mathrm{TV}}+\left\Vert \omega\right\Vert _{\infty}\right)  \left\Vert
\phi\right\Vert _{\infty}}{tu_{m}\check{\pi}_{0,m}}\nonumber\\
&  \leq\frac{C}{t\check{\pi}_{0,m}}\left(  1+\left\Vert \phi\right\Vert
_{\infty}+\frac{1}{u_{m}}\right)  . \label{ExtensionOracleSpeed}%
\end{align}
where $u_{m}=\min_{\tau\in\left\{  a,b\right\}  }\min_{\left\{  j:\mu_{j}%
\neq\tau\right\}  }\left\vert \mu_{j}-\tau\right\vert $ and $\left\Vert
\phi\right\Vert _{1,\infty}=\sup_{\mu\in\left[  a,b\right]  }C_{\mu}\left(
\phi\right)  $. Therefore,
\[
\left\vert \check{\pi}_{0,m}^{-1}\varphi_{m}\left(  t,\boldsymbol{\mu}\right)
-1\right\vert \rightarrow0\text{ if }t_{m}^{-1}\left(  1+u_{m}^{-1}\right)
=o\left(  \check{\pi}_{0,m}\right)  .
\]
So, in view of (\ref{ExtentionErrorConvergingInProb}),
(\ref{ExtensionOracleSpeed}) and the decomposition%
\[
\left\vert \check{\pi}_{0,m}^{-1}\hat{\varphi}_{m}\left(  t,\mathbf{z}\right)
-1\right\vert \leq\check{\pi}_{0,m}^{-1}\left\vert \hat{\varphi}_{m}\left(
t,\mathbf{z}\right)  -\varphi_{m}\left(  t,\boldsymbol{\mu}\right)
\right\vert +\left\vert \check{\pi}_{0,m}^{-1}\varphi_{m}\left(
t_{m},\boldsymbol{\mu}\right)  -1\right\vert ,
\]
we have the desired consistency%
\begin{equation}
\left\{
\begin{array}
[c]{l}%
\left(  2\tau_{m}+R_{1}\left(  \rho_{m}\right)  \right)  ^{-1}m^{\vartheta
-1/2}\gamma_{m}^{-1}\ln\tau_{m}\rightarrow\infty\\
\tau_{m}^{2}m^{2\vartheta}\exp\left(  -2^{-1}\gamma_{m}^{2}\right)
+m^{-2\vartheta}\gamma_{m}^{2}\ln^{-2}\tau_{m}=o\left(  1\right) \\
\tau_{m}\check{\pi}_{0,m}^{-1}\Upsilon\left(  \gamma_{m}m^{-1/2},\tau
_{m}\right)  =o\left(  1\right) \\
\tau_{m}^{-1}\left(  1+u_{m}^{-1}\right)  =o\left(  \check{\pi}_{0,m}\right)
\end{array}
\right.  \Longrightarrow\left\vert \check{\pi}_{0,m}^{-1}\hat{\varphi}%
_{m}\left(  \tau_{m},\mathbf{z}\right)  -1\right\vert \rightsquigarrow0.
\label{ExtensionConsistency}%
\end{equation}
For uniform consistency, we have from (\ref{eq11exh}) and
(\ref{ExtensionConsistency}) a uniform consistency class as%
\begin{equation}
\mathcal{Q}\left(  \mathcal{F}\right)  =\left\{
\begin{array}
[c]{l}%
\left(  2\tau_{m}+R_{1}\left(  \rho_{m}\right)  \right)  ^{-1}m^{\vartheta
-1/2}\gamma_{m}^{-1}\ln\tau_{m}\rightarrow\infty\\
\tau_{m}^{2}m^{2\vartheta}\exp\left(  -2^{-1}\gamma_{m}^{2}\right)
+m^{-2\vartheta}\gamma_{m}^{2}\ln^{-2}\tau_{m}=o\left(  1\right) \\
\tau_{m}\Upsilon\left(  \gamma_{m}m^{-1/2},\tau_{m}\right)  \sup
_{\boldsymbol{\mu}\in\Omega_{m}\left(  \rho_{m}\right)  }\check{\pi}%
_{0,m}^{-1}=o\left(  1\right) \\
\tau_{m}^{-1}\left(  1+u_{m}^{-1}\right)  \sup_{\boldsymbol{\mu}\in\Omega
_{m}\left(  \rho_{m}\right)  }\check{\pi}_{0,m}^{-1}=o\left(  1\right)
\end{array}
\right\}  , \label{OneSidedUCCInProof}%
\end{equation}
which has the exact representation as the $\mathcal{Q}\left(  \mathcal{F}%
\right)  $ derived by \autoref{ThmIVa}.

Now consider $\hat{\varphi}_{1,m}\left(  t,\mathbf{z}\right)  =m^{-1}%
\sum_{i=1}^{m}K_{1}\left(  t,z_{i}\right)  $ with $K_{1}$ in (\ref{eq13a})
that estimates
\[
\tilde{\pi}_{0,m}=\check{\pi}_{0,m}+m^{-1}\sum\nolimits_{\left\{  \mu_{i}%
\in\left\{  a,b\right\}  :1\leq i\leq m\right\}  }2^{-1}\phi\left(  \mu
_{i}\right)  .
\]
Recall $\varphi_{1,m}\left(  t,\mathbf{z}\right)  =\mathbb{E}\left\{
\hat{\varphi}_{1,m}\left(  t,\mathbf{z}\right)  \right\}  $. Then
(\ref{eq16b}) directly gives%
\[
\mathbb{V}\left\{  e_{1,m}\left(  t\right)  \right\}  \leq\frac{t^{2}\left(
b-a\right)  ^{2}\left\Vert \phi\right\Vert _{\infty}^{2}}{\pi^{2}m}%
g^{2}\left(  t,0\right)  ,
\]
and (\ref{eq10dx3}) gives, for all $m$ large enough, with probability at least
$1-\hat{p}_{1,m}^{\ast}$,%
\[
\sup_{\boldsymbol{\mu}\in\mathcal{B}_{1,m}\left(  \rho_{m}\right)  }\sup
_{t\in\left[  0,\tau_{m}\right]  }\left\vert e_{1,m}\left(  t\right)
\right\vert \leq\frac{\left(  b-a\right)  \tau_{m}\left\Vert \phi\right\Vert
_{\infty}}{2\pi}\Upsilon\left(  \tau_{m}m^{-1/2},\tau_{m}\right)  .
\]
On the other hand, \autoref{SpeedOracleExt}, whose $\tilde{\psi}_{1,0}$ has
parameter $\sigma>0$ and whose $\tilde{\psi}_{1,0}$ with $\sigma=1$ reduces to
$\psi_{1,0}$, implies%
\begin{equation}
\left\vert \tilde{\pi}_{0,m}^{-1}\varphi_{1,m}\left(  t,\mathbf{z}\right)
-1\right\vert \leq\frac{C}{t\tilde{\pi}_{0,m}}\left(  1+\left\Vert
\phi\right\Vert _{1,\infty}+\frac{1}{u_{m}}\right)  \text{ when }tu_{m}%
\geq2\text{ and }t\left(  b-a\right)  \geq2. \label{OracleSpeedExtensionA}%
\end{equation}
So, (\ref{ExtensionConsistency}) and (\ref{OneSidedUCCInProof}) hold with
$\tilde{\pi}_{0,m}$ in place of $\check{\pi}_{0,m}$ and $\hat{\varphi}_{1,m}$
in place of $\hat{\varphi}_{m}$.

\subsubsection{Specialization to Gaussian family}

Fourth, we specialize the above result to the Gaussian family with scale
parameter $\sigma>0$. Again the results in the proof of \autoref{ThmIVa} can
be directly used here. Set $\tau_{m}=\sqrt{2\gamma\ln m}$ and $\gamma
_{m}=\sqrt{2q\ln m}$ such that $q>2\vartheta>0$, and set $\vartheta>2^{-1}$
and $R_{1}\left(  \rho_{m}\right)  =O\left(  m^{\vartheta^{\prime}}\right)  $
with $0\leq\vartheta^{\prime}<\vartheta-2^{-1}$. Then, as verified by
(\ref{BndNullProbToZeroGaussian}) and (\ref{TheOrderForShiftRho}), $2\hat
{p}_{2,m}^{\ast}+\hat{p}_{1,m}^{\ast}=o\left(  1\right)  $. Inserting into
(\ref{eq11exh}) the inequality%
\[
{\Upsilon\left(  \gamma_{m}m^{-1/2},\tau_{m}\right)  }\leq\frac{4\sqrt{2q\ln
m}}{2\gamma\ln m\times\sigma^{2}}m^{\sigma^{2}\gamma-1/2}%
\]
provided by (\ref{boundOnUpsilon}), we obtain the analog of
(\ref{ErrorBoundGaussianBndNull}) as%
\[
\sup_{\boldsymbol{\mu}\in\mathcal{B}_{1,m}\left(  \rho_{m}\right)  }\sup
_{t\in\left[  0,\tau_{m}\right]  }\left\vert e_{m}\left(  t\right)
\right\vert \leq C\left\Vert \phi\right\Vert _{\infty}\frac{4\sqrt{q}%
m^{\sigma^{2}\gamma-1/2}}{\sigma^{2}\sqrt{\gamma}}\leq Cm^{\sigma^{2}%
\gamma-1/2}.
\]
So, $\tau_{m}\check{\pi}_{0,m}^{-1}\Upsilon\left(  \gamma_{m}m^{-1/2},\tau
_{m}\right)  =o\left(  1\right)  $ if\ $\check{\pi}_{0,m}^{-1}m^{\sigma
^{2}\gamma-1/2}=o\left(  1\right)  $. Further, (\ref{ExtensionOracleSpeed})
implies $\check{\pi}_{0,m}^{-1}\varphi_{m}\left(  \tau_{m},\boldsymbol{\mu
}\right)  \rightarrow1$ when $\tau_{m}^{-1}\left(  1+u_{m}^{-1}\right)
=o\left(  \check{\pi}_{0,m}\right)  $, which also implies $\tau_{m}%
u_{m}\rightarrow\infty$ and $\check{\pi}_{0,m}^{-1}=O\left(  \sqrt{\ln
m}\right)  $. Therefore, we obtain a uniform consistency class as%
\begin{equation}
\mathcal{Q}\left(  \mathcal{F}\right)  =\left\{
\begin{array}
[c]{l}%
0<\gamma<2^{-1}\sigma^{-2};\tau_{m}=\sqrt{2\gamma\ln m}\\
q>2\vartheta;0\leq\vartheta^{\prime}<\vartheta-2^{-1};R_{1}\left(  \rho
_{m}\right)  =O\left(  m^{\vartheta^{\prime}}\right) \\
\check{\pi}_{0,m}^{-1}=O\left(  \sqrt{\ln m}\right)  ;u_{m}\geq\left(
\tau_{m}\right)  ^{-1}\ln\ln m
\end{array}
\right\}  . \label{ExtensionUnifromConsistencyGaussian}%
\end{equation}
The set $\mathcal{Q}\left(  \mathcal{F}\right)  $ can be interpreted as
follows: let $\tilde{C}>0$ be a constant, and define%
\begin{equation}
\mathcal{U}_{m}=\left\{  \boldsymbol{\mu}\in\mathcal{B}_{1,m}\left(  \rho
_{m}\right)  :%
\begin{array}
[c]{l}%
\vartheta=3/4;q=7/4;\tau_{m}=\sqrt{2^{-1}\sigma^{-2}\ln m}\\
0\leq\vartheta^{\prime}<1/4;R_{1}\left(  \rho_{m}\right)  \leq\tilde
{C}m^{\vartheta^{\prime}}\\
\check{\pi}_{0,m}^{-1}\leq\tilde{C}\sqrt{\ln m};u_{m}\geq\left(  \tau
_{m}\right)  ^{-1}\ln\ln m
\end{array}
\right\}  , \label{UCCGaussianExample}%
\end{equation}
then $\sup\nolimits_{\boldsymbol{\mu}\in\mathcal{U}_{m}}\left(  \left\vert
\check{\pi}_{0,m}^{-1}\sup\nolimits_{t\in\left\{  \tau_{m}\right\}  }%
\hat{\varphi}_{m}\left(  t,\mathbf{z}\right)  -1\right\vert \right)
\rightsquigarrow0$ \ as \ $m\rightarrow\infty.$ Note that $\mathcal{Q}\left(
\mathcal{F}\right)  $ and $\mathcal{U}_{m}$ are the same as those for the
Gaussian family provided by \autoref{ThmIVa}.

Finally, we show the claim about the estimator $\hat{\varphi}_{1,m}\left(
t,\mathbf{z}\right)  =m^{-1}\sum_{i=1}^{m}K_{1}\left(  t,z_{i}\right)  $ with
$K_{1}$ in (\ref{eq13a}), which can be regarded as a special case of the
specialization just given above. Specifically, we can just regard
$e_{m}\left(  t\right)  $ as $e_{1,m}\left(  t\right)  $ and $\check{\pi
}_{0,m}$ as $\tilde{\pi}_{0,m}$ in the specialization given above and then
utilize (\ref{OracleSpeedExtensionA}) while keeping other settings the same,
in order to obtain a uniform consistency class as
(\ref{ExtensionUnifromConsistencyGaussian}) and an example as
(\ref{UCCGaussianExample}) while replacing $\check{\pi}_{0,m}$ by $\tilde{\pi
}_{0,m}$ in these two sets.\qed

\section{Correction to the proof of Theorem 3 of \cite{Chen:2018a}}

\label{SecConsistencyI}

The work of \cite{Chen:2018a} dealt with a point null $\Theta_{0}=\left\{
\mu_{0}\right\}  $. Recall the following 2 definitions from \cite{Chen:2018a}:

\begin{enumerate}
\item Definition 1 of \cite{Chen:2018a}:\ If%
\begin{equation}
\left\{  t\in\mathbb{R}:\hat{F}_{\mu_{0}}\left(  t\right)  =0\right\}
=\varnothing\label{eq6bx}%
\end{equation}
and, for each $\mu\in U\backslash\left\{  \mu_{0}\right\}  $%
\begin{equation}
\sup_{t\in\mathbb{R}}\frac{r_{\mu}\left(  t\right)  }{r_{\mu_{0}}\left(
t\right)  }<\infty\label{eq6x}%
\end{equation}
and%
\begin{equation}
\lim_{t\rightarrow\infty}\frac{1}{t}\int_{\left[  -t,t\right]  }\frac{\hat
{F}_{\mu}\left(  y\right)  }{\hat{F}_{\mu_{0}}\left(  y\right)  }dy=0,
\label{eq6dx}%
\end{equation}
then $\mathcal{\hat{F}}=\left\{  \hat{F}_{\mu}:\mu\in U\right\}  $ is said to
be of \textquotedblleft Riemann-Lebesgue type (RL type)\textquotedblright\ (at
$\mu_{0}$ on $U$).

\item Definition 2 of \cite{Chen:2018a}:\ If a function $\omega:\left[
-1,1\right]  \rightarrow\mathbb{R}$ is non-negative and bounded such that
$\int_{\left[  -1,1\right]  }\omega\left(  s\right)  ds=1$, then it is called
\textquotedblleft admissible\textquotedblright. If additionally $\omega$ is
even on $\left[  -1,1\right]  $ and continuous on $\left(  -1,1\right)  $,
then it is called \textquotedblleft eligible\textquotedblright. If $\omega$ is
eligible and $\omega\left(  t\right)  \leq\tilde{\omega}\left(  1-t\right)  $
for all $t\in\left(  0,1\right)  $ for some convex, super-additive function
$\tilde{\omega}$ over $\left(  0,1\right)  $, then it is called
\textquotedblleft good\textquotedblright.
\end{enumerate}

Recall%
\[
\mathcal{B}_{m}\left(  \rho_{m}\right)  =\left\{  \boldsymbol{\mu}%
\in\mathbb{R}^{m}:m^{-1}\sum\nolimits_{i=1}^{m}\left\vert \mu_{i}\right\vert
\leq\rho_{m}\right\}  \text{ \ for some }\rho_{m}>0
\]
and define $u_{m,0}=\min\left\{  \left\vert \mu_{j}-\mu_{0}\right\vert
:\mu_{j}\neq\mu_{0}\right\}  $.

Theorem 3 of \cite{Chen:2018a} is stated as follows only with a slight
modification to term \textquotedblleft$R_{0}\left(  \rho_{m},\mu_{0}\right)
$\textquotedblright\ and the additional requirement that $\omega$ is good:

\begin{theorem}
[Theorem 3 of \cite{Chen:2018a}]\label{ThmLocationShiftSimple}Assume that
$\omega$ is good, and that $\mathcal{F}$ is a location-shift family for which
(\ref{eq6bx}) holds and $\int\left\vert x\right\vert ^{2}dF_{\mu}\left(
x\right)  <\infty$ for each $\mu\in U$. If $\sup_{y\in\mathbb{R}}\left\vert
\frac{d}{dy}h_{\mu_{0}}\left(  y\right)  \right\vert =C_{\mu_{0}}<\infty$,
then for the estimator $\hat{\varphi}_{m}\left(  t,\mathbf{z}\right)  $ from
Construction I, a uniform consistency class is
\[
\mathcal{Q}\left(  \mathcal{F}\right)  =\left\{
\begin{array}
[c]{c}%
q\gamma^{\prime}>\vartheta>2^{-1},\gamma^{\prime}>0,\gamma^{\prime\prime
}>0,0\leq\vartheta^{\prime}<\vartheta-1/2,\\
R_{0}\left(  \rho_{m},\mu_{0}\right)  =O\left(  m^{\vartheta^{\prime}}\right)
,\tau_{m}\leq\gamma_{m},u_{m,0}\geq\frac{\ln\ln m}{\gamma^{\prime\prime}%
\tau_{m}},\\
t\in\left[  0,\tau_{m}\right]  ,\lim\limits_{m\rightarrow\infty}
\textcolor{blue}{\sup_{\boldsymbol{\mu}\in\mathcal{B}_{m}\left(  \rho_{m}\right)
}}\pi_{1,m}%
^{-1}\tilde{\Upsilon}\left(  q,\tau_{m},\gamma_{m},r_{\mu_{0}}\right)  =0
\end{array}
\right\}  \label{MainClassx}%
\]
where $q$, $\gamma^{\prime}$, $\gamma^{\prime\prime}$, $\vartheta$ and
$\vartheta^{\prime}$ are constants, $R_{0}\left(  \rho_{m},\mu_{0}\right)
=2\int{\left\vert x\right\vert dF_{\mu_{0}}\left(  x\right)  }+ 2\rho
_{m}+2C_{\mu_{0}}+\textcolor{blue}{2}\left\vert \mu_{0}\right\vert $,
$\gamma_{m}=\gamma^{\prime}\ln m$ and%
\[
\tilde{\Upsilon}\left(  q,\tau_{m},\gamma_{m},r_{\mu_{0}}\right)
=\frac{2\left\Vert \omega\right\Vert _{\infty}\sqrt{2q\gamma_{m}}}{\sqrt{m}%
}\sup_{t\in\left[  0,\tau_{m}\right]  }\int_{\left[  0,1\right]  }\frac
{ds}{r_{\mu_{0}}\left(  ts\right)  }.\label{eq12gx}%
\]
Moreover, for all sufficiently large $m$,
\begin{equation}
\sup_{\boldsymbol{\mu}\in\mathcal{B}_{m}\left(  \rho_{m}\right)  }\sup
_{t\in\left[  0,\tau_{m}\right]  }\left\vert \hat{\varphi}_{0,m}\left(
t,\mathbf{z}\right)  -\varphi_{0,m}\left(  t,\boldsymbol{\mu}\right)
\right\vert \leq\tilde{\Upsilon}\left(  q,\tau_{m},\gamma_{m},r_{\mu_{0}%
}\right)  \label{eq11dx}%
\end{equation}
holds with probability at least $1-o\left(  1\right)  $.
\end{theorem}

%

\begin{table}[tbp] \centering
\begin{tabular}
[c]{|l|l|}\hline
\multicolumn{2}{|c|}{Notations with the \textbf{same} definitions}\\\hline
Notations here & Notations in proof Theorem 3 \cite{Chen:2018a}\\\hline
$\hat{\varphi}_{0,m}\left(  t,\mathbf{z}\right)  $ & $\hat{\varphi}_{m}\left(
t,\mathbf{z}\right)  $\\\hline
$\varphi_{0,m}\left(  t,\boldsymbol{\mu}\right)  $ & $\varphi_{m}\left(
t,\boldsymbol{\mu}\right)  $\\\hline
$g\left(  t;\mu_{0}\right)  $ & $a\left(  t;\mu_{0}\right)  $\\\hline
$\mathcal{B}_{m}\left(  \rho_{m}\right)  $ & $\mathcal{B}_{m}\left(
\rho\right)  $\\\hline
$X_{\left(  \mu_{0}\right)  }$ & $X_{1}$\\\hline
$u_{m,0}$ & $u_{m}$\\\hline
$\mathcal{Q}\left(  \mathcal{F}\right)  $ & $\mathcal{Q}_{m}\left(
\boldsymbol{\mu},t;\mathcal{F}\right)  $\\\hline
\multicolumn{2}{|c|}{Notations with the \textbf{different} definitions}%
\\\hline
$A_{0,m}=\Pr\left(  \max\limits_{y\in G_{m}}\left\vert S_{m}\left(  y\right)
\right\vert \geq\frac{\sqrt{2q\gamma_{m}}}{\sqrt{m}}\right)  $ & $B_{0}%
=\Pr\left(  \sup\limits_{\boldsymbol{\mu}\in\mathcal{B}_{m}\left(
\rho\right)  }\max\limits_{y\in G_{m}}\left\vert S_{m}\left(  y\right)
\right\vert \geq\frac{\sqrt{2q\gamma_{m}}}{\sqrt{m}}\right)  $\\\hline
$A_{1,m}=\Pr\left(  \max\limits_{1\leq i\leq l_{\ast}}\left\vert
d_{m}^{\left(  0\right)  }\left(  y_{i}\right)  \right\vert \geq\tilde
{p}_{0,m,1}\right)  $ & $B_{1}=\Pr\left(  \sup\limits_{\boldsymbol{\mu}%
\in\mathcal{B}_{m}\left(  \rho\right)  }\max\limits_{1\leq i\leq l_{\ast}%
}\left\vert d_{m}^{\left(  0\right)  }\left(  y_{i}\right)  \right\vert
\geq\tilde{p}_{0,m,1}\right)  $\\\hline
$A_{2,m}=\Pr\left(  \sup\limits_{\boldsymbol{\mu}\in\mathcal{B}_{m}\left(
\rho\right)  }\sup\limits_{y\in\mathbb{R}}\left\vert \partial_{y}%
d_{m}^{\left(  0\right)  }\left(  y\right)  \right\vert \geq\tilde{p}%
_{0,m,2}\right)  $ & $B_{2}=\Pr\left(  \sup\limits_{\boldsymbol{\mu}%
\in\mathcal{B}_{m}\left(  \rho\right)  }\sup_{y\in\mathbb{R}}\left\vert
\partial_{y}d_{m}^{\left(  0\right)  }\left(  y\right)  \right\vert \geq
\tilde{p}_{0,m,2}\right)  $\\\hline
$R_{0}\left(  \rho_{m},\mu_{0}\right)  $ & $R_{m}\left(  \rho\right)
$\\\hline
\end{tabular}
\caption{Comparison between notations used here and those in the Proof of Theorem 3 of \cite{Chen:2018a}.}\label{NotationComparison}%
\end{table}%

In order to present a correction to the two minor issues in the proof Theorem
3 in \cite{Chen:2018a} in notations that are compatible with this manuscript
and \cite{Chen:2018a}, we provide Table \ref{NotationComparison}. In Table
\ref{NotationComparison}, we have $\gamma_{m}=\gamma^{\prime}\ln m$,
$d_{m}^{\left(  0\right)  }\left(  y\right)  =\hat{s}_{m}\left(  y\right)
-s_{m}\left(  y\right)  $,%
\[
\tilde{p}_{0,m,1}=\frac{\sqrt{2q\gamma_{m}}-\left(  2q\gamma_{m}\right)
^{-1/2}\ln\gamma_{m}}{\sqrt{m}}\text{ \ and \ }\tilde{p}_{0,m,2}=\frac
{\Delta^{-1}\left(  2q\gamma_{m}\right)  ^{-1/2}\ln\gamma_{m}}{\sqrt{m}}.
\]
The notations and quantities that have appeared above and in Table
\ref{NotationComparison} but not in \cite{Chen:2018a} will be specified later.
With Table \ref{NotationComparison}, the two minor issues mentioned above can
be stated as follows:\

\begin{itemize}
\item In Part I of the proof of Theorem 3 in \cite{Chen:2018a}, the author
insufficiently handled the event%
\[
\left\{  \sup\limits_{\boldsymbol{\mu}\in\mathcal{B}_{m}\left(  \rho
_{m}\right)  }\max\limits_{1\leq i\leq l_{\ast}}\left\vert d_{m}^{\left(
0\right)  }\left(  y_{i}\right)  \right\vert \geq\tilde{p}_{0,m,1}\right\}  ,
\]
which led to insufficient handling of the event%
\[
\left\{  \sup_{\boldsymbol{\mu}\in\mathcal{B}_{m}\left(  \rho_{m}\right)
}\max_{y\in G_{m}}\left\vert S_{m}\left(  y\right)  \right\vert \geq
\frac{\sqrt{2q\gamma_{m}}}{\sqrt{m}}\right\}  .
\]
Therefore, a correction to this issue is to just handle the event%
\[
\left\{  \max_{y\in G_{m}}\left\vert S_{m}\left(  y\right)  \right\vert
\geq\frac{\sqrt{2q\gamma_{m}}}{\sqrt{m}}\right\}  .
\]

\item In Part III of the proof of Theorem 3 in \cite{Chen:2018a}, the author
insufficiently handled the ratio%
\[
\sup_{t\in\left[  0,\tau_{m}\right]  }\left\vert \pi_{1,m}^{-1}\varphi
_{0,m}\left(  t,\boldsymbol{\mu}\right)  -1\right\vert ,
\]
where $\lim_{t\rightarrow\infty}\varphi_{0,m}\left(  t,\boldsymbol{\mu
}\right)  =\pi_{1,m}$ for all $\boldsymbol{\mu}$ and $m$ is called the
\textquotedblleft Oracle\textquotedblright.
\end{itemize}

We present the corrected proof below, which differs much from Part III of the
original proof of Theorem 3 in \cite{Chen:2018a}.

\subsection{Proof of \autoref{ThmLocationShiftSimple}}

Set $w_{i}\left(  y\right)  =\cos\left(  yz_{i}-h_{\mu_{0}}\left(  y\right)
\right)  $ for each $i$ and $y\in\mathbb{R}$ and define%
\begin{equation}
S_{m}\left(  y\right)  =\frac{1}{m}\sum_{i=1}^{m}\left(  w_{i}\left(
y\right)  -\mathbb{E}\left[  w_{i}\left(  y\right)  \right]  \right)  .
\label{eq2ex}%
\end{equation}
Let $\hat{s}_{m}\left(  y\right)  =m^{-1}\sum_{i=1}^{m}w_{i}\left(  y\right)
$ and $s_{m}\left(  y\right)  =\mathbb{E}\left[  \hat{s}_{m}\left(  y\right)
\right]  $. Then $S_{m}\left(  y\right)  =\hat{s}_{m}\left(  y\right)
-s_{m}\left(  y\right)  $. For the rest of the proof, we will first assume the
existence of the positive constants $\gamma^{\prime}$, $\gamma^{\prime\prime}%
$, $q$, $\vartheta$ and the non-negative constant $\vartheta^{\prime}$ and
then determine them at the end of the proof. Let $\gamma_{m}=\gamma^{\prime
}\ln m$. The rest of the proof is divided into three parts.

\textbf{Part I}: to show the assertion \textquotedblleft if%
\begin{equation}
\lim_{m\rightarrow\infty}\frac{m^{\vartheta}\ln\gamma_{m}}{R_{0}\left(
\rho_{m},\mu_{0}\right)  \sqrt{m}\sqrt{2q\gamma_{m}}}=\infty\label{eq11bx}%
\end{equation}
where $R_{0}\left(  \rho_{m},\mu_{0}\right)  =2\mathbb{E}\left[  \left\vert
X_{\left(  \mu_{0}\right)  }\right\vert \right]  +2\rho_{m}+2C_{\mu_{0}%
}+2\left\vert \mu_{0}\right\vert $ and $X_{\left(  \mu_{0}\right)  }$ has CDF
$F_{\mu_{0}}$, then, for all large $m$,%
\begin{equation}
\sup_{y\in\left[  \textcolor{blue}{-\gamma_m},\gamma_{m}\right]  }\left\vert
\hat{s}_{m}\left(  y\right)  -\mathbb{E}\left[  \hat{s}_{m}\left(  y\right)
\right]  \right\vert \leq\frac{\sqrt{2q\gamma_{m}}}{\sqrt{m}} \label{eq11jx}%
\end{equation}
holds with probability at least $1-p_{m}\left(  \vartheta,q,h_{\mu_{0}}%
,\gamma_{m}\right)  $, where%
\begin{equation}
p_{m}\left(  \vartheta,q,h_{\mu_{0}},\gamma_{m}\right)  =2 \left(
\textcolor{blue}{2} m^{\vartheta}\gamma_{m}\textcolor{blue}{+1}\right)
\exp\left(  -q\gamma_{m}\right)  +\textcolor{blue}{4}A_{\mu_{0}}q\gamma
_{m}m^{-2\vartheta}\left(  \ln\gamma_{m}\right)  ^{-2} \label{eq11ex}%
\end{equation}
and $A_{\mu_{0}}$ is the variance of $\left\vert X_{\left(  \mu_{0}\right)
}\right\vert $\textquotedblright.

Define the closed interval $G_{m}=\left[  \textcolor{blue}{-\gamma_m},\gamma
_{m}\right]  $. Let $\mathcal{K}=\left\{  y_{1},\ldots,y_{l_{\ast}}\right\}  $
for some $l_{\ast}\in\mathbb{N}_{+}$ with $y_{j}<y_{j+1}$ be a partition of
$G_{m}$ with norm $\Delta=\max_{1\leq j\leq l_{\ast}-1}\left\vert
y_{j+1}-y_{j}\right\vert $ such that $\Delta=m^{-\vartheta}$
\textcolor{blue}{and that $G_m$ is split into $l_{\ast}-1$
sub-intervals of equal lengths $\Delta$}. For each $y\in G_{m}$, pick
$y_{j}\in\mathcal{K}$ that is the closest to $y$. By Lagrange mean value
theorem,%
\begin{align*}
\left\vert \hat{s}_{m}\left(  y\right)  -s_{m}\left(  y\right)  \right\vert
&  \leq\left\vert \hat{s}_{m}\left(  y_{i}\right)  -s_{m}\left(  y_{i}\right)
\right\vert +\left\vert \left(  \hat{s}_{m}\left(  y\right)  -\hat{s}%
_{m}\left(  y_{i}\right)  \right)  -\left(  s_{m}\left(  y\right)
-s_{m}\left(  y_{i}\right)  \right)  \right\vert \\
&  \leq\left\vert \hat{s}_{m}\left(  y_{i}\right)  -s_{m}\left(  y_{i}\right)
\right\vert +\Delta\sup_{y\in\mathbb{R}}\left\vert \partial_{y}\left(  \hat
{s}_{m}\left(  y\right)  -s_{m}\left(  y\right)  \right)  \right\vert ,
\end{align*}
i.e.,%
\[
\left\vert S_{m}\left(  y\right)  \right\vert \leq\left\vert S_{m}\left(
y_{i}\right)  \right\vert +\Delta\sup_{y\in\mathbb{R}}\left\vert \partial
_{y}S_{m}\left(  y\right)  \right\vert ,
\]
where $\partial_{\cdot}$ denotes the derivative with respect to the subscript.
Set%
\[
\tilde{p}_{0,m,1}=\frac{\sqrt{2q\gamma_{m}}-\left(  2q\gamma_{m}\right)
^{-1/2}\ln\gamma_{m}}{\sqrt{m}}\text{ \ and \ \ }\tilde{p}_{0,m,2}%
=\frac{\Delta^{-1}\left(  2q\gamma_{m}\right)  ^{-1/2}\ln\gamma_{m}}{\sqrt{m}%
}.
\]
Then,%
\begin{equation}
A_{0,m}=\Pr\left(  \max_{y\in G_{m}}\left\vert S_{m}\left(  y\right)
\right\vert \geq\frac{\sqrt{2q\gamma_{m}}}{\sqrt{m}}\right)  \leq
A_{1,m}+A_{2,m}, \label{eq10hx}%
\end{equation}
where%
\[
A_{1,m}=\Pr\left(  \max_{1\leq i\leq l_{\ast}}\left\vert S_{m}\left(
y_{i}\right)  \right\vert \geq\tilde{p}_{0,m,1}\right)
\]
and%
\[
A_{2,m}=\Pr\left(  \sup_{y\in\mathbb{R}}\left\vert \partial_{y}S_{m}\left(
y\right)  \right\vert \geq\tilde{p}_{0,m,2}\right)  .
\]

Since
\textcolor{blue}{$\left\vert w_{i}\left(  y\right)  \right\vert \leq1$ uniformly in
$y \in \mathbb{R}  $ for all $i$} and $\left\{  z_{i}\right\}  _{i=1}^{m}$ are
independent, Hoeffding inequality of \cite{Hoeffding:1963} implies%
\begin{equation}
\Pr\left(  \left\vert S_{m}\left(  \textcolor{blue}{y_j}\right)  \right\vert
\geq\frac{\lambda}{\sqrt{m}}\right)  \leq2\exp\left(  -2^{-1}\lambda
^{2}\right)  \text{ for any }\lambda>0 \label{eq2cx}%
\end{equation}
\textcolor{blue}{for all $m$ and each fixed $j$ between $1$ and $l_{\ast}$}.
Applying to $A_{1,m}$ the union bound and Hoeffding inequality gives%
\begin{equation}
A_{1,m}\leq2l_{\ast}\exp\left(  -q\gamma_{m}+\ln\gamma_{m}\right)  \exp\left(
-\frac{\left(  \ln\gamma_{m}\right)  ^{2}}{4q\gamma_{m}}\right)  \leq2 \left(
\textcolor{blue}{2} m^{\vartheta}\gamma_{m}\textcolor{blue}{+1}\right)
\exp\left(  -q\gamma_{m}\right)  . \label{eq10gx}%
\end{equation}

On the other hand, since $r_{\mu_{0}}$ has no real zeros and $\mathcal{F}$ is
a location-shift family, then $r_{\mu}\left(  t\right)  $ has no real zeros
for each $\mu\neq\mu_{0}$ and the argument $h_{\mu}\left(  t\right)  $ is
well-defined and continuous in $t$ on $\mathbb{R}$ for each $\mu\in U$.
\textcolor{blue}{Since $\int\left\vert
x\right\vert ^{2}dF_{\mu}\left(  x\right)  <\infty$ for each $\mu\in U$, we
know that for each fixed $\mu\in U$, the CF $\hat{F}_{\mu}\left(  t\right)
,t\in\mathbb{R}$ is differentiable in $t\in\mathbb{R}$, and hence a branch of
$h_{\mu}\left(  t\right)  $ that is differentiable in $t\in\mathbb{R}$ can be
chosen.} Therefore, $\frac{d}{dy}h_{\mu_{0}}\left(  y\right)  $
\textcolor{blue}{is well-defined}. Further, $\partial_{y}w_{i}\left(
y\right)  =-\left(  z_{i}-\partial_{y}h_{\mu_{0}}\left(  y\right)  \right)
\sin\left(  yz_{i}-h_{\mu_{0}}\left(  y\right)  \right)  $, and%
\[
\partial_{y}\mathbb{E}\left[  w_{i}\left(  y\right)  \right]  =\mathbb{E}%
\left[  \partial_{y}w_{i}\left(  y\right)  \right]  =-\mathbb{E}\left[
\left(  z_{i}-\partial_{y}h_{\mu_{0}}\left(  y\right)  \right)  \sin\left(
yz_{i}-h_{\mu_{0}}\left(  y\right)  \right)  \right]
\]
holds since $\int\left\vert x\right\vert ^{2}dF_{\mu}\left(  x\right)
<\infty$ for each $\mu\in U$ and $\sup_{y\in\mathbb{R}}\left\vert \frac{d}%
{dy}h_{\mu_{0}}\left(  y\right)  \right\vert =C_{\mu_{0}}<\infty$. So,%
\begin{equation}
\sup_{y\in\mathbb{R}}\left\vert \partial_{y}S_{m}\left(  y\right)  \right\vert
\leq\frac{1}{m}\sum_{i=1}^{m}\left\vert z_{i}\right\vert +2C_{\mu_{0}}%
+\frac{1}{m}\sum_{i=1}^{m}\mathbb{E}\left[  \left\vert z_{i}\right\vert
\right]  . \label{eq10dx}%
\end{equation}

Since $\mathcal{F}$ is a location-shift family, there are independent and
identically distributed (i.i.d.) $\left\{  X_{i}\right\}  _{i=1}^{m}$ with
common CDF $F_{\mu_{0}}$ such that $z_{i}=\left(  \mu_{i}-\mu_{0}\right)
+X_{i}$ for $1\leq i\leq m$. Therefore, the upper bound in (\ref{eq10dx})
satisfies%
\begin{align}
&  \frac{1}{m}\sum_{i=1}^{m}\left\vert z_{i}\right\vert +2C_{\mu_{0}}+\frac
{1}{m}\sum_{i=1}^{m}\mathbb{E}\left[  \left\vert z_{i}\right\vert \right]
\nonumber\\
&  \leq\frac{1}{m}\sum_{i=1}^{m}\left\vert X_{i}\right\vert +\frac{2}{m}%
\sum_{i=1}^{m}\left\vert \mu_{i}-\mu_{0}\right\vert +2C_{\mu_{0}}+\frac{1}%
{m}\sum_{i=1}^{m}\mathbb{E}\left[  \left\vert X_{i}\right\vert \right]
\nonumber\\
&  \leq\frac{1}{m}\sum_{i=1}^{m}\left(  \left\vert X_{i}\right\vert
-\mathbb{E}\left[  \left\vert X_{i}\right\vert \right]  \right)  +2\left\vert
\mu_{0}\right\vert +2\mathbb{E}\left[  \left\vert X_{\left(  \mu_{0}\right)
}\right\vert \right]  +2C_{\mu_{0}}+\frac{2}{m}\sum_{i=1}^{m}\left\vert
\mu_{i}\right\vert , \label{eq10dx1}%
\end{align}
where $X_{\left(  \mu_{0}\right)  }$ has CDF $F_{\mu_{0}}$. Set $\tilde{T}%
_{m}=m^{-1}\sum_{i=1}^{m}\left(  \left\vert X_{i}\right\vert -\mathbb{E}%
\left[  \left\vert X_{i}\right\vert \right]  \right)  $. Then (\ref{eq10dx})
and (\ref{eq10dx1}) imply the following: if $\boldsymbol{\mu}\in
\mathcal{B}_{m}\left(  \rho_{m}\right)  $, i.e., $m^{-1}\sum_{i=1}%
^{m}\left\vert \mu_{i}\right\vert \leq\rho_{m}$, then
\[
\sup_{y\in\mathbb{R}}\left\vert \partial_{y}S_{m}\left(  y\right)  \right\vert
\leq\tilde{T}_{m}+R_{0}\left(  \rho_{m},\mu_{0}\right)  ,
\]
and%
\[
\sup_{\boldsymbol{\mu}\in\mathcal{B}_{m}\left(  \rho_{m}\right)  }\sup
_{y\in\mathbb{R}}\left\vert \partial_{y}S_{m}\left(  y\right)  \right\vert
\leq\vert\tilde{T}_{m}+R_{0}\left(  \rho_{m},\mu_{0}\right)  ,
\]
where we recall%
\[
R_{0}\left(  \rho_{m},\mu_{0}\right)  =2\mathbb{E}\left[  \left\vert
X_{\left(  \mu_{0}\right)  }\right\vert \right]  +2\rho_{m}+2C_{\mu_{0}%
}+2\left\vert \mu_{0}\right\vert .
\]
When (\ref{eq11bx}) holds, i.e., $\lim_{m\rightarrow\infty}\tilde{p}%
_{0,m,2}/R_{0}\left(  \rho_{m},\mu_{0}\right)  =\infty$,
\textcolor{blue}{we must have, for all $m$ large enough,\[
\tilde{p}_{0,m,2}-R_{0}\left(  \rho_{m},\mu_{0}\right)  >0\text{ and
\ }1-\frac{R_{0}\left(  \rho_{m},\mu_{0}\right)  }{\tilde{p}_{0,m,2}}\geq
\sqrt{2^{-1}}\text{.}\]
Therefore, for all $m$ large enough, $\tilde{T}_{m}=\left\vert \tilde{T}_{m}\right\vert $ on the event $\left\{  \tilde{T}_{m}\geq\tilde{p}_{0,m,2}-R_{0}\left(  \rho_{m},\mu_{0}\right)  \right\}  $ and}
Chebyshev inequality implies%
\begin{align*}
A_{2,m,1}  &  =\Pr\left(  \textcolor{blue}{\tilde{T}_m} \geq\tilde{p}%
_{0,m,2}-R_{0}\left(  \rho_{m},\mu_{0}\right)  \right)
\textcolor{blue}{\leq\frac{A_{\mu_{0}}}{m\tilde{p}_{0,m,2}^{2}\left( 1-\frac{R_{0}\left( \rho_{m},\mu_{0}\right) }{\tilde {p}_{0,m,2}}\right) ^{2}}}\\
&  \leq\textcolor{blue}{\frac{2A_{\mu_{0}}}{m\tilde{p}_{0,m,2}^{2}}=4}A_{\mu
_{0}}m^{-2\vartheta}q\gamma_{m}\left(  \log\gamma_{m}\right)  ^{-2},
\end{align*}
where $A_{\mu_{0}}$ is the variance of $\left\vert X_{\left(  \mu_{0}\right)
}\right\vert $. Thus, for all $m$ large enough,%
\begin{equation}
\max\left\{  A_{3,m},A_{2,m}\right\}  \leq A_{2,m,1}\leq
\textcolor{blue}{4}A_{\mu_{0}}q\gamma_{m}m^{-2\vartheta}\left(  \ln\gamma
_{m}\right)  ^{-2}, \label{eq10gx1}%
\end{equation}
where%
\[
A_{3,m}=\Pr\left(  \sup_{\boldsymbol{\mu}\in\mathcal{B}_{m}\left(  \rho
_{m}\right)  }\sup_{y\in\mathbb{R}}\left\vert \partial_{y}S_{m}\left(
y\right)  \right\vert \geq\tilde{p}_{0,m,2}\right)  .
\]
Combining (\ref{eq10gx}), (\ref{eq10hx}), (\ref{eq10gx1}) gives%
\[
A_{0,m}=\Pr\left(  \max_{y\in G_{m}}\left\vert S_{m}\left(  y\right)
\right\vert \geq\frac{\sqrt{2q\gamma_{m}}}{\sqrt{m}}\right)  \leq p_{m}\left(
\vartheta,q,h_{\mu_{0}},\gamma_{m}\right)
\]
for all $m$ large enough. This justifies the assertion.

\textcolor{blue}{In fact, we have also achieved a stronger conclusion, explained as follows. When $s\in\left[  -1,1\right]  $, we have $\left\{
ts:s\in\left[  -1,1\right]  ,t\in G_{m}\right\}  =G_{m}=\left[  -\gamma
_{m},\gamma_{m}\right]  $. So, the continuity of $S_{m}\left(  \cdot\right)  $
implies $$\max_{\left(  s,t\right)  \in\left[  -1,1\right]  \times G_{m}}\left\vert S_{m}\left(  ts\right)  \right\vert =\max_{y\in G_{m}}\left\vert
S_{m}\left(  y\right)  \right\vert .$$ Therefore,\[
\Pr\left(  \max_{\left(  s,t\right)  \in\left[  -1,1\right]  \times G_{m}}\left\vert S_{m}\left(  ts\right)  \right\vert \geq\frac{\sqrt{2q\gamma_{m}}}{\sqrt{m}}\right)  \leq p_{m}\left(  \vartheta,q,h_{\mu_{0}},\gamma
_{m}\right)
\]
for all $m$ large enough.}

\textbf{Part II}: to show the uniform bound on $\left\vert \hat{\varphi}%
_{0,m}\left(  t,\mathbf{z}\right)  -\varphi_{0,m}\left(  t,\boldsymbol{\mu
}\right)  \right\vert $. Recall
\[
K\left(  t,x;\mu_{0}\right)  =\int_{\left[  -1,1\right]  }\frac{\omega\left(
s\right)  \cos\left(  tsx-h_{\mu_{0}}\left(  ts\right)  \right)  }{r_{\mu_{0}%
}\left(  ts\right)  }ds,
\]
$\hat{\varphi}_{0,m}\left(  t,\mathbf{z}\right)  =m^{-1}\sum_{i=1}^{m}\left(
1-K\left(  t,z_{i};\mu_{0}\right)  \right)  $ and $\varphi_{0,m}\left(
t,\boldsymbol{\mu}\right)  =\mathbb{E}\left[  \hat{\varphi}_{0,m}\left(
t,\mathbf{z}\right)  \right]  $. Set $e_{0,m}\left(  t,\mathbf{z}%
,\boldsymbol{\mu}\right)  =\hat{\varphi}_{0,m}\left(  t,\mathbf{z}\right)
-\varphi_{0,m}\left(  t,\boldsymbol{\mu}\right)  $. Then%
\[
e_{0,m}\left(  t,\mathbf{z},\boldsymbol{\mu}\right)  =\dfrac{-1}{m}\sum
_{i=1}^{m}\left(  K\left(  t,z_{i};\mu_{0}\right)  -\mathbb{E}\left[  K\left(
t,z_{i};\mu_{0}\right)  \right]  \right)  =-\int_{\left[  -1,1\right]  }%
\frac{\omega\left(  s\right)  }{r_{\mu_{0}}\left(  ts\right)  }S_{m}\left(
ts\right)  ds.
\]
Pick a positive sequence $\left\{  \tau_{m}:m\geq1\right\}  $ such that
$\tau_{m}\leq\gamma_{m}$ for all large $m$ and $\tau_{m}\rightarrow\infty$,
and set
\[
\tilde{\Upsilon}\left(  q,\tau_{m},\gamma_{m},r_{\mu_{0}}\right)
=\frac{2\left\Vert \omega\right\Vert _{\infty}\sqrt{2q\gamma_{m}}}{\sqrt{m}%
}\sup_{t\in\left[  0,\tau_{m}\right]  }\int_{\left[  0,1\right]  }\frac
{ds}{r_{\mu_{0}}\left(  ts\right)  }.
\]
Then, \textbf{Part I} implies that, with probability at least $1-p_{m}\left(
\vartheta,q,h_{\mu_{0}},\gamma_{m}\right)  $,
\begin{align}
&  \sup_{t\in\left[  0,\tau_{m}\right]  }\left\vert \hat{\varphi}_{0,m}\left(
t,\mathbf{z}\right)  -\varphi_{0,m}\left(  t,\boldsymbol{\mu}\right)
\right\vert \nonumber\\
&  \leq\sup_{t\in\left[  0,\tau_{m}\right]  }\int_{\left[  -1,1\right]
}\omega\left(  s\right)  \frac{\sup_{
\textcolor{blue}{\left(s,t\right) \in\left[ -1,1\right] \times \left[ 0,\gamma_m\right] }}%
\left\vert S_{m}\left(  ts\right)  \right\vert }{r_{\mu_{0}}\left(  ts\right)
}ds\leq\tilde{\Upsilon}\left(  q,\tau_{m},\gamma_{m},r_{\mu_{0}}\right)
\label{eq3x1}%
\end{align}
for all sufficiently large $m$ and, since the right-hand side of the above
inequality is independent of $\boldsymbol{\mu}$, we also get, with probability
at least $1-p_{m}\left(  \vartheta,q,h_{\mu_{0}},\gamma_{m}\right)  $,
\begin{align}
&  \sup_{\boldsymbol{\mu}\in\mathcal{B}_{m}\left(  \rho_{m}\right)  }\left(
\sup_{t\in\left[  0,\tau_{m}\right]  }\left\vert \hat{\varphi}_{0,m}\left(
t,\mathbf{z}\right)  -\varphi_{0,m}\left(  t,\boldsymbol{\mu}\right)
\right\vert \right) \nonumber\\
&  \leq\sup_{\boldsymbol{\mu}\in\mathcal{B}_{m}\left(  \rho_{m}\right)  }%
\sup_{t\in\left[  0,\tau_{m}\right]  }\int_{\left[  -1,1\right]  }%
\omega\left(  s\right)  \frac{\sup
_{\textcolor{blue}{\left(s,t\right) \in\left[ -1,1\right] \times \left[ 0,\gamma_m\right] }}%
\left\vert S_{m}\left(  ts\right)  \right\vert }{r_{\mu_{0}}\left(  ts\right)
}ds\leq\tilde{\Upsilon}\left(  q,\tau_{m},\gamma_{m},r_{\mu_{0}}\right)
\label{eq2ix2}%
\end{align}
for all sufficiently large $m$.

\textbf{Part III}: to determine the constants $\gamma^{\prime}$,
$\gamma^{\prime\prime}$, $q$, $\vartheta$ and $\vartheta^{\prime}$ and a
uniform consistency class. Set $\gamma^{\prime}$, $\vartheta$ and $q$ such
that $q\gamma^{\prime}>\vartheta$ and $0\leq\vartheta^{\prime}<\vartheta-1/2$.
Then $p_{m}\left(  \vartheta,q,h_{\mu_{0}},\gamma_{m}\right)  \rightarrow0$,
$A_{0,m}\rightarrow0$ and $m^{\vartheta-1/2}\gamma_{m}^{-1/2}\ln\gamma
_{m}\rightarrow\infty$ as $m\rightarrow\infty$. If additionally $R_{0}\left(
\rho_{m},\mu_{0}\right)  =O\left(  m^{\vartheta^{\prime}}\right)  $ and
$u_{m,0}\geq\left(  \gamma^{\prime\prime}\tau_{m}\right)  ^{-1}\ln\ln m$. Then
$\tau_{m}u_{m,0}\rightarrow\infty$ as $m\rightarrow\infty$ and (\ref{eq11bx}) holds.

Recall%
\[
\psi\left(  t,\mu_{i};\mu_{0}\right)  =\int_{\left[  -1,1\right]  }%
\omega\left(  s\right)  \cos\left(  ts\left(  \mu_{i}-\mu_{0}\right)  \right)
ds.
\]
We obtain
\begin{align*}
\left\vert \frac{\sup_{t\in\left[  0,\tau_{m}\right]  }\hat{\varphi}%
_{0,m}\left(  t,\mathbf{z}\right)  }{\pi_{1,m}}-1\right\vert  &  =\left\vert
\sup_{t\in\left[  0,\tau_{m}\right]  }\left(  \frac{\hat{\varphi}_{0,m}\left(
t,\mathbf{z}\right)  -\varphi_{0,m}\left(  t,\boldsymbol{\mu}\right)  }%
{\pi_{1,m}}+\left(  \frac{\varphi_{0,m}\left(  t,\boldsymbol{\mu}\right)
}{\pi_{1,m}}-1\right)  \right)  \right\vert \\
&  \leq\left\vert \sup_{t\in\left[  0,\tau_{m}\right]  }\left(  \frac
{\hat{\varphi}_{0,m}\left(  t,\mathbf{z}\right)  -\varphi_{0,m}\left(
t,\boldsymbol{\mu}\right)  }{\pi_{1,m}}\right)  \right\vert +\left\vert
\sup_{t\in\left[  0,\tau_{m}\right]  }\left(  \frac{\varphi_{0,m}\left(
t,\boldsymbol{\mu}\right)  }{\pi_{1,m}}-1\right)  \right\vert
\end{align*}
and%
\begin{align}
\left\vert \pi_{1,m}^{-1}\sup\nolimits_{t\in\left[  0,\tau_{m}\right]  }%
\hat{\varphi}_{0,m}\left(  t,\mathbf{z}\right)  -1\right\vert  &
\leq\underbrace{\pi_{1,m}^{-1}\sup_{t\in\left[  0,\tau_{m}\right]  }\left\vert
\hat{\varphi}_{0,m}\left(  t,\mathbf{z}\right)  -\varphi_{0,m}\left(
t,\boldsymbol{\mu}\right)  \right\vert }_{=\epsilon_{0,m,1}\left(  \tau
_{m,}\boldsymbol{\mu}\right)  }\nonumber\\
&  +\underbrace{\left\vert \pi_{1,m}^{-1}\sup\nolimits_{t\in\left[  0,\tau
_{m}\right]  }\varphi_{0,m}\left(  t,\boldsymbol{\mu}\right)  -1\right\vert
}_{=\epsilon_{0,m,2}\left(  \tau_{m,}\boldsymbol{\mu}\right)  }. \label{eq2ix}%
\end{align}
The term $\epsilon_{0,m,1}\left(  \tau_{m,}\boldsymbol{\mu}\right)  $ in
(\ref{eq2ix}) satisfies, for all sufficiently large $m$,%
\begin{equation}
\epsilon_{0,m,1}\left(  \tau_{m,}\boldsymbol{\mu}\right)  \leq\pi_{1,m}%
^{-1}\tilde{\Upsilon}\left(  q,\tau_{m},\gamma_{m},r_{\mu_{0}}\right)
\label{eq2ix6}%
\end{equation}
due to (\ref{eq3x1}) or (\ref{eq2ix2}), and $\epsilon_{0,m,1}\left(  \tau
_{m,}\boldsymbol{\mu}\right)  \rightarrow0$ when $\pi_{1,m}^{-1}%
\tilde{\Upsilon}\left(  q,\tau_{m},\gamma_{m},r_{\mu_{0}}\right)
\rightarrow0$. For the term $\epsilon_{0,m,2}\left(  \tau_{m,}\boldsymbol{\mu
}\right)  $ in (\ref{eq2ix}), we have%
\begin{align*}
\frac{\varphi_{0,m}\left(  t,\boldsymbol{\mu}\right)  }{\pi_{1,m}}-1  &
=\frac{1}{\pi_{1,m}}\frac{1}{m}\sum\nolimits_{\left\{  i:\mu_{i}\neq\mu
_{0}\right\}  }\left(  1-\psi\left(  t,\mu_{i};\mu_{0}\right)  \right)  -1\\
&  =-1+\frac{1}{\pi_{1,m}}\frac{1}{m}\sum_{\left\{  i:\mu_{i}\neq\mu
_{0}\right\}  }1-\frac{1}{\pi_{1,m}}\frac{1}{m}\sum_{\left\{  i:\mu_{i}\neq
\mu_{0}\right\}  }\psi\left(  t,\mu_{i};\mu_{0}\right) \\
&  =-\frac{1}{m_{1}}\sum\nolimits_{\left\{  i:\mu_{i}\neq\mu_{0}\right\}
}\psi\left(  t,\mu_{i};\mu_{0}\right)  ,
\end{align*}
where $m_{1}=\left\vert I_{1,m}\right\vert $, $I_{1,m}=\left\{  1\leq i\leq
m:\mu_{i}\neq\mu_{0}\right\}  $ and $\pi_{1,m}=m_{1}m^{-1}$. Namely,%
\[
\frac{\varphi_{0,m}\left(  t,\boldsymbol{\mu}\right)  }{\pi_{1,m}}-1=-\frac
{1}{m_{1}}\sum\nolimits_{\left\{  i:\mu_{i}\neq\mu_{0}\right\}  }\psi\left(
t,\mu_{i};\mu_{0}\right)
\]
and%
\begin{equation}
\epsilon_{0,m,2}\left(  \tau_{m,}\boldsymbol{\mu}\right)  =\left\vert
\frac{\sup\nolimits_{t\in\left[  0,\tau_{m}\right]  }\varphi_{0,m}\left(
t,\boldsymbol{\mu}\right)  }{\pi_{1,m}}-1\right\vert =\left\vert -\inf
_{t\in\left[  0,\tau_{m}\right]  }\frac{1}{m_{1}}\sum_{\left\{  i:\mu_{i}%
\neq\mu_{0}\right\}  }\psi\left(  t,\mu_{i};\mu_{0}\right)  \right\vert .
\label{eq2ix3}%
\end{equation}
Since $\omega$ is good, the conclusion of Corollary 2 of \cite{Chen:2018a}
implies that $0\leq\psi\left(  t,\mu;\mu_{0}\right)  \leq1$ for all $\mu$ and
$t$. Further, \autoref{lm:OracleSpeed} implies%
\[
\left\vert \psi\left(  t,\mu_{i};\mu_{0}\right)  \right\vert \leq
\frac{4\left(  \left\Vert \omega\right\Vert _{\infty}+\left\Vert
\omega\right\Vert _{\mathrm{TV}}\right)  }{t\left\vert \mu_{i}-\mu
_{0}\right\vert }\text{ \ when }\mu_{i}\neq\mu_{0} \text{\ and } t \ne0.
\]
So, when $\mu_{i}\neq\mu_{0}$ and $\tau_{m} >0$,%
\begin{equation}
0 \leq\inf_{t\in\left[  0,\tau_{m}\right]  }
\textcolor{blue}{\frac{1}{m_{1}}\sum_{\left\{ i:\mu_{i}\neq\mu_{0}\right\} }}
\psi\left(  t,\mu_{i};\mu_{0}\right)  \leq\inf_{t\in\left[  0,\tau_{m}\right]
}\frac{4\left(  \left\Vert \omega\right\Vert _{\infty}+\left\Vert
\omega\right\Vert _{\mathrm{TV}}\right)  }{t \textcolor{blue}{u_{m,0}}
}\nonumber\\
=\frac{4\left(  \left\Vert \omega\right\Vert _{\infty}+\left\Vert
\omega\right\Vert _{\mathrm{TV}}\right)  }{\tau_{m}u_{m,0}}, \label{eq2ix7}%
\end{equation}
where we recall $u_{m,0}=\min\left\{  \left\vert \mu_{j}-\mu_{0}\right\vert
:\mu_{j}\neq\mu_{0}\right\}  $. Therefore, $\epsilon_{0,m,2}\left(  \tau
_{m,}\boldsymbol{\mu}\right)  $ in (\ref{eq2ix3}) satisfies
\begin{equation}
\epsilon_{0,m,2}\left(  \tau_{m,}\boldsymbol{\mu}\right)
=\textcolor{blue}{\inf_{t\in\left[ 0,\tau_{m}\right] }}\frac{1}{m_{1}}%
\sum_{\left\{  i:\mu_{i}\neq\mu_{0}\right\}  }\psi\left(  t,\mu_{i};\mu
_{0}\right)  \leq\frac{4\left(  \left\Vert \omega\right\Vert _{\infty
}+\left\Vert \omega\right\Vert _{\mathrm{TV}}\right)  }{\tau_{m}u_{m,0}},
\label{eq2ix4}%
\end{equation}

If $\tau_{m}u_{m,0}\rightarrow\infty$, then $\epsilon_{0,m,2}\left(  \tau
_{m,}\boldsymbol{\mu}\right)  \rightarrow0$. In summary, combining
(\ref{eq2ix}), (\ref{eq2ix6}) and (\ref{eq2ix4}) gives%
\begin{align}
\sup_{\boldsymbol{\mu}\in\mathcal{B}_{m}\left(  \rho_{m}\right)  }\left\vert
\frac{\sup_{t\in\left[  0,\tau_{m}\right]  }\hat{\varphi}_{0,m}\left(
t,\mathbf{z}\right)  }{\pi_{1,m}}-1\right\vert  &  \leq\sup_{\boldsymbol{\mu
}\in\mathcal{B}_{m}\left(  \rho_{m}\right)  }\left[  \epsilon_{0,m,1}\left(
\tau_{m,}\boldsymbol{\mu}\right)  +\epsilon_{0,m,2}\left(  \tau_{m,}%
\boldsymbol{\mu}\right)  \right] \nonumber\\
&  \leq
\textcolor{blue}{\sup_{\boldsymbol{\mu}\in\mathcal{B}_{m}\left( \rho_{m}\right) }}\left[
\frac{\tilde{\Upsilon}\left(  q,\tau_{m},\gamma_{m},r_{\mu_{0}}\right)  }%
{\pi_{1,m}}+\frac{4\left(  \left\Vert \omega\right\Vert _{\infty}+\left\Vert
\omega\right\Vert _{\mathrm{TV}}\right)  }{\tau_{m}u_{m,0}}\right]  .
\label{eq2ix5}%
\end{align}
Therefore,%
\begin{equation}
\sup\nolimits_{\boldsymbol{\mu}\in\mathcal{B}_{m}\left(  \rho_{m}\right)
}\left(  \left\vert \frac{\sup\nolimits_{t\in\left[  0,\tau_{m}\right]  }%
\hat{\varphi}_{0,m}\left(  t,\mathbf{z}\right)  }{\pi_{1,m}}-1\right\vert
\right)  \rightsquigarrow0\text{ \ as \ }m\rightarrow\infty\label{eq3ix2}%
\end{equation}
if
$\textcolor{blue}{\sup_{\boldsymbol{\mu}\in\mathcal{B}_{m}\left(  \rho_{m}\right)
}}\pi_{1,m}^{-1}\tilde{\Upsilon}\left(  q,\tau_{m},\gamma_{m},r_{\mu_{0}%
}\right)  \rightarrow0$ and $\tau_{m}%
\textcolor{blue}{\inf_{\boldsymbol{\mu}\in\mathcal{B}_{m}\left(  \rho_{m}\right)
}}u_{m,0}\rightarrow\infty$ as $m\rightarrow\infty$. Consequently,%
\[
\mathcal{Q}\left(  \mathcal{F}\right)  =\left\{
\begin{array}
[c]{c}%
q\gamma^{\prime}>\vartheta>2^{-1},\gamma^{\prime}>0,\gamma^{\prime\prime
}>0,0\leq\vartheta^{\prime}<\vartheta-1/2,\\
R_{0}\left(  \rho_{m},\mu_{0}\right)  =O\left(  m^{\vartheta^{\prime}}\right)
,\tau_{m}\leq\gamma_{m},u_{m,0}\geq\frac{\ln\ln m}{\gamma^{\prime\prime}%
\tau_{m}},\\
t\in\left[  0,\tau_{m}\right]  ,\lim\limits_{m\rightarrow\infty}
\textcolor{blue}{\sup_{\boldsymbol{\mu}\in\mathcal{B}_{m}\left(  \rho_{m}\right)
}}\pi_{1,m}%
^{-1}\tilde{\Upsilon}\left(  q,\tau_{m},\gamma_{m},r_{\mu_{0}}\right)  =0
\end{array}
\right\}
\]
is a uniform consistency class. Since \textquotedblleft$\pi_{1,m}^{-1}%
\tilde{\Upsilon}\left(  q,\tau_{m},\gamma_{m},r_{\mu_{0}}\right)
\rightarrow0$\textquotedblright\ controls the uniformity convergence in
(\ref{eq3ix2}) with respect to $\pi_{1,m}$, the convergence in (\ref{eq3ix2})
can also be interpreted as follows. Let $\left\{  a_{m}\right\}  _{m\geq1}$ be
a positive sequence such that $a_{m}\geq a_{m+1}$ and $\lim_{m\rightarrow
\infty}a_{m}=0$, and let $\tilde{C}>0$ be a constant. Define%
\[
\mathcal{U}_{0,m}=\left\{  \boldsymbol{\mu}\in\mathcal{B}_{m}\left(  \rho
_{m}\right)  :%
\begin{array}
[c]{c}%
R_{0}\left(  \rho_{m},\mu_{0}\right)  \leq\tilde{C}m^{\vartheta^{\prime}}\\
\pi_{1,m}^{-1}\tilde{\Upsilon}\left(  q,\tau_{m},\gamma_{m},r_{\mu_{0}%
}\right)  \leq a_{m},u_{m,0}\geq\frac{\ln\ln m}{\gamma^{\prime\prime}\tau_{m}}%
\end{array}
\right\}
\]
with the constants introduced by $\mathcal{Q}\left(  \mathcal{F}\right)  $.
Then (\ref{eq3ix2}) implies%
\[
\sup\nolimits_{\boldsymbol{\mu}\in\mathcal{U}_{0,m}}\left(  \left\vert
\frac{\sup\nolimits_{t\in\left[  0,\tau_{m}\right]  }\hat{\varphi}%
_{0,m}\left(  t,\mathbf{z}\right)  }{\pi_{1,m}}-1\right\vert \right)
\rightsquigarrow0\text{ \ as \ }m\rightarrow\infty.
\]
\qed



\end{document}